\documentclass[11pt,reqno]{amsart}

\usepackage{amssymb}
\usepackage{graphicx}
\usepackage{epsfig}
\usepackage{color}

\catcode`\@=11

\long\def\@savemarbox#1#2{\global\setbox#1\vtop{\hsize\marginparwidth 
  \@parboxrestore\tiny\raggedright #2}}
\newcommand\lref[1]{\ref{#1}%
\@ifundefined{r@DisplaY #1}{}{ (#1)}}

\newcommand\fakelabel[2]{\@bsphack\if@filesw {\let\thepage\relax
   \newcommand\protect{\noexpand\noexpand\noexpand}%
\xdef\@gtempa{\write\@auxout{\string
      \newlabel{#1}{{#2}{\thepage}}}}}\@gtempa
   \if@nobreak \ifvmode\nobreak\fi\fi\fi\@esphack}

\catcode`\@=12

\def\Empty{}
\newcommand\oplabel[1]{
  \def\OpArg{#1} \ifx \OpArg\Empty {} \else
        \label{#1}
  \fi}
\newtheorem{theoremSt}{Theorem}[section]

\newtheorem{exampleSt}[theoremSt]{Example}
\newtheorem{exerciseSt}[theoremSt]{Exercise}

\newcommand\MakeStEnv[1]{
  \newenvironment{#1}[1]{
  \begin{#1St} \oplabel{##1}%
  \global\def\CrntSt{\thetheoremSt}%
}{ 
  \end{#1St} }
  \newenvironment{#1+}[1]{
  \begin{#1St} \label{##1}%
  \label{DisplaY ##1}%
  \global\def\CrntSt{\thetheoremSt}%
  \def\Labl{##1}\ifx\Labl\Empty{} \else {\em (\Labl)\,}\fi%
}{ 
  \end{#1St} }
}
\MakeStEnv{theorem}
\MakeStEnv{corollary}
\MakeStEnv{proposition}
\MakeStEnv{lemma}
\MakeStEnv{definition}
\MakeStEnv{conjecture}
\MakeStEnv{problem}
\MakeStEnv{question}

\newenvironment{example}[1]{
  \begin{exampleSt} \oplabel{#1}%
  \global\def\CrntSt{\thetheoremSt}%
  \rm %
}{ 
  \end{exampleSt} }

\long\def\state#1#2{
\medskip\par\noindent
{\bf #1} 
{\it #2}
\par\medskip
}

\long\def\realfig#1#2{
\begin{figure}[htbp]
\centerline{\includegraphics{#1.pdf}}
\caption[#1]{#2}
\oplabel{#1}
\end{figure}}

\newlength{\saveu}

\newenvironment{pf}{%
 \begin{proof}%
}{ 
 \end{proof}
}
\newenvironment{pf*}[1]{%
 \begin{proof}[#1]%
}{ 
 \end{proof}
}

\newcommand{\finishproof}[1]{ 
  \def\FPArg{#1}
  \ifx\FPArg\Empty
        \newcommand\FPArg{\CrntSt}  \fi
  \smallbreak\noindent\makebox[\textwidth]{\hfill\fbox{\FPArg}}
  \medbreak\noindent
}

\newcommand\AAA{{\mathcal A}}
\newcommand\BB{{\mathcal B}}
\newcommand\CC{{\mathcal C}}
\newcommand\DD{{\mathcal D}}
\newcommand\EE{{\mathcal E}}
\newcommand\FF{{\mathcal F}}
\newcommand\GG{{\mathcal G}}

\newcommand\II{{\mathcal I}}

\newcommand\KK{{\mathcal K}}
\newcommand\LL{{\mathcal L}}
\newcommand\MM{{\mathcal M}}
\newcommand\NN{{\mathcal N}}

\newcommand\PP{{\mathcal P}}

\newcommand\RR{{\mathcal R}}

\newcommand\UU{{\mathcal U}}
\newcommand\VV{{\mathcal V}}
\newcommand\WW{{\mathcal W}}
\newcommand\XX{{\mathcal X}}
\newcommand\YY{{\mathcal Y}}
\newcommand\ZZ{{\mathcal Z}}

\newcommand\PMF{{\PP\kern-2pt\MM\FF}}

\newcommand\PML{{\PP\kern-2pt\MM\LL}}
\newcommand\GL{{\GG\LL}}

\newcommand\half{{\textstyle{\frac12}}}

\newcommand\ep{\epsilon}

\newcommand\hhat{\widehat}

\newcommand\union{\cup}
\newcommand\intersect{\cap}
\newcommand\bbR{{\mathord{\text{I\kern-2pt R}}}}        
\newcommand\bbH{{\mathord{\text{I\kern-2pt H}}}}        

\newcommand\C{{\mathbb C}}

\newcommand\Z{{\mathbb Z}}
\newcommand\R{{\mathbb R}}

\newcommand\Hyp{{\mathbb H}}

\newcommand\PSL[1]{\text{PSL}_{#1}}

\newcommand\bigrightarrow[1]{\hbox to #1{\rightarrowfill}}
\newcommand\bigleftarrow[1]{\hbox to #1{\leftarrowfill}}
\newcommand\Iff{\Longleftrightarrow}
\newcommand\homeo{\cong}
\newcommand\boundary{\partial}
\newcommand\semidir{\mathrel{\hbox{\vrule depth-.03ex height1.1ex\kern-0.15em$\times$}}}

\newcommand\til{\widetilde}

\newcommand{\diam}{\operatorname{diam}}

\renewcommand{\Re}{\operatorname{Re}}
\renewcommand{\Im}{\operatorname{Im}}

\numberwithin{equation}{section}

\catcode`\@=11

\def\subsection{\@startsection{subsection}{2}%
  \z@{.5\linespacing\@plus.7\linespacing}{.5em}%
  {\normalfont\bfseries\centering}}

\def\section{\@startsection{section}{1}%
  \z@{.7\linespacing\@plus\linespacing}{.5\linespacing}%
  {\normalfont\large\bfseries\centering}}

\def\subsubsection{\@startsection{subsubsection}{3}%
  \z@{.5\linespacing\@plus.7\linespacing}{-.5em}%
  {\normalfont\bfseries}}

\catcode`\@=12

\newcommand{\T}{{\mathbf T}}
\newcommand{\inj}{\operatorname{inj}}

\newcommand{\collar}{\operatorname{\mathbf{collar}}}
\newcommand{\bcollar}{\operatorname{\overline{\mathbf{collar}}}}

\newcommand{\I}{{\mathbf I}}
\newcommand{\tprec}{\prec_t}

\newcommand{\pprec}{\prec_p}
\newcommand{\ppreceq}{\preceq_p}
\newcommand{\sprec}{\prec_s}
\newcommand{\cpreceq}{\preceq_c}
\newcommand{\cprec}{\prec_c}
\newcommand{\topprec}{\prec_{\rm top}}

\newcommand{\fsub}{\mathrel{\scriptstyle\searrow}}
\newcommand{\bsub}{\mathrel{\scriptstyle\swarrow}}
\newcommand{\fsubd}{\mathrel{{\scriptstyle\searrow}\kern-1ex^d\kern0.5ex}}
\newcommand{\bsubd}{\mathrel{{\scriptstyle\swarrow}\kern-1.6ex^d\kern0.8ex}}
\newcommand{\fsubeq}{\mathrel{\raise-.7ex\hbox{$\overset{\searrow}{=}$}}}
\newcommand{\bsubeq}{\mathrel{\raise-.7ex\hbox{$\overset{\swarrow}{=}$}}}

\newcommand{\base}{\operatorname{base}}

\newcommand{\bbar}{\overline}

\newcommand{\EL}{\mathcal{EL}}

\newcommand{\tsh}[1]{\left\{\kern-.9ex\left\{#1\right\}\kern-.9ex\right\}}

\newcommand{\vsucc}{\operatorname{succ}}

\newcommand\sbtop{_{\text{top}}}
\newcommand\sbot{_{\text{bot}}}

\newcommand\boundvert{{\boundary_{||}}}
\newcommand\storus[1]{U(#1)}

\newcommand\modl{M_\nu}
\newcommand\MT{{\mathbb T}}
\newcommand\Teich{{\mathcal T}}

\newcommand{\towall}{\mathrel{\dashv}}
\newcommand{\wallto}{\mathrel{\vdash}}
\newcommand{\wall}{\mathrel{|}}
\newcommand{\stprec}{\mathrel{\vee}}
\newcommand{\lands}{\operatorname{\YY}}
\newcommand\add[1]{{\langle #1 \rangle}}

\newcommand\interior{{\rm int}}

\newcommand\sitem{\medskip\item}

\begin{document}

\title[The classification of Kleinian surface groups, II]{The
classification of Kleinian surface groups, II: The Ending Lamination
Conjecture}

\author{Jeffrey F. Brock}
\address{Brown University}
\author{Richard D. Canary}
\address{University of Michigan}
\author{Yair N. Minsky}
\address{Yale University}
\date{\today}
\thanks{Partially supported by NSF grants 
 DMS-0203698, DMS-0203976,   DMS-0354288,
DMS-0906229, DMS-1006298 and DMS-1005973}


\begin{abstract}
  Thurston's Ending Lamination Conjecture states that a hyperbolic
  3-manifold $N$ with finitely generated fundamental group is uniquely
  determined by its topological type and its end invariants.  In this
  paper we prove this conjecture for Kleinian surface groups; the
  general case when $N$ has incompressible ends relative to its cusps
  follows readily.  The main ingredient is a uniformly bilipschitz
  model for the quotient of $\Hyp^3$ by a Kleinian surface group.  The
  first half of the proof appeared in \cite{minsky:ELCI}, and a
  subsequent paper \cite{brock-canary-minsky:ELCIII} will establish
  the Ending Lamination Conjecture in general.
\end{abstract}

\maketitle

\setcounter{tocdepth}{1}
\tableofcontents

\newcommand\epzero{{\ep_0}}
\newcommand\epone{{\ep_1}}
\newcommand\epotal{{\ep_{\rm u}}}
\newcommand\kotal{k_{\rm u}}
\newcommand\Kmodel{K_0}
\newcommand\Kone{K_1}
\newcommand\Ktwo{K_2}
\newcommand\Lone{L_0}
\newcommand\bdry{\partial} 
\newcommand\stab{\operatorname{stab}}
\newcommand\nslices[2]{#2|_{#1}}
\newcommand\ME{M\kern-4pt E}
\newcommand\bME{\overline{M\kern-4pt E}}

\section{The ending lamination conjecture}
\label{intro}

In the late 1970's Thurston formulated a conjectural classification
scheme for all hyperbolic 3-manifolds with finitely generated
fundamental group. The picture proposed by Thurston generalized what
had been previously understood through the work of Ahlfors
\cite{ahlfors-bers}, Bers \cite{bers:spaces}, Kra \cite{kra:spaces},
Marden \cite{marden:geometry}, Maskit \cite{maskit:self}, Mostow
\cite{mostow:hyperbolic}, Prasad \cite{prasad}, Thurston
\cite{wpt:bull} and others about geometrically finite hyperbolic
3-manifolds.

Thurston's scheme proposes {\em end invariants} that encode the
asymptotic geometry of the ends of the manifold, generalizing the role
the Riemann surfaces at infinity play in the geometrically finite
case. More precisely, the following conjecture appears in
\cite{wpt:bull}:

\state{Ending Lamination Conjecture.}{%
A hyperbolic 3-manifold with finitely generated fundamental group is 
uniquely determined by its topological type and its end invariants.
}

This paper is the second in a series of three which will establish the
Ending Lamination Conjecture for all topologically tame hyperbolic
3-manifolds.  For expository material on this conjecture, and on the
proofs in this paper and in \cite{minsky:ELCI}, we direct the reader
to \cite{minsky:knoxville},\cite{minsky:warwick} and
\cite{minsky:cdm}. We also note that 
Bowditch \cite{bowditch:elcgen}, Rees \cite{rees:elc} and Soma
\cite{soma:elc} have meanwhile written alternate proofs of the
conjecture, in which various aspects have been simplified.  

Together with the recent proofs of Marden's Tameness Conjecture by
Agol \cite{agol:tame} and Calegari-Gabai \cite{calegari-gabai:tame},
this gives a complete classification of all hyperbolic 3-manifolds
with finitely-generated fundamental group.

\medskip

In this paper we will focus on the {\em surface group} case.  A
{\em Kleinian surface group} is a discrete, faithful representation
$\rho:\pi_1(S)\to \PSL 2(\C)$ where $S$ is a compact orientable
surface, such that the restriction of $\rho$ to any boundary loop has
parabolic image. These groups arise naturally as restrictions of more
general Kleinian groups to surface subgroups. Bonahon \cite{bonahon}
and Thurston \cite{wpt:notes} showed that the associated hyperbolic
3-manifold $N_\rho=\Hyp^3/\rho(\pi_1(S))$ is homeomorphic to
$\interior(S)\times \R$ and that $\rho$ has a well-defined pair of end
invariants $(\nu_+,\nu_-)$. Typically, each end invariant is either a
point in the Teichm\"uller space of $S$ or a geodesic lamination on
$S$. In the general situation, each end invariant is a geodesic
lamination on some (possibly empty) subsurface of $S$ and a conformal
structure on the complementary surface. We will prove:

\state{Ending Lamination Theorem for Surface Groups.}{%
A Kleinian surface group $\rho$ is uniquely determined, up to
conjugacy in $\PSL 2(\C)$, by its end invariants. 
}

The main technical result which leads to the Ending Lamination Theorem
is the Bilipschitz Model Theorem, which gives a bilipschitz homeomorphism
from a ``model manifold'' $\modl$ to the hyperbolic manifold $N_\rho$
(See \S\ref{bilip model intro} for a precise statement). The model
$\modl$ was 
constructed in Minsky \cite{minsky:ELCI}, and its crucial property is
that it depends only on the end invariants $\nu=(\nu_+,\nu_-)$, and not on
$\rho$ itself. 
(Actually $\modl$ is mapped to the ``augmented convex core'' of $N_\rho$,
but as this is the same as $N_\rho$ in the main case of interest, we
will ignore the distinction for the rest of the introduction.
See \S\ref{bilip model intro} for details).

The proof of the Bilipschitz Model Theorem will be completed in
Section \ref{proofmain}, and the Ending Lamination Conjecture will be
obtained as a consequence of this and Sullivan's rigidity theorem in
Section \ref{ELT proof}.

\medskip

The surface group case bears directly on the more general setting of
hyperbolic 3-manifolds with finitely generated fundamental group and
{\em incompressible ends}, which we now describe. If $N$ is a hyperbolic
3-manifold with finitely generated group, it is natural to excise a standard 
open neighborhood $\PP$ of the
cusps of $N$ to obtain $$N^0 =
N  \setminus \PP.$$
A relative compact core  $K$ for $N^0$ is a compact submanifold whose
inclusion into
$N^0$ is a homotopy equivalence and whose intersection with each
component of $\boundary \PP$ includes by a homotopy equivalence into
that component. Then $P=K\cap\boundary\PP$  is the {\em parabolic locus}
and $\boundary_0K=\partial K\setminus P$ is called the {\em relative boundary} of
$K$. If each component of the relative boundary is incompressible, then $N^0$ is
said to have {\em incompressible ends}.
In this case,
Bonahon's tameness theorem \cite{bonahon} guarantees that
$N^0\setminus K$ is homeomorphic to $\boundary_0K\times (0,\infty)$.
Then each end $\EE$ has fundamental group a Kleinian surface group. One end
of the associated manifold is a homeomorphic lift of $\EE$, and we
associate its end invariant (a lamination or a point in a
Teichm\"uller space) to the corresponding component of $\boundary_0 K$.
The ending lamination theorem for surface groups, together with a
short topological argument, 
gives the following generalization.

\state{Ending Lamination Theorem for Incompressible Ends.}{%
Let $G$ be a finitely generated, torsion-free, non-abelian group.
If $\rho \colon G \to \PSL 2(\C)$ is a discrete faithful representation so that
$N_\rho^0$ has incompressible ends, then $\rho$ is
  determined, up to conjugacy in $\PSL 2(\C)$, by the marked
  homeomorphism type of its relative compact core and the end
  invariants associated to the ends of $N_\rho^0$.  }

The first part of the proof of the ending lamination theorem for
surface groups appeared in \cite{minsky:ELCI}, and we will refer to
that paper for some of the background and notation, although we will
strive to make this paper readable independently.
Section~\ref{general case} provides a discussion of the proof of the
general Ending Lamination Conjecture, which will appear in
\cite{brock-canary-minsky:ELCIII}.

\subsection{Corollaries}

A positive answer to the Ending Lamination Conjecture allows one to
settle a number of fundamental questions about the structure of
Kleinian groups and their deformation spaces.

The Bers-Sullivan-Thurston density conjecture predicts that every
finitely generated Kleinian group is an algebraic limit of
geometrically finite groups.  
In the surface group case, the density
conjecture follows immediately from our main theorem and results of
Thurston \cite{wpt:II} and Ohshika \cite{ohshika:ending-lams}. We
recall that $AH(S)$ is the space of conjugacy classes of Kleinian
surface groups and that a surface group is {\em quasifuchsian} if $N_\rho^0$
has precisely two ends, each of which is geometrically finite.

\state{Density Theorem for surface groups.}{%
The set of quasifuchsian surface groups is dense in $AH(S)$.
}

Marden \cite{marden:geometry} and Sullivan \cite{sullivan:QCDII}
showed that the interior of $AH(S)$ consists exactly of the
quasifuchsian groups.  Bromberg \cite{bromberg:density} and
Brock-Bromberg \cite{brock-bromberg:density} previously showed that
each representation $\rho$ whose image contains no parabolic elements
and for which $N_\rho$ has incompressible ends is an algebraic limit
of geometrically finite representations, using cone-manifold
techniques and the bounded-geometry version of the Ending Lamination
Conjecture in Minsky \cite{minsky:boundgeom}. 

When $M$ has incompressible boundary the set $AH(M)$ of discrete
faithful representations $\rho \colon \pi_1(M) \to \PSL 2(\C)$ plays
the role of the deformation space, and the above theorems of Marden
and Sullivan guarantee its interior consists of geometrically finite
representations such that every parabolic in their image is associate
to a curve in a toroidal boundary component of $M$.
Then we have the
following generalization of the density theorem.  \state{Density
  Theorem for Incompressible Boundary.}{%
  Let $M$ be a compact 3-manifold with incompressible 
  boundary.  Then we have
$$\overline{\interior(AH(M))} = AH(M).$$ }
A more general density theorem holds in the setting of deformation
spaces of pared manifolds with specified parabolic locus. 
We discuss this in section~\ref{corollaries}.
The general version of the Ending Lamination Theorem is a crucial ingredient
in the resolution of the complete Bers-Sullivan-Thurston Density Conjecture
(see  the sequel \cite{brock-canary-minsky:ELCIII}, Namazi-Souto 
\cite{namazi-souto:density} and Ohshika \cite{ohshika:density,ohshika:constructing-limits} for more
details).

\medskip

The Density Theorem has important consequences for the global topology
of $AH(M)$.  If $M$ has incompressible boundary the components of
$AH(M)$ are enumerated by the set ${\mathcal {A}}(M)$ of (marked)
homeomorphism types of (marked) compact 3-manifolds homotopy
equivalent to $M$ (see \cite{canary-mccullough}).  Anderson, Canary
and McCullough \cite{anderson-canary-mccullough:bumping} introduced a
finite-to-one equivalence relation on ${\mathcal{A}}(M)$, called
primitive shuffle equivalence, and proved that the components of
$\overline{\interior(AH(M))}$ are enumerated by the set
$\hat{\mathcal{A}}(M)$ of equivalence classes with respect to this
equivalence relation.  (Roughly, primitive shuffle equivalences are
homotopy equivalences which are allowed to rearrange the order in
which components of the complement of the characteristic submanifold
are attached to certain solid torus components of the characteristic
submanifold.)  It follows from the Density Theorem above that
components of $AH(M)$ are enumerated by $\hat{\mathcal{A}}(M)$. In
particular, applying results from \cite{canary-mccullough}, one sees
that $AH(M)$ has infinitely many components if and only if there is a
thickened torus component $V$ of the characteristic submanifold of $M$
such that $V\cap \boundary M$ has at least 3 components ($M$ has {\em
  double trouble}).

\medskip

We also obtain a quasiconformal rigidity theorem that gives 
a  common generalization of Mostow's \cite{mostow:hyperbolic} and
Sullivan's \cite{sullivan:rigidity} rigidity theorems.

\state{Rigidity Theorem.}{%
Let $G$ be a finitely generated, torsion-free, non-abelian group.
 If $\rho$ and $\rho'$ are two discrete faithful representations of
  $G$ into $\PSL 2(\C)$ that are conjugate by an
  orientation-preserving homeomorphism of $\hhat \C$ and $N_\rho^0$
  has incompressible ends, then $\rho$ and $\rho'$ are
  quasiconformally conjugate.  }

\medskip

Though a central motivation for producing the model manifold lay in
its application to the ending lamination conjecture and other
deformation theoretic questions, the existence of a model manifold for
the ends of $N^0$ guarantees various quantitative geometric features
of independent interest.  As a key example, we establish McMullen's
conjecture that the volume of the thick part of the convex core
of a hyperbolic 3-manifold grows polynomially.

More precisely, if $x$ lies in the thick part of the convex core
$C_N$, then let $B_r^{thick}(x)$ be the set of points in the
$\ep_1$-thick part of $C_N$ which can be joined to $x$ by a path of
length at most $r$ lying entirely in the $\ep_1$-thick part.

Given a compact connected surface $S$ with $\chi(S) <0$, let
$$
d(S) = \begin{cases}
  -\chi(S) &  \text{genus}(S)>0 \\
   -\chi(S)-1   &  \text{genus}(S)=0.
\end{cases}
$$
When $S = R_1 \sqcup \ldots \sqcup R_k$ is disconnected, we
define $d(S) = \max_{i = 1}^k d(R_i)$.

\state{Volume Growth Theorem.}{%
If $N$ is the quotient of a Kleinian surface group $\rho\in\DD(S)$,
then for any $x$ in the $\ep_1$-thick part of the convex core $C_N$
and $r\ge 1$ we have
$${\rm volume}\left( B^{thick}_r(x)\right) \le c_1r^{d(S)},$$
where $c_1$ depends only on the topological type of $S$. 

In general, if $N$ is a complete hyperbolic 3-manifold with 
relative compact core $(K,P)$ so that $N^0$ has incompressible ends, 
we have
$${\rm volume}\left( B^{thick}_r(x)\right) \le c_1r^{d(\bdry_0 K)} + c_2,$$
where $c_1$ depends only on
the topological type of $\bdry_0 K$, and
$c_2$ depends on the hyperbolic structure of $N$.
}

A different proof of the Volume Growth Theorem is given by
Bowditch in \cite{bowditch:bands}. We are grateful to Bowditch for
pointing out an error in our original definition of $d(S)$.

Proofs of these corollaries are given in Section \ref{corollaries}.
Each of them admits generalizations to the setting of all finitely
generated Kleinian groups and these generalizations will be discussed
in \cite{brock-canary-minsky:ELCIII}. 

In Section \ref{corollaries}, we also prove 
the Length Bound Theorem, which gives
estimates on the lengths of short geodesics in a Kleinian surface group
manifold
(see \S\ref{length estimates} for the statement).

We also remark that using the bi-Lipschitz model theorem for surface
groups and the tameness theorem of
\cite{agol:tame,calegari-gabai:tame} Mahan Mj has announced a proof of
local connectivity for limit sets of finitely generated Kleinian
groups, as well as many other related results concerning the existence
and behavior of Cannon-Thurston maps from the boundary of the Kleinian
group to $\hat\C$.

\subsection{Outline of the proof}
\label{outline}

The Lipschitz Model Theorem, from \cite{minsky:ELCI},
provides a degree 1 homotopy equivalence from the model manifold $\modl$ to
the hyperbolic manifold $N_\rho$
(or in general to the augmented core of $N_\rho$, but we ignore the
distinction in this outline), which respects the thick-thin
decompositions of $\modl$ and $N_\rho$ and is Lipschitz on the thick part of
$\modl$ (see \S\ref{bilip model intro}).

Our main task in this paper is to promote this map to a {\em bilipschitz}
homeomorphism between $\modl$ and $N_\rho$ ,
and this is the content of our main result, the
Bilipschitz Model Theorem.
The proof of the Bilipschitz Model Theorem 
converts the Lipschitz model map to a bilipschitz map
incrementally on various subsets of the model. The main ideas of 
the proof can be summarized as follows:

\subsubsection*{Topological order of subsurfaces}  
In section \ref{knotting} we discuss a ``topological order relation'' among
embedded surfaces in a product 3-manifold $S\times\R$. 
This is the intuitive notion that one surface may lie ``below'' another 
in this product, but this relation does not in fact induce a partial
order and hence a number of technical issues arise. 

We introduce an object called a {\em scaffold}, which is a subset of
$S\times\R$ consisting of a union of unknotted 
solid tori and surfaces in $S\times\R$, each isotopic to a level
subsurface, satisfying certain conditions.
The main theorem in this section is the Scaffold Extension Theorem
(\ref{Scaffold Extension}), which states that, under appropriate
conditions (in particular an ``order-preservation'' condition),
embeddings of a scaffold into $S\times\R$ can be 
extended to global homeomorphisms of $S\times\R$.

Much of the rest of the proof is concerned with analyzing this order
in the model manifold, breaking the model up into pieces separated by
scaffolds, and ensuring that the model map satisfies the appropriate
order-preserving condition.

\subsubsection*{Structure of the model: tubes, surfaces and regions}
The structure of the model $\modl$ is organized by the structure of a
{\em hierarchy of geodesics in the complex of curves} as developed in
\cite{masur-minsky:complex2}, applied in \cite{minsky:ELCI} and
summarized here in \S\ref{hierarchy background}.  In particular, such
a hierarchy, which depends only on the end-invariant data $\nu$,
directly produces a combinatorial 3-manifold $\modl$ homeomorphic to
$S\times\R$ containing a collection of unknotted solid tori that
correspond to the Margulis tubes for short geodesics in $N_\rho$.  The
Lipschitz Model Theorem produces a lipschitz map of the complement of
these tubes to the complement of the corresponding Margulis tubes in
$N_\rho$ that extends to a proper map on each tube.  The model also
contains a large family of {\em split-level surfaces,} namely, surfaces
isotopic to level subsurfaces in $S\times\R$ and
bilipschitz-homeomorphic to bounded-geometry hyperbolic surfaces.
These correspond to {\em slices} of the hierarchy.

In Section \ref{po} we discuss {\em cut systems} in this hierarchy.  A
cut system gives rise to a family of split-level surfaces and we show,
in Lemma \ref{Topological Partial Order}, that one can impose spacing
conditions on cut systems so that after a thinning process the
topological order relation restricted to the split-level surfaces
coming from the cut system generates a partial order.

In Section \ref{regions} we will show how the surfaces of a such a cut system
(together with the model tubes) cut the model into regions whose
geometry is controlled. The collection of split-level surfaces and
Margulis tubes bounding such a region form a  scaffold.

\subsubsection*{Uniform embeddings of model surfaces}
The restriction of the model map to a split-level surface is
essentially a Lipschitz map of a bounded-geometry hyperbolic surface
whose boundary components map to Margulis tubes (we call this an {\em
  anchored surface}). These surfaces are not necessarily themselves
embedded, but we will show that they may be deformed in a controlled
way to bilipschitz embeddings.

In general, a Lipschitz anchored surface may be {\em wrapped}
around a deep Margulis tube in $N_\rho$ and any homotopy to
an embedding must pass through the core of this tube.
In Theorem \ref{Relative embeddability}, we show that this
wrapping phenomenon is the only obstruction to a controlled
homotopy. The proof relies on a geometric limiting argument and
techniques of Anderson-Canary \cite{anderson-canary:cores}.
In section \ref{embedding an individual cut}, we check that we may
choose the spacing constants for our cut system, so that
the associated Lipschitz anchored surfaces are not wrapped, and hence
can be uniformly embedded.

\subsubsection*{Preservation of topological order}
In section \ref{thinning the cut system} we show that any cut system
may be ``thinned'' in a controlled way to yield a new cut system, with
uniform spacing constants, so that if two split-level surfaces lie on
the boundary of the same complementary region in $\modl$, then their
associated anchored embeddings in $N_\rho$ (from \S 8.2) are disjoint.
We adjust the model map so that it is a bilipschitz embedding on
collar neighborhoods of these split-level surfaces.

In \ref{preserving order of embeddings}, 
we check that if two anchored surfaces, associated to the thinned cut system,
are disjoint and ordered in the hyperbolic manifold, then their relative ordering agrees
with the ordering of the associated split-level surfaces in the model.
The idea is to locate {\em insulating regions} in the model --
geometrically defined subsets of the manifold that separate the two  surfaces,
and on which there is sufficient control to show that the topological
order between the insulating region and each of the two surfaces is
preserved. A transitivity argument can then be used to show
that the order between the surfaces is preserved as well. 

The insulating regions are of two types. Sometimes there is a model tube between
the associated split-level surfaces in the topological ordering
and it is fairly immediate from properties of the model map
that its image Margulis tube has the correct separation
properties.  When such a tube is not available we show, in
Theorem \ref{product region}, that there exist
certain {\em subsurface product regions}, product interval bundles
over subsurfaces of $S$, that are bilipschitz to subsets of
bounded-geometry surface group manifolds based on lower-complexity 
surfaces. The control over these regions is obtained by a geometric
limit argument.

\subsubsection*{Bilipschitz extension to the regions}
The union of the split-level surfaces and the solid tori divide the
model manifold up into regions bounded by scaffolds.  The Scaffold
Extension Theorem can be used to show that the embeddings on the
split-level surfaces can be extended to embeddings of these
complementary regions. An additional geometric limit argument, given
in section~\ref{region control}, is needed to obtain bilipschitz
bounds on each of these embeddings.  Piecing together the embeddings,
we obtain a bilipschitz embedding of the ``thick part'' of the model
to the thick part of $N_\rho$. A final brief argument, given in
section~\ref{tube control} shows that the map can be extended also to
the model tubes in a uniform way.  This completes the proof.

This outline ignores the case when the convex hull of $N_\rho$ has
nonempty boundary, and in fact most of the proof on a first reading is
improved by ignoring this case. Dealing with the boundary is mostly an
issue of notation and some attention to special cases; nothing
essentially new happens. In Section~\ref{proofmain} most details of
the case with boundary are postponed to~\S\ref{hybrid}.

\subsection{The general case of the Ending Lamination Conjecture}
\label{general case}

In this section we briefly discuss the proof of the Ending Lamination
Conjecture in the general situation. Details of the Ending Lamination
Theorem for incompressible
ends will appear at the conclusion of the paper, and the general
case will appear in \cite{brock-canary-minsky:ELCIII}.

As above, if $N$ is a hyperbolic 3-manifold with finitely generated
fundamental group, we and let $N^0$ be the complement of standard open
neighborhoods of the cusps of $N$, and let $(K,P)$ denote a relative
compact core for $N^0$.

\subsubsection*{Incompressible ends case}
In the setting of incompressible ends, where each component of
$\boundary_0 K = \boundary K\setminus P$ is incompressible, the
derivation of the Ending Lamination Conjecture from the surface group
case is fairly straightforward.  In this case the restriction of
$\pi_1(N)$ to the fundamental group of any component of $R$ of
$\boundary_0 K$ is a Kleinian surface group.  The Bilipschitz Model
Theorem applies to the cover $N_R$ of $N$ associated to $\pi_1(R)$ to
give a model for $N_R$, and one end of $N_R^0$ embeds isometrically
under the covering projection to the end of $N^0$ cut off by $R$.  In
this way we obtain bilipschitz models for each of the ends of $N^0$.

Two homeomorphic hyperbolic manifolds with the same end invariants
must have the same cusps, and a bilipschitz correspondence between
their ends (the end invariant data specify the cusps so that after
removing the cusps the manifolds remain homeomorphic).  Since what
remains is compact, one may extend the bilipschitz homeomorphism on
the ends to a bilipschitz homeomorphism on the non-cuspidal part,
which in turn extends to a global bilipschitz homeomorphism.  One
again applies Sullivan's rigidity theorem \cite{sullivan:rigidity} to
complete the proof.

\subsubsection*{Compressible boundary case}

When some component $R$ of $\boundary K \setminus P$ is compressible,
the subgroup $\pi_1(R)$ is no longer a Kleinian surface group.  Agol
\cite{agol:tame} and Calegari-Gabai \cite{calegari-gabai:tame} proved
that $N$ is homeomorphic to the interior of $K$. Canary
\cite{canary:ends} showed that the ending invariants are well-defined
in this setting.

The first step of the proof in this case is to apply Canary's
branched-cover trick from \cite{canary:ends}. That is, we find a
suitable closed geodesic $\gamma$ in $N$ and a double branched cover
$\pi:\hat N \to N$ over $\gamma$, such that $\pi_1(\hat N)$ is freely
indecomposable.  The singularities on the branching locus can be
smoothed locally to give a pinched negative curvature metric on $\hat
N$. Since $N$ is topologically tame, one may choose a relative compact
core $K$ for $N^0$ containing $\gamma$, so that $\hat K=\pi^{-1}(K)$
is a relative compact core for $\hat N^0$.  Let $P=\boundary N^0\cap
K$ and $\hat P=\boundary \hat N^0\cap \hat K$.  If $R$ is any
component of $\partial K-P$, then $\pi^{-1}(R)$ consists of two
homeomorphic copies $\hat R_1$ and $\hat R_2$ of $R$, each of which is
incompressible.

Given a component $R$ of $\boundary K -P$, we consider the cover
$\hat N_R$ of $\hat N$ associated to $\pi_1(\hat R_1)$. We then apply the
techniques of \cite{minsky:ELCI} and this paper to obtain a bilipschitz
model for some neighborhood of the end $E_R$ of $\hat N_R^0$ cut off by
$\hat R_1$.  In particular, we need to check that the estimates of 
\cite{minsky:ELCI} apply in suitable neighborhoods of $E_R$.
The key tool we will need is a generalization of 
Thurston's Uniform Injectivity Theorem \cite{wpt:I}
for pleated surfaces to this setting (see also
Namazi \cite{namazi:thesis} and
Namazi-Souto \cite{namazi-souto:uniform-inj}). See
Miyachi-Ohshika \cite{miyachi-ohshika:tamebg} for a discussion of
this line of argument in the ``bounded geometry'' case.

Once we obtain a bilipschitz model for some neighborhood of $E_R$ it
projects down to give  a bilipschitz model for a neighborhood of the end
of $N^0$ cut off by $R$. As before, we obtain a bilipschitz model
for the complement of a compact submanifold of $N^0$ and the proof
proceeds as in the incompressible boundary setting.

It is worth noting that this construction does not yield a {\em
  uniform} model for $N$, in the sense that the bilipschitz constants
  depend on the geometry of $N$ and not only on its topological type
  (for example on the details of what happens in the branched covering
  step). 
The model we develop here for the surface group case is uniform, and
  we expect that in the incompressible-boundary case uniformity of the
  model should not be too hard to obtain. Uniformity in general is
  quite an interesting problem, and would be useful for further applications
  of the model manifold.

\subsubsection*{Acknowledgements} The authors would like to thank the
referee for many helpful suggestions.

\section{Background and statements}
\label{background}

In this section we will introduce and discuss notation and background
results, and then in \S\ref{bilip model intro} we will state the main
technical result of this paper, the Bilipschitz Model Theorem. In
\S\ref{length estimates} we will state the Length Bound Theorem.

\subsection{Surfaces, notation and conventions}
\label{isotopy convention}

We denote by $S_{g,n}$ a compact oriented surface of genus $g$ and $n$
boundary components, and define a complexity $\xi(S_{g,n}) = 3g+n$. 
A subsurface $Y\subset X$ is essential if its boundary components do
not bound disks in $X$ and $Y$ is not homotopic into $\boundary
X$. All subsurfaces which occur in this paper are essential. Note that
$\xi(Y)<\xi(X)$ unless $Y$ is isotopic to $X$. 
(This definition of $\xi$ was used in \cite {masur-minsky:complex2}
and \cite{minsky:ELCI},
but we alert the reader that in some related articles, particularly
\cite{behrstock-kleiner-minsky-mosher:qirigid},
a slightly better convention was adopted of $\xi=3g-3+n$. We retain the
older notation for consistency with \cite{minsky:ELCI}.)

As in \cite{minsky:ELCI}, 
it will be convenient to fix standard representatives of each isotopy
class of subsurfaces in a fixed surface $S$.
Let $\hhat S$ denote a separate copy of 
$int(S)$ with a fixed finite area hyperbolic metric $\sigma_0$.
Then if $v$ is a
homotopy class of simple, homotopically nontrivial curves, let
$\gamma_v$ denote the $\sigma_0$-geodesic representative of $v$,
provided $v$ does not represent a loop around a cusp. 
In \cite[Lemma 3.3]{minsky:ELCI} we fix a version of the standard collar
construction to obtain an open annulus $\collar(v)$ (or $\collar(\gamma_v)$)
which is tubular neighborhood of $\gamma_v$ or a horospherical
neighborhood in the cusp case. This collar has the additional property
that the closures of two such collars are disjoint whenever the core
curves have disjoint representatives. If $\Gamma$ is a collection of simple
homotopically distinct and nontrivial disjoint curves we let
$\collar(\Gamma)$ be the union of collars of components. 

Embed $S$ in $\hhat S$ as the complement of
$\collar(\boundary S)$. Similarly for any essential subsurface
$X\subset S$, our 
standard representative will be the component of $\hhat S \setminus
\collar(\boundary X)$ isotopic to $X$
if $X$ is not an annulus, and $\bcollar(\gamma)$ if
$X$ is an annulus with core curve $\gamma$. 
We will from now on assume that any subsurface of $S$ is of this
form. Note that two such subsurfaces intersect if and only if their
intersection is homotopically essential. We will use the term
``overlap'' to indicate homotopically essential intersection (see also
\S\ref{knotting} for the use of this term in three dimensions).

We will denote by $\DD(S)$ the set of discrete, faithful
representations $\rho:\pi_1(S)\to\PSL 2(\C)$ such that any loop
representing a boundary component is taken to a parabolic element --
that is, the set of Kleinian surface groups for the surface $S$.
If $\rho\in\DD(S)$ we denote by $N_\rho$ its quotient manifold
$\Hyp^3/\rho(\pi_1(S))$.

\subsection{Hierarchies and partial orders}
\label{hierarchy background}

We refer to Minsky \cite{minsky:ELCI} for the basic definitions of 
hierarchies of geodesics in the complex of curves of a
surface. These notions were first developed in 
Masur-Minsky \cite{masur-minsky:complex2}. We will recall the
needed terminology and results here.

\subsubsection*{Complexes, subsurfaces and projections}

We denote by $\CC(X)$ the complex of curves of a surface $X$
(originally due to Harvey \cite{harvey:boundary,harvey:modular})
whose $k$-simplices are
$(k+1)$-tuples of nontrivial nonperipheral homotopy classes of
simple closed curves  with disjoint representatives. For $\xi(X) =
4$, we alter the definition slightly, so that $[vw]$ is an edge
whenever $v$ and $w$ have representatives that intersect once (if
$X=S_{1,1}$)  or  twice (if $X=S_{0,4}$).

When $X$ has boundary we define the
``curve and arc complex'' $\AAA(X)$ similarly, where vertices are proper
nontrivial homotopy 
classes of properly embedded simple arcs {\em or} closed curves.
When $X$ is an annulus the homotopies are assumed to fix the endpoints.

If $X\subset S$ we have a natural map $\pi_X:\AAA(S)\to \AAA(X)$
defined using the essential intersections with $X$  of curves in $S$.
When $X$ is an annulus $\pi_X$ is defined using the lift to the
annular cover of $\hhat S$ associated to $\pi_1(X)\subset\pi_1(\hhat S)$.
If $v$ is a vertex of $\CC(S)$, we let $\AAA(v)$ denote the complex
$\AAA(\collar(v))$. (See Section 4 of \cite{minsky:ELCI} for a more careful discussion
of subsurface projection maps.)

We recall from Masur-Minsky \cite{masur-minsky:complex1} that $\CC(X)$
is $\delta$-hyperbolic (see also Bowditch \cite{bowditch:complex} for
a new proof), and from Klarreich \cite{klarreich:boundary} 
(see also Hamenstadt \cite{hamenstadt:boundary} for an alternate proof)
that its Gromov boundary $\boundary\CC(X)$ can be identified with the
set $\EL(X)$ of minimal filling geodesic
laminations on $X$, with the topology inherited from Thurston's space
of measured laminations under the measure-forgetting map.

\subsubsection*{Markings}
A (generalized) marking $\mu$ in $S$ is a geodesic lamination $\base(\mu)$ in
$\GL(S)$, together
with a (possibly) empty list of ``transversals''. A transversal
is a vertex of $\AAA(v)$ where $v$ is a vertex of $\base(\mu)$ (i.e. a
simple closed curve component of the lamination).
A marking is called {\em maximal} if its base is maximal as a
lamination, and if every closed curve component of the base has a
nonempty transversal. 

Given $\alpha\in\CC_0(S)$, a {\em clean transverse curve} for $\alpha$ is
a curve $\beta\in\CC_0(S)$ such that a regular neighborhood of $\alpha\cup\beta$
is either a 1-holed torus or a 4-holed sphere.
A complete {\em clean marking} $\mu$ is a maximal simplex $\base(\mu)$ in $\CC(S)$
together with a clean transverse curve  for any curve $v$  in $\base(\mu)$
which is disjoint from every other curve in $\base(\mu)$.
If $v\in\base(\mu)$ and $\beta$ is a clean transverse curve for $v$, then
we obtain a transversal to $v$ by projecting $\beta$ to $\AAA(v)$.
Therefore, a complete clean marking gives a well-defined maximal marking.
Moreover, a maximal marking (whose base lamination is a pants decomposition)
gives rise to a complete clean marking, which is well-defined up to
bounded ambiguity (see \cite[Lemma 2.4]{masur-minsky:complex2}).

\subsubsection*{Tight geodesics and subordinacy}
A {\em tight sequence} in (a non-annular) surface $X$ is 
a (finite or infinite) sequence $(w_i)$ of simplices in
the complex of curves $\CC(X)$, with the property
that for any vertex $v_i\in w_i$ and $v_j\in w_j$ with $i\ne j$ we
have $d_{\CC(X)}(v_i,v_j) = |i-j|$, and
the additional property that $w_i$ is the boundary of the subsurface
filled by $w_{i-1}\union w_{i+1}$  if $\xi(X)>4$ (this is 
``tightness,'' see Definition 5.1 of \cite{minsky:ELCI}).

If $X\subset S$ is a nonannular subsurface then a {\em tight geodesic}
$g$ in $X$ is a tight sequence $\{v_i\}$ in $\CC(X)$, together with
two generalized markings $\I(g)$ and $\T(g)$ of $X$ such that the following
holds:
If the sequence $\{v_i\}$ has a first element, $v_0$, we require that
$v_0$
is a vertex of $\base(\I(g))$; otherwise by Klarreich's theorem
\cite{klarreich:boundary} $v_i$ converge as $i\to-\infty$ to a unique
lamination in $\EL(S)$. We choose $\base(\I(g))$ to be this lamination
(and $\I(g)$ has no transversals).
A similar condition holds for $\T(g)$ and the forward
direction of $\{v_i\}$. 
We call $X$ the {\em domain of $g$} and write $X=D(g)$. 

When $X$ is an annulus a tight sequence is any finite geodesic sequence such
that the endpoints on $\boundary X$ of all the vertices are contained
in  the set of endpoints of the first and last vertex.
For a tight geodesic we define $\I$ and $\T$ to be simply the first
and last vertices.  We define the {\em successor} ${\rm succ}(v_i)$ of a simplex $v_i$
of $g$ to be $v_{i+1}$ if $v_i$ is not the last simplex, and $\T(g)$ otherwise.
Similarly we define the predecessor ${\rm pred}(v_i)$ to be
$v_{i-1}$ if $v_i$ is not the first simplex, and $\I(g)$ otherwise. 

For convenience we define $\xi(g)$ to be $\xi(D(g))$ for a tight
geodesic $g$. 

Given a tight geodesic $g$ whose domain is not an annulus and a subsurface
$Y\subset D(g)$, we say that $Y$ is a {\em component domain} of $(D(g),v)$
if  $Y$ is a component of $D(g)\setminus
\collar(v)$ or of $\collar(v)$. We say that
a non-annular component domain of $(D(g),v)$ is {\em directly forward subordinate}
to $g$ at $v$, which we write $Y\fsubd g$, if the successor ${\rm succ}(v)$ of $v$ 
intersects $Y$ nontrivially. 
Notice that the simplex $v$ is uniquely determined by $Y$. 
We similarly say that $Y$ is {\em directly backward subordinate} to $g$ at $v$,
which we denote $g\bsubd Y$, if the predecessor ${\rm pred}(v)$ intersects
$Y$ non-trivially. 

We note some special cases: if $Y$ is an annulus and $Y\fsubd g$ at
$v$ (so $v=[\boundary Y]$) then either $\xi(D(g)) = 4$ and $v$ is not the last vertex, or $v$
is the last vertex andf $\T(g)$ has a transversal associated with
$v$.  If $Y$ is a three-holed sphere and $Y\fsubd g$ at $v$ then
$\xi(D(g)) = 4$ and $v$ cannot be the last vertex. Similar statements
hold for $g\bsubd Y$.
(The subordinacy relation as defined in \cite{minsky:ELCI},
and as used here, is slightly more general than the one defined in \cite{masur-minsky:complex2}
in that it allows 3-holed spheres to be directly subordinate to 4-geodesics.)

This relation yields a subordinacy relation among tight geodesics,
namely that $g\fsubd h$ when 
\begin{itemize}
\item  $D(g)\fsubd h$ at $v$, and
\item $\T(g)$ is the restriction to $D(g)$ of ${\rm succ}(v)$.
\end{itemize}
We define $h\bsubd g$ similarly, replacing $\T$ by $\I$ and  ${\rm
  succ}(v)$ by  ${\rm pred}(v)$.
We let $\fsub $ and $\bsub$ denote the
transitive closures of $\fsubd$ and $\bsubd$. 
Note that $Y\fsub g$ makes sense for a domain $Y$ and geodesic $g$. 
We further say that $Y\fsubeq g$ if either $Y\fsub g$ or $Y=D(g)$, and
similarly $g \bsubeq Y$.

\subsubsection*{Hierarchies}

A {\em hierarchy of tight geodesics} (henceforth just ``hierarchy'')
is a collection of tight geodesics in subsurfaces of $S$ meant to
``connect'' two markings. There is a {\em main geodesic $g_H$} 
whose domain is $D(g_H) = S$, and all other geodesics are obtained
by the rule that, if $Y$ is a subsurface such that $b\bsubd Y \fsubd f$
for some $b,f\in H$, then there should be a (unique) geodesic $h\in H$
such that $D(h)=Y$ and $b\bsubd h \fsubd f$ (this determines $\I(h)$
and $\T(h)$ uniquely). The initial and terminal markings $\I(g_H)$ and
$\T(g_H)$ are denoted $\I(H)$ and $\T(H)$ respectively, and we show 
in \cite{masur-minsky:complex2} that these two markings, when they are
finite, determine $H$ up to finitely many choices.
In \cite{minsky:ELCI} (Lemma 5.13) we extend the construction to the case of
generalized markings, and show that a hierarchy exists for any pair
$\I,\T$ of generalized markings such that no two infinite-leaved components of
$\base(\I)$ and $\base(\T)$ are the same.

\subsubsection*{Hierarchy Structure Theorem}
Theorem 4.7 of \cite{masur-minsky:complex2} (and its slight extension
Theorem 5.6 of \cite{minsky:ELCI}) gives the basic structural
properties of the subordinacy relations, and how they organize the
hierarchy. 
In particular, it states that
for $g,h\in H$ we have $g\fsub h$ if and only if $D(g) \subset D(h)$,
and $\T(h)$ intersects $D(g)$ nontrivially (and similarly replacing
$\fsub$ with $\bsub$ and $\T$ with $\I$). We quote the theorem here
(as it appears in \cite{minsky:ELCI})
for the reader's convenience, as we will use it often.
For a subsurface $Y\subseteq S$ and hierarchy $H$ let $\Sigma^+_H(Y)$
denote the set of geodesics $f\in H$ such that $Y\subseteq D(f)$ and
$\T(f)$ intersects $Y$ essentially
(when $Y$ is an annulus this means either that $\base(\T(f))$
has homotopically nontrival intersection with $Y$, or $\base(\T(f))$
has a component equal to the core of $Y$, with a nonempty
transversal).
Similarly define $\Sigma^-_H(Y)$
with $\I(f)$ replacing $\T(f)$. 

\begin{theorem+}{Descent Sequences}
Let $H$ be a hierarchy in $S$, and $Y$ any essential subsurface of $S$.

\begin{enumerate}
\sitem If $\Sigma^+_H(Y)$ is nonempty then it has the form 
$\{f_0,\ldots,f_n\}$ where $n\ge 0$ and 
$$f_0\fsubd\cdots\fsubd f_n=g_H.$$ 
Similarly, 
if $\Sigma^-_H(Y)$ is nonempty then it has the form 
$\{b_0,\ldots,b_m\}$ with $m\ge 0$, where
$$g_H=b_m\bsubd\cdots\bsubd b_0.$$

\sitem If  
$\Sigma^\pm_H(Y)$ are both nonempty and $\xi(Y)\ne 3$, then $b_0 = f_0$, and
$Y$ intersects every simplex of $f_0$ nontrivially.
\sitem If $Y$ is a component domain in any geodesic $k\in H$
$$f\in \Sigma^+_H(Y) \ \ \iff \ \ Y\fsub f,$$
and similarly, 
$$b\in \Sigma^-_H(Y) \ \ \iff \ \ b\bsub Y.$$

If, furthermore, $\Sigma^\pm_H(Y)$ are both nonempty and $\xi(Y) \ne 3$, 
then in fact $Y$ is the support of $b_0=f_0$. 

\sitem Geodesics in $H$ are determined by their supports. That is, if
  $D(h)=D(h')$ for $h,h'\in H$ then $h=h'$.

\end{enumerate}
\end{theorem+}
When there is no chance of confusion, 
we will often denote $\Sigma^\pm_H$ as $\Sigma^\pm$.

\subsubsection*{Slices and resolutions}
A slice of a hierarchy is a combinatorial analogue of a
cross-sectional surface in $S\times\R$. 
Formally, a slice is a collection $\tau$ of pairs $(h,w)$ where $h$ is a
geodesic in $H$ and $w$ is a simplex in $h$, with the following properties.
$\tau$ contains a distinguished
``bottom pair''  $p_\tau = (g_\tau,v_\tau)$.
For each $(k,u)$ in $\tau$ other than the bottom pair,
there is an $(h,v)\in\tau$ such that $D(k)$ is a component domain of
$(D(h),v)$; moreover any geodesic appears at most once among the pairs
of $\tau$.  (See \cite[\S 5.2]{minsky:ELCI} for more details).

We say a slice $\tau $ is {\em saturated} 
if, for every pair $(h,v)\in\tau$ and every 
geodesic $k\in H$ with $D(k)$ 
a component domain of $(D(h),v)$, there is some (hence exactly one) pair
$(k,u)\in\tau$. 
It is easy to see by induction
that for any pair $(h,u)$ there is a saturated slice with bottom pair $(h,u)$.
A slightly stronger condition is that, for every $(h,u)\in\tau$ and
every component domain $Y$ of $(D(h),u)$ with $\xi(Y)\ne 3$ there is,
in fact, a pair 
$(k,u)\in \tau$ with $D(k)=Y$; we then say that $\tau$ is {\em full}. 
It is a
consequence of Theorem \ref{Descent Sequences} that, if $\I(H)$ and
$\T(H)$ are maximal markings, then every saturated slice is full. 

A slice is {\em maximal} if it is full and its bottom geodesic is the
main geodesic $g_H$. 

A {\em non-annular slice} is a slice in which none of the pairs
$(k,u)\in \tau$  
have annulus domains. A non-annular slice is saturated, full, or maximal
if the conditions above hold with the exception of annulus domains. 
In particular, Theorem \ref{Descent Sequences} implies that
if $\base(\I(H))$ and $\base(\T(H))$ are maximal laminations, then
every saturated non-annular slice is a full non-annular slice. 

To a slice $\tau$ is associated a marking $\mu_\tau$, whose base is
simply the union of simplices $w$ over all pairs $(h,w)\in\tau$ with
nonannular domains. 
The transversals of the marking are determined by the pairs in $\tau$
with annular domains, if any (see  \cite[\S 5]{masur-minsky:complex2}).
We also denote $\base(\mu_\tau)$ by $\base(\tau)$.
We note that if $\tau$ is maximal then $\mu_\tau$ is a maximal
marking, and if $\tau $ is maximal non-annular then $\mu_\tau$ is a
pants decomposition. 

We also refer to \cite[\S 5]{masur-minsky:complex2} for the
notion of ``(forward) elementary move'' on a slice, denoted
$\tau \to \tau'$. The main effect of this move is to replace one pair
$(h,v)$ in $\tau$ with a pair $(h,v')$ in $\tau'$ where $v'$ is the
successor of $v$ in $h$. In addition certain pairs in $\tau$ whose
domains lie in $D(h)$ are replaced with other pairs in $\tau'$. 
The underlying curve system 
$\base(\tau)$ stays the same except in the case that $\xi(h)=4$.
When $\xi(h)=4$, $v$ and $v'$ intersect in a minimal way, and all
other curves of $\base(\tau)$ and $\base(\tau')$ agree; if $\tau$ is
maximal this amounts to a standard elementary move on pants
decompositions. 

A {\em resolution} of a hierarchy $H$ is a sequence of elementary
moves $(\cdots \to \tau_n \to \tau_{n+1} \to \cdots)$ (possibly
infinite or biinfinite), where each $\tau_n$ is a saturated slice with
bottom geodesic $g_H$, with the additional property that every 
pair $(h,u)$ (with $h\in H$ and $u$ a simplex of $h$) appears in some
$\tau_n$.  Lemmas 5.7 and 5.8 of \cite{minsky:ELCI} guarantee that
every hierarchy has a resolution with this property.
A resolution is closely related to a ``sweep'' of
$S\times\R$ by cross-sectional surfaces, and this will be exploited
more fully in \S\ref{po}. 

It is actually useful not to involve the annulus geodesics in a
resolution. Thus given a hierarchy $H$ one can delete all annulus
geodesics to obtain a {\em hierarchy without annuli} $H'$ (see
\cite[\S 8]{masur-minsky:complex2}) and a resolution of $H'$ will be
called a non-annular resolution.

\subsubsection*{Partial orders} 
In \cite{masur-minsky:complex2} we introduce several partial orders on
the objects of a hierarchy $H$. In this section we extend the notion
of ``time order'' $\tprec$ on geodesics to a time order on component
domains, and we recall the properties of the partial order $\pprec$ on
pairs.

First, for a subsurface $Y$ of $D(g)$, 
define the  {\em footprint} $\phi_g(Y)$ to be the set of simplices of $g$ which
represent curves  disjoint from $Y$. Tightness
implies that $\phi_g(Y)$ is always an interval of $g$, and the 
triangle inequality in $\CC(D(g))$ implies it has diameter at most 2. 

If $X$ and $Y$ are component domains arising in $H$, 
we say that  $X\tprec Y$ 
whenever there is a ``comparison geodesic'' $m\in H$ such that
$D(m)$ contains $X$ and $Y$ with nonempty footprints, and
\begin{equation}\label{footprints ordered}
\max \phi_m(X) < \min \phi_m(Y).
\end{equation}
(Max and min are with respect to the natural order $v_i < v_{i+1}$ of
the simplices of $m$).
Note that (\ref{footprints ordered}) also implies that
$Y\fsub m$ and $m\bsub X$ in this case, by Theorem \ref{Descent Sequences}.

For geodesics $g$ and $h $ in $H$ we can define $g\tprec h$ if $D(g) \tprec
D(h)$, and this is equivalent to definition 4.16 in
\cite{masur-minsky:complex2}.
(We can similarly define $g\tprec Y$ and $Y\tprec h$.)
In Lemma 4.18 of
\cite{masur-minsky:complex2} it is shown, among other things, that
$\tprec$ is a strict 
partial order on the geodesics in $H$, and it  follows immediately
that it is a strict partial order, with our definition, on all domains
of geodesics in $H$. It is not hard to generalize this and the rest of
that lemma to 
the set of {\em all} component domains in $H$, which for completeness
we do here.
(The main point of this generalization is to deal appropriately with
3-holed spheres, which can be component domains but never support geodesics.)

\begin{lemma}{time order for domains}
Suppose that $H$ is a hierarchy with $\base(\I(H))$ or $\base(\T(H))$ maximal. 
The relation $\tprec$ is a strict partial order on the set of
component domains occuring in $H$. Moreover, if $Y$ and
$Z$ are component domains, then 
\begin{enumerate}
\item \label{nested not ordered}
If $Y\subseteq Z$ then $Y$ and $Z$ are not $\tprec$-ordered.
\item \label{overlap ordered}
Suppose that $Y\intersect Z \ne \emptyset$ and neither domain is
contained in the other. Then $Y$ and $Z$ are $\tprec$-ordered. 
\item \label{lattice diagram}
If $b \bsub Y \fsub f$ then either $b=f$, $b\fsub f$, $b\bsub f$, or
$b\tprec f$. 
\item \label{composition relation}
If $Y\fsub m$ and $m\tprec Z$ then $Y\tprec Z$. Similarly if
$Y \tprec m$ and $m\bsub Z$ then $Y\tprec Z$.
\end{enumerate}
\end{lemma}

\begin{proof}
We follow the proof of Lemma 4.18 in \cite{masur-minsky:complex2},
making adjustments for the fact that the domains may not support
geodesics. Let us first prove the following
slight generalization of Corollary 4.14 of
\cite{masur-minsky:complex2}, in which $Y$ was assumed to be the
domain of a geodesic in $H$.

\begin{lemma}{nonempty footprint}
If $h$ is a geodesic in a hierarchy $H$, $Y$ is a component domain in
$H$,  and $Y\subsetneq D(h)$,
then $\phi_h(Y)$ is nonempty.
\end{lemma}

\begin{proof}
If $Y$ fails to intersect $\T(h)$ then in particular it is disjoint
from the last simplex of $h$, and hence $\phi_h(Y)$ is nonempty. 
If $Y$ does intersect
$\T(h)$ then, by Theorem \ref{Descent Sequences}, we have
$Y\fsub h$. This means that there is some $m$ with
$Y\fsubeq m \fsubd h$ and $\phi_h(D(m))$ is therefore nonempty, and
is contained in $\phi_h(Y)$. 
\end{proof}

Now let us assume that $\base(\I(H))$ is a maximal lamination;
the proof works similarly if $\T(H)$ is maximal (this assumption is used just
once in the proof of part (\ref{overlap ordered}), and we suspect that
the lemma should be true without it).

To prove part (\ref{nested not ordered}), suppose $Y\subseteq Z$. Then
for any geodesic $m$ 
with $Z\subseteq D(m)$, we have $\phi_m(Z) \subseteq \phi_m(Y)$. In
particular the footprints of $Y$ and $Z$ can never be disjoint, and
hence they are not $\tprec$-ordered.

To prove part (\ref{overlap ordered}), let us first establish the
following statement:  

\begin{itemize}
\item[(*)]
If $m\in H$ is a geodesic such that $Y\union Z \subset D(m) $, where
$Y$ and $Z$ are component domains in $H$ which intersect but neither
is contained in the other, and in addition we have either $m\bsub Y$
or $m\bsub Z$, then $Y$ and $Z$ are $\tprec$-ordered.
\end{itemize}

The proof will be by induction on $\xi(m)-\max\{\xi(X),\xi(Y)\}$.
If $\xi(m) $ equals $\xi(Y)$ or $\xi(Z)$ then $D(m)$ equals one of $Y$
or $Z$ and then the other is contained in it, but we have assumed this
is not the case, so the statement holds vacuously.

Now assume that $\xi(m) > \max(\xi(Y), \xi(Z))$. Consider the
footprints $\phi_m(Y)$ and $\phi_m(Z)$ (both nonempty by Lemma
\ref{nonempty footprint}). If the footprints are disjoint then $Y$ and
$Z$ are $\tprec$-ordered with $m$ the comparison geodesic, and we are
done.  If the footprints intersect then since they are intervals the
minimum of one must be contained in the other. Let $v=\min\phi_m(Y)$
and $w=\min\phi_m(Z)$. 

If $v<w$ then in particular $\phi_m(Z)$ does not
include the first simplex of $m$, and so by Theorem \ref{Descent Sequences}
we have $m \bsub Z$. This means that there is some $m'$ with
$m\bsubd m' \bsubeq Z$.
$D(m')$ is a component domain of $(D(m),w)$, so
since $w\in \phi_m(Y)$ and $Y\intersect Z\neq \emptyset$ we find that
$Y\subset D(m')$. Now by induction we may conclude that $Y$ and $Z$
are $\tprec$-ordered. 

If $w<v$ we of course apply the same argument with the roles reversed.
If $w=v$ then we use the hypothesis that $m\bsub Z$ or $m\bsub Y$, and
again repeat the previous argument. This concludes the proof of 
assertion (*).

To show that (\ref{overlap ordered}) follows from (*), it suffices to show the
hypothesis holds for $m=g_H$. 
Suppose that $g_H\bsub Y$ fails to
hold.  This means, by Theorem \ref{Descent Sequences},
that $Y$ does not intersect $\base(\I(H))$, and
since $\base(\I(H))$ is maximal $Y$ must be either a 3-holed sphere with
boundary in $\base(\I(H))$, or an annulus with core in
$\base(\I(H))$. Now since $Y$ and $Z$ intersect, it follows that $Z$
does intersect $\base(\I(H))$ nontrivially. Hence again by 
Theorem \ref{Descent Sequences} we have $g_H \bsub Z$.
Thus we may apply (*) to obtain (\ref{overlap ordered}).

To prove (\ref{lattice diagram}), suppose $b\ne f$. Suppose first that
$D(b)\subset D(f)$. Since $Y\fsub f$, 
Theorem \ref{Descent Sequences} implies that $f\in \Sigma^+(Y)$, 
and since  $Y\subset D(b)$, we must have $f\in \Sigma^+(D(b))$ as
well.  Hence $b\fsub f$ by Theorem \ref{Descent Sequences}.
Similarly if $D(f)\subset D(b)$ we have $b\bsub f$. If neither domain
is contained in the other, since they both contain $Y$ we may apply 
(\ref{overlap ordered}) to conclude that they are
$\tprec$-ordered. Suppose by contradiction $f\tprec b$, and 
let $m$ be the comparison geodesic. Thus $f\fsub m \bsub b$, and 
$\max \phi_m(f)  < \min \phi_m(b)$. 
Since $Y\fsub f\fsub m$ we have by 
Lemma 5.5 of \cite{minsky:ELCI} 
that
$\max \phi_m(Y) = \max\phi_m(f)$. Similarly
$\min \phi_m(Y) = \min\phi_m(b)$. This contradicts
$\max\phi_m(f)  < \min\phi_m(b)$, so we conclude
$b\tprec f$.

To prove (\ref{composition relation}), consider the case $Y\fsub m \tprec Z$ (the other case is
similar). Let $l$ be the comparison geodesic for $m$ and $Z$. Then
$\max\phi_l(D(m)) < \min \phi_l(Z)$. By 
Lemma 5.5 of  \cite{minsky:ELCI}, $Y\fsub m$ implies
$\max\phi_l(Y) = \max\phi_l(D(m))$, and hence $Y\tprec Z$.

Now we may 
finish the proof that $\tprec$ is a strict partial order. Suppose that
$X\tprec Y$ and $Y\tprec Z$. Thus we have geodesics $l$ and $m$ such
that  $X\fsub l \bsub Y \fsub m \bsub Z$, and furthermore
$\max \phi_l(X) < \min \phi_l(Y)$, and 
$\max \phi_m(Y) < \min \phi_m(Z)$. 

Applying (\ref{lattice diagram}) to 
$l \bsub Y \fsub m$ we find that $l=m$, $l\bsub m$, $l\fsub m$, or
$l\tprec m$. 
Suppose first that $l=m$. Then $X\fsub m \bsub Z$, and
$\max\phi_m(X) < \min\phi_m(Y) \le \max\phi_m(Y) < \min
\phi_m(Z)$. Thus $X\tprec Z$.

If $l\fsub m$ then $X\fsub m$ and by  
Lemma 5.5 of \cite{minsky:ELCI} 
we have $\max\phi_m(X) = \max\phi_m(D(l))$. Now since $D(l)$ contains
$Y$ we have $\phi_m(D(l)) \subset \phi_m(Y)$ and it follows that
$\max\phi_m(D(l)) \le \max\phi_m(Y) < 
\min \phi_m(Z)$. Thus again we have $X\tprec
Z$. The case $l\bsub m$ is similar. 

If $l\tprec m$ then $X\fsub l \tprec m \bsub Z$ and applying (\ref{composition relation}) twice
we conclude $X\tprec Z$.

Hence $\tprec$ is transitive, and since by definition $X\tprec X$ can
never hold, it is a strict partial order.
\end{proof}

There is a similar order on {\em pairs} $(h,u)$ where $h$ is a
geodesic in $H$, and $u$ is either a simplex of $h$, or
$u\in\{\I(h),\T(h)\}$. We define a {\em generalized footprint}
$\hhat \phi_m(g,u)$ to be $\phi_m(D(g))$ if $D(g)\subset D(m)$, and
simply $\{u\}$ if $g=m$. We then say that
$$
(g,u) \pprec (h,v)
$$
if there is an $m\in H$ with $D(h)\subseteq D(m)$, $D(g)\subseteq
D(m)$,  and
$$\max\hhat\phi_m(g,u) < \min\hhat\phi_m(h,v)$$
In particular if $g=h$ then $\pprec$ reduces to the natural order on
$\{\I(g),u_0,\ldots,u_N,\T(g)\}$ where $\{u_i\}$ are the simplices of
$g$. This relation is also shown to be a partial order in Lemma 4.18 of
\cite{masur-minsky:complex2}. 

In the proof of Lemma 5.1 in \cite{masur-minsky:complex2}, the
following fact is established which will be used in \S\ref{po} and
\S\ref{regions}. It is somewhat analogous to part (1) of 
Lemma \ref{time order for domains}.

\begin{lemma}{slice unordered}
Any two elements $(h,u)$ and $(k,v)$ of a slice $\tau$ of $H$ are
not $\pprec$-ordered.
\end{lemma}

\subsection{Margulis tubes}
\label{tubes}

\subsubsection*{Tube constants}
Let $N_J$ for $J\subset\R$ denote the region 
$\{x\in N| 2\inj_N(x)\in J\}$.
Thus
$N_{(0,\ep]}$ denotes the $\ep$-thin part of a hyperbolic manifold
  $N$, and  $N_{[\ep,\infty)}$ denotes the $\ep$-thick part. 
Let $\epzero$ be a Margulis constant for $\Hyp^3$, so that for $\ep\le
\epzero$, 
$N_{(0,\ep]}$ for a hyperbolic 3-manifold $N$
is a disjoint union of standard closed
tubular neighborhoods of closed geodesics, or horocyclic cusp
neighborhoods. (See e.g. Benedetti-Petronio \cite{benedetti-petronio}
or Thurston \cite{wpt:textbook}.)

We will call a component of $N_{(0,\ep)}$
an (open)  $\ep$-Margulis tube, and denote it by
$\MT_\ep(\gamma)$, where $\gamma$ is the homotopy class of the core
(or in the rank-2 cusp case, any nontrivial homotopy class in the
tube). If $\Gamma$ is a collection of simple closed curves or homotopy
classes we will denote $\MT_\ep(\Gamma)$ the union of the corresponding
Margulis tubes.

Let $\epone<\epzero$ be chosen as in Minsky \cite{minsky:ELCI}, so
 that the $\epzero$-thick part of an essential pleated surface maps
 into the $\epone$-thick part of the target 3-manifold.
This is
the constant used in the Lipschitz Model Theorem (\S\ref{bilip model
  intro}). It will be our ``default'' Margulis constant and we will
usually denote $\MT_{\epone}$ as just $\MT$. (The only place we use
$\epzero$ will be in the definition of the augmented convex core).

Let $\rho\in\DD(S)$ be a Kleinian surface group, and $N=N_\rho$. 
Then $N$ is homeomorphic to $int(S)\times\R$
(Bonahon \cite{bonahon}). More precisely, Bonahon showed that
$N^1\homeo S\times\R$, where $N^1$ is
$N_\rho \setminus \MT(\boundary S)$ -- the complement of the open cusp
neighborhoods associated to $\boundary S$. 

Thurston
showed that a sufficiently short primitive geodesic in $N$ is
homotopic to a simple loop in $S$. Otal proved the following stronger
theorem in \cite{otal:knotting,otal:knotting2}:

\begin{theorem}{otal}
There exists $\epotal> 0$ depending only on the compact surface $S$
such that, if $\rho\in\DD(S)$ and $\Gamma$ is the set of primitive closed
geodesics in $N=N_\rho$ of length at most $\epotal$, then $\Gamma$ is
{\em unknottted and 
  unlinked}. That is, $N^1$ can be identified with $S\times\R$ in such
a way that $\Gamma$ is identified with a collection of disjoint simple closed
curves of the form $c\times\{t\}$.
\end{theorem}

We remark that 
Otal's proof only explicitly treats the case that $S$ is a closed
surface, but the case with boundary is quite similar. One can also
obtain this result, for any finite subcollection of $\Gamma$, 
by applying a special case of Souto \cite{souto:knotting}.

\subsubsection*{Bounded homotopies into tubes}
The next lemma shows that a bounded-length curve homotopic into a
Margulis tube admits a controlled homotopy into the tube. It will be
used at the end of Section \ref{uniform embeddings}.

\begin{lemma}{bounded annulus}
Let $\UU$ be a collection of $\epone$-Margulis tubes in a hyperbolic 3-manifold
$N$
and let $\gamma$ be an essential curve 
which is homotopic 
within $N\setminus \UU$ to $\boundary\UU$.
 
Then such a homotopy can be found whose diameter is bounded
by a constant $r$ 
depending only on $\epone$ and the length $l_N(\gamma)$.
\end{lemma}

\begin{proof}
The choice of $\epone$ strictly less than the Margulis constant for
$\Hyp^3$ implies that $\UU$ has an embedded collar neighborhood of
definite radius, so  possibly enlarging $\UU$ we may assume that 
the radius of each component is at least $R>0$ (depending on $\epone$).
Let $U$ be a component of $\UU$ with core geodesic $c$. 
Agol has shown in \cite{agol:volume} (generalizing a construction of
Kerckhoff) that
there exists a metric $g$ on $U\setminus c$ such that
\begin{enumerate}
\item $g$ agrees with the metric of $N$
 on a neighborhood of $\boundary U$.
\item $g$ has all sectional curvatures between
  $[-\kappa_1,-\kappa_0]$, where $\kappa_1>\kappa_0>0$ 
  depend on $R$. 
\item On some neighborhood of $c$, $g$ is isometric to a
  rank-2 parabolic cusp.
\end{enumerate}

\medskip

Let $\hhat N$ be the complete, negatively-curved manifold obtained by
deleting the cores of 
components of $\UU$ and replacing the original metric by
the metric $g$ in each one. 
The homotopy from $\gamma$ to $\boundary\UU$ can be deformed to a {\em
  ruled annulus} $A:\gamma\times[0,1]\to \hhat N$, i.e. a map such
that $A(\cdot,0) = id$, $A(\cdot,1)$ has image in $\UU$, and
$A|_{x\times[0,1]}$ is a geodesic. This is possible simply by
straightening the trajectories of the original homotopy, since $\hhat
N$ is complete and negatively curved. Because a ruled surface has
non-positive extrinsic curvature, the
pullback metric on $\gamma\times[0,1]$ must have curvatures bounded
above by $-\kappa_0$.
Furthermore, by pushing $A(\gamma\times\{1\})$ sufficiently far into
the cusps of $\hhat N$, we can ensure that the total boundary length
of the annulus is at most $l(\gamma) + 1$.

The area of the annulus (in the pullback metric) is bounded by $Cl(\boundary A)\le
C(l(\gamma)+1)$, where $C$ depends on the curvature bounds. 
Let $\ep = \epone/2$. 

Let $A'$ be the component of $A^{-1}(\hhat N\setminus\UU)$ containing
the outer boundary $\gamma\times\{0\}$. This is a punctured annulus, 
and the punctures can be filled in by disks in $\gamma\times[0,1]$. Let 
$A''$ denote the union of $A'$ with these disks.
The injectivity radius of $A$ at a point in $A'$ is at least $\ep$, and it follows
that the same holds for $A''$, since any essential loop passing
through one of the added disks must also pass through $A'$. 
Let $A''_\ep$ be the complement in $A''$ of an $\ep$-neighborhood of
the outer boundary (in the induced metric). Any point in $A''_\ep$ is
the center of an embedded disk of area at least $\pi\ep^2$, so the
area bound on $A$ implies that  any component of  $A''_\ep$ has diameter
at most $C'(l(\gamma)+1)$. 
This gives a bound $\diam A'' \le C'(l(\gamma)+1) + l(\gamma)+\ep$.

This bounds how far each of the disks of $A'\setminus A''$ reaches
into the tubes $\UU$, and hence bounds the distortion 
caused by pushing these disks back to $\boundary\UU$. Applying this
deformation to $A''$ yields a homotopy of $\gamma$ into $\boundary
\UU$ with bounded diameter, as desired.

\end{proof}

\subsection{Geometric Limits}
\label{geometric limits}
Let us recall the notion of geometric convergence for {\em hyperbolic manifolds with
baseframes} $(N,\hat x)$, where $N$ denotes a hyperbolic manifold and
$\hat x$ is an orthonormal baseframe for $T_x(N)$  at some point $x\in N$.
We say that $\{(N_i,\hat x_i)\}$ {\em converges geometrically}  to $(N,\hat x)$  if for all $i$ there
exists a diffeomorphic embedding $f_i:X_i\to N_i$ with $x\in
X_i\subset N$ and $df_i(\hat x) = \hat x_i$, 
and for any  $\ep>0$ and $R>0$,
there exists $I$ such that for $i\ge I$ we have
$B_R(x)\subset X_i$ and $f_i$ is $\ep$-close, in $C^2$, to a local isometry.
We call the $f_i$  {\em comparison maps}. 

Equivalently one can state the definition with the comparison maps
going the other direction, $g_i:B_R(x_i) \to N$. It will be convenient
for us to use both definitions.

Notice that if we choose a fixed baseframe $\hat x_0$ for $\Hyp^3$, then a hyperbolic
3-manifold with baseframe $(N,\hat x)$ gives rise to a unique Kleinian group $\Gamma$
such that $N=\Hyp^3/\Gamma$ and $\hat x$ is the projection of $\hat x_0$ to $N$.
We say that a sequence $\{\Gamma_i\}$ of Kleinian groups converges
geometrically to $\Gamma$ if $\{( \Hyp^3/\Gamma_i,\hat x_i)\}$ converges
geometrically to $(\Hyp^3/\Gamma,\hat x)$ (where each baseframe is a projection
of $\hat x_0$). This is equivalent to convergence of $\{\Gamma_i\}$ to $\Gamma$ in
the sense of closed subsets of $\PSL 2(\C)$ (where we give the set of closed
subsets of $\PSL 2(\C)$ the Chabauty topology).
See \cite{ceg} or \cite{benedetti-petronio} for more details. 

If $G$ is any finitely generated, torsion-free group, 
the set $\DD(G)$ of discrete, faithful representations 
$\rho:G\to \PSL 2(\C)$ is given the natural topology of convergence on each element of
$G$, also called the topology of algebraic convergence.

If $H$ is a non-abelian subgroup of $G$ and  $\{\rho_n:G\to \PSL 2(\C)\}$  is
a sequence of representations such that $\{\rho_n|_H\}$ converges, one may
pass to a subsequence so that $\{\rho_n(G)\}$ converges geometrically.
Some important aspects of the
relationship between the sequence of representations and the geometric
limit are described by the following Proposition, which gives relative versions of
Lemma 3.6 and Proposition 3.8 in \cite{jorgensen-marden:convergence} and Lemma 3.6
in \cite{anderson-canary-mccullough:bumping}.

\begin{proposition}{JMrel}
Let $G$ be a torsion-free group and let $H$ be a non-abelian subgroup
of $G$. Let $\{\rho_i\}$ be a sequence in ${\mathcal D}(G)$ 
such that $\{\rho_i|_H\}$ converges to $\rho\in {\mathcal D}(H)$.
Then
\begin{enumerate}
\item
If $\{g_i\}$ is a sequence in $G$ and $\{\rho_i(g_i)\}$ converges
to the identity, then $g_i$ is equal to the identity for all
large enough $i$.
\item
There exists a subsequence $\{\rho_{i_n}\}$ of $\{\rho_i\}$ such
that $\{\rho_{i_n}(G)\}$ converges geometrically to a torsion-free
Kleinian group $\Gamma$ such that $\rho(H)\subset \Gamma$.
\item
Let $K$ be a finite complex, $N=\Hyp^3/\rho(H)$, $\hat N=\Hyp^3/\Gamma$
and $\pi:N\to\hat N$ be the natural covering map.
If $h:K\to N$ is a continuous map, then, for all sufficiently large $i_n$, $\pi(h(K))$
is in the domain of the comparison map $f_{i_n}$ and
$(f_{i_n}\circ\pi \circ h)_*:\pi_1(K)\to \pi_1(N_{i_n})$ is conjugate to
$\rho_n\circ \rho^{-1}\circ h_*$.
\end{enumerate}
\end{proposition}

\begin{proof}
Let $h_1$ and $h_2$ be two non-commuting elements of $H$.
Since $\{\rho_i|_H\}$ is algebraically convergent,
there exists $A>0$ such that $d(\rho_i(h_k)(0),0)\le A$ for all $i$
and $k=1,2$. The Margulis lemma implies that given $A$ 
there exists $\epsilon>0$ such that if  $\alpha,\beta\in{\rm PSL}_2({\bf C})\smallsetminus\{id\}$
and $d(\alpha(0),0)<\epsilon$ and $d(\beta(0),0)\le A$, then either
$\alpha$ and $\beta$ commute or the group they generate is indiscrete
or has torsion.
Thus, if $\gamma\in \rho_i(G\smallsetminus\{id\})$, then $d(0,\gamma(0))\ge \epsilon$,
since $\rho_i(G)$ is discrete  and torsion-free (for all $i$)
and $\gamma$ cannot commute with both $\rho_i(h_1)$ and $\rho_i(h_2)$.
The first fact follows immediately,
while the second fact follows from Theorem 3.1.4 in
Canary-Epstein-Green \cite{ceg}.

The proof of the third fact is virtually identical to the proof of Lemma 7.2
in \cite{anderson-canary:cores} (see also \cite[Prop. 3.3]{canary-minsky:tamelimits}).
The key point is that the comparison maps lift to maps from $\Hyp^3$ to $\Hyp^3$
which are converging to the identity.
\end{proof}

The next lemma will allow us to assume that the comparison maps respect
the thin parts (see Evans  \cite[Prop. 4.3]{evans:tamenesspersists}
for a similar statement).

\begin{lemma}{comparisons respect tubes}
Let $\{(N_i,\hat x_i)\}$ be a sequence of hyperbolic 3-manifolds with
baseframe converging to $(N,\hat x)$.
Let $T$ be a finite collection of
components of $N_{(0,\epone)}$. 

Then, given any $R>0$, we may choose the comparison maps $f_i$
such that, for all large enough $i$, there exists a finite collection $T_i$ of
components  of $(N_i)_{(0,\ep_1)}$ such that
\begin{enumerate}
\item
$f_i(T\cap B_R(x))\subset T_i$,
\item
$f_i(\boundary T\cap B_R(x))\subset\boundary T_i$, and
\item
$f_i(B_R(x)\setminus T)\subset N_i\setminus T_i$.
\end{enumerate}

The corresponding statement holds with the comparison maps going in
the opposite direction. 
\end{lemma}

\begin{pf}
We may assume that $B_R(x)$ is a smooth submanifold of $N$ (if not we may
work with a larger $R$ and restrict).
Choose $\ep'\in (0,\epone)$ so that each component $T_m$ of $T$ contains
a curve $\beta_m$ of length $\ep'$ which is homotopic to the core curve of $T_m$.
Choose $\delta\in (1,\epone/\ep')$.
For sufficiently large $i$ we may assume that $f_i$ is
$\delta$-Lipschitz and its domain $X_i$
contains $B_{R+\epone}(x)$ and $B_{\delta(R+\epone)}\subset f_i(X_i)$ .
Moreover, we may assume that
if $p\in T_m\cap B_R(x)$, for any $m$, then there is a homotopically non-trivial 
curve $\gamma_p$ through $p$
of length at most $\epone$ which is homotopic  within $X_i$ to (a power of)
$\beta_m$.

If $p\in T_m\cap B_R(x)$, then $f_i(\gamma_p)$ has length at most
$\delta\epone$ and is homotopic to (a power of) $f_i(\gamma_m)$ which has length less
than $\epone$. If $f_i(\gamma_p)$ were homotopically trivial then it would
bound a disk of diameter at most $\delta\epone/2$ which would thus
be contained in $f_i(X_i)$. In this case the disk would pull back under $f_i$ to
a disk bounded by $\gamma_p$, so $f_i(\gamma_p)$ must be homotopically
non-trivial.
Therefore,
$f_i(T_m \intersect B_R(x) )$ is contained in the component $(T_m)_i'$ of
$(N_i)_{(0,(1+\delta)\epone)}$ which contains $f_i(\gamma_m)$.
By a  similar argument, we may assume that 
$f_i(B_R(x)\setminus T_m )$ does not intersect  $(T_m)_i'\cap (N_i)_{(0,\epone/\delta)}.$
Let $(T_m)_i$
be the component of $(N_n)_{(0,\epone)}$ contained in $f_i(\gamma_m)$ and let
$T_i$ be the collection of all the $(T_m)_i$.

If $\delta$ is chosen close enough to 1, then the region
$N_{(\epone/\delta,\epone\delta)}\cap (T_m)_i'$
is a bicollar neighborhood of radius $O(\delta)$ of $\boundary (T_m)_i$
and the image of
$\boundary T_m\intersect B_R(x)$ can be represented in
the product structure of the collar as the graph of a nearly
constant function over $\boundary (T_m)_i$.
We can then use the  collar structure to adjust the map in this
neighborhood so that  it satisfies claims (1)--(3), and is still $C^2$-close
to a local isometry.
\end{pf}

\subsection{Ends and ending laminations}
\label{ends}

We recap here the definitions of ends and ending laminations. See
\cite[\S 2]{minsky:ELCI} for more details.

Fix $\rho\in\DD(S)$. Let $N=N_\rho$ and let $C_N$ denote the convex
core of $N$.  
Let  $Q$ denote the union of (open)
$\epone$-Margulis tube neighborhoods of the cusps of $N$, and let $N^0
= N\setminus Q$. 
Let $Q_1 \subset Q$ be the union of
tubes associated to $\boundary S$ (thus $Q_1 = \MT(\boundary S)$, and
$N\setminus Q_1$ is $N^1$, as defined in \S\ref{tubes}).
Let $\boundary_\infty N$ denote the conformal boundary of $N$ at
infinity, obtained as the quotient of the domain of discontinuity of
$\rho(\pi_1(S))$. 

As a consequence of Bonahon's tameness theorem \cite{bonahon},
we can fix an orientation-preserving identification of  $N$ with
$\hhat S\times\R$ in such a way 
that  $Q_1 = \collar(\boundary S)\times\R$, and furthermore so that
$\KK\equiv S\times[-1,1]$ meets the closure of $Q$ in a union of disjoint
essential annuli 
$P=P_1 \union P_+ \union P_-$, where $P_1 = \boundary S \times[-1,1]$ and
$P_\pm\subset int(S)\times\{\pm 1\}$. The pair $(\KK,P)$ is the {\em
  relative compact core} of $N_0$, and the components of $N_0\setminus
\KK$ are neighborhoods of the ends of $N_0$. 

For each component $R$ of $\boundary\KK \setminus P$, there is an
invariant $\nu_R$ defined as follows: If the end associated to $R$ is
geometrically finite then $\nu_R$ is the point in the Teichm\"uller
space $\Teich(R)$ associated to the component of the conformal
boundary $\boundary_\infty N$ that faces $R$. 
If the end is geometrically infinite then 
(again by Bonahon \cite{bonahon} and by Thurston's original definition
in \cite{wpt:notes}) $\nu_R$ is a geodesic lamination in $\EL(R)$,
which is the unique limit (in the measure-forgetting topology) of
sequences of simple closed curves in $R$ whose 
geodesic representatives exit the end associated to $R$.

The top ending invariant $\nu_+$ then has a lamination part $\nu_+^L$ and a Riemann surface
part $\nu_+^T$: $\nu_+^L$ is the union of the core curves $p_+$ of the
annuli $P_+$ (these are the ``top parabolics'') and the laminations
$\nu_R$ for those components $R$ of $S\times\{+1\}\setminus P_+$ which
correspond to simply degenerate ends. $\nu_+^T$ is the union of
$\nu_R\in\Teich(R)$ for those components $R$ of
$S\times\{+1\}\setminus P_+$ which correspond to geometrically finite
ends. We define the bottom ending invariant $\nu_-$ similarly.
We let $\nu$, or $\nu(\rho)$, denote the pair $(\nu_+,\nu_-)$.

Note in particular the special case that there are no parabolics
except for $P_1$, and both the $+$ and $-$ ends are degenerate. 
In this case $\nu_+$ and $\nu_-$ are both filling laminations in
$\EL(S)$. This is called the {\em doubly degenerate case}, and it is
helpful for most of this paper to focus just on this case. 

\subsection{Definition of the model manifold}
\label{model definitions}
We recall here the definition of the model manifold $\modl$ from
\cite{minsky:ELCI}. 

Given a pair $\nu=(\nu_+,\nu_-)$ of end invariants, we construct
as in \cite[\S7.1]{minsky:ELCI} a pair of markings $\mu_\pm$
which encode the geometric information in $\nu$ up to bilipschitz
equivalence. In particular,
when $\nu_+$ is a filling lamination (the $+$ end is simply
degenerate), $\base(\mu_+) 
= \nu_+$.  When $\nu_+$ is a point in Teichm\"uller space, 
$\mu_+$ is a minimal-length marking in the corresponding metric on $S$.
In general $\base(\mu_\pm)$ is a maximal lamination (maximal
among supports of measured laminations) whose infinite-leaf components
are ending laminations for ends of the manifold $N^0$ obtained from
$N$ by removing all cusps. The closed-leaf components of
$\base(\mu_\pm)$ which are not equipped with transversals are exactly
the (non-peripheral) parabolics of $N$. 

Note also that a component can be common to $\base(\mu_+)$  and
$\base(\mu_-)$ only if it is a closed curve, and has a tranversal on
at least one of the two. This is because a non-peripheral parabolic in
$N$ corresponds to a cusp on either one side of the compact core or the other.

We let $H=H_\nu$ be a hierarchy such that $\I(H)=\mu_-$ and
$\T(H)=\mu_+$.  

The model manifold $\modl$ is identified with a subset of $\hhat M
\equiv \hhat S\times\R$, and is partitioned into pieces called {\em  blocks} and 
{\em tubes}. It is also endowed with a (piecewise smooth) path metric.
We make implicit use of this identification of $\modl$ with a subset of $\hhat S\times \R$ throughout
the paper.

\subsubsection*{Doubly degenerate case} We give first the
description of the model when both $\nu_\pm$ are filling
laminations. In this case $N$ has two simply degenerate ends and no
non-peripheral parabolics, and the main geodesic $g_H$ is doubly infinite.

The blocks are associated to {\em 4-edges}, which are edges $e$ of
geodesics $h\in H$ with $\xi(h)=4$. For each such $e$ the block $B(e)$
is identified with a subset of $D(h)\times \R$ which
is isotopic to  $D(h)\times[-1,1]$. 
More precisely, we can identify each $B(e)$ abstractly with
\begin{align*}
B(e) =  (D(e)\times [-1,1]) \setminus & \left(\
\collar(e^-)\times[-1,-1/2)\right. \union\\ 
& \ \ \left. \collar(e^+)\times(1/2,1]\ \right).
\end{align*}
That is, $B(e)$ is a product with solid-torus ``trenches'' dug out of
its top and bottom corresponding to the vertices $e^\pm$. 
This abstract block is embedded in $M$ {\em flatly}, which means
that  each connected 
leaf of the horizontal foliation
$Y\times\{t\}$ is mapped to a level set $Y\times\{s\}$
in the image, with the map on the first factor being the identity. 
(In \cite{minsky:ELCI} we first build the abstract union of blocks and
then prove it can be embedded).

The 3-holed spheres coming from $(D(e)\setminus \collar(e^\pm))\times
\{\pm 1\}$ are called the {\em gluing boundaries} of 
the block. We show in \cite{minsky:ELCI} that every 3-holed sphere $Y$
that arises as a component domain in $H$ appears as a gluing boundary
of exactly two blocks, and these blocks are in fact attached along
these boundaries via the identity map on $Y$. The resulting level surface 
$Y\times\{s\}$ in $\hhat M=\hhat S\times\R$ will always be denoted $F_Y$.

The complement of the blocks in $\hhat M$ is a union of solid tori of
the form $U(v) = \collar(v)\times I_v$, where $v$ is a vertex in $H$ or a
boundary component of $S$, and $I_v$ is an interval. 

If $v$ is a boundary component of $S$ then $I_v=\R$. Otherwise, 
since we are describing the doubly degenerate case, $I_v$ is always a
compact interval. 

\subsubsection*{Geometry and tube coefficients}
The model is endowed with a metric in which the (non-boundary) blocks
fall into a finite number of isometry classes (in fact two, depending
on the topological type), and in which all the annuli in the
boundaries are Euclidean, with circumference $\epone$.
Thus every torus $\boundary U(v)$ is equipped with a Euclidean
metric. 

This allows us to associate to $U(v)$ a coefficient
$\omega(v)\in\Hyp^2$ (in \cite{minsky:ELCI} denoted $\omega_M(v)$),
defined as follows:  $\boundary U(v)$ comes with a preferred marking
$(\alpha,\mu)$ where $\alpha$ is the core curve of any of the annuli
making up $\boundary U(v)$ and $\mu$ is a meridian curve of the solid
torus $U(v)$. This together with the Euclidean structure on $\boundary
U(v)$ determines a point in the Teichm\"uller space of the torus which
is just $\Hyp^2$. 

This information uniquely determines a metric on $U(v)$ (modulo
isotopy fixing the boundary) which makes it isometric to a 
hyperbolic Margulis tube.  The radius of this tube is given by 
\begin{equation}\label{r tube}
r = \sinh^{-1}\left(\frac{\epone|\omega|}{2\pi}\right)
\end{equation}
and the complex translation length of the
element generating this tube is given (modulo $2\pi i$) by
\begin{equation}\label{lambda tube}
\lambda = h_r\left(\frac{2\pi i}{\omega}\right)
\end{equation}
where $h_r(z) = \Re z \tanh r + i\Im z$
(see \S3.2 of \cite{minsky:ELCI}).
Note in particular that $r$ grows logarithmically with $|\omega|$, and
that for large $|\omega|$, $2\pi i/\omega$ becomes a good approximation
for $\lambda$. 

We adopt, for a general loxodromic isometry, the convention that the
imaginary part of the complex translation distance lie in
$(-\pi,\pi]$. Thus 
in the setting of a marked tube
boundary, when $|\omega|$ is sufficiently large, the expression in
(\ref{r tube}) agrees with this convention.

When $v$ corresponds to a boundary component of $S$ (or, in the
general case, to a parabolic component of $\base \mu_\pm $), 
we write $\omega(v) = i\infty$ and we make $U(v)$ isometric to a cusp
associated to a rank 1 parabolic group.

We let $\modl[0]$ denote the union of the blocks, i.e. $\modl$
minus the interiors of the tubes. For any $k\in[0,\infty]$ we let
$\modl[k]$ denote the union of $\modl[0]$ with the tubes $U(v)$ for
which $|\omega(v)| < k$ (in particular note that $\modl[\infty]$
excludes exactly the parabolic tubes).

We let $\UU$ denote the union of all the tubes in the model, and let
$\UU[k]$ denote those tubes with $|\omega| \ge k$. Thus
$\modl[k] = \modl \setminus \UU[k]$.

\subsubsection*{The case with boundary}
When $N$ has geometrically finite ends, 
$\nu_\pm$ are not filling laminations, the main geodesic $g_H$ is not
bi-infinite, and the model manifold has some boundary. The
construction then involves a finite number of ``boundary blocks.''

A boundary block is associated to a geometrically finite end of
$N^0$. Let $R$ be a subsurface of $S$ homotopic to a 
component of $S\times\{1\}\setminus P_+$ which faces a geometrically
finite end, and let 
$\nu_R$ be the associated component of $\nu_+^T$ in $\Teich(R)$.
We construct a block $B\sbtop(\nu_R)$ as follows:
Let $\T_R$ be 
the set of curves of $\base(\T(H_\nu))=\base(\mu_+)$ that are
contained in $R$. Define 
$$
B'\sbtop(\nu_R) = R\times[-1,0] \setminus \left(\collar(\T_R) \times[-1,-1/2)\right)
$$
and let
$$
B\sbtop(\nu_R) = B'\sbtop(\nu_R) \union \boundary R\times[0,\infty).
$$
This is called a {\em top boundary block}. Its {\em outer boundary}
$\boundary_o B\sbtop(\nu_R)$ is
$R\times\{0\}\union \boundary R\times[0,\infty)$, which we note is
homeomorphic to $int(R)$. This will correspond to a boundary component
of $\hhat C_N$.
The gluing boundary of this block lies on its bottom: it is  
$$
\boundary_-B\sbtop(\nu_R) = (R\setminus \collar(\T_R))\times\{-1\}.
$$
Similarly if $R$ is a  component of $S\times\{-1\}\setminus P_-$
associated to a geometrically finite end, we
let $\I_R = \I(H_\nu)\intersect R$ and define 
$$
B'\sbot(\nu_R) = R\times[0,1] \setminus \collar(\I_R) \times(1/2,1].
$$
and the corresponding {\em bottom boundary block}
$$
B\sbot(\nu_R) = B'\sbot(\nu_R) \union \boundary R\times(-\infty,0].
$$
The gluing boundary here is $\boundary_+B\sbot(\nu_R) = 
(R\setminus \collar(\I_R))\times\{1\}.$ 

The vertical annulus boundaries are now 
$\boundvert B\sbtop(\nu_R) = \boundary R \times [-1,\infty)$
and the internal annuli
$\boundary_i^\pm$ are are a union of possibly several
component annuli, one for each component of $\T_R$ or $\I_R$.

To put a metric on a boundary block, we let $\sigma^m$ denote the
conformal rescaling of the Poincar\'e metric on $\boundary_\infty N$
which makes the collars of curves of length less than $\epone$ into
Euclidean cylinders (and is the identity outside the
collars). Identifying the outer boundary of the block with the
appropriate component of $\boundary_\infty N$ we pull back $\sigma^m$,
and then extend using the product structure of the block.
See \S8.3  of \cite{minsky:ELCI} for details.

\subsection{The bilipschitz model theorem}
\label{bilip model intro}

\subsubsection*{Lipschitz model theorem}
We begin by describing the main theorem of \cite{minsky:ELCI}.
Again, $\rho\in\DD(S)$ is a Kleinian surface group with quotient
manifold $N=N_\rho$ and end invariants $\nu$.

If $U$ is a tube of the model manifold, let $\MT(U)$ denote the
$\epone$-Margulis tube (if any) whose homotopy class is the image via
$\rho$ of the homotopy class of $U$. For $k>0$ let $\MT[k]$ denote the
union of $\MT(U)$ over tubes $U$ in $\UU[k]$.

The {\em augmented convex core} of $N$ is 
$\hhat C_N = C^1_N \union N_{(0,\epzero]}$
where $C^1_N$ is the closed $1$-neighborhood of the convex core
$C_N$ of $N$.  
We show in \cite{minsky:ELCI} that this is homeomorphic to $C_N$, and
hence to $\modl$.

\begin{definition}{model map def}
A {\em $(K,k)$ model map} for $\rho$  is a map $f:\modl \to \hhat C_N$
satisfying the following properties: 
\begin{enumerate}
\item $f$ is in the homotopy class determined by $\rho$, is proper and
  has degree 1.
\item $f$ maps $\UU[k]$ to $\MT[k]$, and $\modl[k]$ to
  $N_\rho\setminus\MT[k]$.
\item $f$ is $K$-Lipschitz on $\modl[k]$, with respect to the induced
  path metric. 
\item $f:\boundary\modl \to \boundary \hhat C_{N_\rho}$
is a $K$-Lipschitz homeomorphism on the boundaries. 
\item $f$ restricted to each tube $U$ in $\UU$ with
  $|\omega(U)|<\infty$ is $K'$-Lipschitz, 
  where $K'$ depends only on $K$ and $|\omega(U)|$.
\end{enumerate}
\end{definition}

\state{Lipschitz Model Theorem. \cite{minsky:ELCI}}{%
There exist $K,k>0$ depending only on $S$, so that for any
$\rho\in\DD(S)$ with end invariants $\nu$ there exists a $(K,k)$ model map
$$
f:\modl \to \hhat C_{N_\rho}.
$$
}

\subsection*{The exterior of the augmented core}
In fact what we really want is a model for all of $N$, not just its 
augmented convex core. Thus we need a description of the the exterior
of $\hhat C_N$.  
This was done in Minsky \cite{minsky:ELCI}, by a slight
generalization of the work of Sullivan and  Epstein-Marden in
\cite{epstein-marden}. Let $E_N$ denote the closure of
$N\setminus\hhat C_N$ in $N$, $\boundary_\infty N$ the conformal
boundary at infinity of $N$, and $\bar E_N = E_N\union \boundary_\infty N$.
The metric $\sigma^m$ on
$\boundary_\infty N$ is defined as in \S\ref{model definitions}, 
as the Poincar\'e metric
adjusted conformally so that every thin tube and cusp becomes a
Euclidean annulus. 

Let $E_\nu$ denote a copy of $\boundary_\infty N \times [0,\infty)$,
  endowed with the metric 
$$
e^{2r}\sigma^m + dr^2
$$
where $r$ is a coordinate for the second factor.

The boundary of $\modl$ is naturally identified with
$\boundary_\infty N$, and this enables us to 
form $\ME_\nu$ as the union of $\modl$ with $E_\nu$ 
identifying $\boundary_\infty N\times\{0\}$ with $\boundary
\modl$. We attach a boundary at
infinity  $\boundary_\infty N\times\{\infty\}$
to $E_\nu$, obtaining a manifold with
boundary $\bME_\nu$. We also denote this boundary at infinity as
$\boundary_\infty \ME_\nu$.

In Lemma 3.5 of \cite{minsky:ELCI} we give a uniformly bilipschitz
homeomorphism of $E_\nu$ to $E_N$, which extends to conformal
homeomorphism on the boundaries at infinity, and together with the
Lipschitz Model Theorem gives the following (called the Extended Model
Theorem in \cite{minsky:ELCI}):

\begin{theorem}{extended model}
There exists a proper degree 1 map 
$$f:\ME_\nu \to N,$$
which is a $(K,k)$ model map from 
$\modl$ to $\hhat C_N$,
restricts to 
a $K$-bilipschitz homeomorphism $\varphi:E_\nu \to E_N$,
and extends to a conformal isomorphism from $\boundary_\infty
\ME_\nu$ to $\boundary_\infty N$.
 The constants $K$ and $k$ depend only on the
topology of $N$.
\end{theorem}

The main result of this paper will be the upgrading of this model
map to a bilipschitz map:

\state{Bilipschitz Model Theorem.}{%
There exist $K,k>0$ depending only on $S$, so that for any
Kleinian surface group $\rho\in\DD(S)$ with
 end invariants $\nu=(\nu_+,\nu_-)$
there is an orientation-preserving $K$-bilipschitz
homeomorphism of pairs
$$
F: (\modl,\UU[k]) \to (\hhat C_{N_\rho},\MT[k])
$$
in the homotopy class determined by $\rho$.
Furthermore this map extends to a homeomorphism
$$
\bar F: \bME_\nu \to \bar N
$$
which restricts to  a $K$-bilipschitz homeomorphism from $\ME_\nu$ to
$N$, and a conformal isomorphism from $\boundary_\infty \ME_\nu$ to
$\boundary_\infty N$. 
}

\subsection{Length estimates}\label{length estimates}
Let $\lambda(g)$ be the complex translation length of an isometry $g$,
where we adopt throughout the convention that $\Im\lambda(g) \in
(-\pi,\pi]$. The real part $\ell(g) = \Re \lambda(g)\ge0$ gives the
translation distance of $g$, and we denote $\ell_\rho(v) =
\ell(\rho(v))$ and $\lambda_\rho(v) = \lambda(\rho(v))$ where $v$
denotes a
closed curve  in $S$ or the corresponding conjugacy class in
$\pi_1(S)$.

The length $\ell_\rho(v)$ of a simple closed curve $v$ 
was bounded above using end-invariant data
in Minsky \cite{minsky:kgcc}. Lower bounds for $\ell_\rho$ and for the
complex translation length $\lambda_\rho$ were obtained
in \cite{minsky:ELCI} using the Lipschitz Model Theorem.
The following is a slight 
restatement of the second main theorem of
\cite{minsky:ELCI}, which incorporates this information. 

\state{Short Curve Theorem.}{%
There exist $\bar \ep>0$ and $c>0$,
and a function $\Omega:\R_+\to\R_+$, depending only on $S$, 
such that the following holds: Given a surface group $\rho\in\DD(S)$, 
and any vertex $v\in\CC(S)$, 
\begin{enumerate}
\item \label{SCL initial}
If $\ell_\rho(v) < \bar\ep$ then $v$ appears in the hierarchy
$H_{\nu_\rho}$.
\item \label{SCL lower}
(Lower length bounds) 
If $v$ appears in $H_{\nu_\rho}$ then 
$$
|\lambda_\rho(v)| \ge \frac c {|\omega(v)|}
$$
and 
$$
\ell_\rho(v) \ge \frac c {|\omega(v)|^2}.
$$
\item \label{SCL upper}
(Upper length bounds)
If $v$ appears in $H_{\nu_\rho}$ and $\ep>0$ then 
$$
|\omega(v)| \ge \Omega(\ep) \implies \ell_\rho(v) \le \ep.
$$
\end{enumerate}
}

The quantity $|\omega(v)|$ can be estimated 
from the lengths of the geodesics in the hierarchy whose domains
border $v$. In particular, 
Theorem 9.1 and Proposition 9.7 of \cite{minsky:ELCI} together imply the following:

\begin{lemma}{big h big omega}
There exist positive constants $b_1$ and $b_2$ depending on $S$ such
that, for any hierarchy $H$ and associated model, if $v$ is any vertex
of $\CC(S)$, 
\begin{equation}\label{omega lower bound}
|\omega(v)| \ge  - b_1 + b_2  \sum_{h\in X_v} |h|
\end{equation}
where $X_v$ is the collection of geodesics $h$  in $H$
such that $v$ is homotopic to a component of $\partial D(h)$ and $|h|$
is the length of $h$.
\end{lemma}

Putting the Short Curve Theorem together with Lemma \ref{big h big
  omega}, we obtain:

\begin{lemma}{big h short l}
There is a function $\LL:\R_+\to\R_+$, depending only on $S$, such
that, given $\rho\in\DD(S)$, for any geodesic $h$ in $H_{\nu_\rho}$, 
and any $\ep>0$, 
\begin{equation}\label{big h short l function}
|h| \ge \LL(\ep) \implies \ell_\rho(\boundary D(h)) \le \ep.
\end{equation}
\end{lemma}

The Bilipschitz Model Theorem proved in this paper will allow us to
give the following improvement of the Short Curve Theorem:

\state{Length Bound Theorem.}{%
There exist $\bar \ep>0$ and $c>0$ depending only on $S$, such that
the following holds:

Let $\rho:\pi_1(S)\to \PSL 2(\C)$ be a Kleinian surface group and
$v$ a vertex of $\CC(S)$, and let $H_{\nu_\rho}$ be an associated
hierarchy.
\begin{enumerate}
\item If $\ell_\rho(v) < \bar\ep$ then $v$  appears in $H_{\nu_\rho}$.

\item  
If $v$ appears in $H_{\nu_\rho}$ then
$$
d_{\Hyp^2}\left(\omega(v),\frac{2\pi i}{\lambda_\rho(v)}\right)
\le c
$$
\end{enumerate}
}

The distance estimate in part (2) is natural because $\omega$ 
is a Teichm\"uller parameter for the boundary torus of a Margulis
tube, as is $2\pi i/\lambda$, and $d_{\Hyp^2}$ is the Teichm\"uller
distance. A bound on $d_{\Hyp^2}$ corresponds directly to a model for
the Margulis tube that is correct up to bilipschitz distortion. See
the discussion in Minsky \cite{minsky:torus,minsky:ELCI}.         
The proof of the Length Bound theorem will be given 
in \S\ref{corollaries}.

\subsubsection*{Improved maps of tubes}
The requirements of the
definition of a $(K,k)$ model map $f$ specify a Lipschitz bound for
the restriction of $f$ to each tube $U$ in $\UU$, but we will need a
little more structure, for technical reasons that occur in the proof of
Theorem \ref{product region}. 

\begin{lemma}{radial tube maps}
Given $(K,k)$ there is a proper function $t:[0,\infty)\to
  [-1,\infty)$, with $t(r) \le r-1$, such that given a
$(K,k)$-model map $f:\modl \to \hhat C_N$ there exists a $(K,k)$
model map $f'$ which agrees with $f$ on $\modl[k]$, and which
satisfies the following for each tube $U$  in $\UU[k]$:
\begin{enumerate}
\item If $r$ is the radius of $U$ then the radius of $\MT(U)$ is at
  least $t(r) + 1$.
\item On the radius $t(r)$ collar neighborhood $U'\subset U$ of $\boundary
  U$, $f'$ takes radial lines to radial lines of $\MT(U)$, and
  preserves distance from the boundary. 
\item $f'$ maps $U\setminus U'$ to $\MT(U) \setminus f'(U')$. 
\end{enumerate}
\end{lemma}

Note that in \cite[\S10]{minsky:ELCI}, the maps on tubes are constructed
using a coning argument. This gives the Lipschitz control of
Definition \ref{model map def}, but does not preserve the foliation by
radial lines. 

\begin{proof}
The existence of the proper function $t(r)$ follows directly from the Short
Curve Theorem -- that is, a deep tube in $\UU$ has large $|\omega|$
and hence the corresponding Margulis tube in $N$ is deep as
well. Property (2) uniquely determines $f'$ on the collar $U'$ of the
boundary of a tube. Thus the only thing to check is the
Lipschitz property from Part (5) of Definition \ref{model
  map def}. On $U'$ this follows from the fact that the level sets of
the distance function from the boundary in a Margulis tube are tori such that radial
projection to the boundary is bilipschitz with constant depending only
on the radii and bounded away from infinity as long as distance to the
core is bounded below. Since $U'$ is at least distance 1 from the
core, and the same for its image, $f'$ inherits a Lipschitz bound from
the values of $f$ on $\boundary U$. On the remaining solid torus
$U\setminus U'$ we can extend by the same coning argument as used in
\cite{minsky:ELCI} (\S10, step 6), to obtain a map with controlled Lipschitz
constant. 
\end{proof}

In Section \ref{insulating} we shall assume that the model maps have
been adjusted to satisfy the conclusions of Lemma \ref{radial tube
  maps}. 

\section{Scaffolds and partial order of subsurfaces}
\label{knotting}

\newcommand\Emb{\operatorname{{\mathbf {emb}}}}
\newcommand\Map{\operatorname{{\mathbf {map}}}}
\newcommand\Cgood{C_{\rm good}}

In this section we will study the topological ordering of surfaces in
a product manifold $M=S\times\R$, where $S$ is a compact surface. 
We will define ``scaffolds'' in $M$, which are
collections of embedded surfaces and solid tori satisfying certain
conditions.   Scaffolds arise naturally in the model manifold as unions
of cut surfaces and tubes. 
Our main goal, encapsulated in Theorem \lref{Scaffold Extension},
is to show that two scaffolds embedded
in $M$ and satisfying consistent topological order relations have
homeomorphic complements.  This will
allow us, in the final part of the proof (\S\ref{region control}), to adjust
our model map to be a homeomorphism on selected regions.

In sections \ref{wrapping coeff} and \ref{some ordering lemmas} we use
the technology developed in the proof of Theorem \lref{Scaffold Extension}
to  develop technical lemmas which will be useful later in the paper.

\subsection{Topological order relation}  
\label{topprec defs}

Recall that $M=S\times \R$ and $\hhat M = \hhat S \times\R$ and that $S$ has been
identified with a compact core for the noncompact surface $\hhat S$.
Let $s_t: \hhat S\to \hhat M$ be the map 
$s_t(x) = (x,t)$, and let $\pi:\hhat M\to \hhat S$ be the map $\pi(x,t)=x$.

For $R\subseteq S$ a connected essential non-annular surface, let
$\Map(R,M)$ denote the homotopy class  $[s_0|_R]$.

If $R$ is a closed annulus we want $\Map(R,M)$ to denote a
certain collection of maps of solid tori into $\hhat M$. 
Thus, we consider
proper maps of the form  $f:V\to \hhat M$ 
where $V=R\times J$, $J$ is a closed connected subset of 
$\R$, and for any $t\in J$
$f\circ s_t:R\to M$ is in $ [s_0|_R]$.
If $R$ is a nonperipheral annulus  then $J$ is a finite or
half-infinite interval. If $R$ is peripheral then we require $J=\R$. 
We say that these maps are of ``annulus type''.

Let $\Map(M)$ denote the disjoint union of $\Map(R,M)$ over all
essential  connected subsurfaces $R$. 

We say that $f\in \Map(R_1,M)$ and $g\in \Map(R_2,M)$ 
{\em overlap} if $R_1$ and $R_2$ have essential intersection.

We now define a ``topological order relation'' $\topprec$ on
$\Map(M)$
(which, despite its appellation, does not
extend to a partial order -- see Example \ref{not acyclic} below).
First, we say that $f:R\to M$ is {\em homotopic to
$-\infty$ in the complement of} $X\subset M$ if for some $r$
there is a proper map
$$F:R\times(-\infty,0] \to M\setminus (X\union S\times [r,\infty))$$
such that $F(\cdot,0) = f$. We define {\em homotopic to
$+\infty$ in the complement of} $X$ similarly. (The definition when
$R$ is an annulus is 
similar, where we then consider the map of the whole solid torus
$V=R\times J$.) 

Now, 
given $f\in\Map(R,M)$ and $g\in\Map(Q,M)$ we write
$f \topprec g$ if and only if

\begin{enumerate}
\item $f$ and $g$ have disjoint images.
\item $f$ is homotopic to $-\infty$ in the complement of $g(Q)$, 
but $f$ is not homotopic to $+\infty$ in the complement of $g(Q)$. 
\item $g$ is homotopic to $+\infty$ in the complement of $f(R)$, 
but $g$ is not homotopic to $-\infty$ in the complement of $f(R)$.
\end{enumerate}

The next lemma states some elementary observations about $\topprec$.

\begin{lemma}{basic ordering properties}
Let $R$ and $Q$ be essential subsurfaces of $S$ which intersect
essentially. 
\begin{enumerate}
\item
If $f\in\Map(R,M)$ and $g\in\Map(Q,M)$ have disjoint images and
$f$ is homotopic to $-\infty$ in the complement of $g(Q)$, 
then $f$ cannot be homotopic to $+\infty$ in the complement of
$g(Q)$. 
\item Similarly if
$g$ is homotopic to $+\infty$ in the complement of $f(R)$, 
then $g$ is not homotopic to $-\infty$ in the complement of $f(R)$.
\item
For the level mappings $s_t(x)=(x,t)$, we have
$$
s_t|_R \topprec s_r|_Q
$$
if and only if $t<r$.
\end{enumerate}
\end{lemma}
\begin{proof}
Since $R$ and $Q$ overlap,
there exist curves $\alpha$ in $R$ and $\beta$ in $Q$ that
intersect essentially. If $f$ is homotopic to both $+\infty$ and 
$-\infty$ in the complement of $g$ then we may construct a  map of 
$\alpha\times\R$ to $M$ which is properly homotopic to the inclusion map
and misses $g(\beta)$.  Since $g(\beta)$ is homotopic
to $\beta\times \{ 0\}$, this contradicts the essential intersection
of $\alpha$ and $\beta$. This gives (1), and a similar argument gives (2).

For (3), it is clear that when $t<r$ $s_r|_Q$ is homotopic to
$+\infty$ in the complement of $s_t(R)$, and that $s_t|_R$ is
homotopic to $-\infty$ in the complement of $s_r(Q)$. The rest follows
from (1) and (2).
\end{proof}

\begin{example}{not acyclic}
The relation $\topprec$ does not extend to a partial order on
$\Map(M)$, because it contains cycles. Let $\gamma_1$, $\gamma_2$ and $\gamma_3$ be
three disjoint curves on a surface $S$ of genus 4 such that the components
$A_1$, $A_2$ and $A_3$ of $S\setminus \union_{i=1}^3 \collar(\gamma_i)$ are
all twice-punctured tori. Moreover, we may assume that
$P=A_1\cup \collar(\gamma_1)\cup A_2$, $Q=A_2\cup\collar(\gamma_2)\cup A_3$
and $R=A_3\cup \collar(\gamma_3)\cup A_1$, are all connected.
Let $f:P\to M$ map $A_1$ to $A_1\times \{ 0\}$, $A_2$ to $A_2\times \{ 1\}$
and $\collar(\gamma_1)$ to an annulus in $\collar(\gamma_1)\times [0,1]$.
Similarly, let $g:Q\to M$ map $A_2$ to $A_2\times \{ 0\}$, $A_3$ to
$A_3\times \{ 1\}$ and $\collar(\gamma_2)$ to an annulus in
$\collar(\gamma_2)\times [0,1]$, and let $h:R\to M$ map $A_3$
to $A_3\times \{ 0\}$, $A_1$ to $A_1\times \{ 1\}$
and $\collar(\gamma_3)$ to an annulus in $\collar(\gamma_3)\times [0,1]$.
It is clear that $f\topprec h\topprec g \topprec f$.
\end{example}

\subsubsection*{Ordering disconnected surfaces}
One can extend $\topprec$ to maps $f:R\to M$ where $R$ is a
disconnected subsurface of $S$, with a bit of care. We say that a (possibly)
disconnected subsurface $R$ of $S$ is essential and non-annular, if
each of its components is essential and non-annular. If $R$ and $T$ are
essential non-annular subsurfaces which intersect essentially
and $f:R\to M$ and $g:T\to M$ are in the homotopy class of $s_0|_R$ and $s_0|_T$,
we say that $f\topprec g$ provided that, for any pair $R'$ and $T'$ of
intersecting components of $R$ and $T$, we have
$f|_{R'} \topprec g|_{T'}$. 
(A similar definition can be made which allows for annular components and their corresponding
maps of solid tori.)

It is easy to see, using Lemma \ref{basic ordering properties},
that 

\begin{lemma}{subsurface inherits}
Let $R$ and $T$ be essential, non-annular subsurfaces of $S$.
If $f:R\to M$ and $g:T\to M$ are in the homotopy class of $s_0|_R$ and $s_0|_T$,
then $f\topprec g$ implies that $f|_{R_0}\topprec g|_{T_0}$ whenever $R_0$ and $T_0$
are essential non-annular subsurfaces of $R$ and $T$ which intersect essentially.
\end{lemma}

Note that $R_0$ and $T_0$ are not {\em components} of $R$ and $T$,
just subsurfaces. This will be applied in Section 
\ref{preserving order of embeddings}, where $R_0$
and $T_0$ are equal to $R$ and $T$ minus a union of annuli.

\subsubsection{Embeddings and Scaffolds}
Let $\Emb(R,M)$ be the set of images of 
those maps in $\Map(R,M)$ that are
embeddings. When $R$ is a (closed) annulus 
we also require the map of the solid torus $R\times J$ to be proper
and orientation-preserving. 
Define $\Emb(M) = \union_R \Emb(R,M)$.
Note that an embedding in $\Map(M)$ is determined by its image, up to
reparametrization of the domain of the map by a map homotopic to the
identity, and it follows that we can extend
$\topprec$ to a well-defined relation on $\Emb(M)$. Similarly we can
define the notion of ``overlap'' of members of $\Emb(M) $ to coincide
with overlap of their domains.

A non-annular surface in $\Emb(M)$ is {\em straight} if it is a level
surface, i.e. of the form $R\times\{t\}$. 
A solid torus in $\Emb(M)$ is {\em straight} if it is of the form
$\bcollar(v)\times J$ for some $v$, where $J$ is a closed connected
subset of $\R$. If $v$ is nonperipheral in
$S$ we allow $J$ to be of the form $[a,b]$, $[a,\infty)$ or
$(-\infty,b]$. If $v$ is a component of $\boundary S$ then we require
$J=\R$. 

We are now ready to define scaffolds, which are the primary object of
study in this section.

\begin{definition}{def scaffold}
A {\em scaffold} $\Sigma \subset \hhat M$ is a union of two sets
$\FF_\Sigma$ and $\VV_\Sigma$, where
\begin{enumerate}
\item $\FF_\Sigma$ is a finite disjoint union of elements of $\Emb(M)$,
of non-annulus type.
\item $\VV_\Sigma$ is a finite disjoint union of elements in $\Emb(M)$ of
annulus type (that is, solid tori).
\item \label{straight tori}
$\VV_\Sigma$ is {\em unknotted and unlinked}: it is isotopic in $M$ to
a union of  straight solid tori.
\item $\FF_\Sigma$ only meets $\VV_\Sigma$ along boundary curves of 
surfaces in $\FF_\Sigma$, and conversely for every component $F$ of
$\FF_\Sigma$, $\boundary F$ is embedded in $\boundary \VV_\Sigma$.
\item No two elements of $\VV_\Sigma$ are homotopic.
\end{enumerate}
The components of $\FF_\Sigma$ and $\VV_\Sigma$ are called the
{\em pieces} of $\Sigma$.
\end{definition}

In a straight solid torus $V$ let the {\em level homotopy class}
denote the homotopy class in $\boundary V$ of the curves of
the form $\alpha\times \{t\}$ where $t\in J$ and $\alpha$ is an essential curve in $S$.
If $V$ is isotopic to a straight solid torus we define the
level homotopy class as the one containing the isotopes of the level
curves. The following lemma gives us a common situation 
where elements of $\Map(M)$ take boundary curves to level homotopy classes.

\begin{lemma}{level curves immersed}
Let $\VV\subset M$ be a union of disjoint straight solid tori, no two of which are
homotopic, and let $R$ be an essential non-annular subsurface. If $h\in \Map(R,M)$,
$h(\boundary R)
\subset \boundary \VV$ and $h(R)\intersect int(\VV) = \emptyset$, then
$h$ maps each component of $\boundary R$ to the level homotopy class
of the corresponding component of $\VV$. 
\end{lemma}

\begin{proof}
Let $\gamma$ be a boundary component of $R$ such that $h(\gamma) $ is
contained in a component
$V=\bcollar(v)\times J$ of $\VV$. 
  
If $v$ is peripheral in $M$ then $J=\R$, $\boundary V$ is an
annulus, and there is a unique homotopy class in $\boundary V$
representing the core, the level homotopy class. 
and $h(\gamma)$, being primitive, must be in the 
the level homotopy class. Hence we are done. 

Thus we may assume $v$ is nonperipheral, $J\ne\R$, and without loss of
generality that $b = \sup J < \infty$. Let
$\gamma_b\subset \boundary V$ be the level curve
$\gamma_v\times \{b\}$ where $\gamma_v$ is the standard
representative of $v$ in $S$. 
Since $h(\gamma)$ is homotopic to $\gamma_b$ in $V$, to show that they
are homotopic in $\boundary V$ it suffices to show that the
algebraic intersection number $h(\gamma)\cdot\gamma_b$ vanishes.

Since $R$ is an essential subsurface of $S$, $h(\boundary R)$ meets
$\boundary V$ in either one or two curves. Consider first the
case that $\gamma$ is the only component of $\boundary R$ mapping to
$\boundary V$. 

Compactify $M$ to get $\bbar M =  S\times[-\infty,\infty]$.
Let $\VV_h$ denote the union of components of $\VV$ meeting
$h(\boundary R)$, and let $X = \bbar M \setminus int(\VV_h)$. 
Let $B\subset \bbar M$ be the annulus $\gamma_v \times
[b,\infty]$. Since the components of $\boundary R$ do not overlap, the
solid tori of $\VV_h$ have disjoint projections to $S$ and it follows
that $B$ is contained in $X$. Let $A = \boundary V \union
S\times\{\infty\}$. $B$ determines a class
$[B]\in H_2(X,A)$, and intersection number with $B$
gives a cohomology class (its Lefschetz dual) $i_B\in H^1(X,\partial X-A)$.

If $\alpha$ is a closed curve in $\boundary V$,
then $i_B(\alpha)$ is the algebraic intersection number
of $\alpha$ and $\partial B$ on $\partial V$.
In particular, since
$\boundary B\intersect \boundary V = \gamma_b$, we have
$$
h(\gamma)\cdot \gamma_b = i_B(h(\gamma)).
$$
Now since $h^{-1}(A)\intersect \boundary R = \gamma$, we find that $[h(\gamma)]$
vanishes in $H_1(X,\boundary X-A)$, and it follows that $i_B(h(\gamma)) =
0$. This concludes the proof in this case. 

If $h$ takes two components of $\boundary R$ to $\boundary V$, there is a
double cover of $M$ to which $h$ has two lifts, each of
which has no pair of homotopic boundary components. Picking one
of these lifts and a lift of $V$, we repeat the above argument in this
cover. 
\end{proof}

As an immediate consequence we obtain a statement for scaffolds:
\begin{lemma}{level curves}
Let $\Sigma$ be a scaffold in $\hhat M$. The intersection curves
$\FF_\Sigma \intersect \VV_\Sigma$ are in the level homotopy classes
of the components of $\boundary\VV_\Sigma$.
\end{lemma}

\begin{proof}
By property (\ref{straight tori}) of the definition, and the isotopy
extension theorem \cite{rourke-sanderson}, we may assume that $\VV_\Sigma$
is a union of straight solid tori. We then apply Lemma \ref{level
  curves immersed} to each component of $\FF_\Sigma$. 
\end{proof}

\begin{definition}{def scaffold straight}
A scaffold $\Sigma$ is {\em straight} if every piece of $\Sigma$ is straight.
\end{definition}

Let $\topprec|_\Sigma$ denote the restriction of the $\topprec$
relation to the pieces of $\Sigma$. We can capture the essential
properties of being a straight scaffold with this definition
(see Lemma \ref{straighten scaffold}).

\begin{definition}{def comb straight}
A scaffold $\Sigma$ is
{\em combinatorially straight} provided  $\topprec|_\Sigma$ satisfies
these conditions:
\begin{enumerate}
\item (Overlap condition) Whenever two pieces $p$ and $q$
of $\Sigma$ overlap, either
$p\topprec q$ or $q\topprec p$,
\item (Acyclic condition) The transitive closure of $\topprec|_\Sigma $ is
a partial order. 
\end{enumerate}
\end{definition}

Notice that one may use a construction similar to the one in
Example \ref{not acyclic} to construct scaffolds which satisfy
the overlap condition but not the acyclic condition.

\subsection{Scaffold extensions and isotopies}

Our technology will allow us to study ``good'' maps of scaffolds
into $M$, those which, among other things, have a scaffold as image
and preserve the topological ordering of pieces.

\begin{definition}{def good map}
Let $\Sigma$ be a scaffold in $\hhat M$. 
A map $f:\Sigma \to \hhat M$ is a {\em good scaffold map} if the following
holds:
\begin{enumerate}
\item $f$ is homotopic to the identity.
\item 
$f(\Sigma)$ is a scaffold $\Sigma'$ with
$\VV_{\Sigma'}=f(\VV_\Sigma)$, and 
$\FF_{\Sigma'} = f(\FF_\Sigma)$
\item
For each component $V$ of $\VV_\Sigma$, $f(V)$ is a component of
$\VV_{\Sigma'}$, and 
$f|_V: V\to f(V)$ is proper.
\item $f$ is an embedding on  $\FF_{\Sigma}$.
\item $f$ is order-preserving. That is, for any two pieces $p$ and $q$ of
$\Sigma$ which overlap, 
$$ f(p) \topprec f(q) \Iff  p \topprec q.$$
\end{enumerate}
\end{definition}

The main theorem of this section gives that a well-behaved map $F$ of
$\hhat M$ to itself which restricts to a good scaffold map of a
combinatorially straight scaffold $\Sigma$ is homotopic to a homeomorphism
which agrees with $F$ on $\FF_{\Sigma}$.
In particular this implies that the complements of $\Sigma$ and $F(\Sigma)$
are homeomorphic. 

\begin{theorem+}{Scaffold Extension}
Let $\Sigma\subset \hhat M$ be a
combinatorially straight scaffold, and let 
$F:\hhat M\to \hhat M$ be a proper degree 1 map homotopic to the
identity, such that 
$F|_\Sigma$ is a good scaffold map, and
$F(M\setminus int(\VV_\Sigma)) \subset M\setminus int(F(\VV_\Sigma))$.

Then there exists a homeomorphism $F':\hhat M\to \hhat M$, homotopic
to $F$, such that
\begin{enumerate}
\item $F'|_{\FF_\Sigma} = F|_{\FF_\Sigma}$
\item On each component $V$ of $\VV_\Sigma$, $F'|_V$ is homotopic to
$F|_V$ rel $\FF_\Sigma$, through proper maps $V\to F(V)$. 
\end{enumerate}
\end{theorem+}

We will derive the Scaffold Extension theorem from the Scaffold Isotopy  
theorem, which essentially states that the image of a good scaffold
map can be ambiently isotoped back to the original scaffold.
The proof of the Scaffold Isotopy theorem will be deferred to
section \ref{isotopy proof}.

\begin{theorem+}{Scaffold Isotopy}
Let $\Sigma$ be a straight scaffold, and let $f:\Sigma \to \hhat M$ be a
good scaffold map. There exists an isotopy
$H:\hhat M\times[0,1] \to \hhat M$ 
such that $H_0=id$, $H_1\circ f(\Sigma) = \Sigma$, and
$H_1\circ f$ is the identity on $\FF_\Sigma$.
\end{theorem+}

\subsection{Straightening}

Now assuming Theorem \ref{Scaffold Isotopy}, let us
prove the following corollary, which allows us to treat
combinatorially straight scaffolds as if they were straight.

\begin{lemma}{straighten scaffold}
A scaffold is combinatorially straight if and only if it is ambient
isotopic to a straight scaffold. 
\end{lemma}

\begin{pf}
If $\Sigma$ is straight, then by Lemma \ref{basic ordering properties}
two disjoint pieces are ordered whenever they overlap.
Lemma \ref{basic ordering properties} also implies that 
$\topprec$ is determined
by the ordering of the $\R$ coordinates, and must therefore be
acyclic. Hence 
$\Sigma$ is combinatorially straight. The property of being combinatorially
straight is preserved by isotopy and hence holds for any scaffold
isotopic to a straight scaffold. This concludes the ``if'' direction.

Now suppose that $\Sigma$ is combinatorially straight.
Let $\PP$ denote the set of pieces of $\Sigma$.
We will first construct a straight scaffold $\Sigma_0$ together with a
bijective correspondence $c:\PP\to\PP_0$ taking pieces of $\Sigma$ to
pieces of $\Sigma_0$, such that whenever $p$ and $q$ are overlapping
pieces of $\Sigma$, we have $p\topprec q \iff c(p) \topprec c(q)$.

Let $\topprec'$ denote the transitive closure of $\topprec|_\Sigma$,
which by hypothesis is a partial order. It is then an easy
exercise to show that there is a map $l:\PP\to\Z$ which is order
preserving, i.e. 
$p\topprec' q \implies l(p)<l(q)$.

Now for each component $F$ of $\FF_\Sigma$ let $c(F)$ be the 
level embedding $s_{l(F)}(F)$, and let $\FF_{\Sigma_0}$ be
the union of these level embeddings. 
Two components of $\FF_\Sigma$ have the same $l$ value only if they are
unordered, and since $\topprec$ satisfies the overlap condition, this
implies they have disjoint domains. It follows that these level
embeddings are all disjoint.

We next construct the solid tori in $\VV_{\Sigma_0}$.
Let $V$ be a component of $\VV_\Sigma$, and
let $v$ be its homotopy class in $S$. Recall that $V$ is isotopic to 
$\bcollar(v) \times J$ where $J\subset\R$ is closed and connected. 
The solid torus $c(V)$ will be $\bcollar(v)\times J_0$, 
where $J_0 = [a,b]\intersect \R$, with $a=a(V)$ and $b=b(V)$ defined
as follows.  Let $\beta(V)$ be the set of $l$ values for surfaces
bordering $V$.  
If $\inf J = -\infty$ let $a=-\infty$, and if $\sup J = \infty$ let
$b=\infty$. In all other cases, 
let
$$a(V) = \min \beta(V)\union\{l(V)\} - 1/3$$
 and 
$$b(V) = \max \beta(V)\union\{l(V)\} + 1/3.$$ 
The union of $c(V)$ over $V\in\VV_\Sigma$  gives $\VV_{\Sigma_0}$.

Note that this definition implies that, whenever any component $F$ of
$\FF_\Sigma$ and $V$ of $\VV_\Sigma$ intersect along a boundary 
component $\gamma$ of $F$, the corresponding $c(F)$ and $c(V)$
intersect along the boundary component $\gamma'$ of $c(F)$
corresponding to $\gamma$. 

Once we check that these are the only intersections between the pieces
of $\Sigma_0$, it will follow that $\Sigma_0$ is a (straight) scaffold. The
order-preserving property of the correspondence will follow from the
same argument. We have already observed that all the level embeddings
of pieces of $\FF_\Sigma$ are disjoint, so it remains to consider
intersections involving solid tori. 

First let us establish the following claim about the ordering. 
Suppose that $p$ is a piece of $\Sigma$ which overlaps 
$V\in \VV_\Sigma$ and $F$ is a piece of $\FF_\Sigma$ with a boundary component
$\gamma$ in $\boundary V$ (hence $l(F)\in \beta(V)$). We claim that
\begin{equation}\label{F P ordering 1}
V \topprec p  \implies F \topprec p
\end{equation}
and 
\begin{equation}\label{F P ordering 2}
p \topprec V  \implies p \topprec F.
\end{equation}
Since $V$ and $p$ overlap, $F$ and $p$ must also overlap and hence are
$\topprec$-ordered because $\Sigma$ is combinatorially straight.
Thus we only have to rule out $V\topprec p \topprec F$ and 
$F\topprec p\topprec V$. Assume the former without loss of
generality. $V\topprec p$ implies that $V$ is homotopic to $-\infty$
in the complement of $p$. Since $V$ is homotopic into $F$ by a
homotopy taking $V$ into itself, $p\topprec F$ implies that $V$ is
homotopic to $+\infty$ in the complement of $p$. This is a
contradiction by Lemma \ref{basic ordering properties}.
Together with the corresponding argument for
the case $F\topprec p\topprec V$, this establishes
(\ref{F P ordering 1}) and (\ref{F P ordering 2}).

Now if $p$ is a piece of $\FF_\Sigma$ overlapping $V\in\VV_\Sigma$, suppose
without loss 
of generality that $V \topprec p$. Then $l(V) < l(p)$, and 
for each $F$ with boundary
component on $V$ we have $F \topprec p$ by (\ref{F P ordering 1}), and
hence $l(F) < l(p)$. Furthermore note that it is not possible to have
$b(V)=\infty$ in this case because then $V\topprec p$ could not hold. 
It follows that $b(V) < l(p)$, and therefore
$c(V)$ and $c(p)$ are disjoint and $c(V) \topprec c(p)$. 

If $V'$ is a component of $\VV_\Sigma$ overlapping $V$, and (without loss of
generality) $V\topprec V'$, a similar argument
implies that $b(V) < a(V')$, and again $c(V)$ and $c(V')$ are disjoint,
and $c(V) \topprec c(V')$.

This establishes disjointness for overlapping pieces of
$\Sigma_0$. Disjointness for non-overlapping pieces is immediate from
the definition of $\Sigma_0$. We have also established one direction
of order-preservation, namely $p\topprec q \implies c(p) \topprec
c(q)$ for overlapping pieces $p$ and $q$ of $\Sigma$. For the opposite
direction we need just observe that $c(p)$ and $c(q)$ must be
$\topprec$-ordered, so that the opposite direction follows from the
forward direction with the roles of $p$ and $q$ reversed.

Next, we construct a good scaffold map $h:\Sigma_0 \to \Sigma$: On each
component $F_0=c(F)$ of $\FF_{\Sigma_0}$, by construction, there is a
homeomorphism to $F$, in the homotopy
class of the identity. After defining $h$ on $\FF_{\Sigma_0}$, for each
$V_0=c(V)$ in $\VV_{\Sigma_0}$, $h$ is already defined on the circles
$\boundary V_0 
\intersect \FF_{\Sigma_0}$, and it is easy to extend this to a
proper map of $V_0$ to the corresponding solid torus $V$ (although the
map is not guaranteed to be an embedding). If $V_0$ does not meet any
component of $\FF_{\Sigma_0}$ then we simply define the map $V_0\to V$ to be a
homeomorphism that takes level curves of $\boundary V_0$ to the
homotopy class of level curves on $\boundary V$. 
This gives $h$, which satisfies all the conditions of a good scaffold
map: property (5), order preservation, follows from the
order-preserving property of the correspondence $c$, 
while the other properties are implicit in the construction.

Now apply Theorem \ref{Scaffold Isotopy} to the map $h$, producing an
isotopy $H$ with $H_0=id$ and  $H_1\circ h(\Sigma_0) = \Sigma_0$. 
Thus $H_1(\Sigma) = \Sigma_0$, and we have exhibited the desired
ambient isotopy from $\Sigma$ to a straight scaffold.
\end{pf}

\subsection{Proof of the scaffold extension theorem}
\label{proof extension}

We now give the proof of Theorem \ref{Scaffold Extension}, again
assuming Theorem \ref{Scaffold Isotopy}.

First we note that it suffices to prove the theorem when $\Sigma$ is
straight. For by Lemma \ref{straighten scaffold}, if $\Sigma$ is
combinatorially straight there is a homeomorphism $\Phi:M\to M$ isotopic to the
identity such that $\Phi(\Sigma)$ is straight. We can apply the result
for $\Phi(\Sigma)$ and conjugate the answer by $\Phi$. 

Denote $\FF=\FF_\Sigma$ and $\VV = \VV_\Sigma$. 
Apply the Scaffold Isotopy theorem \ref{Scaffold Isotopy} to the good
scaffold map
$F|_\Sigma$, obtaining an isotopy $H$ of $M$ with $H_0=id$,
$H_1\circ F(\Sigma) = \Sigma$, and $H_1\circ F$ equal to the identity
on $\FF$. Let $F_1 = H_1 \circ F$.

Our desired map will just be $F' = H_1^{-1}$. We have immediately that
 $F'(\Sigma) = F(\Sigma)$, and $F' = F$ on $\FF$. 
It remains to show that $F_1|_V$, for any component $V$ of $\VV$, is homotopic
rel $\FF$ through proper maps of $V$  to the identity. Composing this 
homotopy by $H_1^{-1}$, we will have that 
$F'|_V$ is homotopic  rel $\FF$ through proper maps of $V$ to $F$, as desired.

Fixing a component $V$ of $\VV$, let us first find a
homotopy of $F_1|_{\boundary V}$, through maps of $\boundary V \to
\boundary V$ fixing $\FF\intersect \boundary V$ pointwise,
to the identity on $\boundary V$.

We claim that $F_1|_{\boundary V}$ preserves the level homotopy
class. If $\boundary V$ meets $\FF$, this is immediate since $F_1$
fixes $\FF$ and $\Sigma$ is straight.
If $V$ is disjoint from $\FF$, consider a level curve
$\gamma$ in $\boundary V$. Since $\gamma$ is homotopic to $+\infty$ in
the complement of $int(V)$, and since $F_1$ takes $\hhat M\setminus int(V)$
to itself, it follows that $F_1(\gamma)$ is homotopic to $+\infty$ in
the complement of $int(V)$. Adding
a boundary $\hhat S\times\{\infty\}$ to $\hhat M$, we conclude
that $\gamma$ and $F_1(\gamma)$ are homotopic within $\hhat M
\setminus int(V)$ to curves in this boundary, which are of course
homotopic and hence have vanishing algebraic intersection number. 
Thus the intersection number
$F_1(\gamma)\cdot\gamma$ on $\boundary V$ vanishes as well,
so as in
Lemma \ref{level curves} 
it follows that
$F_1(\gamma)$ and $\gamma$ are homotopic in $\boundary V$. 

$F_1|_{\boundary V}$ also takes meridians to
powers of meridians, since it is a restriction of a self-map of $V$.
Thus it is homotopic to the identity provided $F_1|_{\boundary
V}:\boundary V \to \boundary V$ has degree 1, or 
equivalently that $F_1|_V:V\to V$ has degree 1.  
Since $F_1$ takes $M\setminus int(\VV)$  to $M\setminus
int(\VV)$ by the hypotheses of the theorem, and no two  components of
$\VV$ are homotopic, $F_1$ must take $M\setminus int(V)$ to
$M\setminus int(V)$, and it follows that the degree of $F_1$ on $V$ equals
the degree of $F_1$, which is 1 by hypothesis.

If $\boundary V$ meets no components of $\FF$, this suffices to give
the desired homotopy of $F_1|_{\boundary V}$ to the identity. 
In the general case, $\FF$ may meet $\boundary V$ in one or more
level curves, which break up $\boundary V$ into annuli.
For each such
annulus $A$ we must show that $F_1|_A:A \to \boundary V$ is homotopic
to the identity rel $\boundary A$. 
Let $\gamma\subset S$ be the curve
in the homotopy class of $V$, so that $V = \collar(\gamma) \times J$.

There are two obstructions to this, which we call $d(A)$ and $t(A)$.
Inducing an orientation on $A$ from $\boundary V$,
the 2-chain $A-F_1(A)$ is closed, and determines a homology class in
$H_2(\boundary V)$. If $V$ is compact then $H_2(\boundary V)=\Z$
and we may define $d(A)$ by the equation $[A-F_1(A)] = d(A)[\boundary
V]$. If $V$ is
noncompact then $\boundary V$ is an annulus, $H_2(\boundary
V)=0$, and we define $d(A)=0$.

To define $t(A)$, choose an arc
$\alpha$ connecting the boundary
components of $A$, and note that $\alpha * F_1(\alpha)$ is a closed curve,
which is homotopic in the solid torus to some multiple of
$\gamma$. Let $t(A)$ be this multiple. 

If $d(A)=0$, then $F_1|_A$ is homotopic in $\boundary V$, rel
$\boundary A$, to a homeomorphism from $A$ to itself which is the
identity on the boundary. The number
$t(A)$ then measures the Dehn-twisting of this homeomorphism, and if
$t(A)=0$ then the homeomorphism is homotopic, rel boundary, to the
identity. Thus we must establish that $d(A)=0$ and $t(A)=0$. 

To prove $t(A)=0$, recall that $F_1$ and $F$ are homotopic. Since 
$F$ in turn is homotopic to the identity, there is a homotopy
$G:M\times[0,1] \to M$ with $G_0 = id$ and $G_1=F_1$.
Since $F_1|_{\FF}$ is the identity, the trajectories
$G(x\times[0,1])$ for $x\in\FF$ are closed loops. We claim these loops
are homotopically trivial: Let $F$ be the component of $\FF$
containing $x$. Since $F$ is not an annulus,
$x$ is contained in two loops $\xi$ and
$\eta$ in $F$ that are not commensurable (say that two elements in
$\pi_1(M,x)$ are commensurable if they are powers of a common element.)
$G(\xi\times[0,1])$ is a torus and hence 
homotopically non essential in $M$, so $G(x\times[0,1])$ is commensurable
with $\xi$. Similarly it is commensurable with $\eta$. But since $\xi$
and $\eta$ are not commensurable, $G(x\times[0,1])$ must be
homotopically trivial.

Now place any complete, non-positively curved metric on $M$ for which
all the solid tori are convex
(e.g. put a fixed hyperbolic metric on $S$ and take the product metric
on $S\times\R$)
and deform the trajectories
of $G$ to their unique geodesic representatives. The result is a new
homotopy $G'$ from the identity to $F_1$,  which is constant on $\FF$.
It follows that the arc $\alpha$, whose endpoints are on $\boundary
A$, is homotopic {\em rel endpoints}, inside $V$, to
$F_1(\alpha)$. Hence the loop $\alpha* 
F_1(\alpha)$ bounds a disk in $V$, so that $t(A)=0$.

Next we argue that $d(A)=0$, which breaks down into several cases. 
We may assume that $\boundary V$ is a torus, since $d(A)=0$ by
definition when $V$ is non-compact. 

Suppose that $A$ is of the form $\beta \times [s,t]$ with
$[s,t]\subset J $ (where $V=\collar(\gamma)\times J$) and
$\beta$ is a boundary component of $\collar(\gamma)$. 

If $\gamma$ separates $S$, let $R$ be the component of $S\setminus
\collar(\gamma)$ which is adjacent to $\beta$, and let 
$Q = R\times [s,t]$. Let $B$ be the vertical annulus $\gamma \times
[\max J,\infty]$ in $\bbar M = S\times[-\infty,\infty]$. With the 
natural orientation, the intersection 
$B\intersect F_1(\boundary Q)$ (which we may assume transverse) defines a 
class in $H_1(X)$. This is just the intersection pairing
$i:H_2(X)\times H_2(X,\boundary X)\to H_1(X)$, where 
$X=\bbar M \setminus int(V)$ as in Lemma \ref{level curves}.
In our case since all the 
intersection curves are trivial or homotopic to $\gamma$,
this class is a multiple of $[\gamma]$.
Let $i(F_1(\boundary Q),B)$ denote this multiple. 
As $F_1$ maps $M\setminus \VV$ to $M\setminus \VV$, it follows that
$F_1(Q)\subset X$, so $[F_1(\boundary Q)]=0$ in $H_2(X)$.
Therefore, $i(F_1(\boundary Q),B)=0$.

We will show that $i(F_1(\boundary Q), B) = \pm d(A)$, which implies
that $d(A)=0$.
The components of $\boundary R$ other than $\beta$ are in $\boundary
S$, and hence map to $\boundary S\times \R$, and miss $B$.
Hence $F_1^{-1}(B)\intersect \boundary Q$ is contained in
the surfaces $\boundary_+Q = R\times\{t\}$, $\boundary_-Q =
R\times\{s\}$, and $A$. 
Now $\boundary_+Q$ and $\boundary_-Q$
contain components $Y$ and $Z$ of $\FF$ which meet $V$
at the boundary of $A$. Since $F_1$ is the identity on the (straight)
pieces $Y$ and $Z$, 
$F_1(Y)$ and $F_1(Z)$ are disjoint from $B$. Thus, any curve of
$F_1^{-1}(B)\intersect \boundary_{\pm}Q$ lies in a component
of $\boundary_\pm Q$ which does not contain any curves homotopic to $\gamma$,
so must be homotopically trivial. 
We conclude that $i(F_1(\boundary Q),B) = i(F_1(A),B)$.
But $F_1(A)$ only meets $B$ in its boundary curve
$\gamma\times \{q\}$, and the number of essential intersections,
counted with 
signs, is exactly the degree  with which $F_1(A)$ covers the
complementary annulus $\boundary V \setminus A$. Hence this is (up to sign)
the degree with which $A-F_1(A)$ covers $\boundary V$, 
which is $d(A)$.

Next consider the case that $\gamma$ does not separate $S$. There is a
double cover of $S$ to which $\collar(\gamma)$ lifts to two disjoint
copies $C_1$ and $C_2$, which separate $S$ into $R_1$ and $R_2$. 
In the corresponding double cover of $M$ we have two lifts $V_1$ and
$V_2$ of $V$. Letting $Q = R_1 \times [s,t]$, we can repeat the above
argument to obtain $d(A_1)=0$, where $A_1$ is the lift of $A$ to
$V_1$. Projecting back to $M$, we have $d(A)=0$.

Another type of annulus $A$ is one that contains the bottom annulus
$\collar(\gamma) \times \{\min J\}$, and whose boundaries are
of the form $\beta \times \{s\}$ and $\beta'\times\{s'\}$, where
$\beta'$ is the other boundary component of $\collar(\gamma)$,
and $s,s'\in int J$ are heights where components of $\FF$ meet $V$
on its two sides. Suppose again that $\collar(\gamma)$ separates $S$
into two components $R$ and $R'$ adjacent to $\beta$ and $\beta'$
respectively. 
In this case we let $Q$ be the region of $M$ below
$V \union R\times\{s\}\union R'\times\{s'\}$. Again we can show that 
$d(A) = \pm i(F_1(\boundary Q),B) = 0.$
If $\gamma$ is non-separating we can again use a double cover.
The case where $A$ contains the top annulus
$\collar(\gamma)\times \{\max J\}$ is treated similarly.

Finally it is possible that $\boundary V$ only meets $\FF$ on one
side, say on $\beta\times J$, and hence there may be one annulus
$A$ which is the closure of the complement of $\bar A = \beta\times
[s,t]$ for some $[s,t]\subset int J$ (possibly $s=t$, when $\boundary V$
meets $\FF$ in a unique circle).
In this case, $[\bar A-F_1(\bar A)] =0$ since
$\bar A$ is a
concatenation of annuli for which we have already proved $d=0$
(or a single circle viewed as a singular 2-chain which is fixed by $F_1$).
Since
$[A-F_1(A)]+[\bar A-F_1(\bar A)]=[\boundary V -F_1(\boundary V)]$,
and we have already shown
that  $F_1$ takes $\boundary V$ to itself with degree 1, this is 0
and it follows that $d(A)=0$ as well. 

We conclude, then, that $F_1|_A$ can be deformed to the identity rel
$\boundary A$, for each annulus $A$.
The resulting homotopy of $F_1|_{\boundary V}$ to the identity can
be extended to the interior of the solid torus $V$ by a coning argument, 
yielding a homotopy to the identity through proper maps $V\to V$,
fixing $\FF \intersect \boundary V$. 

Thus, we can now let
$H_1^{-1}$ be the desired map $F'$, and this concludes the proof. 
\qed

\subsection{Intersection patterns}
Before we prove Theorem \ref{Scaffold Isotopy},
we will need to consider carefully the ways in which embedded surfaces
intersect in $M$. 

\subsubsection*{Pockets}

A pair $(Y_1,Y_2)$ of connected, compact incompressible surfaces in $M$
is a {\em parallel pair} if $\boundary Y_1 = \boundary Y_2$,
$\interior(Y_1) \intersect \interior(Y_2) = \emptyset$, and
there is a homotopy $H:Y_1\times[0,1]\to M$ such that $H_0$ is the
inclusion of $Y_1$ into $M$ and $H_1:Y_1\to Y_2$ is a homeomorphism that takes
each boundary component to itself. (Note that, if no two boundary
components of $Y_1$ are homotopic in $M$, the last condition on $H_1$ is
automatic).

\begin{lemma}{pockets}
A parallel pair $(Y_1,Y_2)$ in $M$
is the boundary of a unique compact 
region, which is homeomorphic to 
$$ 
Y_1 \times [0,1] / \sim
$$ 
where $(x,t)\sim (x,t')$ for any $x\in\boundary Y_1$ and
$t,t'\in[0,1]$, by a homeomorphism taking $Y_1\times\{0\}$ to $Y_1$ and
$Y_1\times\{1\}$ to $Y_2$.
\end{lemma}
This region is called a {\em pocket}, and the surfaces $Y_1$ and $Y_2$ 
are its {\em boundary surfaces}. We often denote a pocket by its
boundary surfaces; e.g. an {\em annulus pocket} is a solid torus with
annulus boundary surfaces.

\begin{proof}
Let the map
$
H:Y_1\times[0,1] \to M
$
be as in the definition of parallel pair. 
Proposition 5.4 of Waldhausen \cite{waldhausen} implies that if
$H$ is constant on $\boundary Y_1$, then
the parallel pair bounds a compact region of the desired homeomorphism type.
We will adjust $H$ to obtain a homotopy $H'$ which is constant on
$\boundary Y_1$.

The map $H_1$ takes each
component $\gamma$ of $\boundary Y_1$ to itself. We may assume that
$H_1(x)=H_0(x)=x$ for some point $x\in\gamma$, and let $t$ be
$[H(x\times[0,1])]$ in $\pi_1(M,x)$. 
If $H_1$ were to reverse
orientation on $\gamma$, we would obtain a relation of the
form $tat^{-1} = a^{-1}$ with $a=[\gamma]$,
but this is impossible in the fundamental group of an orientable surface.
Thus we must have $tat^{-1}=a$, and since $a$ is primitive in
$\pi_1(M)$ and $S$ is not a torus, $t=a^m$ for some $m$. 
Hence,
after possibly adjusting $H_1$ by a further twist in the collar of
$\boundary Y_1$, 
we may assume that $H_1$ is the identity on $\gamma$, and
furthermore that $m=0$. Thus 
$H|_{\gamma\times[0,1]}$ is homotopic rel boundary to the map
$(x,s)\mapsto x$. A modification on a collar of
$\boundary Y_1\times[0,1]$ yields a homotopy $H'$ which is
constant on $\boundary Y_1$.
This establishes existence of the pocket. 
Uniqueness follows from the non-compactness of $M$.
\end{proof}

Given two homotopic embedded surfaces with common boundary, one might
hope that the surfaces can be divided into subsurfaces bounding
disjoint pockets. One could then use the pockets to construct a controlled
homotopy which pushed the surfaces off of one another (except at their
common boundary). The following lemma shows that, if the surfaces have
no homotopic boundary components, this is always the case
unless there is one of three specific configurations of disk or annulus
pockets. 

\begin{lemma}{pocket decomposition}
Let $R_1$ and $R_2$ be two homotopic surfaces in $\Emb(M)$
intersecting transversely 
such that $\boundary R_1 = \boundary R_2$.
Suppose also that no two components of $\boundary R_1$ are homotopic
in $M$.
Let $C=R_1\intersect R_2$.
Then there exists a nonempty collection $\XX$ of pockets, such
that each $X\in\XX$ has boundary surfaces
$Y_1\subset R_1$ and $Y_2 \subset R_2$, so that 
$Y_1$ is the closure of a component of $R_1
\setminus C$.
Furthermore at least one of the following holds:
\begin{enumerate}
\item \label{disk}
$\XX$ contains a disk pocket.
\item \label{two annulus}
$\XX$ contains a pair of annulus pockets $X$ and $X'$ in the same, nontrivial,
homotopy class, and their interiors are disjoint from each other, and
from $R_1$ and $R_2$. Furthermore, $X$ and $X'$ are on opposite sides of $R_1$,
as determined by its transverse orientation in $M$. 
\item \label{boundary annulus}
$\XX$ contains an annulus pocket in the homotopy class of a component
of $\boundary R_1$. 
\item \label{pocket decomp}
$\XX$ is a decomposition into pockets: every component of
$R_1\setminus C$ is parallel to some component of  $R_2\setminus C$, 
and the interiors of the resulting pockets are disjoint.
\end{enumerate}
\end{lemma}

\begin{proof}

First, if the intersection locus $C$ has a component that is
homotopically trivial, take such a component $\gamma$ which is
innermost on $R_1$. Thus
$\gamma$
bounds a disk component $Y_1$ of $R_1\setminus C$. On $R_2$, $\gamma$
must also bound a disk $Y_2$, although $int(Y_2)$ may contain
components of $C$. These two disks
must bound a disk pocket since $M$ is irreducible, so we have case
(\ref{disk}). 

From now on we will assume all components of $C$ are homotopically
nontrivial.

Let $H:R_1\times[0,1] \to M$ be a homotopy from the inclusion map $H_0$ of
$R_1$ to a homeomorphism $H_1:R_1\to R_2$.
We note that, since no two boundary components of $R_1$ are homotopic
in $M$,  and $R_1$ is homotopic to an essential subsurface of $S$,
the homotopy class in $M$ of any curve in $R_1$ determines its
homotopy class in $R_1$ (and similarly for $R_2$). Thus, 
for any component $\beta$ of $C$, 
$\beta$ and $H_1(\beta)$
are homotopic in $R_2$. 

\subsubsection*{Parallel internal annuli}
Suppose that a homotopy class $[\beta]$ in $R_1$ that is not
peripheral contains at least 3 elements of $C$. Then we claim 
that annulus pockets as in conclusion (\ref{two annulus}) exist.

The union of all annuli in $R_1$ and $R_2$ bounded by curves
in $[\beta]$ forms a 2-complex in $M$. There is a regular neighborhood of
this complex which is a submanifold $K$ of $M$ all of whose
boundaries are tori and which $R_i$ intersects in a properly embedded
annulus $A_i$ for $i=1,2$. Each component $T$ of $\boundary K$, being a
compressible but not homotopically trivial torus, bounds a unique
solid torus $V_T$ in $M$. If $T$ is the boundary component containing
$\boundary A_i$ then $V_T$ must contain $K$, for otherwise it would
contain $R_i\setminus A_i$, which is
impossible (see Figure \ref{A1A2}).

\realfig{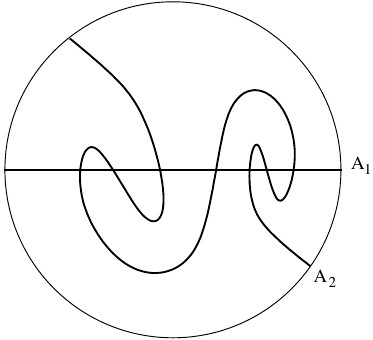}{The solid torus $V_T$ is obtained by crossing this
picture with the circle.}

$A_1$ cuts $V_T$ into two solid tori; Let $V_T^+$ be one of them. 
$A_2$ meets $V_T^+$ in a union of annuli, and
since there are at least 3 intersection curves of $A_1$ and $A_2$, 
these annuli have at least 3 boundary components on $A_1$.
$A_2$ meets $\boundary V_T$ in only 2 circles, so 
there is at least one annulus of $A_2\intersect V_T^+$
whose boundary components are
both in $A_1$. An innermost such annulus in $V_T^+$
yields the desired pocket in
$V_T^+$. Repeating for the other component of $V_T\setminus A_1$, we have
established conclusion (\ref{two annulus}).

\subsubsection*{Peripheral Annuli}
Now suppose that a peripheral homotopy class $[\beta]$ in $R_1$ 
contains at least 2 elements of $C$. There is then an annulus $Y_1$ in
$R_1\setminus C$ whose boundary contains a boundary component of
$R_1$. Again by our assumption that no two boundary components of
$R_1$ are homotopic in $M$, 
the two boundary components of $Y_1$ must also be homotopic in
$R_2$. Thus they  bound an annulus $Y_2$ in $R_2$. The two annuli bound an
annulus pocket by Lemma \ref{pockets}, and this gives 
conclusion (\ref{boundary annulus}). (Note that $Y_2$ is allowed to
have interior intersections with $R_1$).

\subsubsection*{Pocket decomposition}
From now on, we will assume that each nonperipheral homotopy class in $R_1$
contains at most two elements of $C$, and each peripheral homotopy
class in $R_1$ contains exactly one element.

The curves of $C$ define partitions of $R_1$ and $R_2$ whose
components are in one-to-one correspondence by homotopy class. 
In particular, 
if $Y_1$ is the closure of a component of $R_1\setminus C$, its
boundary $\boundary Y_1$ must bound a unique surface $Y_2$ in $R_2$
which is homotopic to $Y_1$.
However we note that $int(Y_2)$ need {\em not} a priori be a
component of $R_2\setminus C$, since two elements of $C$ may be homotopic.
This is the main technical issue we must deal with now. 

Now suppose that no non-annular component of $R_1\setminus C$ has homotopic
boundary components. If $int(Y_1)$ is an annular component of $R_1\setminus C$,
then there is clearly a homotopy from $Y_1$ to $Y_2$ which takes
each boundary to itself. If $Y_1$ is a non-annular component, the
existence of such a homotopy follows from the fact that $Y_1$
has no homotopic boundary components. In either case, $Y_1$ and $Y_2$ form a
parallel pair, and by Lemma \ref{pockets} they bound a pocket
$X_{Y_1}$ (note that $int(Y_1)$ and $int(Y_2)$ are disjoint by choice
of $Y_1$).

Let $\gamma$ be a component of $\boundary Y_1=\boundary Y_2$ which is
nonperipheral in $R_1$. There
are two possible local configurations for $X_{Y_1}$ in a small regular
neighborhood of $\gamma$, shown in Figure \ref{corner}. $X_{Y_1}$ meets the
neighborhood in a solid torus whose boundary contains annuli of $Y_1$
and $Y_2$ adjacent to $\gamma$. $R_1\setminus Y_1$ and
$R_2\setminus Y_2$ meet the neighborhood in two annuli which are
either both outside of $X_{Y_1}$ (case (a)) or inside of $X_{Y_1}$ (case (b)). 
We will rule out case (b). 

\realfig{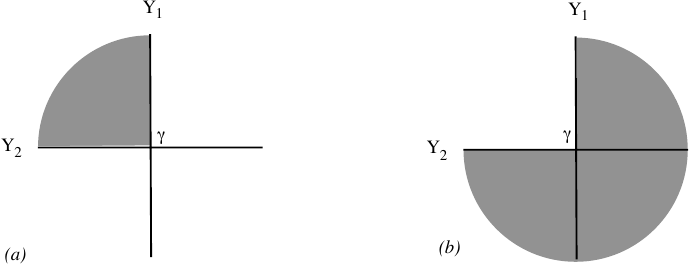}{Local configurations for the corner of a pocket
  $X_{Y_1}$ near a boundary curve $\gamma$. $X_{Y_1}$ is shaded.}

In case (b), there is a component $W$ of $R_2\setminus C$ contained
in $int(X_{Y_1})$, since $W$ cannot intersect $\partial X_{Y_1}=Y_1\cup Y_2$,
whose closure contains  $\gamma$. 
Thus $W$ is homotopic into
$Y_1$, but this can only be if $W$ is an annulus adjacent to
$Y_2$. The second boundary component of $W$ can not be
contained in $Y_2$ or $Y_1$ by definition, nor can it be in
$\boundary Y_1$ since then $Y_1$ and $Y_2$ would also be annuli (and
there is one annulus in $R_2-C$ in each homotopy class). Thus we have a
contradiction and case(b) cannot hold.

Now we can show that, in fact, $int(Y_2)$ cannot meet $C$. For if a
component $\beta$ of $C$ were contained in $int(Y_2)$, there would be
a component $W'$ of $R_1\setminus C$ contained in $int(X_{Y_1})$ with
one boundary component in $Y_2$. As in the previous paragraph
$W'$ has to be an annulus, in the homotopy class of some $\gamma$ in
$\boundary Y_1$. Since there are at most 2 components of $C$ in a
homotopy class, the second boundary component of $W'$ must be $\gamma$
itself, but this is impossible since we have ruled out case (b) of
Figure \ref{corner}. 

From the above, we see that if $int(Y_1)$ and
$int(Y_1')$ are any two components of $R_1\setminus C$, then
the associated pockets $X_{Y_1}$ and $X_{Y'_1}$  are
disjoint except  possibly along common boundary curves of $Y_1$ and
$Y'_1$; any other intersection would lead to the case (b)
configuration or to components of $C$ in $int(Y_2)$.
Thus we obtain 
the desired pocket decomposition of conclusion (\ref{pocket decomp}).

In general, some non-annular components of $R_1\setminus C$
 may have homotopic boundaries. There
is a double cover of $S$ such that each non-annular component
of $R_1\setminus C$ lifts to two
homeomorphic components that do not have homotopic boundaries, and
the surface in the isotopy class of 
$R_1$ has a connected lift. We can repeat the previous arguments in
the corresponding double cover $\til M$ of $M$, and obtain a pocket
decomposition $\til \XX$ there. Every pocket of $\til \XX$ must
embed  under the double cover, 
since there can be no new intersection curves. For the same
reason, pockets downstairs do not intersect except in the expected way
along curves of $C$. Therefore $\XX$ is a pocket decomposition.
\end{proof}

\subsection{Proof of Theorem \ref{Scaffold Isotopy}}
\label{isotopy proof}

We are now prepared to give the proof of Theorem \ref{Scaffold Isotopy}, 
which we restate for the reader's convenience.

\medskip\noindent
{\bf Theorem \ref{Scaffold Isotopy}} (Scaffold Isotopy)
{\em Let $\Sigma$ be a straight scaffold, and let $f:\Sigma \to \hhat M$ be a
good scaffold map. There exists an isotopy
$H:\hhat M\times[0,1] \to \hhat M$
such that $H_0=id$, $H_1\circ f(\Sigma) = \Sigma$, and
$H_1\circ f$ is the identity on $\FF_\Sigma$.}

\medskip

\subsubsection*{Outline} We reduce to the case that
$f(\VV_\Sigma) = \VV_\Sigma$, and then consider the 
intersections of $f(\FF_\Sigma)$
and $\FF_\Sigma$. We wish to progressively simplify these
intersections using the information
provided by Lemma \ref{pocket decomposition}. That lemma describes
embedded surfaces with common boundary, whereas the surfaces of
$\FF_\Sigma$ and $f(\FF_\Sigma)$ have disjoint boundaries that meet solid tori, and are
wrapped around these solid tori in possibly complicated ways. 
Hence as a first step we extend the surfaces into
the solid tori so as to obtain surfaces whose boundaries meet at
the cores. 

We can then use Lemma \ref{pocket decomposition} to analyze the
pockets that arise and simplify them. The main new complication to
keep in mind here is the structure of the intersections of the
pockets with the solid tori.

\begin{pf}
Let $\FF=\FF_\Sigma$ and $\VV=\VV_\Sigma$.
We first reduce to the case that $f(\VV) = \VV$.

By assumption, $\Sigma' = f(\Sigma)$ is a scaffold so there is
a map $\Phi:\hhat M\to \hhat M$ isotopic to the identity, such that 
$\Phi(\VV_{\Sigma'})$ is a union of tori of the form
$\bcollar(v)\times J_v$ with $J_v$ an interval.
Since $\Sigma$  is straight and $f$ is homotopic to the identity,
$\VV$ is also a union of tori of the form $\bcollar(v)\times I_v$, where the set of homotopy classes $v$ is the same. It remains to construct a further isotopy which takes
$\bcollar(v)\times J_v$ to $\bcollar(v)\times I_v$ for each $v$.

Let $h_v:J_v \to I_v$ be an affine orientation preserving
homeomorphism, and let $g_{v,t}(x) = (1-t)x + th_v(x)$. Thus $g_{v,t}$
($t\in[0,1]$)
is an ``affine slide'' of  $J_v$ to $I_v$.  This allows us to define a
family of maps $G_t$ on the tubes of $\Phi(\VV_{\Sigma'})$, so that for each
homotopy class $v$ the restriction of
$G_t$ to $\bcollar(v)\times J_v$ is $G_t(p,x) =
(p,g_{v,t}(x))$. This slides $\bcollar(v)\times J_v$ to 
$\bcollar(v)\times I_v$. 

If $v$ and $w$ are disjoint, so are their collars and the
$G_t$-images of the 
corresponding tubes do not collide.

Whenever $v$ and $w$ overlap, $J_v\intersect J_w = \emptyset$, 
and $I_v\intersect I_w = \emptyset$. 
Supposing
$\max J_v < \min J_w$, the order-preserving property of $f$ implies that
$\max I_v < \min I_w$. Thus it follows that $g_{v,t}(J_v)$
and $g_{w,t}(J_w)$ are disjoint for any $t$, and again the images of
the corresponding tubes do not collide. 

By the isotopy extension theorem \cite{rourke-sanderson},
this isotopy of $\Phi(\VV_{\Sigma'})$ can be extended to an
isotopy $\Psi_t$ of $\hhat M$. 

Thus, after replacing $f$ with $\Psi_1 \circ \Phi \circ f$, 
we may from now on assume that $\VV_{\Sigma'} = \VV$.
We will build an isotopy
of maps of pairs 
$$
H:(\hhat M\times[0,1],\VV\times[0,1]) \to (\hhat M,\VV)
$$
such that $H_0=id$ and $H_1\circ f$ is the identity on $\FF$.
This will be done by induction on the
pieces of $\FF$.

In the inductive step, we may assume that on some union of components
$\EE\subset \FF$, $f$ is already equal to the identity. We let  $R$ be a component
of $\FF\setminus\EE$,
and build an isotopy of pairs $(\hhat M,\VV)$ which fixes
pointwise a neighborhood of $\EE$, and moves $f|_R$ to the identity.

Our first step is to apply an isotopy so that $\boundary R$ and
$f(\boundary R)$ are disjoint.
Let $\gamma$ be a boundary component of $R$, lying in the boundary of
a solid torus $V$ in $\VV$, and let $\gamma'$ be a component of
$f(\boundary R)$ lying in $\boundary V$.

By Lemma \ref{level curves}, since both $\Sigma$ and $f(\Sigma)$ are
scaffolds, both $\gamma$ and $\gamma'$ are in the homotopy class of
level curves and hence homotopic to each other in $\boundary V$. We
claim that they are 
isotopic, within $\boundary V \setminus \EE$, to disjoint curves. 
If they are in different components of $\boundary V \setminus \EE$
then they are already disjoint. If they are in the same component $A$,
then $A$ is either an annulus or a torus in which $\gamma'$ and
$\gamma$ are homotopic, so 
$\gamma'$ is isotopic within $A$ to a curve disjoint from $\gamma$.

This isotopy may be extended to a small neighborhood of $\boundary V$
to have support in the complement of $\EE$, so after applying this
isotopy to $f$ we may
assume that $\boundary R$ and $f(\boundary R)$ are disjoint.

Now, in order to apply Lemma \ref{pocket decomposition}, let us
enlarge $R$ and $f(R)$ to surfaces $R_1$ and $R_2$, as follows.
If $V$
is a  component of $\VV$ meeting just one boundary component $\gamma$
of $R$, let $A_1$ and $A_2$ be embedded annuli in $V$ joining $\gamma$ 
and $f(\gamma)$, respectively, to 
a fixed core curve of $V$, and let $A_1$ and $A_2$ be disjoint except
at the core. 
If $V$ meets two components of $\boundary R$, join them with an
embedded annulus $A_1$, and similarly join the corresponding
pair of components of $f(\boundary R)$ with an embedded $A_2$, so that
$A_1$ and $A_2$ intersect transversely and minimally -- either not at
all, or transversely in one core curve. Figure \ref{addannuli}
illustrates these possibilities. Repeating for each $V$, 
let $R_1$ be the union of $R$ with
all the annuli $A_1$, and let 
$R_2$ be the union of $f(R)$ with
all the annuli $A_2$.

Since $R$ is a level surface, we may choose the annuli $A_1$ to be
level, so that $R_1$ is still a level surface. 

\realfig{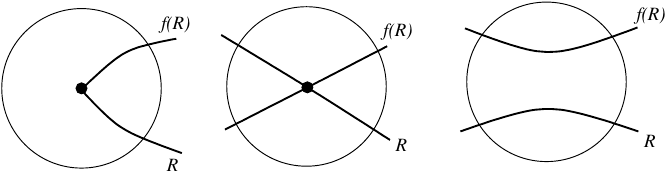}{The three ways in which annuli are added to $R$
and $f(R)$. To obtain the true picture, cross each diagram with $S^1$.}

Because we have joined homotopic pairs of boundary components, $R_1$
and $R_2$ satisfy the conditions of Lemma \ref{pocket
decomposition}. Let $\XX$ be the collection of pockets described in
the lemma. For each of the possible cases we will describe how to
simplify the picture by an isotopy of $(M,\VV)$. The general move,
given a pocket $X$ bounded by $Y_1\subset R_1$ and 
$Y_2\subset R_2$, is to apply an isotopy in a neighborhood of $X$
which pushes $Y_2$ off $R_1$ using the product structure of the
pocket. However we have to be careful to deal correctly with the
possible intersections of $X$ with $\EE$ and $\VV$. 
In particular we will maintain inductively the property that
the intersection of $R_1\union R_2$ with each solid torus is always 
one of the configurations in Figure \ref{addannuli}.

\medskip

{\bf \ref{disk}:}  Suppose $\XX$ contains a disk pocket $X$. Since $X$ is a
ball, no component of $\VV$ or $\EE$ can be contained in it. 
$\EE$ is disjoint from $R_1$ and $R_2$ and hence cannot intersect $X$
at all. Since according to Figure \ref{addannuli} any intersection
curve of $V$ with $R_1\union R_2$ is homotopically nontrivial, 
$\VV$ cannot intersect $\boundary X$, and hence $X$, either. Hence
$R_2$ can be  
isotoped through $X$, reducing the number of intersections, and the
isotopy is the identity on $\VV$ and $\EE$. After a finite number of
such moves we may assume there are no disk pockets. 

{\bf \ref{two annulus}:} Suppose $\XX$ contains two homotopic non-peripheral
annulus pockets $X$ and $X'$, with interiors disjoint from each other and
from $R_1$ and $R_2$, and on opposite sides of $R_1$.

If one of the annulus pockets misses $\VV$ then it cannot
meet any component of $\EE$, and hence  we may
isotope across it to reduce the number of intersection curves. 

If both of the pockets meet $\VV$, it must be in a single solid torus $V$ in
the same homotopy class. By the inductive hypothesis, the
intersections with $V$ must be as in Figure \ref{addannuli}, 
so that $X$ and $X'$ each meet $V$ in a solid
torus. There are two intersection patterns in $V$, depicted in Figure
\ref{twoannuli1}. In case (a), $X$ and $X'$ intersect at the core curve of $V$.
The product structure in $X$ can be
chosen so that $\boundary V\intersect X$ is vertical, and again we can
isotope through $X$, preserving $V$, to reduce the number of
intersection curves. 

\realfig{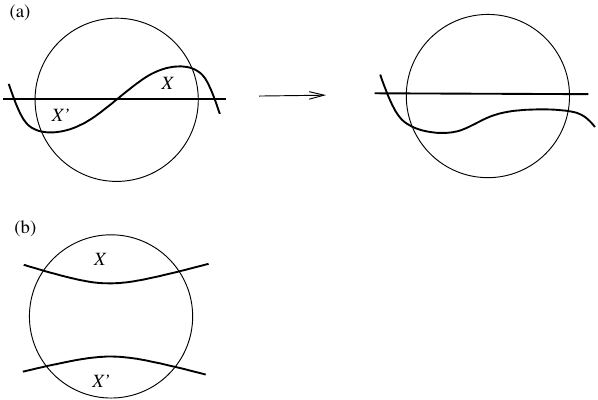}{The two ways that a solid torus $V$ can
intersect two homotopic annulus pockets.}

In case (b), $X$ and $X'$ do not intersect in $V$ and
we have to consider the picture more globally. Since $X$
and $X'$ meet $R_1$ on opposite sides, the possible configurations
are
as in Figure \ref{twoannuli2} (up to orientation).
In each case, $X$ meets $V$ in a region
bounded by a subset of $R_2$, and $X'$ meets $V$ in a region bounded by
a subset of $R_1$. 
Outside $V$ there is a solid torus
$Z$ bounded by three annuli, one in $\boundary V$, one in $\boundary X$
and one in $R_1$, and $int(Z)$ is disjoint from $R_1$. 
In the top
case, $int(Z)$ is disjoint from $R_2$ and we may push the annulus
$\boundary Z \intersect \boundary X$ across the rest of $Z$, thus
reducing the number of intersection curves outside $\VV$ (although we
introduce an intersection in $V$, as shown).

In the bottom case, we can find an innermost annulus of $R_2
\intersect Z$, and push across the resulting pocket.
Note that in these cases the isotopy can be done in the complement of
$\EE$ since no component of $\EE$ can be contained in a solid torus. 

After a finite number of such moves we may assume that $\XX$ does not
contain a pair of homotopic non-peripheral annulus pockets.

\realfig{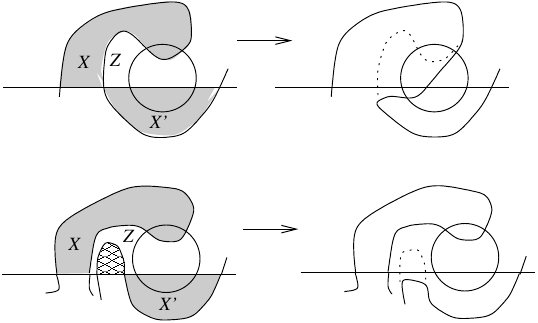}{The two possible moves in case (b)}

{\bf \ref{boundary annulus}:} If $\XX$ contains a peripheral annulus
pocket $X$, let $V$ denote the component of $\VV$ in the same
homotopy class. The intersection of $\boundary V$ with $R_1$ and
$R_2$ must be as in the first picture in Figure \ref{addannuli}.
We can adjust the product
structure of the pocket so that $\boundary V$ is vertical, and
hence we can isotope $Y_2$ off of $Y_1$ while preserving $V$ (see
figure \ref{bdryint}). Again, this can be done in the complement of $\EE$. 

\realfig{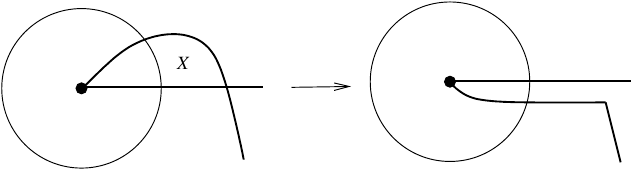}{A peripheral annulus pocket can only intersect $V$
as shown (multiplied by $S^1$). The isotopy move is shown as well.}

{\bf \ref{pocket decomp}:}
If $\XX$ is a pocket decomposition, we first claim that no pocket
contains all of a component of $\VV$ or $\EE$. Suppose that $X\in\XX$
does contain such a component $Z$.  We will obtain a
contradiction to the order-preserving properties of the map $f$.

$Z$ is homotopic into $R_1$, and we claim it 
must also overlap $R$ -- the alternative is that $Z$ is
homotopic to one of the
annuli in $R_1\setminus R$ -- but then it would have to be one of the
solid tori that intersect $R_1$, contradicting the fact that $Z$ lies
within $X$. 

Recall that $R_1$ is a level surface, let $Y_1$ be the 
subsurface of $R_1$ in the boundary of $X$, and assume without loss of
generality that $X$ is adjacent to $R_1$ from below. 
Because $R_1$ is a level surface $Z$ can be pushed to $-\infty$ in the
complement of $R_1$. Since $Z$ overlaps $R$ and $\Sigma$ is
straight,  they are ordered
and so $Z\topprec R$. Since $f$ is a good scaffold map 
we have $f(Z) \topprec f(R)$, but $f(Z)=Z$ so $Z\topprec f(R)$. 

There is an isotopy of $M$ supported on a small neighborhood of the
pockets {\em different from $X$}, which pushes all of them 
outside of the region above
$Y_1$, which we can call $Y_1\times(t,\infty)$. Let $\Psi$ be the end
result of this isotopy.
We can push $Z$ through the product
structure of $X$ to just above $Y_1$, and then to 
$+\infty$ through 
the region $\Psi^{-1}(Y_1\times(t,\infty))$, avoiding $R_2$ and in
particular $f(R)$ (figure \ref{noZ}).	
This contradicts $Z\topprec f(R)$.

\realfig{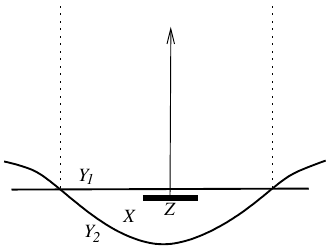}{If $Z$ is contained in a pocket, it contradicts order
preservation.}

The only remaining issue is that a pocket $X$ might intersect, but
not contain, one
of the solid tori $V$. A priori there are six possible intersection
patterns of $V$ with the pockets, as in figure \ref{XUintersect}.

\realfig{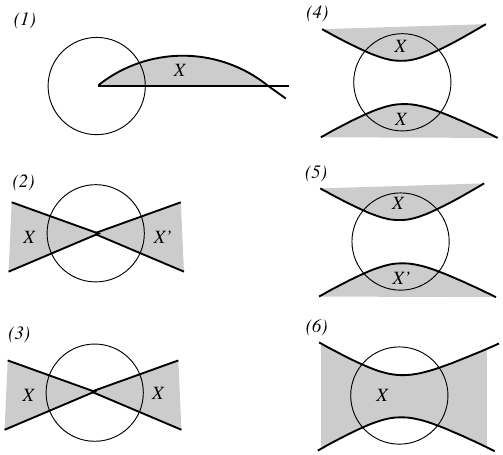}{The six intersection patterns of a pocket $X$ with
a tube $V$}

In case (1), $V$ corresponds to a peripheral homotopy class in $R_1$
and $X$ meets $V$ in a solid torus with the core of $V$ at its
boundary. The product structure of $X$ can be adjusted so that 
the annulus $\boundary V\intersect X$ is vertical.
(This is similar to the peripheral annulus pocket case.)

In case (2), $V$ meets two pockets $X$ and $X'$ and there is an
intersection curve in the core of $V$. Again $\boundary V \intersect
X$ can be made vertical.

In case (3), the local pattern is the same as in case (2) but we
consider the possibility that both intersections are part of $X$. This
case cannot occur: The orientation of $X$
induced from $M$ induces an orientation on each of the two boundary
surfaces $Y_1$ and $Y_2$. However in the local picture in $V$, each $Y_i$
inherits inconsistent orientation from the two sides, since $Y_1$ and $Y_2$
intersect transversely. 

In case (4), $X$ meets $V$ in two disjoint solid tori. This case can
also be ruled out:  Consider an annulus $A$ in $V$ with homotopically
nontrivial boundaries in
the two components of $X\intersect V$. Because $X$ is a pocket, and 
the map of $S$ to $M$ is a homotopy-equivalence, 
the boundaries of $A$ can also be joined by an annulus $A'$ in
$X$. This produces a torus that intersects the level surface
$R_1$ in exactly one, essential, curve $\gamma$.
The curve is in the homotopy
class of $V$, which in this case is non-peripheral in $R_1$. 
But this is impossible: Let
$S_1$ be the full level surface containing $R_1$. Any intersection of
the torus with $S_1 \setminus R_1$ cannot be essential because then it
would be homotopic to $\gamma$, which is nonperipheral in $R_1$.  The
nonessential intersections can be removed by surgeries, yielding a
torus that cuts through $S_1$ in just one curve -- but $S_1$ separates
$M$, a contradiction.

In case (5), $V$ meets $X$ and $X'$ in disjoint solid tori. Hence $V$
must be homotopic to a boundary component of $X$ and of $X'$.  This
gives us a configuration which agrees, in a neighborhood of $V$, 
with the top picture in figure \ref{twoannuli2} -- the only difference is 
that one of the pockets is not an annulus pocket now.
(Notice that the bottom picture in figure \ref{twoannuli2} cannot occur,
since we have assumed that $C= R_1\intersect R_2$ never contains more than 2 homotopic curves.)
The same isotopy move as in figure \ref{twoannuli2} simplifies the
situation by reducing the
number of intersection curves outside of $\VV$ (and changing the local
configuration at $V$ to case (2)).

In  case (6), $V$ intersects $\boundary X$ in two disjoint annuli. The
intersection $V\intersect X$ is a solid torus $Z$, and the 
product structure of $X$  can be adjusted so that the
two annuli of $\boundary Z \intersect \boundary V$ are vertical. (This
is essentially because  in a product $R\times [0,1]$ there is only one
isotopy class of embedded annulus for each isotopy class in $R$.)

In conclusion: cases (3) and (4) do not occur. Case (5) can be
removed locally, yielding a simpler situation. 
Hence we can assume that all pockets only intersect $\VV$
in the patterns of cases (1), (2) and (6).
The product structure of each pocket $X$ can then be adjusted so that all the
annuli of the form $\boundary V \intersect X$ are simultaneously
vertical, and then a single push of $R_2$ through each pocket
yields an isotopy of $(M,\VV)$ that takes $R_2$ to $R_1$, and $f(R)$
to $R$, 
and fixes a neighborhood of $\EE$. A  final isotopy within $R$ takes
$f$ to the identity.  This completes the proof 
of Theorem  \ref{Scaffold Isotopy}.
\end{pf}

\subsection{Wrapping coefficients}
\label{wrapping coeff}

In this section, we develop a criterion to guarantee that a surface
is embedded which will be used in the proof of
Theorem \ref{Relative embeddability}.
Roughly, we show that if a
surface is embedded off of an annulus immersed in the boundary of
a tube, and is homotopic to $+\infty$ or $-\infty$ in the complement of
that tube, then the surface is embedded.

\begin{lemma}{unwrapped is embedded}
Let $R\subseteq S$ 
be an essential subsurface and $V$ a solid torus in $\Emb(M)$
homotopic to a nonperipheral curve $\gamma$ in $R$.
Let $R' = R\setminus \collar(\gamma)$.
Suppose that $h\in \Map(R,M)$ is a map such that $h|_{R'}$ is 
an embedding into $M\setminus int(V)$ with $h^{-1}(\boundary V) =
\bcollar(\gamma)$, $h(\boundary R')$ is unknotted and
unlinked,  
and $h|_{\collar(\gamma)}$ is an immersion of $\collar(\gamma)$ in
$\boundary V$. 

If $h$ is homotopic to
either $-\infty $ or to $+\infty $ in the complement 
of $int(V)$ then $h|_{\collar(\gamma)}$ (and hence $h$ itself) is an embedding.
\end{lemma}

We first define an algebraic measure of wrapping.
Suppose that $V$ is a straight solid torus in $M$, $R\subseteq S$
is an
essential subsurface, and the core of $V$ is homotopic to a
nonperipheral curve in $R$. Consider $f\in\Map(R,M)$ whose image is
disjoint from $int(V)$, and suppose that $f(\boundary R)$ is {\em not linked 
with $V$}, which means that  $f|_{\boundary R}$ is homotopic to
both $-\infty$ and $+\infty$ 
in the complement of $V$. (Note that it is sufficient to assume it is
homotopic to $-\infty$ in the complement of $V$: since $V$ is
nonperipheral in $R$ and straight, $\boundary R\times\R$ 
is disjoint from $V$, so once $f(\boundary R)$ is pushed far enough
below it can be pushed above along $\boundary R\times\R$). 

We then define a ``wrapping
coefficient'' $d_-(f,V)\in \Z$  as 
follows: Let $r\in\R$ be smaller than the minimal level of $V$, 
and let $G:R\times[0,1]\to M$ be a homotopy with $G_0=s_r:R\to
R\times\{r\}$, and $G_1=f$. By the unlinking assumption on
$f(\boundary R)$ we may choose $G$ so that 
$G(\boundary R\times[0,1])$ is disjoint
from $int(V)$. Then the degree of $G$ over $int(V)$ is well-defined, and we
denote it $\deg(G,V)$.

If $G'$ is another such map we claim $\deg(G,V) = \deg(G',V)$. 
Viewing $G$ and $G'$ as 3-chains, the difference $G-G'$ has boundary
equal to $G(\boundary R\times[0,1]) - G'(\boundary R\times[0,1])$, a
union of singular tori. Since $M\setminus V$ is atoroidal and these
tori are not homotopic to $\boundary V$, they must bound a union of
(singular) solid tori $W$ in the complement of $V$, and we can write
$\boundary (G-G'-W) = 0$.  
Since $H_3(M)=\{0\}$ we know $\deg(G-G'-W,V)=0$ and since
the contribution from  $W$ is $0$ we must have $\deg(G,V)=\deg(G',V)$.
Thus we are justified in defining $d_-(f,V) \equiv \deg(G,V)$. 
We also drop $V$ from the notation, writing $d_-(f)$, $deg(G)$ etc, where
convenient. 

It is clear that $d_-(s_r) = 0$ and, more generally, that if $f$ is
deformable to $-\infty$ in the complement of $V$ then
$d_-(f)=0$. Furthermore if $H$ is a homotopy
such that $H(\boundary R\times[0,1])$ is disjoint from $V$, $H_0=f$,
$H_1=g$, and $f(R)$ and $g(R)$ are disjoint from $V$, then 
\begin{equation}\label{degree of difference}
d_-(g)-d_-(f) = \deg(H).
\end{equation}

We can define $d_+(f)$ similarly as $\deg(G)$, where
$G$ is a homotopy such that $G_0=f$
and $G_1 = s_t$, with $t$ larger than the top of $V$. Then $d_+(f)$
must be $0$ for $f$ to be deformable to $+\infty$  in the complement
of $V$.

\begin{proof}[Proof of Lemma \ref{unwrapped is embedded}]
Since $h(\boundary R')$ is unknotted and unlinked we may attach solid tori
to the boundary components of $h(\boundary R)$ so that, together with
$V$ and $h(R')$, we have a scaffold $\Sigma$.
(Notice that since each component of $h(\partial R')\cap \partial V$ is a simple closed
curve which represents an indivisible element of $\pi_1(M)$, it must be isotopic to
the core curve of $V$.) 
Since $\Sigma$ has
no overlapping pieces, it is combinatorially straight.
By Lemma \ref{straighten scaffold} we may assume, after 
re-identifying $M$ with $S\times\R$, that $V$ is a straight solid torus
and $h(R')$ is a level surface. Assume without loss of generality that
$h$ is homotopic to $-\infty$ in the complement of $V$.

Let $A_-$ be the annulus in $\boundary V$
consisting of the part of the boundary below $h(R')$,
and let $h_-:R\to M$ be an extension of $h|_{R'}$
that maps $\collar(\gamma)$ to $A_-$ and is an embedding.

Clearly $h_-$ is deformable to $-\infty$ in the complement of $V$, 
and hence $d_-(h_-) = 0$. Similarly, $d_-(h)=0$.

The difference between $h(\collar(\gamma))$ and
$h_-(\collar(\gamma))$, considered as 2-chains, gives a cycle
$k[\boundary V]$ with $k\in\Z$. 
We immediately have $d_-(h) = d_-(h_-) + k = k$. 
Since $d_-(h)=d_-(h_-)=0$,
we conclude that $k=0$. This implies that, since $h$ and $h_-$
are both immersions on 
$\collar(\gamma)$ which agree on $\boundary\collar(\gamma)$,
they must have the same image. It follows that
$h|_{\collar(\gamma)}$ is an embedding.
\end{proof}

\subsection{Some ordering lemmas}
\label{some ordering lemmas}

We now apply the scaffold machinery to obtain some basic properties of 
the $\topprec$ relation.
All these properties will be used
in the proof of  Lemma \ref{preserve cut order}, which allows
us to choose a cut system whose images under the, suitably
altered, model map are correctly ordered.

The following lemma gives us a
transitivity property for $\topprec$ in some special situations. 

\begin{lemma}{simple transitivity}
Let $R_1$ and $R_2$ be disjoint homotopic surfaces in $\Emb(M)$.
Let $\VV$ be an unlinked unknotted collection of solid tori in $M$
with one component for each homotopy class of component of $\boundary
R_1$, so that 
$\boundary R_1$ and $\boundary R_2$ are embedded in $\boundary \VV$.

\begin{enumerate}
\item $R_1$ and $R_2$ are $\topprec$-ordered. 
\item
Let $Q\in \Emb(M)$ be disjoint from $R_1 \union R_2 \union \VV$, 
so that the domain of $Q$ is contained in the
domain $R$ of $R_1$ and $R_2$.

Suppose that $R_1 \topprec Q$ and $Q\topprec R_2$. Then
$R_1 \topprec R_2$. 
\end{enumerate}
\end{lemma}

\begin{proof}
Let $\Sigma$ denote the scaffold with  $\FF_\Sigma = R_1$ and
$\VV_\Sigma = \VV$. 
$\Sigma$ has no overlapping pieces, so
it is combinatorially
straight. Using the straightening lemma \ref{straighten scaffold},
after an isotopy we may assume that $\Sigma$ is straight. 

Now the Pocket Lemma \ref{pockets} implies that there is a product
region $X$ homeomorphic to $R\times[0,1]$ whose boundary consists of
$R_1$, $R_2$, and annuli in the tubes of $\VV_\Sigma$ associated to
their boundaries. Thus $R_2$ is isotopic to $R_1$ by an isotopy
keeping the boundaries in $\VV$, so that 
we may assume that $X$ is actually equal to $R\times[0,1]$ in
$M=S\times\R$, with $R_1 \union R_2 = R\times\{0,1\}$.
It follows from this that if $R_1=R\times\{0\}$ then $R_1\topprec
R_2$, and if $R_1 = R\times\{1\}$ then $R_2\topprec R_1$. Hence we
have part (1). 

It remains to show that, given the hypothesis on $Q$,  $R_1 =
R\times\{0\}$, from which $R_1\topprec 
R_2$ and hence (2) follows. 

If $Q\subset X$ then we must have $R\times\{0\} \topprec Q$ and
$Q \topprec R\times\{1\}$, so we are done.

Now let us suppose that $Q$ is in the complement of $X$, and obtain a
contradiction. 
The condition $R_1\topprec Q$ implies that there
exists a homotopy $H:Q\times[0,\infty)\to M$ such that
$H(Q\times\{t\})$ goes to $+\infty$ as $t\to+\infty$, and which avoids
$R_1$. We claim that $H$ can be chosen to avoid $R_2$ as well. 
Let $W$ be a neighborhood of $X$ disjoint from $Q$. 
There is a homeomorphism $\phi$ taking $M\setminus R_1$ to $M\setminus
X$, which is the identity outside $W$. Thus $\phi\circ H$ is the
desired homotopy. This contradicts $Q\topprec R_2$, so again we are done.
\end{proof}

The next lemma tells us that for disjoint overlapping non-homotopic
non-annular surfaces  with boundary in an unknotted and unlinked collection
of solid tori, the
$\topprec$-ordering is determined by their boundaries.

\begin{lemma}{boundary ordering}
Let $P,R\in \Emb(M)$ be disjoint overlapping non-homotopic non-annular
surfaces such that
$\boundary R \union 
\boundary P$ is embedded in a collection $\VV$ of unknotted,
unlinked, homotopically distinct solid tori, so that each component of
$\VV$ intersects $\boundary P$ or $\boundary R$, and
$\VV\intersect (int(R)\union int(P)) = \emptyset$.
Suppose that for each component
$\alpha$ of $\boundary R$ that overlaps $P$ we have $\alpha \topprec
P$, and for each component
$\beta$ of $\boundary P$ that overlaps $R$ we have $R \topprec \beta$.
Then $R \topprec P$.
\end{lemma}

\begin{proof}
Assume, possibly renaming $P$ and $R$ and reversing directions,
that the domain of $R$ is not contained in the domain of $P$.

Let $\Sigma$ be the scaffold with $\FF_\Sigma = R$ and $\VV_\Sigma =
\VV$. By hypothesis $\VV$ is isotopic to a union of straight solid
tori, so that $\topprec|_\VV$ is acyclic and satisfies the overlap
condition. Since the hypotheses also 
give us that $R\topprec V$ for each component $V$ of $\VV$ that overlaps
$R$ we conclude that $\topprec|_\Sigma$ is still acyclic and satisfies
the overlap condition. Hence $\Sigma$ is combinatorially straight and
by Lemma \ref{straighten scaffold}, after isotopy we may assume that
$\Sigma$ is straight. 

In particular we may assume that $R = R'\times\{0\}$.
Let $X= R' \times (-\infty,0]$. 

After pushing $\boundary P$ in $\boundary \VV$ (by an isotopy
supported in a small neighborhood of $\VV$ and leaving $\Sigma$ invariant)
we may assume that $\boundary P$ is outside of $X$. 

Now consider $P\intersect \boundary X = P\intersect \boundary R'
\times (-\infty,0]$ (we may assume this intersection is
transverse). All inessential curves bound disks in both $P$ 
and $\boundary X$, and so an isotopy of $P$ will remove them. The
remaining curves are in the homotopy classes of the components of
$\boundary R'$. 
Let $A$ be a component of $P\intersect X$. If $A$ is an annulus
then it and an annulus in $\boundary X$ bound a solid torus (annulus pocket)
in $X$, and
taking an innermost such solid torus we may remove it by an isotopy of $P$
without producing new intersection curves with $\boundary X$. After 
finitely many moves we may assume there are no annular intersections.

Any non-annular $A$ must be a subsurface of $P$, and hence is the
image of a subsurface $A'$ of $S$. On the other hand, $A\subset X$
implies that $A'$ is homotopic into $R'$, and $\boundary A\subset
\boundary X$ implies that $\boundary A'$ is
homotopic into $\boundary R'$. Since $A'$ is not an annulus, the only
way this can happen is if $A'$ is isotopic to $R'$ (see, for example,
Theorem 13.1 in \cite{hempel}.)
But this contradicts the assumption that the domain of $R$ is not
contained in the domain of $P$.  Thus there is no non-annular component $A$. 

We conclude that, after an isotopy that
does not move $R$, $P$ may be assumed to lie outside
of $X$. Thus $P$ is homotopic to  $+\infty$ avoiding
$X$ and hence $R$, and $R$ is homotopic to $-\infty$ (through $X$)
in the complement of $P$. $R$ and $P$ overlap, so we conclude by 
Lemma \ref{basic ordering properties} that $R\topprec P$. 
\end{proof}

The next lemma will allow us to check more easily that tubes
are ordered with respect to non-annular embedded surfaces with
boundary in a collection of straight solid tori.

\begin{lemma}{embedded ordered}
Let $R\in \Emb(M)$ be non-annular, with boundary embedded in a
collection $\VV$ of straight solid tori.
Let $U$ be a straight solid torus disjoint from $R\union \VV$
and overlapping $R$. If $R$ is homotopic to $-\infty$ in the
complement of $U$ then $R\topprec U$, and if $R$ is homotopic to
$+\infty $ in the complement of $U$ then $U\topprec R$.
\end{lemma}

\begin{proof}
Assume without loss of generality that $R$ is homotopic to $-\infty$
in the complement of $U$. 
Let $B=\gamma\times[t,\infty)$
where $\gamma\times\{t\}$ is embedded in $\boundary U$ and
homotopic to the core of $U$. Since $U$ and $\VV$ are straight and
$\VV$ is homotopic to $-\infty $ in the complement of $U$ (since $R$
is), $B$ must be disjoint from $\VV$. We may assume that $B$
intersects $R$ transversely, so that components of $B\intersect R$ are
either homotopically trivial or homotopic to the core of $U$. 

Homotopically trivial components may be removed by isotopy of $B$.
The nontrivial components are signed via the natural orientation of
$B$ and $R$, and the fact that $R$ is homotopic to $-\infty$ in the
complement of $U$ means that the signs sum up to 0. Two components of
opposite signs that are adjacent on $B$ can be removed by an isotopy
of $B$, and we conclude that $B$ can be isotoped away from $R$. Thus
$U$ is homotopic to $+\infty$ in the complement of $R$, and we
conclude (invoking Lemma \ref{basic ordering properties}) that
$R\topprec U$.
\end{proof}

\section{Cut systems and partial orders}
\label{po}

\newcommand{\reals}{\mathbb{R}}

\newcommand{\deep}{5}

\newcommand{\core}{\operatorname{core}}

In this section, we link combinatorial information from the hierarchy
$H$ to topological ordering information of {\em split-level surfaces}
in $\modl$.  A {\em split-level surface} is an embedded surface in the
model manifold associated to a slice of $H$ that is made up of level
subsurfaces arising as the upper and lower boundaries of blocks in the
model.  As these split-level surfaces and their images in $M$ will
play an important role in what follows, we now develop some control
over them and their interactions within the model, aiming in
particular for a consistency result (Proposition~\ref{topprec and
cprec}) comparing topological ordering in $\modl$ and a more
combinatorial ordering we define on corresponding slices in $H$.

This consistency result will not apply generally to all slices in
$H$ and their associated domains, but after a thinning procedure we
arrive at a collection of slices called a {\em cut system} whose
split-level surfaces are well behaved with respect to the topological
partial ordering and divide the model into regions of controlled size
as we will see in \S\ref{regions}.

\subsection{Split-level surfaces associated to slices}
\label{split level}
If $a$ is a slice of $H$, we recall from Section~\ref{hierarchy background}
that $g_a$ denotes its bottom geodesic and $v_a$ is the bottom simplex
of $a$.  Then $p_a = (g_a,v_a)$ is the bottom pair of $a$.
Let $$D(a) = D(p_a) = D(g_a)$$ be the domain of
$g_a$.  If $D(a)$ is not an annulus, let $$\check D(a) = D(a)
\setminus \collar(\base(a))$$ be the complement in $D(a)$ of the
standard annular neighborhoods of the curves in $\base(a)$.  When $a$
is a saturated slice, the subsurface $\check D(a) \subset D(a)$ is a
collection of pairwise disjoint 3-holed spheres.  If $D(a)$ is an
annulus, we let $\check D(a) = D(a)$.

Each slice in $H$ gives rise to a properly embedded surface in
$\modl[0]$ called a {\em split-level surface}.  Given a non-annular
slice $a$ of $H$, each 3-holed sphere $Y \subset \check D(a)$ admits a
natural level embedding $F_Y \subset \modl[0]$.  This embedded copy
$F_Y$ of $Y$ lies in the top boundary and the bottom boundary of the
two blocks that are glued along $F_Y$.
The {\em split-level surface} $F_a$ associated to $a$ is obtained by
letting $$F_a = \bigcup_{Y \subset \check D(a)} F_Y.$$

Given a slice $a$ and $v \in \base(a)$, we say $\gamma_v$ is a {\em
hyperbolic base curve for $a$} if there is a solid torus $\storus{v}$
in $\modl[\infty]$ whose closure is compact; otherwise we say
$\gamma_v$ is a parabolic base curve. 
For each $v \in a$ with $\gamma_v$ a hyperbolic base curve, 
we extend the above embedding of $\check D(a)$ across the
annulus $\collar(v)$ to a map of $\collar(v)$ into $\storus{v}$:
the core $\gamma_v$ is sent to the core of the tube $\storus{v}$ with
its model hyperbolic metric, and the pair of annuli 
$\collar(v) \setminus \gamma_v$ are mapped in such a way that radial
lines in $\collar(v)$ map to radial geodesics in the tube $\storus{v}$.
Given  $v \in \base(a)$ for which $\gamma_v$ is a parabolic base
curve, we extend across the corresponding annulus $\collar(\gamma_v)$
to any embedding of $\collar(\gamma_v)$ into $\storus{v}$.

We remark that given a slice $a$ the only base curves that fail to be
hyperbolic correspond to vertices $a \in \base(a)$ for which $v$ is a
vertex in $\base(I(H))$ without a transversal or $v$ is a vertex in
$\base(T(H))$ without a tranversal.

Extending over each annulus $\collar(v)$ for $v \in \base(a)$ in this
way, we obtain an embedding of $D(a)$ into $\modl$ whose image we
denote by $\hhat F_a$.  For each integer $k \in [0,\infty]$ we denote
by $\hhat F_a[k]$ the intersection $$\hhat F_a[k] = \hhat F_a \cap
\modl[k].$$  
We call the surfaces $\hhat F_a[k]$ {\em extended split-level
surfaces}.

When $a$ is an annular slice, there is a vertex $v$ so that $D(a)
=\overline{\collar(v)}$.  We refer to this vertex $v$ as the {\em core
vertex} of 
$a$, and denote it by $v = \core(a)$.  Then we have the associated
solid torus $\storus{v} \subset \modl$.  In the interest of comparing
all slices in $C$ and their associated topological objects in $\modl$,
we adopt the convention that for each integer $k \in [0,\infty]$ and annular
slice $a$ we
have $$F_a[k] = \hhat F_a[k] = \hhat F_a = \storus{v}.$$

\subsection{Resolution sweeps of the model}
\label{sweep}

A resolution $\{\tau_n\}$ of the hierarchy $H_\nu$ yields a ``sweep'' of
the model manifold by split-level surfaces, which is monotonic with
respect to the $\topprec$ relation. More specifically,
in \S8.2 of \cite{minsky:ELCI}, the embedding of $\modl[0]$ in
$S\times\R$ is constructed inductively using an exhaustion of $\modl[0]$ by
submanifolds $M_i^j$ that are unions of blocks.
Each $F_{\tau_j}$ appears as the ``top'' boundary of $M_i^j$, so
that  $int(M_i^j)$ lies below the cross-sectional surface
$\hhat F_{\tau_j}$
 in the product structure of
$S\times\R$. When $\tau_j\to\tau_{j+1}$ is a move associated to 
to a 4-edge $e_j$, there is a block $B_j=B(e_j)$ 
which is appended above $\hhat F_{\tau_j}$, and $F_{\tau_{j+1}}$
is obtained from $F_{\tau_j}$ by replacing the bottom boundary
of $B_j$ by its top boundary. We say that $B_j$ ``is appended at time
$j$'' in the resolution.
(For other types of elementary moves the
surfaces $\hhat F_{\tau_j}$ and $\hhat F_{\tau_{j+1}}$ are the same.)

The following statement about $\topprec$ is an immediate consequence
of this construction, and Lemma \ref{basic ordering properties}.

\begin{lemma}{sweep order}
Fix a resolution $\{\tau_n\}$ of $H$. If $i<j$ and $W\subset \hhat
F_{\tau_i}$ and $W' \subset \hhat F_{\tau_j}$ are essential
subsurfaces which overlap and are disjoint, then $W\topprec W'$. 

Similarly, if $B$ is appended at time $j>i$ and $W_B$ is the middle
surface of $B$ then $\hhat F_{\tau_i} \topprec W_B$. If $B$ is appended at
time $j<i$ then $W_B\topprec\hhat F_{\tau_i}$.
\end{lemma}

\subsection{Cut systems}

Given a collection $C$ of slices of $H$, we let 
$$
C \vert_h \equiv \{\tau\in C:g_\tau = h\}
$$
denote the slices in $C$ with bottom geodesic $h$.
Let 
$$\deep \le d_1 < d_2\le\infty$$ be fixed elements in
${\mathbb N}\union\{\infty\}$.  Then the collection $C$ is a {\em cut system 
satisfying a $(d_1,d_2)$ spacing condition}
if the following hold:

\begin{enumerate}

\item \label{cut distribution} {\em Distribution of bottom pairs: }
For each $h\in H$ with $\xi(D(h))\ge 4$, 
the set $\{v_\tau:\tau\in C
\vert_h\}$ of bottom vertices on $h$ cuts $h$ into
intervals of length at most $d_2$, and, if non-empty, cuts $h$ into at
least three intervals of size at least $d_1$.  Futhermore, no two
slices have the same bottom pair, and no $v_\tau$ is the first or last
simplex of $h$.

\item \label{cuts initial} {\em Initial pairs:}
For every pair $(h,w)\in\tau\in C$ that is not a bottom pair of $\tau$, $w$ is
the first simplex of $h$.

\item \label{cuts saturated} {\em Saturation:}
Each slice $\tau \in C$ with non-annular bottom geodesic is 
a saturated non-annular slice.

\item \label{annulus cuts} {\em Annular cut slices:}  For any annular
geodesic $g$ there is at most one slice $\tau \in C$ with $g_\tau =
g$.

\end{enumerate}

Note that the spacing condition (\ref{cut distribution}) puts no restriction on annular slices
in $C$. An annular slice consists of an annular geodesic and a choice
of vertex, but in fact the vertex plays no role in the rest of the
argument and is only included for notational consistency. 


The following lemma  will allow us to exploit
the standing assumption that $d_1\ge 5$ in the
definition of a cut system.

\begin{lemma}{deep placement}
Let $H$ be a hierarchy, and let $a$ be a slice of $H$ with bottom
pair $p_a = (g_a,v_a)$.  If $v_a$ has distance at least $3$ from
${\bf I}(g_a)$ and ${\bf T}(g_a)$ along $g_a$ then for any pair $p =
(h,v) \in a$
we have $g_a \bsub h$, $h \fsub g_a$, and
$$(g_a , {\bf I}(g_a)) \pprec p \pprec (g_a, {\bf T}(g_a)).$$  

If $a$ is an annular slice, then for any pair $p =
(g_a,v)$ we have
$$(g_a , {\bf I}(g_a)) \ppreceq p \ppreceq (g_a, {\bf T}(g_a)).$$  
\end{lemma}

\begin{pf}
We first assume that $a$ is non-annular.
Given a pair $p = (h,v)$ with $p \in a$, the
footprint $\phi_{g_a}(D(h))$ has diameter at most $2$ and contains
$v_a$.  If $v_a$ lies distance at least $3$ from ${\bf I}(g_a)$ and
${\bf T}(g_a)$, we have 
$$\max
\phi_{g_a}(D(h)) < {\bf T}(g_a) \ \ \ \text{and} \ \ \  {\bf I}(g_a) < \min
\phi_{g_a}(D(h)).$$ 
Therefore, 
by the definition of $\pprec$ in \S\ref{hierarchy background},
$$(g_a , {\bf I}(g_a)) \pprec p \pprec (g_a, {\bf T}(g_a)).$$
Moreover,
it follows that both ${\bf I}(g_a)$ and ${\bf T}(g_a)$ intersect $D(h)$, so
Theorem~\ref{Descent Sequences} implies that
$h \fsub g_a \bsub h$.

If $a$ is annular then $\pprec$-order is just linear order on the 
pairs with bottom geodesic $g_a$ so the second statement follows
immediately.
\end{pf}

We next show that cut systems exist.

\begin{lemma}{cut systems exist}
Given positive integers $d_1 \ge \deep$ and $d_2 \in [3d_1,\infty]$
there is a cut system $C$ satisfying a $(d_1,d_2)$-spacing condition.
\end{lemma}

\begin{proof} Given a geodesic $g \in H$ with non-annular domain, so
that $|g| > d_2$ we may choose a non-empty collection of pairs along
$g$ satisfying condition~(\ref{cut distribution}).  As in
section~\ref{hierarchy background}, for each such pair $(g,u)$ there
is a choice of a saturated non-annular slice $\tau$ of $H$ with $(g,u)$ as its bottom
pair.  The choice of $\tau$ proceeds inductively by at each stage choosing a pair
from a geodesic in $H$ whose domains arise as a component domain of
$(D(k),v)$ where $(k,v)$ is a pair already in $\tau$.  If $h$ is a geodesic in $H$ such that $D(h)$
is a component domain of $(D(k),v)$, then we can obtain a new slice by adding any
pair $(h,w)$. To satisfy condition (\ref{cuts initial}) of the definition we need to have
$w$ be the first simplex of $h$, so we must check that $h$ has a first
simplex. Lemma \ref{deep placement} implies that
$g \bsub h$, so there is a backward sequence $g \bsubd h_m \cdots h_1
= h$. The footprint of $h_m$ contains $u$, so by the 
spacing condition it is not adjacent to the initial marking of $g$. It
follows that $\I(h_m)$ is the restriction to $D(h_m)$ of one of the
simplices of $g$ and so is a simplex (and not an arational
lamination). By induction  each $\I(h_i)$ is a simplex, and in
particular $h$ has an initial simplex.
After filling in every non-annular component domain which arises,
we obtain a non-annular slice satisfying conditions~(\ref{cuts initial}) and (\ref{cuts saturated}).

In general any choice of a collection of slices on annular geodesics 
will satisfy condition~(\ref{annulus cuts}) provided there is at most
one slice for each annular geodesic, so we may make any such choice to
conclude the proof of the lemma.
\end{proof}




%



\subsection{Hierarchy partial order and split level-surfaces}

Given a cut system $C$ we have the topological ordering relation
$\topprec$ on its associated extended split-level surfaces $\hhat F_a$, $a \in
C$.  Because the surfaces $\hhat F_a$ are not themselves level
surfaces, however, it does not immediately follow from the preceding
section that the transitive closure of $\topprec$ on these split-level
surfaces is a partial order: it could conceivably have cycles.  We
devote the remainder of this section to establishing that this is not
the case.

To show there are no cycles, we employ the ordering properties
inherent in the hierarchy $H$ to construct an order relation on the
slices in a cut system whose transitive closure is a partial order.  We
prove that this {\em cut ordering} is consistent with topological
ordering of overlapping associated split level surfaces
(Proposition~\ref{topprec and cprec}) from which it follows directly that
the transitive closure of $\topprec$ on the split-level surfaces
associated to $C$ is a partial order.

\medskip

Let $C$ be a cut system, and  begin with the following relation on
slices. Given slices $a,b\in C$, let 
$$
a \cprec' b
$$
hold whenver there exist pairs $p\in a$ and $p'\in b$ such that $p
\pprec p'$. 

We will then define $\cprec$ to be the transitive closure of
$\cprec'$. 

In order to analyze $\cprec$, we break it down  into various
possibilities. Given slices $a$ and $b$ in $C$, we make the following definitions:
\begin{enumerate}
\item $\displaystyle a \stprec b $
means that for all $p\in a$ and all $p'\in b$ we have
$$p \pprec p'.$$
\item $\displaystyle a \towall b$ 
means that $D(a)\subset D(b)$ and 
for some $p'\in b$ and {\em all} $p\in a$  we have
$$p \pprec p'.$$ 
\item $\displaystyle a \wallto b$
means that $D(a) \supset D(b)$ and for some $p\in a$
and {\em all} $p'\in b$ we have
$$ p \pprec p'.$$
\item $\displaystyle a  \wall b$
means that there exists a third slice $x$ (called a comparison slice) 
such that
$$ a \towall x \wallto b.$$
\end{enumerate}
\noindent {\bf Remark:} These possibilities need not be mutually exclusive.

\medskip

Our main lemma is the following.
\begin{lemma}{hierarchy partial order}
The transitive closure $\cprec$ of $\cprec'$ defines a (strict)
partial order on $C$. Furthermore, 
we have $a\cprec b$  if and only if at least one of the following
holds
\begin{itemize}
\item $a \stprec b$,
\item $a \towall b$, 
\item $a \wallto b$, or 
\item  $a  \wall b$.
\end{itemize}
\end{lemma}

\begin{pf}
We begin by proving a consistency lemma, which ensures that slices
in $C$ comparable via these relations unless $\check D(a)$ and $\check D(b)$ do not
overlap and $g_a$ and $g_b$ are not $\tprec$-ordered

\begin{lemma}{consistency}
For any two distinct slices $a$ and $b$ in a cut system $C$, 
\begin{enumerate}

\item \label{equal} If $D(a) = D(b)$ then $g_a=g_b$, 
we have $a \stprec b$ or $b \stprec a$, and the
$\stprec$-ordering is consistent with the order of bottom simplices
$v_a$ and $v_b$ along $g_a = g_b$.

\item \label{non-nested} If $D(a) \not= D(b)$, $D(a)$ and $D(b)$
overlap, and neither is strictly contained in the other, then
$a\stprec b$ or $b \stprec a$.

More generally  if $g_a \tprec g_b$, even without
intersection of the domains, then $a\stprec b$.

\item \label{nested} If $D(a) \subset D(b)$ and $D(a) \not\subseteq
\overline{\collar}(\base(b))$ then 
we have $a\towall b$ or $b \wallto a$.  
Furthermore if $a\towall b$ or $a \wallto b$ then for no pairs $p\in
a$ and $p'\in b$ do we have $p'\pprec p$.

\item \label{disjoint} If $\check D(a)$ and $\check D(b)$ do not
overlap, and $g_a$ and
$g_b$ are not $\tprec$-ordered, then for all $p\in a$ and $p'\in b$,
$p$ and $p'$ are {\em not} $\pprec$-ordered.
\end{enumerate}
\end{lemma}

A corollary of this lemma is: 

\begin{corollary}{consistency cor}
The relation $a\cprec' b$ holds if and only if one of 
$a \stprec b$, $a \towall b$, or $a \wallto b$ holds. 
\end{corollary}
\begin{pf}[Proof of Corollary]
If $a=b$ then by Lemma \ref{slice unordered} no pairs of $a$ and $b$
can be $\pprec$-ordered, hence none of the relations hold. 

If $a\ne b$ then exactly one of the four cases of Lemma
\ref{consistency} holds. Assuming $a\cprec' b$ there exists $p\in a$
and $p'\in b$ such that $p\pprec p'$, and this rules out case (4). The
first three cases give us $a\stprec b$, $a\wallto b$ or $a\towall b$
depending on the domains $D(a)$ and $D(b)$, and the ordering of $g_a$
and $g_b$. Conversely, $a\stprec b$, $a\towall b$, and $a\wallto b$
each imply $a\cprec' b$ by definition.

\end{pf}

\begin{pf}[Proof of Lemma \ref{consistency}] 
Proof of Part~(\ref{equal}). If  $D(a) = D(b)$ then $g_a = g_b$. 
Moreover, $D(a) = D(b)$ is not an annulus, since $C$ contains
at most a single slice on any annular geodesic and $a$ and $b$ are
assumed distinct.  Since non-annular slices of a cut system satisfy the
$(d_1,d_2)$ spacing condition with $d_1 \ge \deep$, the
bottom simplices $v_a$ and $v_b$ have distance at least $\deep$ on
the geodesic $g_a = g_b$.  Assume that $v_a < v_b$.

Given any pairs $p \in a$ and $p' \in b$, we wish to show that $p \pprec
p'$.  Since $\diam(\hat \phi_{g_a}(p)) \le 2$ and $\diam(\hat
\phi_{g_a}(p')) \le 2$, and since $v_a \in \hat \phi_{g_a}(p)$ and
$v_b \in \hat \phi_{g_a}(p')$, we have $ \max \hat \phi_{g_a}(p) <
\min \hat \phi_{g_a}(p')$.  It follows that $p \pprec p'$, and we
conclude that $a \stprec b$.

\medskip

Proof of Part~(\ref{non-nested}).  
Assume that $D(a)
\not= D(b)$.  If $g_a
\tprec g_b$, then by definition there is a geodesic $m \in H$ so that
$g_a \fsub m \bsub g_b$ and $\max \phi_m(D(g_a)) < \min
\phi_m(D(g_b))$.  In particular, we have $(g_a, {\bf T}(g_a)) \pprec
(g_b, {\bf I}(g_b))$.  Let $p \in a$ and $p' \in b$ be pairs in the
slices $a$ and $b$.  
Applying Lemma~\ref{deep placement} we have
$$p \ppreceq (g_a, {\bf T}(g_a)) \pprec (g_b, {\bf I}(g_b)) \ppreceq
p'.$$ By transitivity of $\pprec$ we conclude that $p \pprec p'$.

If $D(a)\intersect D(b)\ne \emptyset$ and neither domain is strictly contained
in the other, then $g_a$ and $g_b$ are $\tprec$-ordered by
Lemma~\ref{time order for domains}.
It follows that either $a \stprec b$ or $b \stprec a$, and if $g_a
\tprec g_b$ then $a \stprec b$.

\medskip

Proof of Part~(\ref{nested}). Suppose $D(a)\subset D(b)$; note that in
particular this guarantees that $b$ is non-annular.  Let $(h,v)$ be a
pair in $b$ with $D(a)\subset D(h)$ (the bottom pair $p_b$ has this
property). Then either $v$ intersects $D(a)$, or $D(a)$ is contained
in one of the component domains of $(D(h),v)$.  Since $b$ is a
saturated non-annular slice, either this component domain supports a
pair $(h',v')\in b$ or the component domain is an annulus in
$\collar(\base(b))$ and we have $$D(a) \subseteq \overline{\collar}(\base(b)).$$

Thus, provided $D(a) \not\subseteq \overline{\collar}(\base(b))$, we may
begin with $p_b$ and proceed inductively to arrive at a unique
$(h,v)\in b$ such that $D(a)\subseteq D(h)$ and
$v\notin\phi_h(D(a))$.  Since $\phi_h(D(a))$ is non-empty by
Lemma~\ref{nonempty footprint} we may assume without loss of
generality that $\max\phi_h(D(a))< v$ which guarantees that
$(g_a,\T(g_a)) \pprec (h,v)$.

Applying Lemma~\ref{deep placement}, for any pair $p\in a$ we have $p
\ppreceq (g_a,\T(g_a))$ so we may conclude that $$ p\pprec (h,v) $$ for
all $p\in a$. Now if there were some $p'\in b$ and $p\in a$ such that
$p'\pprec p$, then $p'\pprec (h,v)$, contradicting Lemma~\ref{slice
unordered}, which guarantees the
non-orderability of pairs in a single slice.  This proves the
second paragraph of Part~(\ref{nested}).

\medskip

Proof of Part~(\ref{disjoint}).  Assume $g_a$ and $g_b$ are not
$\tprec$-ordered, and that $\check D(a)$ and $\check D(b)$ do not overlap.
This implies that either
\begin{enumerate}
\item $D(a)$ and $D(b)$ are disjoint domains, 
\item $D(a)$ is an annulus with 
$D(a) \subseteq \overline{\collar}(\base(b))$, or 
\item $D(b)$ is an annulus with
$D(b) \subseteq \overline{\collar}(\base(a))$.
\end{enumerate}

Suppose that $p\pprec p'$ for some $p = (h,v) \in a$
and  $p' = (h',v') \in b$. 
Then there is a comparison geodesic $m$,
with $h \fsubeq m \bsubeq h'$ and $\max \hat \phi_m(p) < \min
\hat \phi_m(p')$.
Assume first that $D(a)$ and $D(b)$ do not overlap.  We note that
$D(m)$ contains both $D(h)$ and $D(h')$, while $D(h) \subseteq D(a)$
and $D(h')
\subseteq D(b)$.  Thus, disjointness of $D(a)$ and $D(b)$
implies that $D(m)$ is not contained in either $D(a)$ or $D(b)$.  The
proof of Lemma~\ref{deep placement} guarantees $h
\fsubeq g_a$ and $g_b \bsubeq h'$.  It follows that $g_a$ and $m$ 
lie in the forward sequence $\Sigma^+(D(h))$ 
while $g_b$ and $m$ 
lie in the backward sequence $\Sigma^-(D(h'))$ (see
Theorem~\ref{Descent Sequences}).
Since the domains in $\Sigma^+(D(h))$ are nested and likewise for $\Sigma^-(D(h'))$, 
we conclude that $D(m)$ is either equal to, contained
in, or contains each of $D(a)$ and $D(b)$. 

It follows that $D(a)
\subset D(m)$ and $D(b) \subset D(m)$.  Since $ \phi_m(D(a)) \subseteq
\phi_m(D(h))$ and $\phi_m(D(b)) \subseteq \phi_m(D(h'))$, then, it
follows that $$\max \phi_m(D(a)) < \min \phi_m(D(b))$$ and hence 
$g_a \tprec g_b$,  contradicting our assumption.

Without loss of generality, the final case is that $D(a)$ is an
annulus with $D(a) \subseteq \overline{\collar}(\base(b))$.  Then
$D(a)$ is an annulus component domain of a pair $(h,v) \in b$. Since
$b$ is a saturated non-annular slice, $D(a)$ supports no pair in $b$.
Thus $a$ can be added to $b$ to form a slice $\tau$ in the hierarchy
$H$.  The non-orderability of pairs in a slice (Lemma~\ref{slice
unordered}) implies that for each $p$ in $a$ (here $p = p_a$) and each
pair $p' \in b$, we have $p$ and $p'$ are not $\pprec$ ordered
contrary to our assumption.  This completes the proof of
Part~(\ref{disjoint}).

\end{pf}

With Lemma \ref{consistency} and its corollary in hand, 
we can reduce expressions of length 3  in the relations $\cprec'$ and $\wall$ 
to length 2: 

\subsubsection*{Claim 1:}  $a\cprec' b \cprec' c$ implies  $a\cprec' c$
or $a \wall c$, where the former occurs unless $a\towall b \wallto c$. 

Applying Corollary \ref{consistency cor}, we have $a\wallto b$,
$a\stprec b$, or $a\towall b$; and similarly for $b$ and $c$. 

If $a\wallto b$ or $a\stprec b$, then there exists $p\in a$
such that $p\pprec p'$ for all $p'\in b$. Since there exists $p'\in b$
and $p''\in c$ with $p'\pprec p''$, transitivity of $\pprec$ implies
that $p\pprec p''$, and so $a \cprec' c$. 

If $b\towall c$ or $b\stprec c$, the same argument holds with the
names changed. 

The remaining possibility is $a\towall b \wallto c$, and 
in this case $a\wall c$ by definition.

\subsubsection*{Claim 2:} Expressions of the form $a \cprec' b \wall c$
and $a\wall b \cprec' c$ can be reduced to expressions of length 2. 

Proof: $a\cprec' b \wall c$ means there exists $x\in C$ such that 
$a \cprec' b \towall x \wallto c$. Now Claim 1 reduces $a\cprec' b
\towall x$ to $a\cprec' x$, and a second
application of Claim 1 completes the job. The other case is similar. 

\medskip

We therefore obtain the latter part of the lemma: the transitive closure of 
$\cprec'$, which is the same as that of 
$\towall$, $\wallto$ and $\stprec$,  is obtained just by adjoining the 
relation $\wall$.

The fact that $\cprec$ is a partial order now reduces to checking that
$a \cprec a$ is never true. For $a\cprec' a$, 
this follows from Lemma~\ref{slice unordered}. 
If $a\wall a$,
then for some $x$ we have $a\towall x \wallto a$. But this means $x
\wallto a \towall x$ which by Claim 1 means $x\cprec' x$, which again
cannot occur.
\end{pf}

We deduce the following as a corollary.
\begin{corollary}{pprec order consistency}
If $a\cprec b$ then, for every $p\in a$ and
$p'\in b$ that are $\pprec$-ordered we have $p \pprec p'$. 
\end{corollary}

\begin{proof} 
If there exists $p\in a$ and
$p'\in b$ such that $p'\pprec p$, then by definition $b\cprec' a$. But
then $a \cprec b \cprec a$, which contradicts the fact that $\cprec$
is a strict partial order. 
\end{proof}

\subsection{Topological partial order}
Given a cut system $C$, each slice $a \in C$ determines either a split-level
surface $F_a$ as the disjoint union of the 3-holed spheres $\check
D(a)$ in $M_\nu$, or, if $a$ is annular, $a$ determines a solid torus
$\storus v$ where $v = \core(a)$.  We will now relate the
$\cprec$-order on the slices of a cut system to the $\topprec$-order on
their associated split-level surfaces and solid tori in $\modl$.
(We remind the reader that the ordering $\topprec$ is defined on
disconnected subsurfaces of $S$ in Section~\ref{topprec defs}, and
therefore $\topprec$ applies to the split-level surfaces $F_a$).

To begin with, we relate $\tprec$-ordering properties in the hierarchy
$H_\nu$ of the 3-holed spheres $Y$ arising as component domains in
$H_\nu$ and the annulus geodesics $k_v$ arising for each vertex $v \in
H_\nu$ to the topological ordering $\topprec$ applied to the level
surfaces $F_Y$ and solid tori $\storus{v}$ in $\modl$.  We will then
use these ordering relations to relate the $\cprec$-ordering to
$\topprec$-order on the surfaces $\hhat F_a$, for slices $a$ in a cut
system $C$.

In order to discuss this relationship, fix a resolution $\{\tau_i\}$
for the hierarchy $H$.  The following definitions allow us to keep
track of the parts of the resolution sequence where certain objects
appear. Let $v$ denote a simplex whose vertices appear in $H$, $Y$ a
3-holed sphere in $S$, $k$ a geodesic in $H$, and $p=(h,w)$ a
geodesic-simplex pair in $H$. Let $h'$ be a subsegment of the geodesic
$h$. Define:
\begin{align*}
J(v) &= \{i: v\subset \base(\tau_i)\}. \\
J(p) &= \{i:p\in \tau_i\}.\\
J(Y) &= \{i: Y\subset S \setminus \collar(\base(\tau_i))\}.\\
J(h) &= \{i: \exists v\in h, (h,v)\in\tau_i\}.\\
J(h') &= \{i:\exists v'\in h', (h,v/) \in\tau_i\}.
\end{align*}

To relate the appearance of these intervals in $\mathbb{Z}$ to partial
orderings in the hierarchy $H$, we record the following consequence of
the {\em slice order} $\sprec$ in \cite[\S5]{masur-minsky:complex2}.

\begin{lemma}{monotonicity}
Let $\tau_i$ and $\tau_j$ be slices in a resolution $\{\tau_n\}$ of
$H$ with $i < j$.  Then if $p \in \tau_i$ and $p' \in \tau_j$ are
$\pprec$-ordered, we have $$p \pprec p'.$$
\end{lemma}

\begin{proof} The slices in the resolution $\{ \tau_n \}$ are ordered
with respect to the order $\sprec$ on complete slices with bottom
geodesic the main geodesic of $H$.  By \cite[Lemma
5.3]{masur-minsky:complex2}, we have $\tau_i \sprec \tau_j$ and
therefore that each pair $q \in \tau_i$ either also lies in $\tau_j$
or there is a $q' \in \tau_j$ with $q \pprec q'$ (by the definition of
$\sprec$).

Assume that $p' \pprec p$.  Then $p$ and $p'$ cannot both lie in
$\tau_i$ by Lemma~\ref{slice unordered}, so it follows that there is a
$p''$ in $\tau_j$ with $p \pprec p''$.  By transitivity of $\pprec$ we
have $p' \pprec p''$, with $p'$ and $p''$ both in $\tau_j$, which
contradicts Lemma~\ref{slice unordered} applied to $\tau_j$.
\end{proof}

Clearly $J(Y)=J([\boundary Y])$, and if $p=(h,w)$ and $v\subset w$
then $J(p)\subset J(v)$. We also have:
\begin{lemma}{J interval}
Let $Y$, $p=(h,w)$,  $v \subset w$, $h$ and $h'$ be as above.  
Then $J(Y)$, $J(v)$, $J(p)$, $J(h)$ and $J(h')$  are all intervals.
\end{lemma}

\begin{proof}
The conclusion for $J(v)$  was proven in \cite[Lemma 5.16]{minsky:ELCI}.
For any simplex $w$, $J(w)=\cap_{v\in w} J(v)$ so
$J(w)$ is an interval too. For a 3-holed sphere $Y$, then, we use
the fact that $J(Y)= J([\boundary Y])$.

Suppose that $J(p)$ is not an interval for some pair $p=(h,v)$.
Then there exists $i<j<k$ such that $i,k\in J(p)$, but $j$ does
not lie in $J(p)$. Therefore, by  \cite[Lemma 5.3]{masur-minsky:complex2},
there exists $q\in \tau_j$ which is $\pprec$-orderable with
respect to $p$. But then Lemma \ref{monotonicity} implies
that both $p\pprec q$ and $q\pprec p$ which contradicts
the fact that $\pprec$ is a partial order. It follows that $J(p)$ is an interval.

We now consider a geodesic $h$. If $h$ has only one simplex $v$, then
$J(h)=J(h,v)$ is an interval by the previous paragraph, so we may assume
that $h$ has at least two simplices.
In any elementary move $\tau_i\to
\tau_{i+1}$, some pair $(h,v)\in\tau_i$ might be ``advanced'' to
$(h,\vsucc(v))$, some pairs $(k,u)$ where $u$ is the last simplex
might be erased, and some pairs $(k',u')$ where $u'$ is the first simplex,
might be created. Thus a geodesic can only appear at its beginning, advance
monotonically, and disappear at its end. If a geodesic $h$ makes two
appearances and $v$ is a simplex in $h$, then $J(h,v)$ will not be
an interval, contradicting the result of the previous paragraph.
This shows that $J(h)$ is an interval, and the same argument
applies to any subsegment of $h$. 
%
\end{proof}

Now we consider ordering relations for vertices and their associated
solid tori in $\modl$.  Given a vertex $v$, note that the interval
$J(v)$ is precisely the interval $J(k_v)$ corresponding to the annulus
geodesic $k_v$ associated to $v$.
\begin{lemma}{vertex ordering}
Let $v,v'$ be vertices appearing in $H$, whose corresponding curves
$\gamma_v$ and $\gamma_{v'}$ intersect non-trivially in $S$. The
following are equivalent:
\begin{enumerate}
\item \label{v order top} $\storus v \topprec \storus{v'}$.
\item \label{v order J} $\max J(v) < \min J(v')$
\item \label{v order k} $k_v \tprec k_{v'}$
\item \label{v order p}
For any pairs $p = (h,w)$ and $p' = (h',w')$,
if $v\in w$ or $h = k_v$,  and 
if $v'\in w'$ or $h' = k_{v'}$
we have $p \pprec p'$.
\end{enumerate}
Furthermore, either these relations or their opposites hold.
\end{lemma}

\begin{proof}
We note first that since $\gamma_v$ and $\gamma_{v'}$ intersect
non-trivially, the tori $\storus v$ and $\storus{v'}$ must be
$\topprec$-ordered.  Further, since $v$ and $v'$ cannot appear
simultaneously in any slice in the resolution, the intervals $J(v)$
and $J(v')$ must be disjoint.  Likewise, since the domains of the
geodesics $k_v$ and $k_v'$ overlap, Lemma~\ref{time order for domains}
guarantees that $k_v$ and $k_{v'}$ must be $\tprec$-ordered. By
Lemma~\ref{sweep order} the resolution $\{\tau_i\}$ yields a
sweep through the model $\modl$ by split-level surfaces which is
monotonic in the sense that if two overlapping level surfaces $W$ and
$W'$ appear in the sweep with $W$ occurring first, then $W\topprec
W'$.  This applies both to level surfaces and to solid tori with
respect to their associated annular domains.  Hence (\ref{v order
top}) and (\ref{v order J}) are equivalent.

Since every slice in the resolution is saturated, $J(v)=J(k_v)$,
so for $i \in J(v)$ we must have a
pair of the form $(k_v,u)$ in $\tau_i$. For $j\in J(v')$ we must have
some $(k_{v'},u')$ in $\tau_j$. Since $\gamma_v$ and $\gamma_{v'}$ intersect,
$k_v$ and $k_v'$ are $\tprec$-ordered, by \cite[Lemma 4.18]{masur-minsky:complex2},
so $(k_v,u)$ and $(k_{v'},u')$ are $\pprec$-ordered and
$(k_v,u)\pprec (k_{v'},u')$ if and only if $k_v\tprec k_{v'}$.
Lemma~\ref{monotonicity}
implies that $i<j$ if and only if $(k_v,u)\pprec (k_{v'},u')$.
Therefore, 
$i< j$  for all $i\in J(v)$ and $j\in J(v')$ if and only if $k_v\tprec k_{v'}$, so we see
that (\ref{v order J}) and (\ref{v order k}) are equivalent.

Now assume that  (\ref{v order J}) holds, i.e. that $\max J(v) < \min J(v')$.
Let $p=(h,w)$ with $v \in w$ or  $h=k_v$  and $p'=(h',w')$
with $v' \in w'$ or $h' = k_{v'}$, 
note that $J(p) \subseteq J(v)$ 
and $J(p') \subseteq J(v')$. 
It follows immediately that  $\max J(p)< \min J(p')$. 
Lemma~\ref{monotonicity} implies that if $p$ and $p'$
are $\pprec$-ordered, then $p\pprec p'$.

Since $\gamma_v$ and $\gamma_{v'}$ intersect
nontrivially, the 
domains $D(h)$ and $D(h')$ intersect. If one is not inside the other
then $h$ and $h'$ are $\tprec$-ordered 
by \cite[Lemma 4.18]{masur-minsky:complex2} and hence $p$ and $p'$
are $\pprec$-ordered. If $D(h) = D(h')$ are equal then, since $w\ne w'$, either
$w<w'$ or $w>w'$ and again they are $\pprec$-ordered. If, say,
$D(h)\subset D(h')$ then, since $\gamma_v$ intersects $\gamma_{v'}$,
we have $w'\notin \phi_{h'}(D(h))$. Lemma \ref{nonempty footprint}
implies that $\phi_{h'}(D(h))$ is non-empty, so it follows from the
definition of $\pprec$ that $p$ and $p'$ are $\pprec$-ordered.

We hence conclude that (\ref{v order J}) $\implies$  (\ref{v order
p}). Again for the opposite direction reverse the roles of $v$ and
$v'$.
\end{proof}

For 3-holed spheres, we have the following:

\begin{lemma}{Y ordering}
Let $Y$, $Y'$ be 3-holed spheres appearing as component domains in
$H$, and suppose that $Y$ and $Y'$ intersect essentially. Then the
following are equivalent:
\begin{enumerate}
\item
$F_Y\topprec F_Y' $
\item 
$ \max J(Y) < \min J(Y')$
\item
$Y \tprec Y'$.
\end{enumerate}
\end{lemma}

\begin{proof}
The equivalence of (1) and (2) follows from the same argument as in
Lemma \ref{vertex ordering}.

We now show that (3) implies (2). Suppose $Y\tprec Y'$.  Recall from
section~\ref{hierarchy background} that this implies there exists a
geodesic $m$ in $H$ such that $$ Y \fsub m \bsub Y' $$ and that $$\max
\phi_m(Y) < \min \phi_m(Y').$$

It follows, from Theorem \ref{Descent  Sequences}
that there exist geodesics
$f, f'$ (possibly the same) such that 
$$ Y \fsubd f \fsubeq m \bsubeq f' \bsubd Y'.$$ 
Since $Y$ is a 3-holed sphere, $\xi(D(f))=\xi(D(f'))=4$.
Let $v$ be the vertex of $f$ such that $Y$ is a
component domain of $(D(f),v)$, and let $v'$ be the vertex of $f'$
such that $Y'$ is a component domain of $(D(f'),v')$.

Since $Y$ is a three-holed sphere $v$ cannot be the last simplex
of $f$ and $v'$ cannot be the first simplex of $f'$.
There is exactly one elementary move in the resolution that replaces
$(f,v)$ with $(f,\vsucc(v))$. Before this move $Y$ is a complementary
domain of the slice marking, and afterwards it is not, since
$\vsucc(v)$ intersects $Y$. Hence, 
$$
\max J(Y) = \max J(v).
$$
The same logic gives us 
$$
\min J(Y') = \min J(v').
$$

We claim that $k_v \tprec k_{v'}$.  
The annulus $D(k_v)$ is a component domain of $(D(f),v)$ and
likewise $D(k_{v'})$ is a component domain of $(D(f'),v')$.  It follows that
$k_v \fsubd f$ and $f' \bsubd k_{v'}$, and thus $$ k_v \fsubd f
\fsubeq m \bsubeq f' \bsubd k_{v'}.$$
The claim follows provided the footprints of $D(k_v)$ and $D(k_{v'})$ on $m$
are disjoint and correctly ordered.
 
In the case that $m = f$ we note that since $f$ is a $4$-geodesic the
vertex $(f,\vsucc(v))$ intersects $v$ and $Y$.  Thus we have 
$$v = \max \phi_m(Y) = \max \phi_m(D(k_v)),$$
and likewise when $m= f'$ then we
have $v' = \min \phi_m(Y') = \min \phi_m(D(k_{v'}))$.

When $f \fsub m$, 
Lemma 5.5 of \cite{minsky:ELCI} implies that 
$$\max \phi_m(D(k_v)) = \max \phi_m(Y) = \max
\phi_m(D(f)).$$ Likewise if $m \bsub f'$ then we similarly conclude
that $$\min
\phi_m(D(k_{v'})) = \min \phi_m(Y') = \min \phi_m(D(f'))$$ We also
know, in particular, that all these footprints are non-empty.

Since we have $\max \phi_m(Y) < \min \phi_m(Y')$, it follows that
$\max \phi_m(D(k_v)) < \min \phi_m(D(k_{v'}))$ and thus that $k_v
\tprec k_{v'}$.

Applying Lemma~\ref{vertex ordering}, then, we have
$$
\max J(v) < \min J(v'),
$$ so we may conclude that $\max J(Y) < \min J(Y')$ and this
establishes (2).  

To show that (2) implies (3), note that since $Y$
and $Y'$ intersect they must be $\tprec$-ordered.  Hence, if (3) is false
we have $Y'\tprec Y$ and we apply the above argument to reach a
contradiction.
\end{proof}

\begin{lemma}{tori-Y order}
Let $Y$ be a 3-holed sphere and $v'$ a vertex appearing in $H$ such
that $\gamma_{v'}$ and $Y$ overlap.  Then the following are equivalent:
\begin{enumerate}
\item \label{Y<v' topprec} $F_Y \topprec U(v')$,
\item \label{Y<v' interval} $\max J(Y) < \min J(v')$, and
\item \label{Y<v' time} $Y \tprec k_{v'}$.
\end{enumerate}
Symmetric conditions hold for $U(v) \topprec F_{Y'}$,
and either these relations or their opposites hold.

\end{lemma}

\begin{proof}
We first show that (\ref{Y<v' topprec}) implies (\ref{Y<v' interval}). We recall that $F_Y$ is a straight surface
and that $U(v)$ is a straight solid torus in $\modl$.
If $F_Y \topprec U(v')$ then there
is a $v \in [\bdry Y]$ so that $\gamma_{v}$ and $\gamma_{v'}$ intersect
nontrivially.  Thus, as in the proof of Lemma~\ref{vertex ordering}, the straight solid tori
$\storus{v}$ and $\storus{v'}$ must be $\topprec$-ordered
and it follows that the height intervals they determine in the vertical
coordinate of $\modl$ are disjoint.  
Likewise, since $F_Y \topprec \storus{v'}$,
the height of $F_Y$ in $\modl$ must be less than the
minimum height of $\storus{v'}$ in $\modl$.  Since 
we have $F_Y \cap \bdry \storus{v} \not= \emptyset$, we conclude that 
$\storus{v} \topprec \storus{v'}$.
Lemma~\ref{vertex ordering} guarantees that  $\max J(v) < \min J(v')$, 
which implies that $\max J(Y) < \min J(v')$ since $J(Y) \subset
J(v)$.  Hence (\ref{Y<v' interval}) follows from (\ref{Y<v'  topprec}). 

Now assume that $\max J(Y) < \min J(v')$.  
Again, choose $v \in [\bdry Y]$ so that $\gamma_{v}$ and $\gamma_{v'}$ intersect
nontrivially. Since $J(v)$ and $J(v')$ are disjoint, $J(v)$ is an interval and $J(Y)
\subset J(v)$, we have $\max J(v) < \min J(v')$.  Applying
Lemma~\ref{vertex ordering} we have $\storus{v} \topprec \storus{v'}$
and therefore that $F_Y\topprec \storus{v'}$ (since $F_Y$ is a level
surface abutting $U(v)$). Hence 
(\ref{Y<v' topprec}) follows from (\ref{Y<v' interval}). 

Finally we show the equivalence of (\ref{Y<v' interval}) and (\ref{Y<v' time}).   
Arguing as in the proof of Lemma \ref{Y ordering}, if $Y \tprec k_{v'}$, then there is a
geodesic $m$ so that $$ Y \fsubd f \fsubeq m \bsubeq f' \bsubd k_{v'}.
$$ and $\max \phi_m(Y) < \min \phi_m(D(k_{v'})).$ Again, let $v \in
[\bdry Y]$ be such that $\gamma_v$ is a component of $\bdry Y$ that
intersects $\gamma_{v'}$.  Then, since $\gamma_v$ and $\gamma_{v'}$ intersect, the
geodesics $k_{v}$ and $k_{v'}$ are $\tprec$-ordered, and since
$\phi_m(Y) \subset\phi_m(D(k_v))$, we have $k_v\tprec k_{v'}$

Thus, by Lemma~\ref{vertex ordering}, we have $\max J(v) < \min J(v')$
and so since $J(Y) \subset J(v)$, we have 
$$\max J(Y) < \min J(v').$$
Therefore,  (\ref{Y<v' time}) implies (\ref{Y<v' interval}).

As before, to show the converse, we simply observe that since
$\gamma_{v'}$ and $Y$ overlap, we have that $Y$ and $k_{v'}$ are
$\tprec$-ordered.  If $k_{v'} \tprec Y$, we apply the above argument
to conclude $\max J(v') < \min J(Y)$, which is a contradiction.  This
shows the equivalence of (\ref{Y<v' interval}) and (\ref{Y<v' time}),
concluding the proof.
\end{proof}

Now let us go back to considering slices in a cut system $C$.  Let $a$
and $b$ be two distinct slices in the cut system $C$ whose domain
surfaces $D(a)$ and $D(b)$ overlap. We would like to ensure that the
surfaces $\hhat F_a$ and $\hhat F_b$ in the model are
$\topprec$-ordered if and only if $a$ and $b$ are consistently
$\cprec$-ordered.

Before commencing the proof, we argue that distinct non-annular slices $a$ and
$b$ in a cut system have no underlying curves in common.

\begin{lemma}{cuts dont meet}
If $a$ and $b$ are two distinct non-annular slices in a cut system
$C$, then $\base(a)$ and $\base(b)$ have no 
vertices in common. 
\end{lemma}

\begin{proof}

Suppose by way of
contradiction that there is a vertex $v$ common to $\base(a)$ and
$\base(b)$. 
If $D(a)=D(b)$ then $g_a=g_b$ and $v_a$ and $v_b$ are simplices on
$g_a$ spaced at least $\deep$ apart. However since $v$ is distance at
most 1
from both in $\CC(D(a))$, this is a contradiction. From now on we
assume that $D(a)\ne D(b)$. Thus, since $\check D(a)$ and $\check D(b)$ overlap,  
Lemma \ref{consistency} implies that $a$ and
$b$ are $\cprec$-ordered. Without loss of generality we may assume
$a \cprec b$, and moreover one of $a\towall b$, $a\wallto b$ or
$a\stprec b$ must occur. 

Let $p_1=(h_1,u_1)$ and $p_2=(h_2,u_2)$ be the pairs of $a$ and $b$,
respectively, such that $u_1$ and $u_2$ contain the vertex $v$. We
claim:
\begin{enumerate}
\item[(*)] There is a pair $q=(k,w)$ in the hierarchy such that 
$$p_1 \pprec q \pprec p_2 $$ 
and $\gamma_w$ intersects $\gamma_v$
non-trivially.
\end{enumerate}

To see this, suppose first that $a\towall b$. 
As in part (3) of the proof of Lemma \ref{consistency}, there exists
a pair $p'=(h',x')\in b$ such that
\begin{equation}\label{all ga before p'}
(g_a,{\bf T}(g_a)) \pprec p',
\end{equation}
and hence for any simplex $u$ of $g_a$ we have $(g_a,u)\pprec p'$. 
Let $k=g_a$,
and let $w$ be a simplex in $k$ such that $v_a < w$ and
$d_{D(a)}(v_a,w) > 3$. This is possible because of the lower spacing
constant $d_1 \ge \deep$, and $q=(k,w)$ will be our desired pair: Since
$d_{D(a)}(v_a,v)\le 1$, we have that $\gamma_v$ must intersect
$\gamma_w$. Since the footprint $\phi_k(D(h_1))$ contains $v_a$ and
has diameter at most 2, we also have $\max \phi_k(D(h_1)) < w$, and so
$p_1 \pprec q$. Because $\gamma_w$ and $\gamma_v$ have non-trivial
intersection, $q$ and $p_2$ must be $\pprec$-ordered, by Lemma
\ref{vertex ordering}. 
Thus, to show that $q\pprec p_2$ it suffices to rule out
$p_2\pprec q$. But since $q\pprec p'$, if $p_2\pprec q$ then
$p_2\pprec p'$, and this contradicts the non-orderability of different
pairs in a slice (Lemma \ref{slice unordered}).  We conclude that $q$
satisfies Claim (*). 
If $a\wallto b$ then a symmetric argument again yields $q$.

Finally if
$a\stprec b$ then (again by Lemma \ref{consistency}) $g_a \tprec g_b$
and (\ref{all ga before p'}) holds for $p' = (g_b,v_b)$, the bottom
pair of $b$. Hence the same argument yields Claim (*). 

\medskip

Now fixing a resolution for the hierarchy, let $J(q)$ and $ J(p_i)$ be
defined as before. Since $\gamma_w$ intersects $\gamma_v$ we must have
that $J(q)$ is disjoint from both $J(p_1)$ and $J(p_2)$. Since
$p_1\pprec q \pprec p_2$, we can apply Lemma
\ref{vertex ordering} to obtain $$ \max J(p_1) < \min J(q) \le \max J(q)
< \min J(p_2).  $$ On the other hand, both $J(p_1)$ and $J(p_2)$ are
contained in $J(v)$, which is an interval disjoint from $J(q)$. This
is a contradiction, and Lemma \ref{cuts dont meet} is established.
\end{proof}
We remark that as a consequence of Lemma~\ref{cuts dont meet}, if $a$
and $b$ are each non-annular slices in a cut system $C$ and $Y \subset
\check D(a)$ and $Y' \subset \check D(b)$ are 3-holed spheres,
then $Y$ and $Y'$ are distinct; otherwise there would be some common vertex
in $\base(a)$ and $\base(b)$ in their common boundary.

\subsection{Comparing topological and cut ordering}
To relate $\cprec$-order of $a$ and $b$ to topological ordering in the
model, we relate the order properties we have obtained for the
constituent annuli and $3$-holed spheres making up $D(a)$ and $D(b)$
to the corresponding level 3-holed spheres and annuli in $\hhat F_a$
and $\hhat F_b$ or the corresponding solid tori when either $a$ or $b$
is annular.  We call level $3$-holed spheres and annuli in $\hhat F_a$
{\em pieces} of $\hhat F_a$, when $a$ is non-annular, and likewise for
$\hhat F_b$.  Then our first task will be to demonstrate that whenever
pieces of $\hhat F_a$ and $\hhat F_b$ overlap (or the solid tori, when
$a$ or $b$ is annular) they are topologically ordered consistently
with the cut ordering on the slices $a$ and $b$.

\begin{lemma}{pieces ordered}
Let $a$ and $b$ be non-annular slices in a cut system $C$ such that $a
\cprec b$ 
and $D(a)$ and $D(b)$ overlap.  Let $v \in \base(a)$ and $v' \in
\base(b)$ be vertices of $H$, and let $Y \subset \check D(a)$
and $Y' \subset \check D(b)$ be 3-holed spheres.  Then the following
holds:
\begin{enumerate}
\item if  $\gamma_v \cap \gamma_{v'} \not= \emptyset$ then $\storus{v}
\topprec \storus{v'}$,
\label{vertices}
\item if $\gamma_v \cap Y' \not= \emptyset$ then $\storus{v} \topprec
F_{Y'}$,
\label{mixed1}
\item if $Y \cap \gamma_{v'} \not= \emptyset$ then $F_Y \topprec
\gamma_{v'}$, and 
\label{mixed2}
\item if $Y \cap Y' \not= \emptyset$ then $F_Y \topprec F_{Y'}$.
\label{Ys}
\end{enumerate}
If $a$ is annular take  $v = \core(a)$, and if $b$ is annular take 
 $v' = \core(b)$. With this notation, 
$(\ref{vertices})$ holds in all cases, 
$(\ref{mixed1})$ holds if just $a$ is annular,
and $(\ref{mixed2})$ holds if just $b$ is annular. 
\end{lemma}

\begin{pf}
First assume neither $a$ nor $b$ is annular.

Suppose that $v\in\base(a)$ and $v'\in\base(b)$ and that 
$\gamma_v\cap\gamma_{v'}\not=\emptyset$.
Let $q=(h,w)\in a$ be a pair such that $v\in w$ and let
$q'=(h',w')\in b$ be a pair such that $v'\in w'$. Then
$q$ and $q'$ are $\pprec$-ordered, since $\gamma_v$ and $\gamma_{v'}$
intersect.
Corollary \ref{pprec order consistency} then gives
that $q\pprec q'$. Lemma \ref{monotonicity} implies that
$\max J(q) < \min J(q')$. Thus, since $J(q)\subset J(v)$,
$J(q')\subset J(v')$ and $J(v)$ and $J(v')$ are disjoint
intervals, we have
$$\max J(v) < \min J(v').$$ 
Lemma \ref{vertex ordering} then implies that
$$\storus{v} \topprec \storus{v'}.$$
So, we have established (1).

If $v$ lies in $[\bdry D(a)]$, we claim that for each $p \in a$ we
have $$J(p) \subseteq J(v).$$ To see this, we note that since $v$
represents a curve in the boundary of $D(a)$, the vertex $v$ is
present in the base of any complete slice containing a pair with $g_a$
as its geodesic; in particular we have that $$J(g_a)
\subseteq J(v).$$ 
By Lemma~\ref{deep placement}, we have $$(g_a , {\bf I}(g_a)) \ppreceq p
\ppreceq (g_a, {\bf T}(g_a)).$$ 
Applying
Lemma~\ref{monotonicity} we have 
$$\min J(g_a,\I(g_a)) \le \min J(p) \ \ \text{and} \ \ \max J(p) \le
\max J(g_a,\T(g_a))$$ and so we conclude that
$$\min J(v)  \le \min J(g_a) \le \min J(p) 
\le \max J(p) \le \max J(g_a) \le \max J(v)$$
and thus $J(p) \subseteq J(v)$, since $J(v)$ is an interval by
Lemma~\ref{J interval}. 

Similarly, we see that if $v'$ lies in $[\bdry D(b)]$, then
$J(p')\subset J(v')$ for all $p'\in b$.

We will need the following generalization of (1) in the proofs of (2)
and (3). 

{\em \begin{enumerate}
\item[($1+$)] Assume $a$ and $b$ are non-annular and $\gamma_v\intersect \gamma_{v'} \ne \emptyset$.
If  either 
\begin{enumerate}
\item[(i)] $D(a)$ is not contained in $D(b)$, $v\in [\bdry D(a)]$ and
$v'\in\base(b)$, 
\item[(ii)]
$D(b)$ is not contained in $D(a)$, $v\in \base(a)$ and
$v'\in[\bdry D(b)]$, or
\item[(iii)]
$D(a)$ and $D(b)$ are not nested, $v\in[\bdry D(a)]$ and
$v'\in[\bdry D(b)]$,
\end{enumerate}
then 
$$\max J(v) < \min J(v')\ \ \ \ {\rm and}\ \ \ \ \ \storus{v} \topprec
\storus{v'}.$$
\end{enumerate}
}

\begin{proof}{}
First assume that $D(a)$ is not contained in $D(b)$, $v\in [\bdry D(a)]$ and $v'\in\base(b)$.
Lemma \ref{consistency} implies that either $a\stprec b$ or $a\wallto b$. Therefore,
we may choose $q\in a$ such that $q\pprec q'$ for all $q'\in b$. In particular,
if we choose $q'=(h,w')$ such that $v'\in w'$, then $J(q)\subset J(v)$
and $J(q')\subset J(v')$. Again, since $J(v)$ and $J(v')$ are disjoint
intervals and, by Lemma \ref{monotonicity}, $\max J(q) < \min J(q')$, we see
that $\max J(v) < \min J(v')$ and, applying Lemma \ref{vertex ordering},
$\storus{v} \topprec \storus{v'}.$ So, we have established ($1+$) in
case (i).
A very similar argument handles cases (ii) and (iii).
\end{proof}


We now establish part~(\ref{mixed1}) of Lemma \ref{pieces ordered}.  If $v \in \base(a)$ and there
is a $Y' \subset \check D(b)$ with $\gamma_v
\cap Y' \not= \emptyset$, then since $\base(a)$ and $\base(b)$ share no
vertices, by Lemma~\ref{cuts dont meet}, either
\begin{enumerate}
\item[{ (A)}] $D(a) \subset D(b)$ and there is
a $v' \in \base(b) \cap [\bdry Y']$ with $\gamma_v \cap \gamma_{v'} \not= \emptyset$,  or
\item[{ (B)}] $D(a)$ is not contained in  $D(b)$ and there is a
$v' \in [\bdry Y']$ for
which $\gamma_v \cap \gamma_{v'} \not= \emptyset$.   
\end{enumerate}
(Note that in case~(B) we must allow the possibility that the vertex $v'$ lies
in $[\bdry D(b)]$).
Applying part (1)  in case (A) and
part (i) of ($1+$)
in case (B),
we have
$$\max J(v) < \min J(v') \le \min J(Y'),$$
and therefore 
$$\storus{v} \topprec F_{Y'}$$ by Lemma~\ref{tori-Y order}.  The
argument for part~(\ref{mixed2}) is symmetrical.

Finally, for part~(\ref{Ys}), if there are 3-holed spheres $Y \subset
\check D(a)$ and $Y'
\subset \check D(b)$ that overlap then once again
by Lemma~\ref{cuts dont meet} either
\begin{enumerate}
\item[{\bf (i)}] $D(a)=D(b)$ and there exists $v \in \base(a)\cap [\bdry Y]$
and $v' \in \base(b) \cap [\bdry Y']$ with $\gamma_v \cap \gamma_{v'} \not=
\emptyset$,
\item[{\bf (ii)}]  $D(a) \subset D(b)$ and there exists $v \in [\bdry Y]$
and $v' \in \base(b) \cap [\bdry Y']$ with $\gamma_v \cap \gamma_{v'} \not=
\emptyset$,
\item[{\bf (iii)}]  $D(b) \subset D(a)$ and there exists $v \in\base(a)\cap [\bdry Y]$
and $v' \in [\bdry Y']$ with $\gamma_v \cap \gamma_{v'} \not=
\emptyset$,
or 
\item[{\bf (iv)}]  $D(a)$ and $D(b)$ are non-nested and there are
vertices $v \in [\bdry Y]$ and $v' \in [\bdry Y']$ with $\gamma_v \cap
\gamma_{v'} \not= \emptyset$.   
\end{enumerate}
In each of these cases, part (1)  or ($1+$) implies that we
have $$\max J(Y) \le \max J(v) < \min J(v') \le \min J(Y'),$$ so we
may apply Lemma~\ref{Y ordering} to conclude that $$F_Y
\topprec F_{Y'}$$ as desired. 

We now consider the cases when either $a$, $b$ or both $a$ and $b$ are annular
If $a$ is annular, choose $q=(k_v,w)$ to be the unique pair in the slice $a$.
Similarly, if $b$ is annular, choose $q'=(k_{v'},w')$ to be the unique pair in the
slice $b$. The proof of (1) now goes through verbatim in all cases.

Now suppose that $a$ is annular and $Y$ is a 3-holed sphere in $\check D(b)$.
Since $\gamma_v$ and $Y$ overlap, there exists $v'\in [\bdry Y]$ such that
$\gamma_v$ and $\gamma_{v'}$ intersect. One may argue, as in the proof
of ($1+$),
that $\max J(v) < \min J(v')$ and $\storus{v} \topprec \storus{v'}.$
The proof of (2) in this case then proceeds as in the non-annular case.
The proof of (3) when $Y$ is a 3-holed sphere in $\check D(a)$ and $b$ is
annular proceeds similarly.


\end{pf}

We are now ready to prove the following.
\begin{proposition}{topprec and cprec}
If $a$ and $b$ are slices in a cut system $C$ with overlapping domains, then
$$
\hhat F_a \topprec \hhat F_b \iff a \cprec b.
$$
\end{proposition}

\begin{proof}
Assume that $a \cprec b$.  
Since $D(a)$ and $D(b)$ overlap, we know
that $\check D(a)$ and $\check D(b)$ also overlap, since otherwise
either $$D(a) \subset \overline{\collar}(\base(b)) \ \ \ \text{or} \ \
\ D(b) \subset \overline{\collar}(\base(a))$$ and in this case $a$ and $b$ are
not $\cprec$-ordered by Lemma~\ref{consistency}.  It follows that all
vertices of $\base(a)$ and $\base(b)$, or the vertices corresponding
to the cores of $D(a)$ or $D(b)$ if either is annular, satisfy the
hypotheses of Lemma~\ref{pieces ordered}.  Thus, whenever pieces or
the solid tori making up $\hhat F_a$ and $\hhat F_b$ overlap they are
consistently topologically ordered in $\modl$.

Given $t \in \reals$, let $T_t \colon S\times \reals \to S\times
\reals$ be the translation 
$T_t(x,s) = (x,s+t)$ in the vertical ($\reals$) direction, and
consider the embeddings of $\hhat F_a$ and $\hhat F_b$ into $S \times
\reals$ as subsets of $\modl$ (see section~\ref{model definitions}).
Then the consistent topological ordering guarantees that for each
positive $s$ and $t$ we have $$T_{-s}(\hhat F_a) \cap T_{t} (\hhat
F_b) =
\emptyset$$
(recall that when $a$ is not annular, each annular piece of $\hhat
F_a$
is contained in a solid torus $\storus{v}$ for $v \in \base(a)$).

Thus, these translations provide homotopies of $\hhat F_a$ to
$-\infty$ in the complement of $\hhat F_b$ and of $\hhat
F_b$ to $+\infty$ in the complement of $\hhat F_a$.
Applying Lemma~\ref{basic ordering properties}, it follows that
$$\hhat F_a \topprec \hhat F_b.$$

Conversely, we assume that $\hhat F_a \topprec \hhat F_b$.  Since
$D(a)$ and $D(b)$ are overlapping domains, Lemma~\ref{consistency}
guarantees that either 
\begin{itemize}
\item $a$ and $b$ are $\cprec$-ordered, or 
\item we have $D(a) \subset \overline{\collar}(\base(b))$ or $D(b)
\subset \overline{\collar}(\base(a)).$   
\end{itemize}
In the latter case, if $D(a) \subset \overline{\collar}(\base(b))$
then $a$ is an annular slice, and $D(a) = \overline{\collar}(v)$ for
some $v \in \base(b)$.  Then $\hhat F_b$ intersects $\storus{v} =
\hhat F_a$ in an
annulus, and so $\hhat F_b$ and $\hhat F_a$ are not
$\topprec$-ordered, which is a contradiction.  The symmetric argument
rules out $D(b) \subset \overline{\collar}(\base(a))$.

Thus $a$ and $b$ are $\cprec$
ordered.  If $b \cprec a$, then the previous argument guarantees that
$$\hhat F_b \topprec \hhat F_a$$ 
contradicting the hypothesis.  This completes the proof.

\end{proof}

We conclude with the following consequence, guaranteeing that
topological order on the extended split-level surfaces arising from a
cut system is a partial order.

\begin{proposition+}{Topological Partial Order} The relation
$\topprec$ on the components of $\{\hhat F_a: a\in C\}$ has no
cycles, and hence its transitive closure is a partial order.
\end{proposition+}

\begin{proof}
We have shown that $\topprec$ is equivalent to the relation $\cprec$
restricted to surfaces $\{\hhat F_a : a \in C \}$  whose domains
overlap.  

Thus the transitive closure of
$\topprec$ over all the cut surfaces is a subrelation of $\cprec$
(which was already transitive). Since $\cprec$ is a partial order,
$\topprec$ has no cycles.
\end{proof}

\section{Regions and addresses}
\label{regions}

In this section we will explore the way in which a cut system
divides the model manifold into complementary regions, whose size and
geometry are bounded in terms of the spacing constants of the cut system.

For the remainder of the section
we fix a cut system $C$. The split-level surfaces $\{F_\tau:
\tau\in C\}$  divide $\modl[0]$ into components which we call
{\em complementary regions} of $C$ (or just {\em regions}).

In \S\ref{define address}, we will define the {\em address} of a block in
$\modl[0]$ in terms of the way the block is nested among the
split-level surfaces of $C$.
In \S\ref{address region structure} we will then
describe the structure of each subset $\XX(\alpha)\subset
\modl[0]$ consisting of blocks with address $\alpha$. In particular
Lemma \ref{product minus products} will show that, roughly
speaking, $\XX(\alpha)$ can be described as a product region bounded
between two split-level surfaces, minus a union of smaller product
regions (and tubes).
We will also prove Lemma \ref{region address}, which shows that each
complementary region of $C$
lies in a unique $\XX(\alpha)$. 

In \S\ref{region size} we will bound the size (i.e. number of blocks)
of each $\XX(\alpha)$, and hence of each complementary region of $C$.

In \S\ref{filled regions} we will extend the discussion to the filled
model $\modl[k]$ with $k\in[0,\infty]$. The filled cut surfaces $\hhat
F_\tau[k]$ cut $\modl[k]$ into connected components, and we shall show in
Proposition \ref{modl[k] regions} that, under appropriate
assumptions on the spacing 
constants of $C$, these components correspond in a simple way to the 
components in $\modl[0]$.
 
In the rest of the section, for an internal block $B$ let $W_B$ denote
the ``halfway surface'' $D(B)\times\{0\}$ in the parametrization of
$B$ as a subset of $D(B)\times[-1,1]$. If $B$ is a boundary block let
$W_B$ denote its outer boundary.

\subsection{More ordering lemmas}

Before we get started let us prove three lemmas 
involving slice surfaces and $\topprec$.

The first is another ``transitivity'' lemma.

\begin{lemma}{WB transitive}
Let $c$ and $d$ be two slices in a cut system $C$, and let $B$ be a
block with $D(B) \subset D(c)\intersect D(d)$. If the halfway surface
$W_B$ satisfies
$$\hhat F_c \topprec W_B \qquad\text{and}\qquad W_B \topprec \hhat F_d$$
then
$$\hhat F_c \topprec \hhat F_d.$$
\end{lemma}

Although this statement does not seem surprising we note that, 
since $\topprec$ is not in general transitive, and $W_B$ is
not itself a cut surface, we must take care in the proof.

\begin{proof}
Assume first that $B$ is an internal block.
A cut surface $\hhat F_\tau$ is a union of 
level surfaces (3-holed sphere gluing surfaces) and 
annuli embedded in straight solid tori. Call the level surfaces and the
solid tori the ``pieces'' associated to $\hhat F_\tau$.
$W_B$ 
is also a level surface, and moreover avoids (the interiors of) all
solid tori and gluing 
surfaces in $\modl$. It is therefore 
$\topprec$-ordered with any piece 
which it overlaps. For overlapping pieces $x,y,z$ it is easy
to see that $x\topprec y$ and $y\topprec z$ implies $x\topprec z$.
Now let $x$ and $y$ be pieces associated with $c$ and $d$,
respectively, which overlap each other and $W_B$. These exist 
since $D(B)\subset D(c) \intersect D(d)$, and the projections of the
pieces of $c$ and $d$ to $D(c)$ and $D(d)$, respectively,
decompose them into essential subsurfaces. These subsurfaces cover all
of $D(B)$, and hence must intersect each other there.

From the hypotheses of the lemma we conclude that
$x\topprec W_B$, and $W_B \topprec y$, and therefore
$x\topprec y$. 

Now since $c$ and $d$ have overlapping domains they are
$\cprec$-ordered by Lemma \ref{consistency}, 
and by Lemma \ref{topprec and cprec} we may conclude that 
either
$\hhat F_c\topprec \hhat F_d$ or
$\hhat F_d\topprec \hhat F_c$. The latter implies $y\topprec x$, which
contradicts $x\topprec y$. We conclude 
that $\hhat F_c\topprec \hhat F_d$.

If $B$ is a boundary block then the theorem is vacuous, since $W_B$ is
part of the boundary of $\modl$, and is embedded in $\hhat S\times\R$
in such a way that nothing in $\modl$ lies above it (if it is in the
top boundary) or below it (if it is in the bottom).
\end{proof}

The following lemma tells us that we can compare blocks and cut
surfaces, whenever they overlap. 

\begin{lemma}{block slice order}
Let $B$ be any block  and let $\tau$ be
any saturated nonannular slice. If $W_B$ and $\hhat F_\tau$ overlap, 
then they are $\topprec$-ordered. 
\end{lemma}

\begin{proof}
Again the lemma is immediate for boundary blocks so we may assume that
$B$ is internal.
The first step is to extend $\tau$ to a maximal nonannular slice. 
Note that since $\base(\I(H))$ and $\base(\T(H))$ are maximal
laminations, any saturated slice $\tau$ is full (see 
\S\ref{hierarchy background}). 
Hence if the bottom geodesic $g_\tau$ is $g_H$, we are done. If not, there
is some $h$ such that $g_\tau \fsubd h$, and a simplex $w$ in $h$ 
such that $D(g_\tau)$ is a component domain of $(D(h),w)$. Add $(h,w)$
to $\tau$, and successively fill in the components of $D(h)\setminus
D(g_\tau)$ to obtain a saturated nonannular slice $\tau'$ with $g_{\tau'}  =
h$. Repeat this inductively until we get a saturated nonannular slice
$\tau_0$ with bottom geodesic $g_H$, hence a maximal nonannular slice.

Now by Lemma 5.7 of \cite{minsky:ELCI}, there exists a (nonannular) resolution
with $\tau_0$ as one of its slices. 
If we now consider the sweep 
through $\modl$ determined by this resolution (see \S\ref{sweep}), we see that
there is some moment when the block $B$ is appended.
Applying Lemma \ref{sweep order}, for any slice
$\tau_i$ that occurs in the resolution before this moment we 
have $\hhat F_{\tau_i}\topprec W_B$, and for any $\tau_i$ that occurs after, we
have $W_B \topprec \hhat F_{\tau_i}$. Since $\hhat F_\tau$ is an
essential subsurface of 
$\hhat F_{\tau_0}$, 
and $\hhat F_\tau$ and $W_B$ overlap, this implies (using 
Lemma \ref{basic ordering properties})
that they are $\topprec$-ordered. 
\end{proof}

The next lemma allows us to compare tubes and cut surfaces, and will be 
used to prove the ``unwrapping property'' at the end of the proof in
Section \ref{proofmain}. It shows, in particular, that a slice surface
$\hhat F_c$ and a disjoint tube $U$ can be moved to $-\infty$ and
$+\infty$, respectively (or vice versa) without intersecting each
other. In Section \ref{proofmain} we will apply this to their images
in a hyperbolic 3-manifold $N$ to conclude that certain surfaces cannot be
wrapped around ``deep enough'' Margulis tubes, and this will allow us to construct 
controlled embedded surfaces in $N$.

\begin{lemma}{unwrapping}
Let $\tau$ be any saturated nonannular
slice in $H_\nu$ and let $w$ be a vertex of $H_\nu$,
such that $\collar(w)$ and $\check D(\tau)$ have non-trivial
intersection. 

Then either $\hhat F_\tau \topprec U(w)$ or
$U(w) \topprec \hhat F_\tau$.
\end{lemma}

\begin{proof}
As in Lemma \ref{block slice order}, we extend $\tau$ to a maximal
slice $\tau_0$, and fix a resolution of $H$ that includes $\tau_0$. 
The assumption that $\collar(w)$ and $\check D(\tau)$ intersect
implies that $\hhat F_{\tau_0}$ does not intersect $U(w)$.
Thus in the sweep through $\modl$ defined by the resolution, $\hhat
F_{\tau_0}$
is reached either before or after $U(w)$, and it follows as in
Lemma \ref{block slice order}
that they are $\topprec$-ordered. 
\end{proof}

\subsection{Definition of addresses}
\label{define address}

An {\em  address pair} for a block $B$ in $\modl$
is a pair of cuts $(c,c')$ with $D(B)\subset D(c)=D(c')$, such that
$$
\hhat F_{c} \topprec W_B \topprec \hhat F_{c'}.
$$
We say that an address pair $(c,c')$ is {\em nested within} a
different address pair 
$(d,d')$ if $d\cpreceq c$ and $c'\cpreceq d'$. 
We say that an address pair is 
{\em innermost} if it is minimal with respect to the relation
of nesting among address pairs for $B$.

\begin{lemma}{Addresses unique}
If $B$ has at least one address pair then 
it has a unique innermost address pair $(c,c')$ and furthermore
$(c,c')$ is nested within every other address pair for $B$. 
\end{lemma}

\begin{proof}
Let $(c,c')$ and $(d,d')$ be address pairs for $B$. 
We first claim that one of $D(c)$ and $D(d)$ must be contained in the other,
and that if $D(c) \subsetneq D(d)$ then $(c,c')$ is nested within $(d,d')$.

Since $D(B) \subseteq D(c)$ and $D(B) \subseteq D(d)$, the domains
$D(c)$ and $D(d)$ intersect. 
First assume that neither $D(c)$ nor $D(d)$ is contained in
the other.  In this case the bottom geodesics $g_c$ and $g_d$ are
$\tprec$-ordered (by Lemma \ref{time order for domains}), and without loss
of generality we may assume $g_c \tprec g_d$.
Note also that $g_{c'}=g_c$.

By Lemma~\ref{consistency}, we have $c' \cprec d$ which implies that
$\hhat F_{c'} \topprec \hhat F_d$ by Lemma~\ref{topprec and cprec}.  
On the other hand, by definition of address pairs, we have
$$
 \hhat F_d \topprec W_B 
\qquad\text{and}\qquad
W_B \topprec \hhat F_{c'},
$$
which by Lemma \ref{WB transitive} then implies 
$\hhat F_d \topprec \hhat F_{c'}$. But 
this contradicts $\hhat F_{c'} \topprec \hhat F_d$, by definition of $\topprec$.
We conclude that one of the domains is contained in the other.

Suppose that $D(c) \subsetneq D(d)$.  We claim that in this case
we must have
\begin{equation}
  \label{c d nesting}
d \cprec c \cprec c' \cprec d'  
\end{equation}
so that $(c,c')$ is nested within $(d,d')$.

To see this, note that by Lemma~\ref{consistency} we have that either
$$ c \towall d  \ \ \ \text{or}  \ \ \ d \wallto c $$ in the partial
order on cuts.  Suppose first that 
$c \towall d$. Then 
there is some $p\in d$ such that
for the bottom pair $p_c$ of $c$, 
$p_c \pprec p$, and in fact
the proof of
Lemma~\ref{consistency} shows that 
$(g_c,\T(g_c)) \pprec p$. 
Lemma \ref{deep placement} then shows that 
for  any pair $q\in c'$,
$q \pprec p$. This implies that
$$c' \towall d.$$
By Lemma~\ref{topprec and cprec} we have
$$\hhat F_{c'} \topprec \hhat F_d,$$
and since $(c,c')$ and $(d,d')$ are address pairs we have both
$$
\hhat F_d \topprec W_B
\qquad \text{and} \qquad
W_B \topprec \hhat F_{c'}.
$$ 
By Lemma \ref{WB transitive}, this implies $\hhat F_d \topprec \hhat
F_{c'}$, again a contradiction.

Thus we have ruled out $c\towall d$, and 
it follows that $d \wallto c$. By the same  argument with directions
reversed, we may also conclude that $c' 
\towall d'$.  This establishes the nesting claim (\ref{c d nesting}). 

Now suppose that $D(c) = D(d)$. 
The relation $\topprec$  on surfaces $\{\hhat F_\tau: \tau\in C, \
g_\tau = g_c\}$ 
is the same as the linear order on their bottom simplices
$\{v_\tau\}$. Thus by Lemma \ref{WB transitive}
the sets $\{\tau\in C: g_\tau=g_c, W_B\topprec \hhat F_\tau\}$
and $\{\tau\in C: g_\tau=g_c, \hhat F_\tau\topprec W_B\}$ form
disjoint intervals in this order,
and there is a unique innermost
pair.

Since we have shown that the domains of address pairs are linearly
ordered by inclusion, there is a unique domain of minimal complexity,
and among the pairs with that domain there is a unique innermost
one. This is the desired address pair. 
\end{proof}

We are now justified in making the following definition.
\begin{definition}{address}
If $(c,c')$ is the innermost address pair for $B$ then we say $B$ has {\em
address} $\add{ c,c'}$.  If $B$ has no address pairs, 
we say that $B$ has the {\em empty
address} denoted $\add{ \varnothing}$.
\end{definition}

We let $D(\add{c,c'})$ denote $D(c)=D(c')$ and let $D(\add\varnothing)=S$.
Note that if $\add{c,c'}$ is an address then $c$ and $c'$ are
successive in the $\cprec$ order on $C|_{g_c}$. 

\subsection{Structure of address regions}
\label{address region structure}
Having shown that each block has a well defined address, let
$\XX(\alpha)$ denote the union of blocks
with address $\add\alpha$.
We will now describe the structure of $\XX(\alpha)$
as, roughly speaking, a product region minus a union of smaller
product regions.

If $(c,c')$ is any address pair, 
note (e.g. by Lemma \ref{topprec and cprec})
that $\hhat F_c$ and $\hhat F_{c'}$ are disjoint unknotted properly
embedded
surfaces in $D(c)\times \R \subset S\times\R$, which are isotopic to
level surfaces, and transverse to the $\R$ direction. 
Hence they cut off from $D(c)\times\R$ a region $\BB=\BB(c,c')$
homeomorphic to 
$D(c)\times[0,1]$ containing (the closure of)
all points above $\hhat F_{c}$ and below $\hhat F_{c'}$ (in the $\R$ coordinate).
Define also $\BB(c,\cdot)$ to be the set of all points in
$D(c)\times\R$ that are above $\hhat F_c$, and define $\BB(\cdot,c')$
similarly (these are useful for considering infinite geodesics).

The boundary of $\BB$ in $S\times\R$ is therefore the union of 
$\hhat F_c \union \hhat F_{c'}$ with annuli
in $\boundary D(c)\times\R$, one for each component of $\boundary
D(c)$. 
Indeed these annuli lie in $\boundary U(\boundary D(c))$, and their
boundaries are the circles $\boundary \hhat F_c$ and $\boundary \hhat
F_{c'}$.

It is clear from this description that
a block $B$ is contained in $\BB(c,c')$ if and only if $(c,c') $ is an
address pair for $B$. If we define
$\BB(\varnothing)$ to be all of $\modl$,
we can generally say that 
$\XX(\alpha) \subset \BB(\alpha)$.

Furthermore if $(d,d')$ is any address pair that is
nested within $(c,c')$  
then $\BB(d,d')$ has interior disjoint from
$\XX(\alpha)$. Similarly for any $(d,d')$, $\BB(d,d')$ has interior
disjoint from $\XX(\varnothing)$.

In fact it follows from the definitions that $\XX(\alpha)$
is obtained by deleting from the product region $\BB(\alpha)$  all 
(interiors of) such product regions $\BB(d,d')$, as well as the tubes
$\UU$.  Recall that $\UU$ is the collection of all tubes in the model manifold.

If $g$ is a geodesic with $C|_g$ nonempty, 
let $a_g$ and $z_g$ be the first and last slices of $C|_g$, if they
exist ($g$ may be infinite in either direction)
and define
$\BB(g) = \BB(a_g,z_g)$ if $g$ is finite, $\BB(g) = \BB(a_g,\cdot)$ if 
$a_g$ exists but not $z_g$, and $\BB(g) = \BB(\cdot,z_g)$ if $z_g$
exists but not $a_g$. 

We call a geodesic $h$ an {\em inner boundary geodesic} for $\alpha$
if $D(h)\subsetneq D(\alpha)$,  $h$ supports slices $d,d'\in C|_h$
which are nested 
within $\alpha$, and $D(g)$ is maximal by inclusion among such geodesics.
For $\alpha=\add\varnothing$, the same
definition holds with the convention that every pair $(d,d')$ is said
to be nested in $\varnothing$.

The following lemma describes the region  $\XX(\alpha)$.

\begin{lemma}{product minus products}
If $\alpha$ is an address  for $C$, then:
\begin{enumerate}
\item If $h$ is an inner boundary geodesic for $\alpha$, then $\BB(h)\subset \BB(\alpha)$
and $int(\BB(h))\cap \XX(\alpha)=\emptyset$.
\item If $h,h'$ are inner boundary geodesics for $\alpha$, then $\BB(h)\cap\BB(h')=\emptyset$.
\item
$$
\XX(\alpha)=\BB(\alpha) \setminus \bigl(\UU \union \bigcup_h int(\BB(h))\bigr)
$$
where the union is over all
inner boundary geodesics $h$ for $\alpha$.
\end{enumerate}
Moreover, if $h\in H$ and $C|_h$ is  non-empty, then
$h$ is an inner boundary geodesic for exactly one address $\alpha$.
\end{lemma}

When $\alpha = \add{c,c'}$ 
we call $F_c$ and $F_{c'}$ the outer boundaries of $\XX(\alpha)$. 
The surfaces $F_{a_h}$ and $F_{z_h}$ for any inner boundary geodesic
$h$ are called inner boundary subsurfaces.
When $\alpha=\add{\varnothing}$, the outer boundaries of $\XX(\alpha)$
are the outer boundaries of $\modl$.
The boundary of $\XX(\alpha)$ consists of 
these inner and outer boundary surfaces together with annuli and tori in
$\boundary\UU$.

\begin{proof}
We first note that if $h\in H$, $d,d'\in C|_h$ and $(c,c')$ is a pair
such that $c\wallto d \wallto c'$, then the argument in the proof of Lemma \ref{Addresses unique} implies that $c \wallto d' \towall c'$ as well. 

It follows that if $h$ is an inner boundary geodesic for $\alpha$ and
$B$ is a block in $\BB(h)$, then $\alpha$ is an address pair for $B$
but it is not an innermost address pair. Therefore,  $\BB(h)\subset \BB(\alpha)$
and $int(\BB(h))\cap \XX(\alpha)=\emptyset$, establishing (1).

If $h$ and $h'$ are inner boundary geodesics for
$\alpha$, a nonempty 
intersection of $\BB(h)$ and $\BB(g)$ implies that, for some pairs
$a,a'\in C|_h$ and $b,b'\in C|_g$, there
is a block for which $(a,a')$ and $(b,b')$ are both address
pairs. However as in the proof of lemma \ref{Addresses unique} this
implies that 
one of $D(g)$ and $D(h)$ must be strictly contained in the other, and
this then implies that one of $(a,a')$ and $(b,b')$ is nested in the
other, which contradicts the definition of an inner boundary
geodesic. This establishes (2).

(3) then follows from (1) and the definition of $\XX(\alpha)$.

To show the last statement, notice that if $h\in H$
and $C|_h$ is non-empty, then either there is no $(c,c')$ 
such that $c\wallto d \towall c'$ for any
$d\in C|_h$, 
or there is a nonempty collection $\DD$ of pairs $(c,c')$ for which
$c\wallto d \towall c'$ for {\em all} $d\in C|_h$. 
In the first case, $h$ must be an inner boundary geodesic for
$\alpha=\add\varnothing$. 
In the second case it follows, 
as in the proof of Lemma \ref{Addresses unique}, that
there is a unique innermost $(c,c')$ in $\DD$ and this must be
the unique address $\alpha$ for which $h$ is an inner boundary geodesic. 
\end{proof}

With this picture in mind, we can relate addresses to connected
components of the complement of the surfaces of $C$:

\begin{lemma}{region address}
All blocks in a complementary region of $C$ have the same address.
\end{lemma}

\begin{proof}
Any two blocks $B$ and $B'$ in the same connected component
are connected by
a chain $B=B_0,\ldots,B_n=B'$ where $B_i$ and $B_{i+1}$ are adjacent
along a gluing surface which does not lie in any of the cut surfaces
$\{ F_c:c\in C\}$. It thus suffices
to consider the case that $B $ and $B'$ are adjacent along gluing
surfaces that are not in the cuts.

Let us show that any address pair 
$(c,c')$ for $B$ is also an address pair for $B'$ (and, by symmetry,
vice versa). This will imply that 
the innermost pairs, and hence the addresses, are the same.

The region $\BB(c,c')$ contains $B$. Since $\boundary\BB(c,c')$
consists of $F_c$, $F_c'$ and portions of the boundaries of tubes, the
gluing surface connecting $B$ 
to $B'$ is not in this boundary. It follows that $B'$ is also
contained in $\BB(c,c')$, and hence $(c,c')$ is an address pair for
$B'$.
This completes the proof.
\end{proof}

\subsection{Sizes of regions}
\label{region size}

Our next lemma will bound the
number of blocks in any $\XX(\alpha)$. As an immediate
consequence of Lemma \ref{region address}, we also get a bound on the
size of any complementary region of the cut system.

\begin{lemma}{Bound regions}
The number of blocks in $\XX(\alpha)$ for any address $\alpha$ is
bounded by a constant $K$ depending only on $S$ and $d_2$.
\end{lemma}

\begin{proof}
Fix an address $\alpha$. If $\alpha=\add{c,c'}$,
let $g_\alpha=g_c=g_{c'}$ be the
bottom geodesic for $c$ and $c'$. If $\alpha=\add\varnothing$,
let $g_\alpha = g_H$.

Let $\ZZ=\ZZ_\alpha$ be the set of all three-holed spheres $Y$ such
that $F_Y$ is a component of $\boundary_- B$ for
some internal block $B$ in $\XX(\alpha)$.
Since every internal block $B$ has nonempty $\boundary_- B$
and every $F_Y$ is in the $\boundary_-$ gluing boundary of at
most one block, it follows that the number of internal blocks in
$\XX(\alpha)$ is at most $|\ZZ|$. Since there is a bound on the
number of boundary blocks depending only on $S$, 
it will suffice to find a bound on $|\ZZ|$ that depends only
on $S$ and $d_2$. 

For a geodesic $h$, 
let $\lands(h)$ be the set of all three-holed spheres $Y$ 
for which $Y \fsub h$.
For each geodesic $h \fsubeq g_\alpha$, we define
$$J_\ZZ(h) = \{ v\in h \; : \; v = \max \phi_h(Y) \ \text{for some} \
Y \in \ZZ, Y \subset D(h) \},$$ 
the set of landing points  on $h$ of forward sequences for $Y \in \ZZ$.

\newcommand\Jbound{m}

The bound on $|\ZZ|$ will follow from the following four claims:
\begin{enumerate}
\item $\ZZ \subset \lands(g_\alpha)$
\item If $h\fsubeq g_\alpha$ and $\xi(h)>4$ then
$$
\ZZ \cap \lands(h) = \bigsqcup_{\stackrel{k \fsubd h}{\max \phi_h(D(k)) \in
J_\ZZ(h)}} 
\left( \ZZ \cap \lands(k) \right)
$$
\item If $h\fsubeq g_\alpha$ and $\xi(h)=4$ then
$$
|\ZZ\cap \lands(h)| \le 2|J_\ZZ(h)|
$$
\item For any $h\fsubeq g_\alpha$,
$$|J_\ZZ(h)| \le \Jbound,$$
where $\Jbound$ is a constant depending only on $d_2$.
\end{enumerate}
Assuming these claims hold, we can
prove a bound  $|\ZZ\intersect\lands(h)|\le
K_{\xi(h)}$ by 
induction on $\xi(h)$. For $\xi(h)=4$, we obtain the bound
$K_4=2\Jbound$ by claims (3) and (4). In the induction step suppose we
already have a bound $K_{\xi(h)-1}$. For any interior simplex
$v\in h$ there is a unique
geodesic $k\fsubd h$ with $\max \phi_h(D(k))=v$ ($D(k)$ can only be
the component domain of $(D(h),v)$ that intersects the successor of $v$).
If $v$ is the first or last simplex of $h$ then it is a vertex
so there are at most two non-annular component domains of $(D(h),v)$
and hence at most two non-annular $k\fsubd h$ with $\max\phi_h(D(k)) = v$.
Thus the union in claim (2) has at most
$|J_\ZZ(h)|+2$  terms. Together with 
claim (4) we obtain
a bound of $K_{\xi(h)} = K_{\xi(h)-1}(\Jbound+2)$. Claim (1) then
gives us our desired bound $|\ZZ| \le K_{\xi(g_\alpha)}$.

Before proving the claims, we note the following facts: 
If $Y$ and $Y'$ are three-holed spheres in the hierarchy and are both
contained in $D(h)$ for some $h$, then
\begin{equation}
\label{footprints tprec}
\max\phi_h(Y) < \min \phi_h(Y')  \implies
Y \tprec Y'.
\end{equation}
This follows directly from the definition of $\tprec$, noting that 
if $\max\phi_h(Y) < \min \phi_h(Y') $ then in particular the last
vertex of $h$ is not in $\phi_h(Y)$ and the first is not in
$\phi_h(Y')$, so that $Y\fsub h \bsub Y'$.

From the contrapositive, with $Y$ and $Y'$ interchanged, we obtain
\begin{equation}\label{tprec footprints}
Y\tprec Y' \implies \min\phi_h(Y) \le \max\phi_h(Y').
\end{equation}

Now we prove claim (1). If 
$\alpha=\add\varnothing$ then $g_\alpha=g_H$ and (1) is
immediate since $\YY(g_H)$ contains all three-holed spheres
which are component domains of the
hierarchy except those that are component domains of $\T(H)$,
and those are excluded from $\ZZ$ (they correspond to $\boundary_-$ gluing
surfaces of boundary blocks).

Assume $\alpha=\add{c,c'}$. Let $B$ be a block in $\XX(\alpha)$
and $F_Y$ a component of $\boundary_- B$.
Since $W_B \topprec \hhat F_{c'}$ and the interior of
$B$ is disjoint from $\hhat F_{c'}$, we must have
$F_Y\topprec \hhat F_{c'}$.

Let 
$Y'$ be a component domain of $(D(g_\alpha),\base(c'))$ which overlaps $Y$.
It follows (applying Lemma \ref{basic ordering properties}) that
$F_Y\topprec F_{Y'}$.   Lemma \ref{Y ordering} tells us that
$$ 
Y \tprec Y'.
$$
Thus by (\ref{tprec footprints}) we have
$\min\phi_{g_\alpha}(Y) \le \max\phi_{g_\alpha}(Y')$. 
(These footprints are nonempty by Lemma \ref{nonempty footprint}.)
Letting
$(g_\alpha,v')$ denote the bottom pair of $c'$, 
since $Y'$ is a complementary component of $\base(c')$ it follows 
that $\phi_{g_\alpha}(Y')$ contains $v'$.
Since the lower
spacing bound $d_1$ for the cut system is at least 5, and footprints
have diameters at most 2, it follows that $\max\phi_{g_\alpha}(Y')$ is at least
3 away from the last simplex of $g_\alpha$, and hence
$\max\phi_{g_\alpha}(Y)$ is at least 1 away. Thus $Y\fsub g_\alpha$,
or $Y\in\lands(g_\alpha)$. This establishes claim (1).

This discussion also proves claim (4) for $h=g_\alpha$ when
$\alpha=\add{c,c'}$: if 
$(g_\alpha,v)$ is the bottom pair of $c$ then $v$ and $v'$ are at most
$d_2$ apart. The above argument shows that $\max\phi_{g_\alpha}(Y)$ is
at most 4 forward of $v'$, and the same argument run in the opposite
order (with $c$ replacing $c'$) shows that $\max\phi_{g_\alpha}(Y)$
occurs no further back than 2 behind $v$. This restricts it to an
interval of diameter $d_2+6$, which gives us claim (4) for $g_\alpha$
provided $\Jbound > d_2+6$.

Now consider claim (4) for $h\fsub g_\alpha$, or for $h= g_\alpha$
when $\alpha=\add\varnothing$. We claim that if $Y\fsub
h$ and $Y\in\ZZ$ then $\max \phi_h(Y)$ occurs within $d_2 + d_1/2 + 3$ of the
endpoints of $h$. Suppose this is not the case and let us look for a
contradiction. The length of $h$ is then greater than $2d_2
+ d_1$ (possibly it is infinite), which means that there must be at
least two slices of $C$ based on $h$.  
There exist slices $d,d'\in C|_h$  whose bottom simplices $v,v'$
satisfy $v< \max\phi_h(Y) < v'$ and are at least $d_1/2+3$ away from
$\max\phi_h(Y)$: These can be the first and last slices of $h$, if
these exist, since they are within $d_2$ of the endpoints of $h$; or 
if $h$ is infinite in the backward or forward direction a sufficiently
far away slice will do for $d$ or $d'$ respectively. Note that $d_1/2+3>5$.
For any component $Y'$ of $D(h)\setminus\collar(\base(d'))$,
$\phi_h(Y')$ contains $v'$ and it follows that $\max \phi_h(Y) <
\min\phi_h(Y')$ so that $Y\tprec Y'$ by (\ref{footprints
tprec}) and hence $F_Y \topprec F_{Y'}$ by Lemma \ref{Y ordering}.
It follows, as in the proof of Proposition \ref{topprec and cprec},
that $F_Y \topprec \hhat  F_{d'}$. 
A similar argument yields $\hhat F_d \topprec F_Y$.

Now let $B$ be a block in $\XX(\alpha)$ with $F_Y\subset \boundary_\pm B$.
By Lemma \ref{block slice order}, $W_B$ and $\hhat F_d$ must be
$\topprec$-ordered, and similarly for $W_B$ and $\hhat F_{d'}$. 
Since 
the interior of $B$ does not meet $\hhat F_d$ or $\hhat F_{d'}$,
the ordering we've established for $F_Y$ implies that 
 $\hhat F_d \topprec W_B \topprec F_{d'}$.

We also note that $D(B)\subset D(h)$, as follows. The block $B$ is
associated to an edge in a geodesic $k$ with $D(k)=D(B)$.  Assume
without loss of generality that $Y\subset \boundary_-B$. 
Thus if $e$ is the edge of $k$ defining $B$ we have that $Y$ is
component domain of $(D(h),e^-)$ which intersects $e^+$, and
in particular $Y\fsubd k$. We also have
$Y\fsub h$, since $\phi_h(Y)$ is far from the ends of $h$. Hence $h$
and $k$ are in the forward sequence $\Sigma^+(Y)$, so one is contained
in the other. Since $\xi(k)=4$, we must have $D(k)\subset D(h)$.

We can therefore conclude that $(d,d')$ is an address pair for $B$. 
If $\alpha=\add\varnothing$ then this is already a contradiction.
If not, then since the domain of $d$ is strictly smaller than that of
$c$ we must have $(d,d')$ nested within $(c,c')$ by (the proof of)
Lemma \ref{Addresses unique}, a contradiction to the assumption that
$(c,c')$ is the innermost address pair for $B$. 

This contradiction establishes our claim, so that $\max\phi_h(Y)$ is
confined to a pair of intervals of total length $2d_2 + d_1 + 6\le
3d_2+6$. This gives the desired bound on $|J_\ZZ(h)|$ for Claim (4).

For Claim (2), let $h\fsubeq g_\alpha$ with $\xi(h) >4$.
Suppose that $Y\in\ZZ\intersect \lands(h)$. Then $Y\fsub h$ but 
we cannot have $Y\fsubd h$ since $\xi(h)>4$. Thus
there is a $k\in\Sigma^+(Y)$
such that $k\fsubd h$, and hence
(by Lemma 5.5 of \cite{minsky:ELCI})
$\max\phi_h(D(k)) = \max\phi_h(Y)$.
In particular $\max\phi_h(D(k)) \in J_\ZZ(h)$. Since $Y\fsub k$ 
we also have $Y\subset \ZZ\cap \lands(k)$. Note that $Y$
cannot be in $\ZZ\cap \lands(k')$ for a different $k'\fsubd h$
by the uniqueness of the forward sequence $\Sigma^+(Y)$.
Thus we obtain the partition of $\ZZ\cap \lands(h)$ described in Claim
(2).

For Claim (3), let  $h\fsubeq g_\alpha$ with $\xi(h)=4$.
Now $Y\fsub h$ exactly if $Y\fsubd h$, and this occurs when $Y$ is a
complementary component of $D(h)\setminus\collar(v)$ for
$v=\max\phi_h(Y)$. There is one such component for each $v$ when
$D(h)$ is a one-holed torus, and two when $D(h)$ is a four-holed
sphere. The inequality of claim (3) follows. 

Thus we have established the bound $|\ZZ|\le K_{\xi(S)}$, where $K_{\xi(S)}$
depends only on $S$ and $d_2$.
\end{proof}

\subsection{Filled regions}
\label{filled regions}

We will also need to consider regions determined by a cut system
$C$ in the {\em filled model} $\modl[k]$ for some constant $k >0$.    
If $C$ is a cut system, then the surfaces $\{ \hhat F_\tau[k] : \tau \in
 C\}$ again decompose the model $\modl[k]$ into regions.  We
wish to verify that if the lower space bound is chosen large enough, then
these regions in a filled model differ from the
regions determined by $\{ F_\tau : \tau \in  C\}$ only by
filling in certain tubes whose boundaries lie entirely in a given
region. More precisely,
let 
$$\WW_i = \modl[i] \setminus \bigcup_{c \in C} \hhat F_c[i].$$
Thus the components of $\WW_0$ are the the complementary regions in
$\modl[0]$ of the cut system, which we have been considering 
up til now. 

\begin{proposition}{modl[k] regions}
Given $k>0$ there is a constant $d_1\ge 5$ such that, 
if $C$ is a cut system with a spacing lower bound of at least $d_1$, 
then the connected components of $\WW_0$ are
are precisely the connected components of  $\WW_k$
minus the tubes of size $|\omega|<k$. 

In particular, all blocks in a connected component of 
$\WW_k$
have the same address. 
\end{proposition}

The main step in the proof of Proposition \ref{modl[k] regions} is the
following lemma, which shows that if $d_1$ is chosen large enough and $U$ is a tube in $\modl[k]$,
then $U$ meets at most one split-level surface associated to a non-annular
slice in $C$.

\begin{lemma}{tube filling}
Given $k$ there exists $d_1\ge5$ so that for any cut
system $C$ with spacing lower bound of at least $d_1$ and each tube $\storus v$
in $\modl[k]$, 
there is at most one nonannular $a\in C$ such that 
$\boundary \storus v$ meets $F_a$.
\end{lemma}

\begin{proof}
Let $v$ be a vertex in $H$ so that $|\omega(v)| <k$, and hence
$\storus v \subset \modl[k]$.
Suppose $\boundary \storus  v$ meets a cut surface $F_a$ for some
$a\in C$. 
This implies that $v\in[\boundary\check D(a)]$,
so either
\begin{enumerate}
\item $v \in [\bdry D(a)]$, or
\item $v \in \base(a)$.
\end{enumerate}

The lower spacing bound on $C$ means that  the
bottom geodesic $h_a$ has length at least $3d_1$, 
so if $v \in [\bdry D(a)]$ this yields a lower bound on $|\omega(v)|$.
In particular letting $b_1$ and $b_2$ be the constants in Lemma 
\ref{big h big omega}, 
if we have chosen $d_1\ge (k + b_1)/3b_2$, Lemma 
\ref{big h big omega} would imply $|\omega(v)|\ge k$
which is a contradiction, 
and hence case (1) cannot occur.

Now suppose that there are two  slices $a,b\in C$ such that
$\boundary \storus v$ meets $F_a$ and $F_b$, and hence that
$v\in \base(a)$ and $v\in \base (b)$.  This 
possibility is ruled out by Lemma \ref{cuts dont meet},
and this completes the proof of Lemma \ref{tube filling}. 
\end{proof}

We can now complete the proof of the Proposition:

\begin{proof}[Proof of Proposition \ref{modl[k] regions}]
Choose $d_1$ to be the constant  given by Lemma \ref{tube filling}.
$\WW_k$ is obtained from
$\WW_0$ by attaching, for each tube $U$ with $|\omega(U)| < k$, the
set 
$$ U \setminus \bigcup_{c\in C} \hhat F_c.$$
By Lemma~\ref{tube filling}, 
$U$ meets at most one surface $\hhat F_c$ with $c\in C$ non-annular, and if it does
so then the intersection is a single annulus. Thus each component of
$ U \setminus \bigcup_{c\in C} \hhat F_c$ is a solid torus which 
either meets $\boundary U$  in the entire boundary,
or in a single annulus. In either case each component meets
$\boundary U$ in a connected set. This means that the components
of $U \setminus \bigcup_{c\in C} \hhat F_c$ cannot connect different
components of $\WW_0$.

It follows that a connected component of $\WW_k$ is equal to a
connected component of $\WW_0$ union the adjacent pieces of tubes.

The final statement of the proposition is an immediate consequence of
Lemma \ref{region address}.

\end{proof}

\section{Uniform embeddings of Lipschitz surfaces}
\label{uniform embeddings}

The main result of this section is 
Theorem \ref{Relative embeddability}, which proves
that a Lipschitz map of a surface with
bounded geometry into the manifold $N_\rho$ can be deformed to an
embedding in a controlled way, provided it satisfy a number of
conditions, the most important being an ``unwrapping condition'' that
rules out the possibility that the homotopy will be forced to go through
a deep Margulis tube.

We begin by introducing a series of definitions which allow us to
describe the type of surfaces we allow and to express 
what it means to deform to an embedding in a controlled way. 

A compact hyperbolic surface $X$ (possibly disconnected) with
geodesic boundary is said to be {\em $L$-bounded} (or has a
$L$-bounded metric) if
no homotopically non-trivial curve in $X$ has length
less than $1/L$ and  every boundary
component has length in  $[1/L,L]$.
A map $f:X\to N$ is $L$-bounded if $X$ is $L$-bounded and
$f$ is $L$-Lipschitz.

An {\em anchored surface} (or map) is a map of pairs
$$
f:(X,\boundary X)\to (N\setminus \MT({\boundary X}),\boundary
\MT({\boundary X}))
$$
where $X\subseteq S$ is an essential subsurface and
$f$ is in the homotopy class determined by $\rho$.
An anchored surface is {\em $\epsilon$-anchored} if
$\ell_\rho(\gamma)< \ep$ for each component
$\gamma$ of $\boundary X$.

\medskip

If $X$ has a hyperbolic metric,
an anchored surface $f:X\to N$ is {\em $(K,\hat\ep)$-uniformly embeddable}
if there exists a homotopy, called a {\em $(K,\hat \epsilon)$-uniform homotopy},
\begin{equation}\label{nice homotopy}
H:(X\times [0,1],\boundary X\times [0,1]) \to
(N\setminus \MT({\boundary X}),\boundary
\MT({\boundary X}))
\end{equation}
with $H(\cdot,0) = f$ such that
\begin{itemize}
\item $H$ is $K$-Lipschitz
\item $H$ restricted to $X\times[1/2,1]$ is a $K$-bilipschitz
$C^2$ embedding with the norm of the second derivatives bounded by $K$.
\item
For all $t\in [1/2,1]$, $H(\boundary X\times \{t\})$ is a collection of 
geodesic circles in $\boundary\MT(\boundary X)$.
\item $H(X\times[1/2,1])$ avoids all $\epone$-Margulis tubes
with core length less than $\hat\ep$.
\item $H(X\times[0,1])$ avoids all $\hat\ep$-Margulis tubes.
\end{itemize}

\medskip

An $L$-bounded map which is homotopic to an embedding may not be
uniformly embeddable (with constants depending only on $L$), due to an
obstruction called wrapping (see also \S\ref{wrapping coeff}): 
In Anderson-Canary \cite{anderson-canary:pages} and McMullen \cite{mcmullen:grafting},
there is a construction, 
with fixed $L$, of a sequence of manifolds $N_{\rho_n}$
with Margulis tubes $\MT(\alpha,n)$ of depth going to $\infty$, and
immersed $L$-bounded surfaces ``wrapped'' around these tubes. Each such
surface is homotopic to an embedding but the homotopy is forced to
pass through the core of $\MT(\alpha,n)$, and hence there is no fixed
$(K,\ep)$ such that these surfaces are $(K,\ep)$-uniformly
embeddable. 

In view of this obstruction, 
we say that a map $f:X \to N^1$ is {\em unwrapped} with respect to a curve
$\alpha\in\CC_0(X)$ if $\ell_\rho(\alpha)<\ep_1$ and $f$ is homotopic to either $+\infty$ or
$-\infty$ in $N^1\setminus \MT(\alpha)$.  

We recall that $\ep_1$ is our fixed choice of Margulis constant, and  $\MT(\alpha)$ denotes
a $\ep_1$-Margulis tube. Moreover, $\ep_0$ and $\epotal$ are also Margulis
constants, $\ep_0>\ep_1$ and $\epotal$ is chosen so that the collection of
geodesics in $N$ of length at most $\epotal$ is unknotted and unlinked
(see section \ref{tubes}). The manifold
$N_\rho^1=N_\rho\setminus\MT(\boundary S)$,
$C_{N_\rho}$ is the convex core of $N_\rho$ and $C_{N_\rho}^r$ is its closed neighborhood
of radius $r$.

When a surface is unwrapped we can show that it {\em is }
uniformly embeddable, provided it is anchored on sufficiently deep
tubes. More precisely:

\begin{proposition}{Easy embeddability}
Let $S$ be an oriented compact surface. Given $L>1$ and $\delta<\epotal$ there
exist $\ep, \hat\ep>0$ and $K>1$ such that the following holds: 

Let $\rho\in \DD(S)$, $R\subset S$ an essential subsurface. Suppose
that
$$
f:R \to C_{N_\rho}\intersect N^1_\rho
$$
is $\ep$-anchored, and is unwrapped with respect to $\alpha\in \CC(R)$
whenever $\ell_\rho(\alpha) < \delta$. 

Then $f$ is $(K,\hat\ep)$-uniformly embeddable.
\end{proposition}

This proposition is a special case of Theorem \ref{Relative
  embeddability} which we will prove below. For motivation, let us
sketch the proof of the special case.

Supposing that the statement is false, we may fix $L$ and $\delta$ and
find a sequence of $L$-bounded $\ep_n$-anchored maps with $\ep_n\to 0$ which are
unwrapped but not uniformly embeddable. Restricting to a subsequence
and remarking the maps we may assume the domain $R$ is fixed, 
and extract a geometric limit $N$ of the target manifolds and a
limiting map $f:R\to N$. The boundary curves of $R$ map to cusps in
$N$, but there may be additional curves in $R$ whose images are
parabolic. Let $P$ be a maximal collection of disjoint curves in $R$
mapping to parabolics. Using the techniques
of Anderson-Canary \cite{anderson-canary:cores} (see Section \ref{embedding machine}), the restriction of $f$
to $R\setminus \collar(P)$ can be homotoped to an anchored embedding $h$ (where the
anchoring tubes are cusps). 

By adding annuli in the boundaries of the cusps associated to $P$, we can extend $h$ to an
embedding of all of $R$. The unwrapping condition with respect to the
curves of $P$ implies that this can be done in such a way
that the result is homotopic to $f$ in $N$ (\S\ref{resewing along parabolics}). This homotopy now pulls
back to give a uniform sequence of homotopies in the approximates,
yielding a contradiction(\S\ref{contradiction}).

\medskip

In the proof of the Bilipschitz Model Theorem, we will need to prove
uniform embeddability for the Lipschitz model map restricted to 
certain extended split-level surfaces $\hat F_\tau[k]$
associated to cuts $\tau$ in
the model, for some uniform choice of $k$ (see Section \ref{split
  level}). One can ensure that the (images of) boundary curves 
of the base surface $\hat F_\tau$ are short, by requiring that
the base geodesic of $\tau$ be long. 

Choosing $k$ large would guarantee that the 
boundary curves of $\hat F_\tau[k]$ which are interior to $\hat
F_\tau$ are shorter than any desired $\ep$, but at the price of
including curves of length close to $\ep$  in the interior, thus
degrading the boundedness of the surface. We therefore
cannot use Proposition \ref{Easy embeddability} directly to uniformly
embed such maps. We will need to establish the following more
complicated statement:

\begin{theorem}{Relative embeddability}
Let $S$ be an oriented compact surface. 
Given $\delta, L>0$, 
there exist $\ep,\hat\ep\in (0,\epotal)$ 
and $K>1$, so that the following holds.  

Let $\rho\in\DD(S)$, $R\subset S$ an essential subsurface, and
$\Gamma$ a simplex in $\CC(R)$, and let 
$X = R \setminus \collar(\Gamma)$.
Suppose that 
there is an $\epotal$-anchored, $L$-bounded
surface
$$f : X \to C_{N_\rho}\intersect N_\rho^1$$
in the homotopy class determined by $\rho$, and there exists
an extension
\newline
$\bbar f:R\to N_\rho$  of $f$ such that
\begin{enumerate}
\item
$\bbar f$ takes $\collar(\Gamma)$ into
$\MT(\Gamma)$,
\item 
$\bbar f$ is unwrapped with respect to 
any $\alpha\in\CC_0(R)\setminus \Gamma$ for which $\ell_\rho(\alpha)<\delta$.
\item $\bbar f$ is $\ep$-anchored.
\end{enumerate}
Then $f$ is $(K,\hat\ep)$-uniformly embeddable and the uniform homotopy
$H$ has image in $C^{1/2}_{N_\rho}\union (N_{\rho})_{(0,\ep_2]}$ where
$\ep_2=(\ep_0+\ep_1)/2$.
\end{theorem}

Note that the length bound $\epotal$ for the internal curves $\Gamma$ in $R$ is
fixed ahead of time together with $L$, whereas the bound $\ep$ for the
boundary curves of $R$ depends on $L$. A key new difficulty is 
that, since the output is a uniform
embedding for the map $f$, not the map $\bbar f$,  we must
anchor on the internal tubes $\MT(\Gamma)$ which need  not be extremely deep. 

In the argument by contradiction, such tubes may not become cusps in
the geometric limit, and the machinery of
Section \ref{embedding machine} will yield an anchored embedding of
$R$ and not $X$.  Thus we will need to ``re-anchor'' the embeddings on $\Gamma$.

After developing the machinery for embeddings in geometric limits, we
will return in \S \ref{embedding proof} to the proof of Theorem
\ref{Relative embeddability}. At that point we will also give a more
detailed outline of the rest of the proof. 

\subsection{Embedding in geometric limits}
\label{embedding machine}

In this section we show that given a sequence $\{\rho_n\}$
of representations which converge on a subsurface $F$ of $S$ so that
$\{\ell_{\rho_n}(\partial F)\}$ converges to 0 and the limits
have no non-peripheral parabolics, we can produce an anchored 
embedding of $F$ into the geometric limit of $\{N_{\rho_n}\}$.

Let $F$ be a (possibly disconnected) subsurface of $S$ which has no
annulus components.
Note that given a component $F_i$ of $F$ there is a family of homomorphisms
$\sigma:\pi_1(F_i)\to \pi_1(S)$ consistent with the inclusion map
(depending on choice of basepoints and arcs connecting them), and
any two of these differ by conjugation in $\pi_1(S)$.

\begin{definition}{convergence on subsurface}
A sequence $\{\rho_n\}$ in $\DD(S)$ is
{\em convergent on $F$} if,  for each $i$ 
there is a sequence $\{\sigma_n^i:\pi_1(F_i)\to\pi_1(S)\}$ consistent
with the inclusion map so that the sequence of representations
$$
\rho_n^i   = \rho_n\circ\sigma_n^i
$$
converges to a representation $\rho^i:\pi_1(F_i) \to \PSL 2(\C)$.
\end{definition}

We call the $\rho^i$ {\em limit representations on $F$} of
$\{\rho_n\}$ (but note that they depend on the choice of $\sigma_n^i$).

\begin{proposition}{general subsurface embed}
Let $S$ be an orientable surface and let $F$ be an essential
subsurface with components $\{F_i\}$, none of which are annuli.
Suppose that $\{\rho_n\}$ is a sequence in ${\mathcal D}(\pi_1(S))$
such that 
\begin{enumerate}
\item
$\{\rho_n\}$ is convergent on $F$ with limit representations
$\rho^i\in\DD(\pi_1(F_i))$
\item
$\rho^i(g)$ is parabolic if and only if $g$ is peripheral in
$\pi_1(F_i)$, and
\item
$\{\rho_n(\pi_1(S)\}$ converges geometrically to $\Gamma$.
\end{enumerate}
Then letting  $\hhat N=\Hyp^3/\Gamma$, there exists 
an anchored embedding
$$
h:F \to \hhat N
$$
such that $h|_{F_i}$ is in the homotopy class determined by $\rho^i$ for each $F_i$.
\end{proposition}

In the proof of Proposition \ref{general subsurface embed}
we will need to consider separately
the components $F_i$ for which $\rho^i(\pi_1(F_i))$ is geometrically
finite, and those for which it is geometrically infinite. 
The geometrically finite subsurfaces will be handled  using (relative
versions of) machinery
developed by Anderson-Canary-Culler-Shalen \cite{anderson-canary-culler-shalen}
and Anderson-Canary \cite{anderson-canary:cores}, while the geometrically
infinite subsurfaces will be handled using Thurston's covering theorem.

\subsubsection{The limit set machine}

We first establish Proposition \ref{general subsurface embed} 
when the algebraic limits are quasifuchsian:

\begin{proposition}{limit embed quasifuchsian}
Let $S$ be an orientable surface and let $F$ be an essential
subsurface with components $\{F_i\}$, none of which are annuli.
Suppose that $\{\rho_n\}$ is a sequence in ${\mathcal D}(\pi_1(S))$
such that 
\begin{enumerate}
\item
$\{\rho_n\}$ is convergent on $F$ with limit representations
$\rho^i\in\DD(\pi_1(F_i))$
\item
$\rho^i$ is a quasifuchsian representation of $F_i$ for all $i$, and
\item
$\{\rho_n(\pi_1(S)\}$ converges geometrically to $\Gamma$.
\end{enumerate}
Then letting  $\hhat N=\Hyp^3/\Gamma$, there exists 
an anchored embedding
$$
h: F \to \hhat N
$$
such that $h|_{F_i}$ is in the homotopy class determined by $\rho^i$ for each $F_i$.
\end{proposition}

Let us give an outline of the proof of
Proposition \ref{limit embed quasifuchsian}. 
(The actual proof proceeds in the opposite order to the outline.)
Since limit sets of quasifuchsian groups are Jordan curves,
any essential intersection of the maps of $F_i$ and $F_j$ into the
geometric limits  associated to the representations $\rho^i$ and $\rho^j$
(or essential self-intersection
of the map of $F_i$) would result in  limit sets of
conjugates of $\rho^i(\pi_1(F_i))$ and $\rho^j(\pi_1(F_j))$  which
cross (see Lemma \ref{precisely embed}). A result of
Susskind \cite{susskind:limitsetints}, see Theorem \ref{limitsetints},
implies that the intersection of the limit sets of two geometrically
finite subgroups $\Phi_1$ and $\Phi_2$ of a Kleinian group
consists of the limit set of their intersection $\Phi_1\cap \Phi_2$ along
with certain parabolic fixed points $P(\Phi_1,\Phi_2)$. Therefore,
it suffices to prove that the intersection of $\rho^i(\pi_1(F_i))$ and
a distinct conjugate of $\rho^j(\pi_1(F_j))$ has at most one point
in its limit set (Lemma \ref{purelypar}) and that there are
no problematic parabolic fixed points (see Proposition \ref{onepoint}.)

We first recall Susskind's result which describes the intersections of
the limit sets of two geometrically finite subgroups of a Kleinian group.
(Soma \cite{soma:function} and Anderson \cite{anderson:intersection}
have generalized Susskind's result to allow the subgroups to be
topologically tame.)
Given a pair $\Theta$ and $\Theta'$ of subgroups
of a Kleinian group $\Gamma$, let $P(\Theta,\Theta')$ be the set of
points $x\in\Lambda(\Gamma)$ such that the stabilizers of $x$
in $\Theta$ and $\Theta'$ are rank one parabolic subgroups which
generate a rank two parabolic subgroup of $\Gamma$.

\begin{theorem}{limitsetints}{(Susskind \cite{susskind:limitsetints})}{}
Let $\Gamma$ be a Kleinian group and let $\Phi_1$ and $\Phi_2$ be
nonelementary, geometrically finite subgroups of $\Gamma$.  Then,
$$ \Lambda(\Phi_1)\cap \Lambda(\Phi_2) =
\Lambda(\Phi_1\cap\Phi_2)\union P(\Phi_1,\Phi_2). $$
\end{theorem}

We next show that the intersection of $\rho^i(\pi_1(F_i))$ and
a distinct conjugate of $\rho^j(\pi_1(F_j))$ has at most one point
in its limit set.  This generalizes Lemma 2.4 from
\cite{anderson-canary-culler-shalen}, in the setting of surface groups.

\begin{lemma}{purelypar}
Let $\{\rho_n\}$ be a sequence in $\DD(S)$ which is convergent on an
essential subsurface $F$, with nonannular components $F_i$ and
limit representations $\rho^i$. Suppose that $\{\rho_n(\pi_1(S))\}$
converges geometrically to $\Gamma$, and that  
$\{\ell_{\rho_n}(\boundary F)\}$ converges to $0$.

If $\gamma\in \Gamma$ and either $i\ne j$ or
$i=j$ and $\gamma\notin \rho^i(\pi_1(F_i))$, then
$$\gamma \rho^i(\pi_1(F_i)) \gamma^{-1} \intersect \rho^j(\pi_1(F_j))$$ 
is purely parabolic.
\end{lemma}

\begin{proof}
Let $\{\sigma_n^i\}$ be the sequences of maps,
as in Definition \ref{convergence on subsurface}, such that
$\{\rho_n^i\}=\{\rho_n\circ\sigma_n^i\}$  converges to $\rho^i$ for each $i$.

Let $\gamma \in \Gamma$ and suppose that $\{\rho_n(h_n)\}$
converges to $\gamma$. Suppose that
$\alpha \in \gamma \rho^i(\pi_1(F_i)) \gamma^{-1} \intersect \rho^j(\pi_1(F_j))$
is non-trivial. Then there exist nontrivial $a\in\pi_1(F_i)$
and $b\in \pi_1(F_j)$ such that
\begin{equation}\label{conjugacy in limit}
\rho^j(b) = \gamma\rho^i(a)\gamma^{-1} = \alpha.
\end{equation}
Our goal is to prove that $\alpha$ must be parabolic.

Since $\{\rho_n(\sigma_n^j(b))\}$ and 
$\{\rho_n(h_n \sigma_n^i(a) h_n^{-1})\}$ both converge to $\alpha$,
Proposition \ref{JMrel} (part 1) implies that 
$$ 
\sigma_n^j(b) = h_n \sigma_n^i(a) h_n^{-1}
$$
for all sufficiently large $n$.

In particular $a$ and $b$ represent the same free homotopy class in
$S$. If $i\ne j$, 
$a$ and $b$ must represent boundary components of $F_i$ and $F_j$ that
are freely homotopic to each other, and since we have assumed 
$\{\ell_{\rho_n}(\boundary F)\}$ converge to $0$, we conclude that 
$\alpha$ is parabolic.

If $i=j$, we may re-mark  the sequence $\{\rho_n|_{\pi_1(F_i)}\}$
by precomposing with $\sigma_n^i$, 
so that from now on we may fix an inclusion of $\pi_1(F_i)$ in
$\pi_1(S)$, and set
$\sigma_n^i = id$. After dropping finitely many terms from the
sequence, we have
\begin{equation}\label{nth conjugacy}
b = h_n a h_n^{-1}
\end{equation}
for all $n$. If $h_n\notin \pi_1(F_i)$ for some $n$, then
$a$ and $b$ must represent boundary components of $F_i$
that are homotopic in the complement of $F_i$. Again,
since $\{\ell_{\rho_n}(\boundary F)\}$ converges to $0$, 
we may conclude that $\alpha$ is parabolic.

Thus, suppose that $h_n\in\pi_1(F_i)$ for all $n$.
We will show that, if $\alpha$ is not parabolic, then
$h_n$ is eventually constant. Equation (\ref{nth
conjugacy}) implies that $h_m  h_n^{-1}$ centralizes $b$
for all $m,n$.  Letting $m=1$, applying $\rho_n$, 
and taking a limit as $n\to \infty$, 
we find that 
$\rho^i(h_1)\gamma^{-1} $ centralizes $\rho^i(b)$. 
Since we are assuming that $\rho^i(b)$ is hyperbolic, 
its centralizer in $\Gamma$ is infinite cyclic, 
so there exist non-zero integers
$k$ and $l$ such that
$$(\rho^i(h_1)\gamma^{-1})^k =\rho^i(b)^l.$$
Again Proposition \ref{JMrel} (part 1) implies that 
$(h_1h_n^{-1})^k=b^l$ for all large enough $n$. Since elements of
torsion-free
Kleinian groups have unique roots, we conclude that $\{h_1h_n^{-1}\}$,
and hence $\{ h_n\}$, is eventually constant.
Therefore $\gamma = \lim \rho_n(h_n)$ lies in $\rho^i(\pi_1(F_i))$,
which contradicts the hypotheses of the lemma. We conclude that
$\alpha$ must be parabolic. 
\end{proof}

In order to show that the limit set of $\rho^i(\pi_1(F_i))$ and a
distinct conjugate of $\rho^j(\pi_1(F_j))$ do not cross
it remains to check that there are no problematic parabolic fixed points.
Our proof generalizes the argument of Proposition 2.7 of 
\cite{anderson-canary-culler-shalen}.

\begin{proposition}{onepoint}
Let $\{\rho_n\}$ be a sequence in $\DD(S)$ which is convergent on an
essential subsurface $F$, with nonannular components $F_i$ and
quasifuchsian limit representations $\rho^i$. Suppose also that
$\{\rho_n(\pi_1(S))\}$ converges geometrically to $\Gamma$.

If $\gamma\in \Gamma$ and either $i\ne j$ or
$i=j$ and $\gamma\notin \rho^i(\pi_1(F_i))$, then
the intersection of limit sets
$$\Lambda(\gamma \rho^i(\pi_1(F_i)) \gamma^{-1}) \intersect
\Lambda(\rho^j(\pi_1(F_j)))$$ 
contains at most one point.
\end{proposition}

If one makes use of  Soma and Anderson's generalization of Susskind's 
result and Bonahon's tameness theorem,  one may replace the
assumption in Proposition \ref{onepoint}
that the limit representations are quasifuchsian with the
weaker assumption that $\{\ell_{\rho_n}(\partial F)\}$ converges to 0.

\begin{proof}
The hypothesis that the $\rho^i$ are quasifuchsian representations of
$F_i$ implies in
particular that the lengths $\{\ell_{\rho_n}(\boundary F)\}$ converge
to 0, and hence we may apply  Lemma \ref{purelypar}.

Fixing $\gamma$, $i$ and $j$, let 
$\Phi_1 = \gamma \rho^i(\pi_1(F_i)) \gamma^{-1}$, and
$\Phi_2 = \rho^j(\pi_1(F_j))$.
Lemma \ref{purelypar} implies that $\Phi_1\intersect \Phi_2$ is a purely
parabolic subgroup, and so has at most 1 limit point (0 if it is
trivial). Thus, the proposition follows from Theorem
\ref{limitsetints} once we establish that $P(\Phi_1,\Phi_2) =
\emptyset$.

Let $\{ \rho_n(h_n) \}$ be a sequence converging to $\gamma$, and
suppose there is a point $x\in P(\Phi_1,\Phi_2)$.
The stabilizer $\stab_{\Phi_2}(x)$ is generated by some
$\rho^j(b)$, and
$\stab_{\Phi_1}(x)$ is generated by
$\gamma\rho^i(a)\gamma^{-1}$ where $a$ and $b$ are primitive elements
of $\pi_1(F_i)$ and $\pi_1(F_j)$, respectively.
Since these two elements must commute, Proposition \ref{JMrel} implies
that $h_n\sigma_n^i(a)h_n^{-1}$ commutes with
$\sigma_n^j(b)$ for sufficiently large $n$.
(Here $\sigma_n^i$ are as in the proof of Lemma \ref{purelypar}).
Since $a$ and $b$ are primitive and all abelian subgroups of
$\pi_1(S)$ are cyclic,
$$ h_n\sigma_n^i(a)h_n^{-1} = (\sigma_n^j(b))^{\pm 1} $$
for sufficiently large $n$. Applying $\rho_n$, and taking a
limit we conclude that
$$
\gamma\rho^i(a)\gamma^{-1} = \rho^j(b)^{\pm 1}
$$
but this contradicts the assumption that $\gamma\rho^i(a)\gamma^{-1}$ and
$\rho^j(b)$ generate a rank 2 group. Thus $P(\Phi_1,\Phi_2)$
must be empty and the proposition follows.
\end{proof}

In order to convert these conclusions about limit sets to conclusions
about embedded surfaces, let us
recall from \cite{anderson-canary:cores} that a collection
$\Phi_1,\ldots,\Phi_n$ of nonconjugate quasifuchsian subgroups of a
Kleinian group $\Gamma$ is called {\em precisely embedded} if 
$\stab_\Gamma(\Lambda(\Phi_i)) = \Phi_i$ for each $i$, and 
if every translate of 
$\Lambda(\Phi_i)$ by an element of $\Gamma$ is contained in the closure of a 
component of $\hhat\C \setminus \Lambda(\Phi_j)$, for each $i$ and $j$.

A system of {\em spanning disks} $\{D_1,\ldots,D_n\}$ for $\{\Phi_i\}$
are disks properly embedded in $\Hyp^3\union\hhat\C$ such that
$\boundary D_i = \Lambda(\Phi_i)$, $\stab_\Gamma(D_i)=\Phi_i$ and
$\gamma(D_i)$ is disjoint from $D_j$ unless $i=j$ and $\gamma\in\Phi_i$.
Thus, such disks
would project in $\Hyp^3/\Gamma$ to embedded, disjoint surfaces
$D_i/\Phi_i$.

Anderson and Canary observe, in Lemma 6.3 of
\cite{anderson-canary:cores} and the remark that follows (p. 766),
that

\begin{lemma}{precisely embed}
Any precisely embedded system
$\{\Phi_1,\ldots,\Phi_n\}$ of quasifuchsian subgroups of a Kleinian
group $\Gamma$ admits a system of spanning disks
$\{D_1,\ldots,D_n\}$. Furthermore, one may choose the spanning disks
so that there exists $\epsilon>0$ such that each component of
the intersection of any embedded surfaces $D_i/\Phi_i$ with a
non-compact component of $\MT_\epsilon$ is a properly embedded,
totally geodesic half-open annulus. 
\end{lemma}

We are now ready to complete the proof of the embedding theorem in the
quasifuchsian case.

\subsubsection*{Proof of Proposition \ref{limit embed quasifuchsian}}
Let $S$,  $F$, and $\{\rho_n\}$ be as in the statement of the proposition.
Let $\Gamma $ be the geometric limit  of $\{\rho_n(\pi_1(S))\}$ and
consider the quasifuchsian limit representations $\rho^i:\pi_1(F_i) \to \Gamma$.

Proposition \ref{onepoint} implies that the limit sets
of any two distinct conjugates of $\rho^i(\pi_1(F_i))$ and
$\rho^j(\pi_1(F_j))$ for components $F_i$ and $F_j$ are disjoint,
or intersect in exactly one point. These limit sets are all Jordan
curves since the groups are quasifuchsian, and thus any one of them 
is contained in the closure of a complementary disk of any other.
The conclusion of Proposition \ref{onepoint}, applied to the
conjugates of a single group $\rho^i(\pi_1(F_i))$, imply that 
$\stab_{\Gamma}(\Lambda(\rho^i(\pi_1(F_i))))$ must be 
$\rho^i(\pi_1(F_i))$ itself. Thus $\{\rho^i(\pi_1(F_i))\}$ form a
precisely embedded system of quasifuchsian groups in $\Gamma$, and
we can apply Lemma \ref{precisely embed} to obtain a system of
spanning disks $D_i$ for these groups. 
Note that the quotients $D_i/\rho^i(\pi_1(F_i))$ may be identified
with ${\rm int}(F_i)$, so that 
the resulting embeddings $h_i:{\rm int}(F_i) \to \hhat N$ are disjoint, are in
the homotopy classes determined by $\rho^i$, and so that there exists 
$\ep>0$ so that intersection of $h_i(F_i)$ with each component of
$\MT_\ep({h_i}_*([\boundary F_i]))$ is a properly embedded totally geodesic
half-open annulus.  Therefore, truncating the maps at 
$\MT_\ep({h_i}_*([\boundary F_i]))$ yields embeddings that are anchored with
respect to the $\ep$-Margulis tubes, whose domains can be identified with the
compact surfaces $F_i$.

Finally, composing with a diffeomorphism of $\hhat N$ that takes
the $\ep$-Margulis tubes $\MT_\ep({h_i}_*([\boundary F_i]))$ to
the $\epone$-Margulis tubes
$\MT({h_i}_*([\boundary F_i]))$ and is homotopic to the identity,
we obtain the desired anchored embedding $h$. This concludes
the proof of Proposition
\ref{limit embed quasifuchsian}.
\qed

\subsubsection{Using the covering theorem}

We now consider the case where the algebraic limits are geometrically
infinite. The main statement we need is the following, whose proof is
an application of Thurston's Covering Theorem.

\begin{proposition}{limit embed degenerate}
Let $S$ be an orientable surface and let $R$ be a connected essential
non-annular subsurface.  Let $\{\rho_n\} $ be a sequence  in $\DD(S)$
that is convergent on $R$ with limit representation
$\hat\rho:\pi_1(R)\to \PSL 2(\C)$,
and suppose that $\{\rho_n(\pi_1(S))\}$ converges geometrically to
$\Gamma$.  

Suppose that  
\begin{enumerate}
\item
$\hat\rho(g)$ is parabolic if and only if $g$ represents
a boundary component of $R$, and
\item 
$\hat\rho(\pi_1(R))$ is geometrically infinite.
\end{enumerate}

If $K$ is a compact subset of $\hhat N = \Hyp^3/\Gamma$, 
then there exists an anchored embedding $h:R\to N$ in the homotopy
class determined by $\hat\rho$, whose image does not intersect $K$. 
\end{proposition}

\begin{proof}
The following generalization of Thurston's covering theorem
\cite{wpt:notes} is established in \cite{canary:covering}.
 
\begin{theorem}{covering-theorem}
Let $N$ be  a topologically tame hyperbolic 3-manifold  which covers
an infinite volume hyperbolic 3-manifold $\hhat N$ by a local isometry
$\pi : N \to \hhat N$.
If $E$ is a geometrically infinite end of $N^0$, then 
$E$ has a neighborhood $U$ such that $\pi$ is finite-to-one on $U$.
\end{theorem}
(Here we recall that $N^0$ is obtained from $N$ by removing all its cuspical
$\epone$-Margulis tubes).

Since $R$ has only one component, after remarking $\{\rho_n\}$ by a
sequence of inner automorphisms we may assume that
$\{\rho_n|_{\pi_1(R)}\}$ converges to $\hat\rho$.

Let $N = \Hyp^3/\hat\rho(\pi_1(R))$. By the assumptions and Bonahon's
theorem, $N^0$ may be identified with $R\times\R$, 
and has a geometrically infinite end $E$. Let $\pi:N\to \hhat N$ be
the covering map associated to the inclusion
$\hat\rho(\pi_1(R))\subset \Gamma$. Note that $\hhat N$ has infinite
volume since it is the geometric limit of infinite-volume hyperbolic
3-manifolds. Theorem 
\ref{covering-theorem} then implies that there exists a neighborhood 
of $E$, on which $\pi$ is finite-to-one. 

Suppose that there does not exist a neighborhood
of $E$ on which $\pi$ is one-to-one.
An argument in Proposition 5.2 of \cite{anderson-canary:cores}
then implies that there exists a primitive element $\alpha \in
\hat\rho(\pi_1(R))$ which is a $k$-th power of some $\gamma\in\Gamma$,
with $k>1$. Let $\alpha = \hat\rho(a)$ and $\gamma = \lim \rho_n(g_n)$
for $\{g_n\in\pi_1(S)\}$. By Lemma \ref{JMrel}, for large enough $n$
we must have $a = g_n^k$. However, a primitive element of $\pi_1(R)$
must also be primitive in $\pi_1(S)$, and this is a contradiction.

Thus there is a neighborhood $U$ of $E$ on which $\pi$ is an
embedding, and hence there is a $t\in\R$ such that
$R\times\{t\}\subset U$, and $\pi(R\times\{t\})$ avoids $K$. This
gives our desired anchored embedding.
\end{proof}

\subsubsection{Proof of Proposition \ref{general subsurface embed}}
The proof of the limit embedding theorem in the general case now
follows from Propositions \ref{limit embed quasifuchsian}
and \ref{limit embed degenerate}. Let $F$ be the surface on which 
$\{\rho_n\}$ converges in the sense of Definition \ref{convergence on
subsurface}. The assumption that
$\rho^i(g)$ is parabolic if and only if $g$ is peripheral in
$\pi_1(F_i)$ implies that
$\rho^i(\pi_1(F_i))$ is either quasifuchsian or geometrically
infinite with no non-peripheral parabolics. Let $F'\subset F$ be the
union of the quasifuchsian components. Proposition \ref{limit embed
quasifuchsian} gives us an anchored embedding $h':F'\to\hhat N$. 
Enumerate the components of $F\setminus F'$ as $F_1,\ldots,F_k$.
Let $K_0 = h'(F')$. Applying Proposition \ref{limit embed degenerate}
gives us an anchored embedding $h_1:F_1 \to \hhat N$ avoiding $K_0$.
Now inductively define $K_i = K_0 \union \bigcup_{j \le i} h_j(F_j)$,
and apply Proposition \ref{limit embed degenerate} to obtain
$h_{i+1}$ avoiding $K_i$. The union of maps $h',h_1,\ldots,h_k$
is the desired anchored embedding of $F$.
\qed

\subsection{Proof of Theorem \ref{Relative embeddability}}
\label{embedding proof}

\subsubsection*{Outline}
As in the sketch following Proposition \ref{Easy embeddability}, the
basic strategy is to assume the theorem fails and consider a sequence
of counterexamples in which the anchoring constants $\ep_n$ go to 0
and uniform embeddability fails. After extracting a subsequence and
remarking, we may assume that the domains are a fixed surface $R$, the
curve systems are a fixed $\Gamma$, and that $\Gamma$ can be
partitioned as $\Gamma' \union \Delta$, where the lengths of the
curves in $\Gamma'$ go to 0, and the lengths of the curves in $\Delta$
are bounded from below (\S\ref{embeddability proof setup}). Let $D = R\setminus \collar(\Gamma') = X
\union \collar(\Delta)$. We may further assume that the metric on $D$ is fixed.
Let $N_n$ be the target manifolds, and $f_n:X\to N_n$ be the
maps in the hypothesis. We can extend $f_n$ to  $\hhat f_n:D \to N_n$
in a bounded way because the tubes of $\Delta$ in $N_n$ are not getting too
deep. 

As before, we wish to construct an anchored embedding of  $X$
into a geometric limit of $\{ N_n\}$
and then pull back to $N_n$ to obtain a contradiction for large values
of $n$. We begin by working with the larger surface $D$. 
Since the components of $\hhat f_n(D)$ may be
pulling apart, we may actually need to consider a collection of geometric limits.
Each component of $D$ is
contained in a maximal collection $D^J$ of components of $D$
such that a subsequence of $\hhat f_n|_{D^J}$ converges 
to a limiting map $\hhat f^J:D^J\to N^J$ into
an appropriate geometric limit of $\{N_n\}$ (\S\ref{setup geometric limits}).

We'd like to apply the embedding result from the previous section
to $\hhat f^J$, but
we can't guarantee that there are no unexpected parabolic elements.
Let $P^J$ be a maximal collection of disjoint non-peripheral curves on $D^J$
which are homotopic into cusps of $N^J$. If we let $F^J=D^J-\collar(P^J)$,
then Proposition \ref{general subsurface embed}
guarantees the existence of an anchored embedding
$$\bbar h^J:(F^J,\partial F^J)\to (N^J\setminus \MT(\boundary F^J),
\boundary \MT(\boundary F^J))$$
which is homotopic to $\hhat f^J|_{F^J}$ (\S\ref{anchoring
  on parabolics}).
We then use the unwrapping property and Lemma \ref{bounded annulus}
to extend $\bbar h^J$ to an embedding $\hhat h^J$
defined on all of $D^J$ and homotopic to $\hhat f^J$ 
(\S\ref{resewing along parabolics}).

However, what we want is an anchored embedding of
$X^J=D^J\cap X=D^J-\collar(\Delta)$,
so we must ``reanchor'' on $\Delta$.
We apply a result of Freedman, Hass, and Scott \cite{freedman-hass-scott:area} to produce
an embedding $\hhat g^J:D^J\to N^J$ whose image is contained in a regular
neighborhood 
of $\hhat f^J(D^J)$ and misses $\MT(\Delta)$. Finally, we apply the
Annulus Theorem to produce an embedded annulus joining
$\hhat g^J(D^J)$ to each component of $\MT(\Delta)$. The usual
surgery argument produces the desired anchored embedding, $g^J:X^J\to
N^J$ (\S\ref{reanchoring on Delta}).

This anchored
embedding has a bilipschitz collar neighborhood and is homotopic to
$f^J = \hhat f^J|_{X^J}$, and a further topological argument (Lemma
\ref{anchoring homotopy}) yields a homotopy
through anchored surfaces. 
We can pull back this homotopy to find,
for all large $n$,
a uniform homotopy of $f_n|_{X^J}$ to an anchored embedding
of $X^J$ into $N_n$.
Since we can do this for each collection $D^J$,
we can combine the resulting uniform homotopies to obtain a contradiction
for large values of $n$ (\S\ref{contradiction}).

\bigskip

\subsubsection{Proof: setting up the notation}\label{embeddability proof setup}
Fix $S,L$ and $\delta$, and
suppose by way of contradiction that the theorem fails. Then there is
a sequence $\{(\rho_n,R_n,f_n,\Gamma_n,\ep_n,\hat\ep_n,K_n)\}$ with
$\ep_n\to 0$, $\hat\ep_n\to 0$, and 
$K_n\to \infty$ for which the hypotheses of the theorem hold, but for
which the conclusions fail.

Possibly precomposing $\rho_n$ and $f_n$ with a sequence of
homeomorphisms of $S$ and passing to a subsequence, we may assume that
all $R_n$ are equal to a fixed surface $R$, and $\Gamma_n$ are equal
to a fixed curve system $\Gamma$. 
The surface $X=R\setminus\collar(\Gamma)$ is equipped with a sequence
of $L$-bounded metrics, which we may assume (again after remarking
and passing to a subsequence) converge to an
$L$-bounded metric $\nu$.  Then, possibly adjusting $L$ slightly,
we may assume that $f_n$ is $L$-bounded with respect to $\nu$ for
each $n$. For each $n$ we have an extension $\bbar f_n:R\to N_n$ with
the properties given in the statement of Theorem \ref{Relative
  embeddability}, notably the unwrapping condition.
The failure of the conclusions means that
there is no homotopy
$$H_n:(X\times[0,1],\partial X\times [0,1])\to
(N_n-\MT(\partial X,n),
\partial\MT(\partial X,n))$$
with $H_n(\cdot, 0)=f_n$,
such that
\begin{enumerate}
\item
$H_n$ is $K_n$-Lipschitz,
\item
$H_n$ is a smooth $K_n$-bilipschitz embedding of $X\times [1/2,1]$ with
the norm of the second derivatives bounded by $K_n$, 
\item
For all $t\in [1/2,1]$, $H_n(\boundary X\times \{t\})$ is a collection of
geodesic circles in the intrinsic metric of $\boundary\MT(\boundary X,n)$.
\item
$H_n(X\times [1/2,1])$ does not intersect any $\epsilon_1$-Margulis
tubes with core length less than $\hat \epsilon_n$, and
\item
$H_n(X\times[0,1])$ avoids the $\hat\ep_n$-thin part of $N_n$, and
\item
$H_n(X\times [0,1])$ lies in $C^{1/2}_{N_n}\union (N_n)_{(0,\ep_2]}$.
\end{enumerate}
Here by $\MT(B,n)$ we denote the $\ep_1$-Margulis tubes in $N_n$
associated to $\bbar f_n(B)$ for a curve system $B$ in $R$. 

Possibly restricting to a subsequence
again, we may assume that for each $\gamma\in\Gamma$ the
lengths $\{\ell_{\rho_n}(\gamma)\}$ converge.
Let $\Gamma'$ be the set of $\gamma\in\Gamma$ whose lengths
$\{\ell_{\rho_n}(\gamma)\}$ converge to 0, and let $\Delta = \Gamma\setminus
\Gamma'$. 
It will be convenient to suppose that in the metric $\nu$, all
boundary components of $X$ have the same length. This can be done by
changing $\nu$ by a bilipschitz distortion, and then altering the constant
$L$ appropriately.

Let $D = X\union \collar(\Delta)= R \setminus \collar(\Gamma')$.
The metric on $X$ can be
extended across $\collar(\Delta)$ to a metric which 
makes each component isometric to $S^1\times[0,1]$ (with $S^1$
isometric to a component of $\boundary X$).
Extend each $f_n$ to a map $\hhat f_n$ which, on each component of
$\collar(\Delta)$, takes the intervals
$\{p\}\times[0,1]$ to geodesics whose maximal length is shortest among
all such maps. This length is uniformly bounded above because the
$\rho_n$-lengths of components of $\Delta$ converge
to positive constants, and hence the distance from the cores of their
Margulis tubes to their boundaries is bounded. Notice that
$\hhat f_n (\collar(\Delta))\subset \MT(\Delta,n)$,
$\hhat f_n(D)\subset C_{N_n}\intersect N_n^1$. Moreover note that
$\hhat f_n$ and $\bbar f_n$ agree on $X$, and on $\collar(\Delta)$
they are connected by a homotopy rel boundary whose image lies in $\MT(\Delta,n)$,
followed possibly by a reparameterization of the domain by twists in
$\collar(\Delta)$. 
In particular, $\hhat f_n$
satisfies the same unwrapping condition as $\bbar f_n$.

After remarking the $\rho_n$ by Dehn
twists supported on $\collar(\Delta)$, we may assume
the extended maps $\hhat f_n$ are in the homotopy classes determined by
$\rho_n$.  We now have a fixed metric on $D$ and a sequence of
$L'$-Lipschitz maps $\hhat f_n$ for some constant $L'$.

\bigskip

\subsubsection{Geometric limits}\label{setup geometric limits}
Fix a basepoint $x^j$ for each component $D^j$ of $D$, and let
$y^j_n = f_n(x^j)$. Let $\hat y_n^j$ be a baseframe at $y_n^j$.
Since the metric on $D^j$ is fixed, 
each $\hhat f_n$ is $L'$-Lipschitz, and $\pi_1(D^j)$ is non-abelian,
a standard application of the Margulis lemma gives
an uniform lower bound on the injectivity radius of $N_n$
at $y^j_n$ for all $n$ and $j$. 
We may therefore pass to a subsequence such that, for any fixed $j$, $(N_n,\hat y_n^j)$ converges to
a hyperbolic manifold with baseframe $(N^j,\hat y^j)$. 
Furthermore,
$\{(C_{N_n},\hat y_n^j)\}$ converges to
$(C^j,\hat y^j)$ where $C_{N^j}\subset C^j$. 
(In fact, in our setting one
can  further show that $C^j=C_{N_j}$.)
Moreover, if $\{c_n^j:X_n\to N_n\}$ are the comparison maps, we can assume that 
the sequence of maps
$\{(c_n^j)^{-1}\circ\hhat f_n|_{D^j} : (D^j,x^j) \to  (N^j,y^j)\}$ (which make sense for all
large enough values of $n$)
converges to a map
$\hhat f^j:(D^j,x^j) \to (N^j,y^j)$.

After further restriction to a subsequence we may assume that, for
each pair $(j,j')$, the distances $\{ d(y^j_n,y^{j'}_n)\}$ 
converge to some $d_{jj'}\in [0,\infty]$. The relation $d_{jj'}<\infty$ is
an equivalence relation; fix an equivalence class $J$.
For $j,j'\in J$, we may identify $N^j $ with $ N^{j'}$, naming it $N^J$.
Notice that $d_{N^J}(y^j,y^{j'})=d_{jj'}$.  Let $D^J=\union_{j\in J} D^j$
and let $\hhat f^J:D^J\to N^J$ denote the union of the maps $\hhat f^j$
for all $j\in J$. Since $\hhat f_n(D^J)\subset C_{N_n}$ for all $n$,
$\hhat f^J(D^J)\subset C_{N^J}$.

Our goal now is to apply Proposition \ref{general subsurface embed}
to deform $\hhat f^J$ to an embedding in $N^J$. We must first obtain an
algebraically convergent sequence in the sense of Definition
\ref{convergence on subsurface}.

For each $J$, choose $j_0\in J$ and let $y^J=y^{j_0}$,  $\hat y^J=\hat y^{j_0}$,
$y^J_n=y^{j_0}_n$, $\hat y^J_n=\hat y^{j_0}_n$,
and $c_n^J=c_n^{j_0}$.
Let $\Theta^J$ be an embedded tree in $N^J$
formed by joining each $y^j$ to $y^J$ with an arc.
For all large enough $n$ the pullback $c_n^J(\Theta^J)$ of $\Theta^J$ to $N_n$ 
may be deformed slightly to give an embedded tree $\Theta^J_n$ all of whose
edges join the pullback  of $y^J$ to $y^j_n$ for some $n$.
(We must deform slightly since the endpoints of $c_n^J(\Theta^J)$
are only guaranteed to be near to the $y^j_n$.)
Therefore the total length of $\Theta^J_n$ is bounded for all large $n$.
Using paths in the tree we obtain, for all $j\in J$, homomorphisms
$$
\rho_n^j: \pi_1(D^j,x^j) \to \pi_1(N_n,y^J_n)
$$
consistent with the maps $\hat f_n|_{D^j}$. (More explicitly, if
$e^j_n$ is the edge in $\Theta^J_n$ joining $y_n^J$ to $y_n^j$,
then $\rho_n^j$ takes $[\alpha ]$ to
$[e^j_n * f(\alpha) * \overline{e^j_n} ]$.)
After conjugating $\rho_n$ in
$\PSL 2(\C)$, we may assume that the fixed baseframe $\hat x_0$ for $\Hyp^3$ descends
to $\hat y^J_n$ and 
hence consider $\rho_n$ as an isomorphism
from $\pi_1(S)$ to $\pi_1(N_n,y^J_n)$.
Thus, for each $j\in J$, we can define
$$
\sigma_n^j: \pi_1(D^j,x^j) \to \pi_1(S)
$$
by $\sigma_n^j = \rho_n^{-1} \circ \rho_n^j$. 

Now, $\{\rho^j_n\}$ converges, after restricting to a
subsequence, because each  $\hat f_n|_{D^j}$ is $L'$-lipschitz and
$\Theta_n$ has bounded total length, so the images of any fixed element of
$\pi_1(D^j)$ are represented by loops of uniformly bounded length,
and hence move the origin in $\Hyp^3$ a uniformly bounded amount. 

Thus, after restricting to an appropriate subsequence,
$\{\rho_n\}$ converges on the subsurface $D^J$, in the 
sense of Definition \ref{convergence on subsurface}, using the maps
$\sigma_n^j$ as defined above.
The limiting representation
$\rho^j$, for a component $D^j$ where $j\in J$, 
corresponds to the homotopy class of the limiting map $\hhat f^J|_{D^j}$.

\subsubsection{Anchoring on parabolics}\label{anchoring on parabolics}
These limiting representations may have non-peripheral parabolics.
Let $P^J$ denote a maximal set of disjoint homotopically distinct simple 
closed nonperipheral curves in $D^J$ whose images under the limiting
representations are parabolic.
Let $F^J = D^J \setminus \collar(P^J)$, and for each component $F^i$ of
$F^J$ contained in a component $D^j$ of $D^J$,  fix an injection
$\pi_1(F^i) \to \pi_1(D^j)$ consistent with the 
inclusion map. Then, with the same $\{\sigma_n^j\}$ as before,
restricted to $\pi_1(F^i)$, we have convergence
of $\{\rho_n\}$ on $F^J$, and the limiting representations $\hat\rho^i$
(which are the restrictions of $\rho^j$ to $\pi_1(F^i)$)
have no nonperipheral parabolics. 

We can now apply Proposition \ref{general subsurface embed} to $F^J$,
obtaining an anchored embedding
$$
\bbar h^J:(F^J,\boundary F^J) \to (N^J,\MT(\boundary F^J))
$$
such that, for each component $F^i$ of $F^J$, 
$\bbar h^J|_{F^i}$ is in the homotopy class determined by $\hat\rho^i$,
which is the same as the homotopy class of $\hhat f^J|_{F^i}$. 
Since each component of $\MT(\boundary F^J)$ is a cusp, it is
easy to see that $\bbar h^J|_{F^i}$ is properly homotopic
to $\hhat f^J|_{F^i}$ within $N^J\setminus \MT(\boundary F^J)$.

\subsubsection{Resewing along parabolics}\label{resewing along parabolics}
We next want, for each component $\alpha$ of $P^J$, to add 
an embedded annulus on $\boundary\MT(\alpha)$ to
the image of $\bbar h^J$, 
thus obtaining an anchored embedding of $D^J$ which is
homotopic to $\hhat f^J$. The unwrapping property of each
$\bbar f_n$, and hence each $\hhat f_n$,
will guarantee the existence of such an annulus.

Since all the curves in $P^J$ are homotopic into cusps
of $N^J$, for large enough $n$ they all have $\rho_n$-length less than
$\delta$. The unwrapping condition then implies in particular that, for all large
enough $n$, the image $\bbar f_n(R)$, and hence also $\hhat f_n(D)$, does not
intersect $\MT(P^J,n)$.
Therefore, $\hhat f^J(D^J)$
does not intersect $\MT(P^J)$.  Let
$$\bbar H^J:(F^J\times[0,1],\boundary F^J\times [0,1]) \to
(N^J\setminus \MT(\boundary F^J),\boundary\MT(\boundary F^J)),$$
be a homotopy with $\bbar H^J_0 =\hhat f^J |_{F^J}$ and $\bbar H^J_1 = \bbar h^J$.
 
We now explain how to  use $\bbar H^J$ to extend $\bbar h^J$
across the annuli $\collar(P^J)$ 
to obtain a map $\hhat h^J$ which takes these annuli to $\boundary \MT(P^J)$. 

In the union of  solid tori $\collar(P^J)\times [0,1]$, let 
$\varphi$ be a bilipschitz homeomorphism from the top annuli
$\collar(P^J)\times \{1\}$ to the remainder of the boundary, 
$(\collar(P^J)\times\{0\}) \union (\boundary\collar(P^J)\times[0,1])$,
which is the identity on the intersection curves
$\boundary\collar(P^J) \times \{1\}$ and is homotopic rel boundary
to the identity map.
Extend $\bbar H^J$ to all of $D^J\times\{0\}$ to be equal to $\hhat f^J$, and
then consider the map $\bbar H^J\circ \varphi$ on the annuli 
$\collar(P^J)\times\{1\}$ which maps into the complement of $\MT(P^J)$. 
On each annulus component $A$, this map is homotopic rel boundary, in
the exterior of $\MT(P^J)$, to a unique ``straight'' map to
$\boundary\MT(P^J)$. Here ``straight'' means that
geodesics orthogonal to the core of the annulus are taken to straight
lines in the Euclidean metric of $\boundary\MT(P^J)$. In particular
the map is an immersion. 
Define $\hhat h^J$ to be the extension of $\bbar h^J$ to $\collar(P^J)$ by
this straight map.
Use the homotopy between $\hhat h^J$ and $\hhat f^J$ on $\collar(P^J)$
to extend $\bbar H^J$ across the solid tori to a proper
(in $N^J\setminus \MT(\boundary D^J)$)
homotopy $\hhat H^J$
between $\hhat h^J$ and $\hhat f^J$ which is defined
on $D^J\times [0,1]$ and avoids $\MT(P^J)$.

We recall that there exist comparison maps
$c_n^J:X_n \to N_n$
such that $\{ X_n\}$ exhausts $N^J$ and $c_n^J$ is increasingly $C^2$-close to a local
isometry.
We may choose $R>0$ so that the image of $\hhat H^J$ is contained in $B_R(y^J)$.
For all large enough $n$,
let $\hhat H^J_n=c_n^J\circ \hhat H^J$ be the pullback of
$\hhat H^J$ and let $\hhat h^J_n=c_n^J\circ \hhat h^J$ be the pull back of
$\hhat h^J$.
Lemma \ref{comparisons respect tubes}
allows us to further assume that, for all large enough $n$, there exists a collection
$T_n$ of components of $(N_n)_{(0,\epone)}$ such that
\begin{enumerate}
\item
$c_n^J(\MT(\partial F^J)\cap B_R(y^J))\subset T_n$,
\item
$c_n^J( \boundary\MT(\partial F^J)\cap B_R(y^J))\subset\boundary T_n$, and
\item
$c_n^J(B_R(y^J)-\MT(\partial F^J))\subset N_n-T_n$.
\end{enumerate}
Proposition \ref{JMrel}(3) implies that, again for large enough $n$, $T_n=\MT(\partial F^J,n)$.
Therefore, $\hhat h_n^J$ is an anchored surface.

We will now apply Lemma \ref{unwrapped is embedded} to show
that $\hhat h^J$ is an embedding on $\collar(\alpha)$ for each component $\alpha$ of $P^J$.
Let $\hhat H^J_n=c_n^J\circ \hhat H^J$ be the pullback of
$\hhat H^J$ and let $\hhat h^J_n=c_n^J\circ \hhat h^J$ be the pull back of
$\hhat h^J$.
The unwrapping property implies that $\hhat f_n(D^J)$ can be pushed
to $+\infty $ or $-\infty$ disjointly from $\MT(\alpha,n)$.
Since $\{(c_n^J)^{-1}\circ\hhat f_n|_{D^J}\}$ converges  to $\hhat f^J$, for large enough
$n$ there is a very short homotopy from $\hhat f_n|_{D^J}$ to
$\hhat H^J_n|_{D^J\times\{0\}}$. Since the image of $\hhat H^J_n$
does not intersect $\MT(P^J,n)$, 
we conclude that $\hhat h^J_n$ can also be pushed to either $+\infty$ or
$-\infty$ in $N_n^0\setminus \MT(\alpha,n)$.

Since the $\rho_n$-lengths of $\boundary F^J$ converge to 0 they are
eventually less than $\epotal$, and then 
Otal's Theorem \ref{otal}
implies that $\MT(\boundary F^J,n)$ is unknotted and unlinked in
$N_n^0$.  Therefore, 
we can apply Lemma \ref{unwrapped is embedded} to
the restriction of $\hhat h^J_n$ to the union of $\collar(\alpha)$
with the components of $F^J$ which are adjacent to it, 
concluding that $\hhat h^J_n$ restricted to $\collar(\alpha)$ is an
embedding into $\boundary\MT(\alpha,n)$.
Thus in the geometric limit
$\hhat h^J|_{\collar(\alpha)}$ is an embedding into
$\boundary \MT(\alpha)$. Applying this to all
components of $P^J$ we conclude that the map $\hhat h^J$ is an
embedding into $N^J$.

\subsubsection{Reanchoring on $\Delta$}\label{reanchoring on Delta}
The embedding $\hhat h^J$ is defined on $D^J$, whereas we need an anchored
embedding whose domain is 
$$
X^J = X\intersect D^J = D^J \setminus\collar(\Delta).
$$
Restricting $\hhat h^J$ to
$X^J$ is not sufficient, since it would not be anchored on
$\MT(\Delta)$.

Thus, consider again the map $\hhat f^J$, which meets
$\MT(\Delta)$ only in the embedded annuli
$\hhat f^J(\collar(\Delta))$. Deform these intersection
annuli to the boundary of a small regular neighborhood of
$\MT(\Delta)$, obtaining a map $\underline{f}^J: D^J\to N^J$
which misses $\MT(\Delta)$ and is properly homotopic to $\hhat h^J$ within
$N^J\setminus\MT(\boundary D^J)$.

Let $Y$  be a neighborhood of $\underline{f}^J(D^J)$ within
$N^J\setminus \MT(\boundary X^J)$. 
Let $Z$ be a compact, irreducible submanifold of
$N^J\setminus \MT(\boundary D^J)$ which contains $Y$, 
such that $\underline{f}^J$ is homotopic to  $\hhat h^J$  within $Z$.
Moreover, possibly adjusting $\hhat h^J$ by an isotopy supported in a
neighborhood of $\boundary\MT(\boundary D^J)$, we may assume that the homotopy
between $\hhat h^J$ and $\underline{f}^J$ fixes the
boundary pointwise. 

We are now in a position to use the following result of
Freedman-Hass-Scott \cite{freedman-hass-scott:area}:

\begin{theorem}{FHS}
Let $Z$ be compact and irreducible and $f:(D,\boundary D)\to
(Z,\boundary Z)$ a $\pi_1$-injective map that admits a homotopy, fixing the
boundary pointwise, to an embedding. Then for any neighborhood $Y$ of $f(D)$,
there is a homotopy of $f$, fixing the boundary pointwise, to an embedding that lies in $Y$.
\end{theorem}

\subsubsection*{Remark} Bonahon was the first to observe that Theorem \ref{FHS}
follows quickly from Freedman-Hass-Scott \cite{freedman-hass-scott:area}.
A proof of Theorem \ref{FHS} in the case that $F$ is a closed surface is
given in remark at the end of section 2 of Canary-Minsky \cite{canary-minsky:tamelimits}.
In order to establish the relative version stated above, one simply replaces the use of
Theorem 5.1 in \cite{freedman-hass-scott:area} with the relative version derived
in section 7 of \cite{freedman-hass-scott:area} (see also Jaco-Rubinstein \cite{jaco-rubinstein}).

\medskip

Applying Theorem \ref{FHS}, 
we obtain an anchored embedding
$\hhat g^J:D^J\to N^J$ which is properly homotopic to $\underline{f}^J$
and whose image lies in $Y$.
In particular, $\hhat g^J(D^J)$ misses $\MT(\Delta)$.
(Notice that we can't simply obtain $\hhat g^J$ by 
naively pushing $\hhat h^J$ off of $\MT(\Delta)$, since we
have no a priori control over the intersection of $\hhat h^J(D^J)$
with $\MT(\Delta)$.)

Let $\hhat g^J_n=c_n^J\circ \hhat g^J$ be the pullback of $\hhat g^J$ by the comparison maps to $N_n$
(defined for $n$ sufficiently large). Notice that, again by  Lemma \ref{comparisons respect tubes}
and Proposition \ref{JMrel},
we may assume that $\hhat g^j_n$ is anchored and is in the homotopy class determined
by $\rho_n$ for all large enough $n$.

We claim that for each component $\beta$ of $\Delta$
there is a homotopy from $\hhat g^J_n(\beta)$ to an embedded longitudinal curve on
$\boundary\MT(\beta,n)$  that avoids
$\MT(\Delta\union \partial D^J,n)$.
Notice first, that since $f_n$ is $\epotal$-anchored,
Otal's Theorem \ref{otal} implies that
$\MT(\partial X,n)=\MT(\Gamma \union \boundary R,n)$ is
unknotted and unlinked in $N_n^1$.
Hence $N_n^0$ can be
identified with $S\times\R$ in such a way that the geodesic representative
$\beta^*_n$ is $\beta\times\{0\}$, $\MT(\beta,n)=\collar(\beta)\times [a,b]$,
and $B=\beta\times\R$ is disjoint from 
$\MT(\Delta-\{\beta\},n)\union \MT(\boundary D^J,n)$, and in particular from 
$\hhat g^J_n(\boundary D^J)$ since $\hhat g^J_n$ is anchored.  
Since $\beta $ is non-peripheral in the essential surface $D^J$ and
$\hhat g^J_n$ is in the homotopy class determined by $\rho_n$, $\hhat g^J_n(D^J)$
must intersect $B$.  After proper isotopy of $B$ we may assume the
intersections of $B$ with $\hhat g^J_n(D^J)$ are essential circles, and 
so the closest one to $\beta\times\{0\}$ yields the desired homotopy.

In order to use $c_n^J$ to transport this homotopy to $N^J$ we must
first bound its diameter.  As the bilipschitz constants of $c_n^J$
converge to 1, $\ell_{N_n}(\hhat g^J_n(\beta))\le C$ for some uniform
constant $C$.  By Lemma \ref{bounded annulus}, there is a homotopy
from $\hhat g^J_n(\beta)$ to $\boundary \MT(\beta,n)$ which 
avoids  $\MT(\Delta\union\boundary D^J,n)$ and 
lies in an $r(C)$-neighborhood of $\hhat g^J_n(D^J)$. 
Applying Lemma \ref{comparisons
respect tubes}, we can, for large enough $n$, use $c_n^J$ to pull this homotopy back
to obtain a  homotopy $Q_\beta$ from $\beta$ to $\boundary\MT(\beta)$ which
avoids $\MT(\Delta\union\boundary D^J)$.
We will next apply a version of the Annulus Theorem 
to conclude that there is an {\em embedded} annulus $\hhat Q_\beta$ in
the complement  of $\MT(\Delta\union \boundary D^J)$
joining $\hhat g^J(\beta)$ to $\boundary \MT(\beta)$, whose interior is
disjoint from $\hhat g^J(D^J)$. 

Since $\hhat g^J(D^J)$ is embedded, $Q_\beta^{-1}(\hhat g^J(D^J))$ is a union of
embedded curves in the domain annulus. The inessential ones may be
removed by a homotopy, and the remainder are isotopic to the
boundary.
Hence by restricting to a complementary component  of the remaining
curves of intersection we obtain
a new homotopy which meets $\hhat g^J(D^J)$ only on a boundary
curve. The image of this curve may not be embedded in $\hhat g^J(D^J)$, but
since it is homotopic to $\beta$ we may deform it in $\hhat g^J(D^J)$
to a simple curve.
Shifting this deformation slightly away from $\hhat g^J(D^J)$ we
obtain a new homotopy $Q'_\beta$ which meets $\hhat g^J(D^J)$ in a simple 
curve. 

Let $Z'$ be a compact, irreducible submanifold of $N^J \setminus
\MT(\Delta\union \boundary D^J)$ which contains the
$2r(C)$ neighborhood of $\hhat g^J(D^J)$. Remove from $Z'$ a regular
neighborhood $Y'$ of $\hhat g^J(D^J)$,  to obtain a
Haken manifold $W$.   If $Y$ is chosen small enough then
$Q'_\beta\intersect W$ is a proper
singular  annulus with one boundary embedded in $\boundary Y'$
and the other in $\boundary \MT(\beta)$. 
Now we may apply the Annulus Theorem (see \cite{jaco-shalen} and
\cite[Thm VIII.13]{jaco}) in $Z'-Y'$
to conclude that there is an embedded annulus $\hhat Q_\beta$ joining
$\hhat g^J(\beta)$ to $\boundary\MT(\beta)$ whose interior avoids 
$\hhat g^J(D^J)$ and $\MT(\Delta\union \boundary D^J)$.

Repeat this for every component of $\Delta$. The resulting embedded
annuli may intersect but only in inessential curves, so after surgery
we obtain a union $\hhat Q_\Delta$ of embedded annuli.

A surgery using a regular neighborhood of
$\hhat Q_\Delta$ then yields a smooth embedding $g^J:X^J\to N^J$
which is anchored on $\MT(\boundary X_J)$, and is homotopic to
$f^J=\hhat f^J|_{X^J}$.

This is not quite enough for us since we need a
{\em homotopy of pairs}
$$H^J:(X^J\times [0,1],\boundary X^J\times [0,1])
\to (N^J\setminus \MT(\boundary X^J),\boundary \MT(\boundary X^J)).$$
The issue here is that, although both $f^J$ and $g^J$ are
anchored on all the tubes of $\MT(\boundary X^J)$, including those of
$\Delta$ which are not cusps in the geometric limit, the homotopy that
we have constructed contains steps in which the anchoring on $\Delta$ disappears,
and moreover the homotopy is not constrained to stay away from those
tubes. We can resolve this issue with the following lemma, whose proof
we postpone to \S\ref{S: anchoring homotopy}.

\begin{lemma}{anchoring homotopy}
Let $M$ be an oriented irreducible 3-manifold with boundary,
$D$ a compact oriented surface of negative Euler characteristic, and suppose that
$g:(D,\boundary D) \to
(M,\boundary M)$ is an incompressible embedding and
$f:(D,\boundary D) \to (M,\boundary M)$ is homotopic to $g$.
Let $\Delta$ be an essential curve system in $D$ and $T$ an open regular
neighborhood of $g(\Delta)$ in $M$ such that $A = g^{-1}(T) = f^{-1}(T)$ is a regular
neighborhood of $\Delta$, and  $f|_A = g|_A$. 
Let  $X = \overline{D\setminus A}$.

We further suppose that
\begin{enumerate}
\item $M\setminus T$ is acylindrical
\item A simple curve
$\alpha\subset X$, such that $g(\alpha)$ is homotopic in $M$ to a power
of a component of $g(\boundary X)$, is homotopic within $X$ to a component of
$\boundary X$. 
\end{enumerate}

Then,  $f|_X$ and $g|_X$ are homotopic as maps of
pairs from $(X,\boundary X)$ to $(M\setminus T,\boundary T\union
\boundary M)$. 
\end{lemma}

In order to apply this lemma, with $M$ being $N^J$ minus its cusps, $D
= D^J$ and $\Delta = \Delta^J$, we first note that, by 
Lemma \ref{level curves immersed}, both $f^J$ and $g^J$ take
$\boundary X^J$ to curves on $\MT(\Delta^J)$ in the level homotopy class.
Hence after a homotopy supported near $\MT(\Delta^J)$, we can assume that $f^J$ and $g^J$
agree on $\boundary X^J$. We can then extend both maps to
$\collar(\Delta^J)$ so that they agree there and map $\Delta^J$ to the
geodesic cores of the tubes in $N^J$. We let these extensions to $D$ be $f$
and $g$, respectively. 

The incompressibility of $g$ is clear. 
One quick way to verify assumption (1), acylindricity of $M\setminus T$,
is to recall that, since $g(\Delta)$ is a union of geodesics in a
hyperbolic 3-manifold,  $N^J-g(\Delta)$ admits a metric of pinched negative curvature
such that each component of $T-g(\Delta)$ is a toroidal cusp in this metric
(see Agol \cite{agol:volume} and the discussion in the proof of Lemma \ref{bounded annulus}.)
It remains to verify assumption (2). If $\alpha$ is a simple
closed curve in $X$ and $b$ a boundary component of $X$ such that
$g(\alpha)$ and $g(b^k)$ are homotopic, then the homotopy between
them can be pulled back to an approximate $N_n$, in
which the manifold is homotopy-equivalent to a surface and $X$
corresponds to a subsurface. There, such a homotopy gives rise to a
homotopy within $X$ of $\alpha$ to $b^k$,  which implies that $k=\pm
1$ since $\alpha$ is simple.

We may now apply Lemma \ref{anchoring homotopy}
to obtain the desired homotopy of pairs. 
Since $N^J$ is homeomorphic to ${\rm int}(C^{1/2}_{N^J})\union N^J_{(0,\ep_2)}$,
by a homeomorphism which is the identity on
$C^{1/4}_{N^J}\union N^J_{(0,\ep_1]}$, we may assume
that both $g^J$ and the homotopy $H^J$ to $f^J$ lie entirely
in ${\rm int}(C^{1/2}_{N^J})\union N^J_{(0,\ep_2)}$.
We may further assume that the restriction of $H^J$ to
$X^J\times [1/2,1]$ is an $C^2$-embedding and that
for all $t\in [1/2,1]$, $H^J(\boundary X^J\times \{t\})$ is a collection of
geodesic circles in $\boundary\MT(\boundary X^J)$.

\subsubsection{Obtaining the contradiction}\label{contradiction}
As its image is compact, the homotopy $H^J$
between $g^J$ and $f^J$ avoids the $\hat\ep^J$-thin part
for some $\hat\ep^J>0$. Let $\{\MT_1,\ldots,\MT_n\}$ be the components
of the $\epone$-thin part which $g^J(X^J)$ intersects and which
are either cusps or have core curves of length less than $\hat\ep^J$.
Notice that no $\MT_i$ is a component of $\MT(\boundary X^J)$.
For each $i$, there is a regular neighborhood $\mathcal{U}_i$ of $\MT_i$
which is contained in $N^J_{(0,\ep_2)}$
and a diffeomorphism
$\Upsilon_i:\mathcal{U}_i\setminus (\MT_i\cap N^J_{(0,\hat\ep^J)})
\to \mathcal{U}_i\setminus\MT_i$
which is the identity on $\boundary \mathcal{U}_i$. 
Extend the collection of  $\Upsilon_i$, via the identity outside 
$\union \MT_i$, to an embedding $\Upsilon:N_{[\hat\ep^J,\infty)} \to N$. 
We may replace $H^J$ with 
$\Upsilon\circ H^J$ which has the additional property of avoiding
$\epsilon_1$-Margulis tubes with core length less than $\hat\ep^J$.

Pulling $H^J$ back to $N_n$, after applying Lemma \ref{comparisons respect tubes}, 
we obtain, for
all large enough $n$, a $(2K^J,{\hat\ep^J\over 2})$-uniform homotopy
$H_n^J$ from $c_n^J\circ f^J$ to an anchored embedding.
Since $\{(c_n^J)^{-1}\circ \hhat f_n^J\}$ converges  to $\hhat f^J$, we may, for large enough $n$,
deform $H_n^J$ slightly so that it is a 
$(3K^J,{\hat\ep^J\over 3})$-uniform homotopy between an anchored embedding and
$f^J_n$.
Condition (3) in the definition, that the images of the boundary circles 
$\boundary X\times \{t\}$ are geodesic in $\boundary\MT(\boundary X)$
for $t\in[1/2,1]$, may be obtained by noticing that, because the $C^2$
bounds on the comparison  maps converge to 0, for large enough $n$ the
images are nearly 
geodesic circles and hence can be expressed as graphs of nearly
constant functions over the geodesic circles to which they are
homotopic. Hence the map can be adjusted to satisfy condition (3).

Since, $\{ (C_{N_n},\hat y_n^J)\}$ converges geometrically to $(C^j,\hat y^J)$,
which contains $C_{N^J}$,
the resulting homotopy $H_n^J$  lies within
\hbox{$(C^{1/2}_{N_n}\union N^n_{(0,\ep_2]})\setminus\MT(\partial X^J,n)$}
for all large enough $n$. 

For each equivalence class $J'$,
we obtain a sequence of homotopies $H_n^{J'}$ in the same way.
Since the distance $d_{jj'}$ between basepoints converges to $\infty$ if
$j\in J$, $j'\in J'$, these maps eventually have disjoint images. 
Combining we obtain for all large $n$ a $(K,\hat\ep)$-uniform
homotopy $H_n$ which shows that $f_n$ is $(K,\hat \ep)$-uniformly embeddable,
where $K=\max\{ 3K^J\}$ and $\hat\ep=\min\{ {\hat\ep^J\over 3}\}.$
If $n$ is chosen large enough that $K_n>K$ and
$\hat\ep_n <\hat\ep$, then we have obtained a contradiction.
\qed

\subsubsection{Anchoring the homotopy}
\label{S: anchoring homotopy}
In this final subsection we supply the proof of Lemma \ref{anchoring
  homotopy}. 

\begin{proof}
Since
$f|_A = g|_A$ we may apply a small deformation to $f$ so that  $f|_A$
is parallel to but disjoint from $g(D)$. Then 
using the homotopy from $f$ to $g$, restricted
to $X$, we may obtain a homotopy
$H:X\times[0,1] \to M$, such that $H(\cdot,0) = f|_X$,  $H(X,1)$ is disjoint from
$g(X)$ and parallel to it in a collar neighborhood, and $H(\boundary X
\times\{0\})$ is disjoint from $g(D)$. 

The key point is to adjust $H$ so
that it misses the curves of $g(\Delta)$. 
Thus consider $\Sigma = H^{-1}(g(D))$. Assuming general position,
this is a properly embedded surface in $X\times[0,1]$, and disjoint
from $X\times\{1\}$ and $\boundary X \times\{0\}$. After the usual
surgery operation we may assume that $\Sigma$ is
$\pi_1$-injective in $X\times[0,1]$. See, for example the proof of Lemma 6.5 in Hempel \cite{hempel}, and note that we
need to use the fact that $g(D)$ is incompressible, and that
$M\setminus g(D)$ is irreducible, which follows from our
assumptions. Moreover, again by irreducibility of $M\setminus g(D)$,
we can remove all disk components of $\Sigma$, and thus we
can assume that $\Sigma$ meets $\boundary X\times[0,1]$ only in
essential curves. Since $\Sigma$ does not meet $X\times\{1\}$ at all,
it is properly isotopic to a $\pi_1$-injective surface whose boundary
lies in $X\times\{0\}$, and we can apply 
Waldhausen \cite{waldhausen}, Proposition
3.1 and Corollary 3.2, to conclude that each
component of $\Sigma$ is parallel to a subsurface of 
$X\times\{0\}$.

Now let $\alpha\subset \Sigma$ be a component of $H^{-1}(g(\Delta))$,
which by general position is a loop. Let $\alpha'$ be a projection of
$\alpha$ to $X$. 

If $\alpha'$ is homotopically trivial then $\alpha$ is homotopically
trivial in $\Sigma$ (by $\pi_1$-injectivity), and hence bounds a disk
$E\subset\Sigma$. A regular neighborhood $B$ of $E$ has boundary
sphere mapping to the complement of $g(\Delta)$. Now $M\setminus
g(\Delta)$ is irreducible, since our hypothesis implies that the
universal cover of $M$ is $\R^3$, and hence $H$ can be redefined on
$B$ to a map that misses $g(\Delta)$.

If $\alpha'$ is homotopically nontrivial, then $g(\alpha')$ is homotopic to a power of
$\gamma$, so, by assumption (3), 
$\alpha'$ is actually homotopic to a corresponding boundary component
$\beta$ in $X$. Thus, since each component of $\Sigma$ is parallel to $X\times\{0\}$,
$\alpha$ can be joined to $\beta\times[0,1]$ by 
an embedded annulus $C_\alpha$. Moreover this can be done simultaneously
for all $\alpha$'s so that the $C_\alpha$'s are disjoint from each
other and from $X\times\{0\}$. After deleting a regular
neighborhood of these annuli from $X\times[0,1]$ we therefore obtain a
submanifold $P$ which is still homeomorphic to $X\times[0,1]$, and
contains $X\times\{0,1\}$. The restriction $H|_P$ therefore gives a
homotopy missing $g(\Delta)$, which we rename $H$. 

Now we can construct the desired
homotopy of maps of pairs: First, after retracting $T\setminus
g(\Delta)$ to $\boundary T$ we can assume the image of $H$ is contained in
$M\setminus T$. Then
on each annulus  $\beta\times[0,1]$,
by the acylindricity of $M\setminus T$,  $H$ can be deformed rel $\beta\times\{0,1\}$
to a map with image in $\boundary T$. 
After realizing this homotopy in a
collar neighborhood, we obtain a homotopy that takes $\beta\times[0,1]$ into
$\boundary T$. This completes the proof. 
\end{proof}

\section{Insulating regions}
\label{insulating}

In this section we will establish the existence of long `bounded-geometry
product regions' in the hyperbolic manifold $N$ when the hierarchy satisfies
certain conditions. Roughly, if $H$ contains a very long
geodesic $h$ supported in some non-annular domain $R$, and if there are no very
long geodesics subordinate to $h$, then there is a big region in $N$
isotopic to $R\times[0,1]$, whose geometry is prescribed by $h$.
Furthermore the model map can be adjusted to be an embedding on this 
region, without disturbing too much the structure outside of it. 
In order to quantify this more carefully, let us make the following
definition: 

A segment $\gamma$ of a geodesic $h\in H$ is said to be
{\em $(k_1,k_2)$-thick}, where $0<k_1<k_2$, provided:
\begin{enumerate}
\item $|\gamma|>k_2$
\item For any $m\in H$ with $D(m) \subset D(h)$ and $\phi_h(D(m))
\intersect \gamma \ne\emptyset$, 
$|m|<k_1$.
\end{enumerate}

Let $\tau_1$ and $\tau_2$ be two (full) slices with the same bottom
geodesic $h$,  
and suppose that the bottom simplices $v_{\tau_1}$ and $v_{\tau_2}$
are spaced by at least $5$,
and $v_{\tau_1}<v_{\tau_2}$. As in \S\ref{address region structure}, 
there is a region
$\BB(\tau_1,\tau_2)\subset \modl$, homeomorphic to $D(h)\times[0,1]$
and bounded by $\hhat F_{\tau_1}$ and $\hhat F_{\tau_2}$ and the tori
$\storus {\boundary D(h)}$. It is the geometry of the model map on
such regions that we will control.

\begin{theorem}{product region}
Fix a surface $S$. Given positive constants $K$, $k$, $k_1$ and $Q$,
there exist $k_2$ and $L$ such
that,  if $f:\modl \to N$ is a $(K,k)$ model map,
$\gamma$ is a $(k_1,k_2)$-thick segment of $h\in H_\nu$
and $\xi(h)\ge 4$, then there exist slices
$
\tau_{-2},\tau_{-1},\tau_0,\tau_1,\tau_2
$
with bottom
geodesic $h$ and bottom simplices $v_{\tau_i}$ in $\gamma$ satisfying
$$v_{\tau_{-2}} < v_{\tau_{-1}} < v_{\tau_{0}} < v_{\tau_{1}} <
v_{\tau_{2}},$$ 
with spacing of at least 5 between successive
simplices, and so that 
\begin{enumerate}
\item $f$ can be deformed, by a homotopy supported on the union of
$\BB_2 = \BB(\tau_{-2},\tau_2)$ and the tubes $U(\boundary D(h))$, 
to an $L$-Lipschitz
map $f'$ which is 
an orientation-preserving embedding on $\BB_1 =
\BB(\tau_{-1},\tau_1)$, and
\item $f'$ takes $\modl \setminus \BB_1$ to 
$N\setminus f'(\BB_1)$, and
\item
the distance from $f'(F_{\tau_0})$ to 
$f'(F_{\tau_{-1}})$ and $f'(F_{\tau_{1}})$
is at least $Q$.
\end{enumerate}
\end{theorem}

\begin{proof}
Suppose, by way of contradiction, that the theorem fails. Then
there exist $K,k,Q$ and $k_1$ and sequences $k_n\to \infty$ and $L_n\to\infty$,
representations $\rho_n\in\DD(S)$ with associated hierarchies $H_n$,
$(K,k)$ model maps $f_n:M_{\nu_n} \to N_n$, 
and $(k_1,k_n)$-thick segments $\gamma_n \subset h_n\in H_n$,
but for which the model maps $f_n$ cannot be deformed to an $L_n$-Lipschitz
map satisfying the conclusions of the theorem.

We will extract and study a certain geometric limit in order to obtain
a contradiction.

\subsubsection*{Convergence of hierarchies} 
Let us discuss briefly a natural sense of convergence for a sequence
of hierarchies, which is a mild generalization of the notions used in
 \cite[\S6.5]{masur-minsky:complex2} and \cite[\S5.5]{minsky:ELCI}.

Fix a subsurface $R\subseteq S$ and a basepoint $v_0\in \CC(R)$, and
for $E>0$ let $B_E(v_0)$ denote the $E$-ball around $v_0$ in $\CC(R)$.
If $H_n$ is a sequence of hierarchies containing geodesics $h_n$ with
$D(h_n) = R$, and $H_\infty$ is a hierarchy whose {\em main} geodesic
$h_\infty$ is biinfinite and has domain $R$, we say that $H_n$ {\em converge to $H_\infty$
relative to $R$} if the following holds: 
\begin{enumerate}
\item For any $E>0$, for all sufficiently large $n$,   $h_n\intersect
  B_E(v_0) = h_\infty \intersect B_E(v_0)$. 
\item For any $E>0$, for all sufficiently large $n$, the set of tight geodesics
$$\beta(H_n,E) \equiv \{m\in H_n: D(m)\subsetneq R, [\boundary D(m)]\subset B_E(v_0)\}$$
is equal to the set
$$\beta(H_\infty,E) \equiv \{m\in H_\infty: D(m)\subsetneq R, [\boundary D(m)]\subset B_E(v_0)\}.$$
\end{enumerate}
(Note that convergence is independent of the choice of $v_0$).

Returning now to the sequence $H_n$ from our argument by
contradiction, we can obtain this type of convergence after suitable
adjustments. First note that 
after passing to a subsequence and remarking, we may assume that the 
domains $D(h_n)$ are a constant surface $R$. We then have:
\begin{lemma}{hierarchies converge}
After remarking the sequence $H_n$, one can choose a basepoint in $\CC(R)$
which is within the middle third of each $\gamma_n$, such that a subsequence
of the  $H_n$ converges relative to $R$ to a hierarchy $H_\infty$. 
\end{lemma}

\begin{proof}[Proof of Lemma \ref{hierarchies converge}]
Note that it suffices to show that, 
for each $E$, the subsets $h_n\intersect B_E(v_0)$ and $\beta(H_n,E)$
eventually stabilize. This is because the set $H_\infty$ of tight
geodesics obtained in this way naturally inherits the subordinacy
relations from the $H_n$'s, and hence has the structure of a hierarchy.
(The same argument is made in \cite[\S6.5]{masur-minsky:complex2}).

We start by showing that (after suitable remarking and restriction to
a subsequence) the geodesics $h_n$, and in particular their
subsegments $\gamma_n$, converge in the sense of part (1) of the definition.

Choose a basepoint
$v_{n,0}$ for $\gamma_n$ which is distance at least $k_n/3$ from each
endpoint of $\gamma_n$. Let  $\tau_{n,0}$ be a maximal slice of $H_n$
containing  
the pair $(h_n,v_{n,0})$, and let $\mu_{n,0}$ be its associated clean marking.
Since there are only finitely many clean markings in $S$ up to homeomorphism,
we may assume after remarking and extracting a
subsequence that all the $\mu_{n,0}$ are equal to a
fixed $\mu_0$, and $v_{n,0}\equiv v_0$.

Fix $E>0$ and suppose that $n$ is large enough that
$E < k_n/6$. We claim that there is a finite set of possibilities,
{\em independent of $n$}, for the simplices of $h_n$ that are within
distance $E$ of $v_0$.  To see this, let $w$ be such a simplex in
$h_n$. By Lemmas  5.7 and 5.8 of \cite{minsky:ELCI}, 
there is a resolution of $H_n$ containing
$\tau_{n,0}$ and passing through some slice $\tau$ containing the pair
$(h_n,w)$. 
Now by the monotonicity property of resolutions
(see Lemma \ref{J interval})
every slice in the resolution between $\tau_{n,0}$ and $\tau$
contains a pair $(h_n,u)$ with $v_0 \le u \le w$. 
Therefore, any geodesic appearing in this part of the resolution
and supported in $R$ must have footprint in $h_n$ that intersects the
interval $[v_0,w]$. Because of 
the $(k_1,k_n)$-thick property, all these geodesics have length
bounded by $k_1$. 

The sum of the lengths of all these geodesics can then be bounded by
$O(Ek_1^\alpha)$ where $\alpha\le\xi(S)$, using an inductive counting argument
similar to the one in Section \ref{po}: First, the segment
$[v_0,w]$ has length bounded by $E$. Thus there are at most $O(E)$
geodesics directly subordinate 
to $h_n$ with footprint intersecting this interval. Each of these has length
bounded by $k_1$, so there are $O(Ek_1)$ geodesics directly
subordinate to these geodesics. 
We continue inductively, and note that the complexity $\xi$ decreases
with each step. 
Since every geodesic with footprint in $h_n$ is obtained 
in this way (by the definition of a hierarchy), this gives us the
bound we wanted. 

Each elementary move in the resolution that takes place in $R$
corresponds to an edge in one of these geodesics, so we conclude that
the markings $\mu_0|_R$ and
$\mu_\tau|_R$ are separated by $O(Ek_1^\alpha)$ elementary moves.
(Note that we do not obtain or need a bound on the number of moves
that occur outside of $R$). Lemma 5.5 from \cite{masur-minsky:complex2}
then implies that the associated complete clean markings are separated
by  $O(Ek_1^\alpha)$ elementary moves.
Since the number of possible elementary moves on
any given complete clean marking of $R$ is finite (see \cite[\S 2.5]{masur-minsky:complex2}), this gives a finite
set, independent of $n$, containing all possible simplices of $h_n$
within distance $E$ of $v_0$. 

We conclude that there are only finitely many possibilities for the
segment of length $E$ of $h_n$ around $v_0$, and after extracting a
subsequence these can be assumed constant. Enlarging $E$ and
diagonalizing, we may assume that $h_n$ converges to a bi-infinite
geodesic $h_\infty$.

Now fix $E$ again and consider $\beta(H_n,E)$.
For large enough $n$, all $m\in \beta(H_n,E)$ are forward and
backward subordinate to $h_n$, since 
$D(m)$ intersects $\I(h_n)$ and $\T(h_n)$ by the triangle inequality
in $\CC(R)$. In particular,
if $h_n \bsubd m\fsubd h_n$ then $\I(m)$ and $\T(m)$ are simplices of
$h_n$, and they determine $m$ up to a finite number of choices by
Corollary 6.14 of \cite{masur-minsky:complex2}. Proceeding
inductively we see that
there are only finitely many possibilities for the elements of 
$\beta(H_n,E)$. Thus 
by the same diagonalization argument as the previous paragraph
we may assume that
$\beta(H_n,E)$ eventually
stabilizes. This gives us a limiting collection $H_\infty$ of
tight geodesics supported in subsurfaces of $R$, and as mentioned above the
subordinacy relations of $H_n$ among such geodesics survive to give
$H_\infty$ the structure of a hierarchy. Thus we have the desired
convergence. 
Note that $H_\infty$ has a
biinfinite main geodesic $h_\infty$, and every other geodesic has
length at most $k_1$. 
\end{proof}

\subsubsection*{Convergence of models}
From now on we assume that we have remarked and restricted to a
subsequence, and let $H_\infty$ be the limit hierarchy provided by
Lemma \ref{hierarchies converge}.
The hierarchy $H_\infty$ has an associated model manifold
$M_\infty\homeo \hhat R\times \R$ (where, as in \S\ref{isotopy convention}, 
$\hhat R$ is an open surface containing $R$ such that $R=\hhat R
\setminus \collar(\boundary R)$).
We claim that this is a geometric limit
of the models $M_n$ for $H_n$, with appropriate baseframes. In Section
\ref{geometric limits} we discussed the notion of geometric limits for
hyperbolic manifolds with basepoints. The same idea applies here,
except (since the models are only piecewise smooth) that the
comparison maps can only be required to be bilipschitz and piecewise
smooth. In fact as we shall see the comparison maps we will obtain
shall be isometries on every block in the comparison region. 

After restricting and remarking the sequence as above, we will select
baseframes $\hat x_n$ with basepoints $x_n$ in the split-level surfaces $F_{\tau_{n,0}}$, and a
baseframe $\hat x_\infty$ at a point $x_\infty\in M_\infty$, and show that
$\{(M_n,\hat x_n)\}$ converges geometrically to $(M_\infty,\hat x_\infty)$.
(Here, in the definition of geometric convergence we only assume that
the comparison maps are piecewise smooth, but they eventually map each
block by an isometry.)

\begin{lemma}{geom lim of models}
The models $M_n$ with baseframes $\hat x_n$ 
converge geometrically to $(M_\infty,\hat x_\infty)$. In every compact
subset of $M_\infty$, the comparison maps are eventually
block-preserving and map blocks to blocks by isometries. 
\end{lemma}

Let us first establish the following fact about the models $M_n$, whose purpose is to show
that if two blocks are glued along a 3-holed sphere whose projection
is a subsurface of $R$ whose boundary intersects the endpoints of
$h_n$, then the domains of both blocks are in $R$ too. 

\begin{lemma}{block in R}
Let
$B$ be a block in $M_n$ and $F_Y$ a  gluing surface for $B$, such that
the corresponding 3-holed sphere $Y$ is contained in $R$, and
$\boundary Y$ intersects $\base(\I(h_n))$ and $\base(\T(h_n))$.
Then $D(B)\subseteq R$.
\end{lemma}

\begin{proof}
Let $f$ be the $4$-geodesic
containing the edge $e$ corresponding to $B$. Then $D(f)=D(B)$, and
so $Y\fsubd f$ or $f\bsubd Y$ (often both). Suppose the former. 
The assumption on $\boundary Y$ implies that
$h_n\in \Sigma^+(Y)$, so by Theorem
\ref{Descent Sequences} $Y\fsubd f \fsubeq h_n$. In particular
$D(B)\subseteq D(h_n)=R$.
\end{proof}

\begin{proof}[Proof of  Lemma \ref{geom lim of models}]
Recall that we have already passed to a subsequence such that
$\tau_{n,0}$ is constant. We choose a component $Y$ of
$\base(\tau_{n,0})$ which projects into $R$ and, for each $n$, a point
$x_n\in M_n$ identified with a  fixed point $x\in F_Y$.
We similarly choose an orthonormal baseframe for $T_{x_n}F_Y$, identified
with a fixed orthonormal baseframe, and
extend it to an orthonormal baseframe  $\hat x_n$ for $T_{x_n}M_n$ by adding a unit
vector normal to $F_Y$ and pointing upward (in the natural product structure on $M_n$).
One may similarly choose $x_\infty\in M_\infty$ and a baseframe $\hat x_\infty$.

Let $B_n$ be
one of the two blocks for which $F_Y$ is a gluing surface. By
Lemma \ref{block in R}, $D(B_n)\subseteq R$. 

Let $M_{H_n,h_n}$ denote the union of blocks and tubes in $M_n$ whose
associated forward or backward sequences pass through $h_n$ -- in
particular this is contained in $R\times\R$. The previous paragraph
shows that $B_n\subset M_{H_n,h_n}$.

We claim that, for any fixed $r$, the $r$-neighborhood $\NN_r(x_n)$ in $M_n$
is contained in $M_{H_n,h_n}\union\UU_n(\boundary R)$ for sufficiently
large $n$.  

First, after deforming paths in tubes to tube boundaries, we note that
any point in 
$\NN_r(x_n)\setminus \UU_n$ is reachable from $x_n$ through a sequence
of $s=O(e^r)$ blocks.

Suppose that $B$ and $B'$ are adjacent blocks such that $D(B)\subset R$
and $\phi_{h_n}(D(B))$ (nonempty by lemma \ref{nonempty footprint}) is
distance at least 4 from $\base(\I(h_n))$ and $\base(\T(h_n))$.
Then $B$ and $B'$ share a gluing surface $F_Y$ such that $\boundary Y$
must intersect the curves associated to $\base(\I(h_n))$ and
$\base(\T(h_n))$ (or laminations if $h_n$ is infinite).
Lemma \ref{block in R} implies that $D(B')\subseteq R$ as well. 
Moreover, if $\xi(R)>4$, then $\phi_{h_n}(D(B)) $ and $\phi_{h_n}(D(B'))$ are both in
$\phi_{h_n}(Y)$, which has diameter at most 2.   If $\xi(R)=4$, we conclude that
$B$ and $B'$ are associated to adjacent edges of $h_n$.

We now assume that  $\xi(R)>4$ and complete the argument in this case.
(The case when $\xi(R)=4$ will be handled afterwards.)
If $n$ is large
enough that $v_0$ is more than $2s+8$ from the ends of $h_n$, then any
block $B$ that is reachable in $s$ steps from $B_n$ is still in
$M_{H_n,h_n}$, and moreover $d_{\CC(R)}(v_0,\boundary D(B)) \le 2s$.
We conclude that any block meeting $\NN_r(x_n)$ is in $M_{H_n,h_n}$,
and moreover the boundary of its associated domain is
a bounded distance from $v_0$
in $\CC(R)$, and hence for large enough $n$ its associated 4-geodesic
is equal to a geodesic in $H_\infty$. 

Now let $U$ be a tube in $\UU_n$ meeting
$\NN_r(x_n)$, such that $U$ is not a component of $\UU_n(\boundary
R)$.  
Since $U$ meets $\NN_r(x_n)$ it is adjacent to at least one block $B$
with $D(B)\subset R$ and $d_{\CC(R)}(v_0,\boundary D(B))\le
2s$. Thus $core(U)$ is contained in $R$ and
$d_{\CC(R)}(v_0,core(U))\le 2s+1$. We claim that in fact {\em all} 
blocks $B'$ adjacent to $U$ have $D(B')\subset R$
and $d_{\CC(R)}(v_0,\boundary D(B'))\le 2s+2$. To see this, note first 
that any block $B'$ adjacent to $U$ must either have a boundary component
in the homotopy class of $U$, or contain $core(U)$ in $D(B')$, and
hence if $D(B')\subset R$ we have
$d_{\CC(R)}(v_0,\boundary D(B'))\le 2s+2$. 
Now if $B_1$ is a block
adjacent to $U$ with $D(B_1)\subset R$, and $B_2$ is adjacent to $U$
and to $B_1$, then (for large enough $n$) we can again apply Lemma
\ref{block in R} to conclude that $D(B_2)\subset R$ as well.
It follows by connectivity of $\boundary U$ that, in fact, all blocks
adjacent to it have domain contained in $R$. 

We can use this to show that, for high enough $n$, $\omega(U)$ is
bounded. 
Recall from \cite[\S\S 8.3,9.3]{minsky:ELCI}
that $\Im \omega(U)$, or the ``height'' of $\boundary U$, is
$\epone$ times the number of annuli in $\boundary U$, and this is
estimated up to bounded ratio by the total number of blocks adjacent to $U$. 
The footprint $\phi_{h_n}(D(B))$ for any such block is contained in
$\phi_{h_n}(core(U))$, and so if $D(B) \fsubeq m\fsubd h_n$ there are
a finite number of possibilities for $m$, independently of $n$. 
The length of $m$ is bounded by $k_1$ by the $(k_1,k_n)$-thick condition.
By the same inductive counting argument as used above in the discussion of 
limits of hierarchies,
we therefore know the total number of such blocks 
is bounded by $O(k_1^\alpha)$.
We conclude that $\Im\omega(U)$ is uniformly bounded. 

The magnitude $|\Re\omega(U)|$ is estimated by the length $|l_U|$ of the
annulus geodesic $l_U$ associated to $U$; more precisely
$||\Re\omega(U)| - |l_U|| \le C|\Im\omega(U)|$ (see (9.6) and
(9.17) of \cite{minsky:ELCI}).
Since the footprint of this annulus domain is also at most $2s$ from
$v_0$, for high enough $n$ the $(k_1,k_n)$ condition implies that 
$|l_U|\le k_1$. We conclude that
$|\omega(U)|$ is bounded, and hence so is the diameter of $U$.

Thus, fixing $r$ and letting $n$ grow, we find that $\NN_r(x_n)$ is
contained in $M_{H_n,h_n}\union \UU_n(\boundary R)$, and that 
the geometry of the blocks and tubes (other than
$\UU_n(\boundary R)$)  eventually stabilize. 
It follows that the geometric limit of $(M_{H_n,h_n},\hat x_n)$
is  $(M_\infty,\hat x_\infty)$, minus the parabolic tubes associated to
$\boundary R$. Moreover the comparison maps
can be taken to preserve the block structure, and to be isometries on
each block and each tube that is not parabolic in $M_\infty$.

Furthermore, for each 
$\gamma$ in $\boundary R$ we have
$\Im\omega_n(\gamma)\to \infty$, because
$|h_n|> k_n \to \infty$. Thus $U_n(\gamma)$ converge geometrically to a rank-1
parabolic tube. Thus in fact the geometric limit of $(M_n,\hat x_n)$ is $(M_\infty,\hat x_\infty)$.

It remains to consider the case 
when $\xi(R)=4$. Now the blocks of $M_{H_n,h_n}$ are all associated with
edges of $h_n$, and hence are organized in a linear sequence with each
one glued to its successor. Each tube $U$
in $M_{H_n,h_n}$ is adjacent to exactly two blocks, so that
$\Im\omega(U)$ is uniformly bounded, 
and $\Re\omega(U)$ is bounded by 
the lengths of the associated annulus geodesic, as before. 
Hence the geometric limit is an infinite
sequence of blocks and tubes, and the rest of the conclusions follow
easily in this case too. 
\end{proof}

\subsubsection*{Limit model map}
Let $\tau_0$ be the slice of $H_\infty$ obtained as the limit of the
restrictions of $\tau_{n,0}$ to $R$. We may conjugate $\rho_n$ so that
$f_n(x_n)\in N_n=\Hyp^3/\rho_n(\pi_1(S))$ is the projection of the origin in $\Hyp^3$.
Since the model maps are uniformly
Lipschitz,  we may
extract from $\rho_n|_{\pi_1(R)}$ a convergent
subsequence. Denote the limit by $\rho_\infty$, and 
let $G_\infty = \rho_\infty(\pi_1(R))$. 

Since the boundary curves of $R$ have $|\omega_n|\to \infty$,
their lengths in $N_n$ go to zero (by the Short Curve Theorem of
\cite{minsky:ELCI})
so they must be parabolic in the
limit. So $G_\infty$ is a Kleinian surface group.

After restricting to a further subsequence we may assume 
(Lemma \ref{JMrel}) that  $\{\rho_n(\pi_1(S))\}$
converges geometrically to a group
$\Gamma_\infty$.

We can assume, and will do so for the remainder of this section, that
the model maps $f_n$ satisfy the conclusions of Lemma \ref{radial tube
  maps} -- in particular, for each tube $U$ of the model $M_n$, the
restriction of $f_n$ to 
a $t(r)$-collar neighborhood of $\boundary U$ in $U$ (where $r$ is the
depth of $U$ and $t(r)$ is a proper function)
takes radial lines to radial lines and preserves distance to the tube
boundary.

\begin{lemma}{limit model map}
The group $G_\infty$ is doubly degenerate, and is equal to the
geometric limit $\Gamma_\infty$.
After possibly restricting again to a subsequence,
the model maps $f_n:M_n\to N_n$ converge
geometrically to a model map $f_\infty : M_\infty \to N_\infty$, where
$N_\infty = \Hyp^3/G_\infty$. 
\end{lemma}

\begin{proof}
Let $\lambda_\pm\in\EL(R)$ be the endpoints of $h_\infty$ (by
Klarreich's theorem, see \S\ref{hierarchy background}).
The vertices $v_i$ of $h_\infty$ converging to
$\lambda_\pm$ (as $i\to \pm\infty$)
all have bounded length in $G_\infty$, so their geodesic
representatives leave every compact set in the quotient.
We conclude that
$\lambda_\pm$ are the ending laminations, and $G_\infty$ is doubly
degenerate.

The proof that $\Gamma_\infty = G_\infty$  is 
similar to an argument made by Thurston
\cite{wpt:notes} in a slightly different context
(see also \cite{canary:covering}). Since both ends of $G_\infty$ are
degenerate, 
Thurston's covering theorem (see Theorem \ref{covering-theorem}) tells
us that the covering map  
$\Hyp^3/G_\infty \to \Hyp^3/\Gamma_\infty$ is finite-to-one, 
and hence $[\Gamma_\infty,G_\infty]<\infty$.
If $\gamma\in\Gamma_\infty \setminus G_\infty$ then for some finite $k$
we have $\gamma^k\in G_\infty$. Let $\gamma = \lim \rho_n(g_n)$ 
with $g_n\in\pi_1(S)$, and let $\gamma^k = \rho_\infty (h)$ with
$h\in\pi_1(R)$. By Lemma \ref{JMrel}, since
$\rho_n(h^{-1}g_n^k)$ converges to the identity, we must have
$h = g_n^k$ for large enough $n$. Since $k$-th roots are unique in
$\pi_1(S)$, we find that $g_n$ is eventually constant, and since
$\pi_1(R)$ contains all of its roots in $\pi_1(S)$, $g_n$ is
eventually contained 
in $\pi_1(R)$. Thus $\gamma\in G_\infty$ after all.

Let $N_\infty = \Hyp^3/G_\infty = \Hyp^3/\Gamma_\infty$.
Again restricting to a subsequence, the model maps $f_n$ converge 
on $M_\infty$ minus the tubes to a $K$-Lipschitz map into $N_\infty$.
On the tubes, we can also obtain Lipschitz control: The 
non-peripheral tubes of $M_\infty$ have bounded $|\omega|$ and hence 
(property (5) of Definition \ref{model map def} of model maps) the
maps $f_n$ are uniformly Lipschitz on these tubes. On the
peripheral tubes, those associated to $\boundary R$, the depths go to infinity
but the conclusions of Lemma \ref{radial tube maps} tell us that on
increasingly large subsets of the tube the map is just a radial
extension of its values on the tube boundary. This together with the Lipschitz
control in the complement of the tubes 
guarantees that, for a fixed $K'$, the maps are eventually $K'$-Lipschitz
on each compact subset. Hence 
the maps in fact converge everywhere to a map
$f_\infty:M_\infty \to N_\infty$ which is a
homotopy-equivalence (since it induces an isomorphism on $\pi_1$)

We can see that $f_\infty$
is proper as follows: Any block $B$ in $M_\infty$ 
meets some slice surface $\hhat F_\tau$, and each $\hhat F_\tau$ 
meets a representative $\gamma_u\subset M_\infty$ of some vertex $u$
of the geodesic $h_\infty$ such that $\gamma_u$ has
bounded length. If a sequence $B_i$ of blocks in $M_\infty$ leaves
every compact set, the corresponding vertices $u_i$ go to $\infty$ in
$\CC(R)$, and so the images of $\gamma_{u_i}$ in $N_\infty$, whose
lengths remain bounded, must leave every compact set. Since the
surfaces $\hhat F_\tau$ have bounded diameter (because there is a
uniform bound on $|\omega|$ for all tubes in $M_\infty$), this means that
the sequence $f_\infty(B_i)$ also
leaves every compact set. Each non-peripheral tube
lies in a bounded neighborhood of some block, so the images of these
tubes are
properly mapped as well. On $\UU(\boundary R)$, the 
conclusions of Lemma \ref{radial tube maps} imply that the limiting
map isometrically takes radial lines of these rank-1 cusps to radial
lines, and so it is proper because its
restriction to the cusp boundary is proper. Thus $f_\infty $ is
proper. 

In the remainder of the proof we will use the following lemma several
times. It is essentially a uniform properness property for the
sequence of model maps. 
Let $y_n = f_n(x_n)\in N_n$, and for the limiting
basepoint $x_\infty\in M_\infty$ let $y_\infty=f_\infty(x_\infty)$.
Let $\varphi_n:Y_n\to N_n$ be a sequence of comparison maps where
$\{ Y_n\}$ is a nested exhaustion of $N_\infty$ by compact subsets.

\begin{lemma}{uniform properness}
For each $r>0$ there exist $n(r)$ and $d(r)$ such that, for
$n\ge n(r)$, 
$f_n^{-1}(\NN_r(y_n))$ is contained in 
the $d(r)$-neighborhood of $x_n$ in $M_{H_n,h_n}\union \UU_n(\boundary
R)$.  
\end{lemma}

\begin{proof}
Suppose by way of contradiction that the lemma is false. Then there
exists $r>0$ such that, after possibly restricting again to a
subsequence, there is a sequence $z_n\in M_n$ such that 
$f_n(z_n) \in \NN_r(y_n)$, but $d(z_n,x_n) \to \infty$. 

The tricky point here is that, a priori, the geometric limiting
process only controls the 
maps $f_n$ on large neighborhoods of $x_n$ in $M_{H_n,h_n}$, so 
we have to rule out the possibility that $z_n$ is in an entirely
different part of the model $M_n$. 

Assume without loss of generality that $y_\infty$ is in the $\epone$-thick
part of $N_\infty$.  Suppose first that $z_n$ is contained in a block
$B_n$ for each $n$.

Since the maps $f_n$ are uniformly Lipschitz, 
the image $f_n(B_n)$ remains a bounded distance from $y_n$. For large
enough $n$, its image must be in the compact set $\varphi_n(Y_n)$. 
Identifying $N_\infty$ with $R\times\R$,  we find that
$\varphi_n^{-1}(f_n(B_n))$ is 
homotopic into $R\times\{0\}$. For large enough $n$ this homotopy is
contained in the comparison region $Y_n$ and can be pulled back to $N_n$. 
Since $f_n$ is a homotopy equivalence, we conclude that $D(B_n)$ is
a subsurface of $R$.

Let $e_n$ be the 4-edge associated to $B_n$. 
Now we claim that
the distances $d_{\CC(R)}(e_n,v_0)$ are unbounded. Otherwise
there is a bounded subsequence and
as in the discussion on convergence of hierarchies, for
$n$ in the subsequence we can reach $B_n$
from $x_n$  in a bounded distance, using elementary moves from the
initial marking $\mu_0$ to a marking containing the vertex $e_n^-$
(recall from \S\ref{model definitions} that $e_n^\pm$ are the vertices of $e_n$).
This contradicts the assumption that $d(z_n,x_n)\to \infty$.

The sequence $\{e_n^-\}$, being unbounded in $\CC(R)$, contains infinitely
many distinct elements. However all of these are vertices of the model
and hence the $f_n$-images of the corresponding curves in $B_n$ have
uniformly bounded length in  $N_n$. 
The comparison maps take these 
curves to curves $\{\alpha_n\}$ of bounded length
in a compact subset of $N_\infty$, and this means they fall into
finitely many homotopy classes in $N_\infty$. However, if $\alpha_m$
and $\alpha_{m'}$ are homotopic in $N_\infty$, then for large enough
$n$ the homotopy pulls back and their preimages are homotopic in
$N_n$, and hence in $S$. This
contradicts the fact that that there are infinitely many distinct
$\{e_n^-\}$, and hence rules out  $d(z_n,x_n)\to\infty$ if the $z_n$ are all contained
in blocks $B_n$.

If (restricting to a subsequence) every $z_n$ is contained in a tube
$U_n$, then we can assume that $d(z_n,\partial U_n)\to\infty$,
since otherwise $z_n$ remain a bounded distance from some blocks and
the previous argument can be applied. Hence,  since the depth of $U_n$ goes to $\infty$, 
for all
large $n$ the model map
$f_n$ takes $U_n$ with
degree 1 onto the corresponding tube
$\MT_n\subset N_n$. The
conclusions of Lemma \ref{radial  tube maps} imply that
$d(f_n(z_n),\partial\MT_n)\to\infty$ and hence that
$d(f_n(z_n),y_n)\to \infty$,
giving the desired contradiction.
\end{proof}

As a consequence of this lemma we can show that
$f_\infty$ has degree 1. Indeed, let us 
show that $\deg f_\infty =\deg f_n$ for sufficiently high $n$. 
Let $\psi_n:(X_n,x_\infty)\to (M_n,x_n)$ be the sequence of comparison
maps for the geometric convergence of the model manifolds where
$\{X_n\}$ is a nested exhaustion of $M_\infty$ by compact sets.

If $W\subset M_\infty$ is a compact submanifold 
containing $f_\infty^{-1}(y_\infty)$ (which is compact since
$f_\infty$ is proper), the degree of $f_\infty|_W$ over $y_\infty$
is equal to $\deg f_\infty$.

Now let $d(0)$ and $n(0)$ be given by Lemma \ref{uniform properness},
so that $f_n^{-1}(y_n) \subset \NN_{d(0)}(x_n)$ for all $n>n(0)$.
We may choose $W$ large enough that, for large enough $n$, 
$\psi_n(W)$ contains $\NN_{d(0)+1}(x_n)$. Thus
the degree of $f_n|_{\psi_n^{-1}(W)}$ over $y_n$ is equal
to $\deg f_n$.

Now choosing $W$ according to the previous two paragraphs, 
we know by definition of geometric limits that the maps
$\varphi_n^{-1}\circ f_n \circ \psi_n$ are eventually defined on $W$
and converge to $f_\infty|_W$, so that 
for large enough $n$ the degree of 
$f_n|_{\psi_n^{-1}(W)}$ over $y_n$ equals the degree of
$f_\infty|_W$ over $y_\infty$. Hence $\deg f_\infty = \deg f_n$, and since
$\deg f_n = 1$, we have $\deg f_\infty = 1$ as desired. 

The remaining model map properties in Definition \ref{model map def} -- that
$f_\infty$ takes the tubes of $\UU[k]$ to the corresponding Margulis
tubes, and their complement to the complement of the tubes, and the
$\omega$-dependent Lipschitz bounds within the tubes -- are all
inherited from the properties of the maps $f_n$, via the geometric
convergence of both models and targets. 
This completes the proof of Lemma \ref{limit model map}.
\end{proof}

\subsubsection*{Product regions in the limit}
In order to finish the proof of Theorem \ref{product region}, we need 
a topological lemma about deforming proper homotopy equivalences of pairs.
Let $V$ be the 3-manifold $R\times\R$,
with $\boundary V = 
\boundary R\times\R$. Let $C_s = R\times[-s,s]$, which we note is a
relative compact core for $(V,\boundary V)$. 

\begin{lemma}{embed core}
Suppose that a map of pairs $f:(V,\boundary V) \to (V,\boundary V)$ is 
a proper, degree 1 map homotopic to the identity.
Then there exists
a homotopy of $f$ through maps of pairs to a map $f'$, such that
\begin{enumerate}
\item The homotopy is compactly supported,
\item $f'|_{C_1}$  is the identity,
\item $f'(V\setminus C_1) \subset V\setminus C_1$. 
\end{enumerate}
\end{lemma}

The proof of this lemma is fairly standard and we omit it.

\medskip

Now to apply this to our situation let $M'_\infty = M_\infty \setminus
\UU(\boundary R)$ and $N'_\infty = N_\infty \setminus \MT(\boundary R)$ be the
complements of the peripheral model tubes and Margulis tubes,
respectively, and note that there are orientation-preserving
identifications
$$
\Phi_M : M'_\infty \to V
$$
and 
$$\Phi_N:
N'_\infty \to V 
$$
so that
the map $F = \Phi_N \circ f_\infty \circ \Phi_M^{-1}$
satisfies the conditions of Lemma \ref{embed core}.
We will need to choose these identifications a bit more carefully.

First note that 
in $M_\infty$ every surface $\hhat F_\tau$ for a slice $\tau$
is  isotopic to a level surface.
Constructing a cut system as in \S\ref{po} and 
using Proposition \ref{topprec and cprec}, 
we may choose an ordered sequence of slices
$\{c_i\}_{i\in\Z}$ whose base simplices $v_i$ are separated by at
least 5 in $h_\infty$,  and adjust $\Phi_M$ so that 
$\Phi_M(\hhat F_{c_i}) = R\times\{i\}$.

We may choose the identification $\Phi_N$ so that
$\Phi_N^{-1}(C_1)$ contains a $(Q+1)$-neighborhood
of $\Phi_N^{-1}(R\times\{0\})$, where $Q$ is the constant in part (3)
of the theorem.

Now let $F':(V,\boundary V)\to (V,\boundary V)$ be the map homotopic
to $F$ given by
Lemma \ref{embed core} and let 
$f'_\infty = \Phi_N^{-1} \circ F' \circ \Phi_M$. 
A remaining minor step is to 
extend $f'_\infty$ to a map (still called $f'_\infty$)
which is defined on all of $M_\infty$, and homotopic to $f_\infty$ by
a homotopy of pairs
$(M_\infty,\UU(\boundary R)) \to (N_\infty, \MT(\boundary R))$ which
is supported on a compact set.  
Choose a positive integer $s$ so that the homotopy from $F$ to $F'$ is
supported in the interior of
$C_s=R\times[-s,s]$, pull the annuli $\boundary C_s \intersect
\boundary V$ back 
to annuli in $\boundary \UU(\boundary R)$ via $\Phi_M^{-1}$, 
and pick collar neighborhoods in $\UU(\boundary R)$ of these annuli. 
The extension of the homotopy to one supported in the union of 
$C_s$ and these collar neighborhoods is elementary. 

Let $G$ denote the final homotopy from $f_\infty$ to $f'_\infty$.
Choose slices $\tau_i$ so that
$\tau_0 = c_0$, $\tau_{\pm 1} = c_{\pm 1}$, and $\tau_{\pm 2} = c_{\pm
  s}$. Hence $\BB_2  = \BB(\tau_{-2},\tau_{2}) = \Phi_M^{-1}(C_s)$
together with $\UU(\boundary R)$ contain the support of the homotopy $G$, 
and $f'_\infty$ has all the topological properties described in the conclusions of the
theorem. It remains to pull this picture back to the approximating manifolds. 

Let $r$ be such that 
$\Phi_N^{-1}(C_1)=f_\infty'(\BB_1)$ is contained in
$\NN_{r-1}(y_\infty)$, and let $n(r)$ and $d(r)$ be the constants given by
Lemma \ref{uniform properness}. We may assume that 
$\NN_{2d(r)}(x_\infty)\subset \BB_2\union\UU(\boundary R)$. (If not, we may choose $s$ larger above, so that the inclusion does hold.)

Let $Z$ be a compact manifold contained in the interior of $\BB_2\union\UU(\boundary R)$
which contains $\NN_{2d(r)}(x_\infty)$ and
the support of the homotopy $G$. Let $V$ be a collar neighborhood
of $\boundary Z$ within $\left(\BB_2\union\UU(\boundary R)\right)-Z$.
Since $Z\union V$ is compact, for large enough $n$,
we can define $f_n'=\phi_n\circ f_\infty'\circ \psi_n^{-1}$ on $\psi_n(Z)$.
We define $f_n'=f_n$ on $M_n-\psi_n(Z\union V)$ and use the product structure
on $\psi_n(V)$ to extend $f_n'$ over $\psi_n(V)$.
Since $\phi_n^{-1}\circ f_n\circ \psi_n$ converges to
$f_\infty$, we can extend so that $\max\{ d(f_n(x),f_n'(x))\ |\ x\in \psi_n(V)\}\to 0$.

The slices $\tau_{\pm i}$ ($i=0,1,2$) give, for large enough $n$,
slices in $H_n$ that 
define regions $\BB_i(n)$ which converge, under the comparison maps,
to $\BB_i$. Recalling that the comparison maps preserve the block structure,
we see that $f'_n$ 
are orientation preserving embeddings on
$\BB_1(n)$, and admit homotopies to $f_n$ which are supported in $\BB_2(n)\union
\UU_n(\boundary R)$.
All that remains is to show that, for $n$ large enough,
$f'_n(M_n\setminus \BB_1(n))$ is disjoint from $f'_n(\BB_1(n))$. 

If $n$ is chosen greater than $n(r)$ and also sufficiently large that
$Z\cup V\subset X_n$ and $\psi_n$ is 2-biLipschitz on $Z\cup V$,
then $f_n'$ is defined and $\NN_{d(r)}(x_n)\subset \psi_n(Z)$.
We also note that, for large enough $n$, $f_n'(\BB_1(n))\subset\NN_r(y_n)$.
Lemma \ref{uniform properness} guarantees that 
$f_n(M_n\setminus \psi_n(Z\cup V)),$ which equals $f'_n(M_n\setminus \psi_n(Z\cup V))$, is
disjoint from $\NN_r(y_n)$, and hence from $f'_n(\BB_1(n))$. Similarly,
$f_n(V)$ is disjoint from $\NN_r(y_n)$ and, since $\max\{ d(f_n(x),f_n'(x))\ |\ x\in \psi_n(V)\}\to 0$,
$f_n'(V)$ is disjoint from $f_n'(\BB_1(n))$ for all large enough $n$.
The definitions of $f_n'$ and $f_\infty'$
guarantee that $f'_n(\psi_n(Z)\setminus \BB_1(n))$ is disjoint
from $f'_n(\BB_1(n))$. Thus $f'_n(M_n\setminus \BB_1(n))$ is disjoint from $f'_n(\BB_1(n))$
as desired.

Finally, we note that the entire construction can be performed so that
all maps  are Lipschitz with some constant $L$, simply by using
piecewise-smooth maps in the geometric limit. Thus the sequence of
maps $f_n$ does admit $f'_n$ which satisfy the conclusions of the
theorem with Lipschitz constant $L$, and as soon as $L_n>L$ this
contradicts our original choice of sequence. This contradiction
establishes the theorem. 
\end{proof}

\section{Proof of the bilipschitz model theorem}
\label{proofmain}

We are now ready to  put together the ingredients of the previous
sections and complete the proof of our main technical theorem, which
we restate here:

\state{Bilipschitz Model Theorem.}{%
There exist $K',k'>0$ depending only on $S$, so that for any
Kleinian surface group $\rho\in\DD(S)$ with
 end invariants $\nu=(\nu_+,\nu_-)$
there is an orientation-preserving $K'$-bilipschitz
homeomorphism of pairs
$$
F: (\modl,\UU[k']) \to (\hhat C_{N_\rho},\MT[k'])
$$
in the homotopy class determined by $\rho$.
Furthermore this map extends to a homeomorphism
$$
\bar F: \bME_\nu \to \bar N
$$
which restricts to  a $K'$-bilipschitz homeomorphism from $\ME_\nu$ to
$N$, and a conformal isomorphism from $\boundary_\infty \ME_\nu$ to
$\boundary_\infty N$. 
}

In Section \ref{embedding an individual cut} we will apply the
results of \S\ref{uniform embeddings} to obtain embeddings of the cut
surfaces in a cut system. In Section \ref{thinning the cut   system}
we will show how an appropriately thinned-out cut system results in
cut surfaces whose images are disjoint, and in Section 
\ref{preserving order of embeddings} we will show that, once these
adjustments are made our map will preserve topological order on the
cuts. In \S\ref{region control} and \S\ref{tube control} we will
extend the embedding to the complement of the cut surfaces and extend
control to Margulis tubes, thus 
finishing the proof in the special case where 
$\nu_\pm$ are both laminations in $\EL(S)$ (the doubly degenerate
case). The remaining cases, in which the convex core has nonempty
boundary,  will be treated in Section \ref{hybrid}.

\subsection{Embedding an individual cut}
\label{embedding an individual cut}

Let $f:\modl \to \hhat C_N$ be the
$(K,k)$ model map provided by the Lipschitz Model
Theorem (see \S\ref{bilip model intro}).

Recall the Otal constant $\epotal$ from  Theorem \ref{otal},
and the function $\Omega$ from the Short Curve
Theorem in \S\ref{bilip model intro}.
Let $\kotal= \max(k,\Omega(\epotal))$. The Short Curve Theorem guarantees that all the
model tubes with $|\omega|\ge \kotal$ map to Margulis tubes
with length at most $\epotal$. Theorem \ref{otal} guarantees that these
image Margulis tubes are unknotted.

Each surface $\hhat F_{\tau}[\kotal]$ associated to a saturated non-annular
slice $\tau$ is composed of standard
3-holed spheres attached to bounded-geometry annuli. Hence it admits
an $r$-bounded hyperbolic metric $\sigma_\tau$ with geodesic boundary which is
$r$-bilipschitz equivalent to its original metric, for some $r$
depending on $\kotal$. We will henceforth
consider these surfaces with these adjusted metrics. This together
with the $K$-Lipschitz bounds on the model map $f$ tells us that
$f|_{\hhat F_{\tau}[\kotal]}$ is an $\Lone$-bounded map (as in
\S\ref{uniform embeddings}), where $\Lone$ depends on $\kotal$ and $K$.

\begin{lemma}{embed one cut}
There exist $d_0$, $\Kone$ and $\hat\ep$ (depending only on $S$) such that, 
if $\tau$ is a saturated non-annular slice in $H_\nu$ such that the length
$|g_\tau|$ of 
its base geodesic is at least $d_0$, then
$f|_{\hhat  F_{\tau}[\kotal]}$ is $\epotal$-anchored, and
$(\Kone,\hat\ep)$-uniformly embeddable (with respect to the metric
$\sigma_\tau$ on $\hhat  F_{\tau}[\kotal]$).
\end{lemma}

\begin{proof}
We will check that the conditions of
Theorem \ref{Relative embeddability} hold, and thereby obtain the uniform embeddability. 

Let $\delta>0 $ be sufficiently small that, if $v\in
\CC(S)$ has $\ell_\rho(v) < \delta$, then $v$ is in the hierarchy
$H_\nu$ and $|\omega(v)| > k$. Such a $\delta$ is guaranteed  to exist by the
Short Curve Theorem, parts (\ref{SCL initial}) and (\ref{SCL lower})
(see \S\ref{length estimates}).
Let $\ep\in(0,\epotal)$ be the constant provided by Theorem \ref{Relative
embeddability} for this value of $\delta$ and
$L=\Lone$, and let $\Kone>0$ and $\hat\ep\in(0,\epotal)$ be the uniform embeddability
constants provided by Theorem \ref{Relative embeddability}.
Let $d_0=\LL(\ep)$ where $\LL$ is the function from Lemma \ref{big h short l},
so that $|g_\tau|\ge d_0$ implies that $\ell_\rho(\boundary D(\tau))\le \ep$.

Let $R=D(\tau)$ and let $\Gamma$ be the set of vertices $v$ in $\tau$ with
$|\omega(v)|\ge \kotal$.  By the definition of $\kotal$, $\ell_\rho(v)\le \epotal$
for all $v\in \Gamma$.
The subsurface $X=R\setminus\collar(\Gamma)$ can 
be identified with $\hhat F_\tau[\kotal]$, and this gives it an
$\Lone$-bounded metric. Since $\ep\le\epotal$ and  $\ell_\rho(\boundary D(\tau))\le \ep$,
the map $f|_{\hhat F_{\tau}[\kotal]}$, or $f|_X$,  is $\epotal$-anchored.
Moreover,
$\bar f = f|_{\hhat F_\tau}$
is $\ep$-anchored. Therefore,  $f|_X$  satisfies condition (3) of
Theorem \ref{Relative embeddability}. 
Condition (1) of Theorem \ref{Relative embeddability} follows from
the properties of the $(K,k)$ model map, and the choice of $\kotal$
and $\Gamma$. 

Next we establish the 
unwrapping condition (2) of Theorem \ref{Relative embeddability}. 
If $w$ is any vertex in $\CC_0(R)\setminus\Gamma$ such that $\ell_\rho(w)< \delta$, 
we have $|\omega(w)|>k$ and hence $f$ takes $U(w)$ to $\MT(w)$
and $\modl\setminus U(w)$ to $N\setminus \MT(w)$.
Applying 
Lemma \ref{unwrapping} to $U(w)$ and $\hhat F_\tau$, 
we have in particular that
$\hhat F_\tau$
is homotopic to either $+\infty $ or $-\infty$ in the complement of
$U(w)$, so $f|_{\hhat F_\tau}$ is homotopic to either $+\infty$ or
$-\infty$ in the complement of $\MT(w)$. This is exactly
the unwrapping condition (2).

(When $\modl$ has nonempty boundary, we interpret ``homotopic to $\pm
\infty$  in the model'' by considering $M' = 
\modl \setminus (\boundary \modl \union \UU(\boundary S))$, which is
homeomorphic to $S\times\R$.  
In $N$, we consider $C'=\hhat C_N\setminus (\boundary\hhat C_N \union
\MT(\boundary S))$. Since the model map takes $M'$ properly to $C'$, 
we can make the same arguments.)

Having verified that the conditions of 
Theorem \ref{Relative embeddability} hold, we obtain the desired
$(\Kone,\hat\ep)$-uniform embeddability of $f|_{\hhat F[\kotal]}$. 
\end{proof}

For a nonannular cut $c$ in a cut system $C$ with spacing lower bound at
least $d_0$,  let $G_c:\hhat F_{c}[\kotal] \times[0,1] \to N$
be the homotopy provided by Lemma \ref{embed one cut}, which is a
bilipschitz embedding on 
$\hhat F_{c}[\kotal] \times[1/2,1]$. 
Let $f_c : \hhat F_{c}[\kotal] \to N$ be the embedding defined by 
\begin{equation}
\label{define embedding fc}
f_c(x) = G_c(x,3/4).
\end{equation}

We can extend this map to 
the annuli of $\hhat F_c \setminus \hhat F_c[\kotal]$, if we are
willing to drop the Lipschitz bounds:

\begin{corollary}{fc on annuli}
The homotopy $G_c$ provided by Lemma \ref{embed one cut} can be
extended to a map $\bar G_c:\hhat F_c\times[0,1] \to N$ so that, on each
annulus $A$ in $\hhat F_c\setminus \hhat F_c[\kotal]$ we have
$\bar G_c(A\times[0,1])$ contained in the corresponding Margulis tube $\MT(A)$,
and $\bar f_c(x) = G_c(x,3/4)$ is still an embedding. 
\end{corollary}
\begin{proof}
Note first that, by our choice of $\kotal$, each $A$ indeed has core curve
whose length in $N$ is sufficiently short that $\MT(A)$ is non-empty.
The model map $f$ already takes $A$ into $\MT(A)$.
Since $G_c$ is a homotopy through anchored embeddings, it takes
$\boundary A\times[0,1]$ to $\boundary \MT(A)$, and 
so the existence of $\bar G_c$ is a simple
fact about mappings of annuli into solid tori. 
\end{proof}

\subsubsection*{Annular cuts} If $c$ is an annular cut let $\omega(c)$
denote the meridian coefficient of the corresponding tube $U(c)$, and 
if $|\omega(c)|>\kotal$ let $\MT(c)$ denote the Margulis tube
associated to the homotopy class of the annulus. 
For notational consistency, we let 
\begin{equation}\label{fc for annulus}
f_c  = \MT(c)
\end{equation}
when $|\omega(c)|> \kotal$. As in \S\ref{knotting} we are blurring the
distinction between a map and its image here. 

\subsection{Thinning the cut system}
\label{thinning the cut system}

Lemma \ref{embed one cut} allows us, after bounded homotopy of the
model map, to embed individual slices of a
cut system, but the images of these embeddings may intersect in
unpredictable ways. We will now show that, by thinning out a cut
system in a controlled way we can obtain one for which the cuts that
border any one complementary region have disjoint $f_c$-images.

Because the model manifold is built out of standard pieces, 
for any nonannular slice $c$ there is a paired bicollar 
neighborhood $E^0_c$ for
$(\hhat F_c[\kotal],\boundary\hhat F_c[\kotal])$  in
$(\modl[\kotal],\boundary\UU[\kotal])$ which is
uniformly bilipschitz equivalent to a standard 
product. That is, there is a bilipschitz piecewise smooth
orientation-preserving homeomorphism 
$$ 
\varphi_c:E^0_c \to \hhat F_c[\kotal]\times[-1,1] 
$$
which restricts to $x\mapsto (x,0)$
on  $\hhat F_c[\kotal]$, 
where the bilipschitz constant depends only on the surface
$S$. We are here taking $\hhat F_c[\kotal]\times[-1,1]$ with the
product metric $\sigma_c\times dt$, where $\sigma_c$ is the hyperbolic
metric defined in \S\ref{embedding an individual cut}.
The relative boundary $\varphi_c^{-1}(\boundary\hhat
F_c[\kotal]\times[-1,1])$ is the intersection of $E^0_c$ with
$\boundary\UU[\kotal]$.
Moreover we may choose the collars so that,
if $\hhat F_c[\kotal]$ and $\hhat F_{c'}[\kotal]$ are
disjoint, so are $E^0_c$ and $E^0_{c'}$. Let $E_c$ denote the
subcollars $\varphi_c^{-1}(\hhat F_c[\kotal]\times[-\half,\half])$; our
final map will be an 
embedding on each $E_c$, for an appropriate set of slices $\{c\}$. 

\begin{lemma}{Thin for disjointness}
Given a cut system $C$ with spacing lower bound  $d_1\ge d_0$ and upper bound
$3d_1$,  there exists a
cut system $C'\subset C$ with spacing upper bound $d_2 = d_2(S,d_1)$, and
$\Ktwo = \Ktwo(S)$, such that there is a map of pairs
$$
f': (\modl[\kotal],\boundary \UU[\kotal]) \to 
        (\hhat C_{N_\rho} \setminus \MT[\kotal],\boundary\MT[\kotal])
$$
which is homotopic through maps of pairs to the restriction of the model map 
$f|_{\modl[\kotal]}$, such that
\begin{enumerate}
\item $f'$ is a $(\Ktwo,\kotal)$ model map, and its homotopy to $f$
is supported on the union of collars $\union E^0_c$ over nonannular
$c\in C'$.

\item Inside each subcollar $E_c$ for nonannular $c\in C'$, 
$f'$ is an
orientation-preserving 
$\Ktwo$-bilipschitz embedding, and 
$$
f'|_{\hhat F_c[\kotal]} = f_c.
$$
Furthermore $f'\circ \varphi_c^{-1}$ restricted to 
$\varphi_c(E_c)$ has norms of second derivatives 
bounded by $\Ktwo$, with respect to the metric $\sigma_c\times dt$.
\item 
For each complementary region $W\subset \modl[\kotal]$ of $C'$, the
subcollars $E_c$
of nonannular cut surfaces on $\boundary W$ have disjoint $f'$-images.
\end{enumerate}
\end{lemma}

\begin{proof}
For each $c\in C$ we will use the homotopy
to an embedding provided by Lemma \ref{embed one cut} (via
Theorem \ref{Relative embeddability}) to redefine $f$ in $E^0_c$. 
The resulting map will
immediately satisfy the conclusions of Lemma \ref{Thin for disjointness}
except possibly for the disjointness condition (3).
In order to satisfy (3) we will have to ``thin'' the cut system.

To define the map in each collar, 
fix a nonannular $c\in C$, and 
let $G_c:\hhat F_c[\kotal]\times [0,1] \to N\setminus\MT[\kotal]$
denote the proper homotopy given by  
Lemma \ref{embed one cut} and Theorem \ref{Relative embeddability}, 
where we recall that
$G_c$ restricted to 
$\hhat F_c[\kotal]\times [1/2,1] $ is a $K_1$-bilipschitz embedding
with norm of second derivatives bounded by $K_1$, 
and $f_c(x) = G_c(x,3/4)$.

Define $\sigma: [-1,1]\to [-1,1]$ so that it is affine in the
complement of the ordered 6-tuple
$(-1,-3/4,-1/2,1/2,3/4,1)$, and takes the points of 
the 6-tuple, in order, to 
$(-1,0,1/2,1,0,1)$ (see Figure \ref{sigma}). Note
that $\sigma$ is orientation-reversing in $(1/2,3/4)$ and
orientation-preserving otherwise. 

\realfig{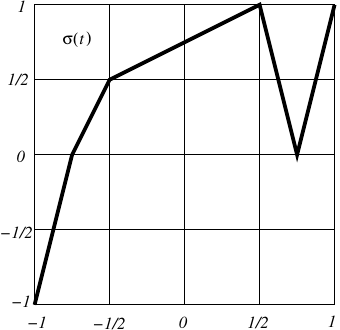}{The graph of the reparametrization function $\sigma(t)$}

Define a map $g_c:\hhat F_c[\kotal]\times[-1,1] \to N$ via
$$
g_c(x,t) = \begin{cases}
  G_c(x,\sigma(t)) &   |t|\le 3/4 \\
  f(\varphi_c^{-1}(x,\sigma(t))) & |t| \ge 3/4.
\end{cases}
$$
Now given a subset $C'\subset C$, let $f'$ restricted to $E^0_c$ 
for $c\in C'$
be $g_c\circ \varphi_c$, and let $f'=f$ on the complement of
$\union_{c\in C'} E^0_c$. 
Since the collars are pairwise disjoint for all $c\in C$, this
definition makes sense. 

One easily verifies that $f'$ satisfies conclusions (1) and (2) of the
lemma. Note in particular that, since 
$\sigma$ takes $[-1/2,1/2]$ to $[1/2,1]$ by an orientation-preserving
homeomorphism, 
$f'$ on $E_c$ is a reparametrization of the embedded part of $G_c$,
and  that since $\sigma(0) = 3/4$ we have 
$f'|_{\hhat F_c[\kotal]} = f_c$.

We will now explain how to choose $C'$ so that
$f'$ will satisfy the disjointness condition (3) as
well.

\medskip

The Margulis lemma gives a bound $n(L)$ on the 
number of loops of length at most $L$ through any one point in a
hyperbolic 3-manifold, if the loops represent
distinct primitive homotopy classes.
This, together with
the Lipschitz bound on $f'$ and the fact that
all vertices in different slices are distinct,
gives a bound $\beta(r)$ on the number of 
$f'(E_c)$ that can touch any given
$r$-ball in $N$.

Let $r_0$ be an upper bound for the diameter of any embedded collar
$f'(E_c)$ (following from the Lipschitz bound on $f'$).

The pigeonhole principle then yields
the following observation, which will be used repeatedly below:

\begin{itemize}
\item[(*)]
 Given a set
$Z\subset C$ of at most $k$ slices, and a set $Y\subset C$ of at least
$k\beta(r_0)+1$ slices, there exists $c\in Y$ such that $f'(E_c)$ is
disjoint from $f'(E_{c'})$ for all $c'\in Z$.
\end{itemize}

We will describe a nested sequence $C=C_0 \supset C_1 \supset
\cdots \supset C_{\xi(S)-3}$ of cut systems, so that $C_j$ satisfies
an upper spacing bound 
$d_2(j)$, and the following condition, 
where for an address $\alpha$ of $C_j$ 
we let $\XX_j(\alpha)$ denote the union of
blocks with address $\alpha$. 
\begin{itemize}
\item[(**)]
For any two slices $c,c'\in C_j$ whose cut surfaces are 
on the boundary of $\XX_j(\alpha)$, and whose
complexities are greater than $\xi(S)-j$,
$f'(E_c)$ and $f'(E_{c'})$ are disjoint. 
\end{itemize}
Thus $C_{\xi(S)-3}$ will
be the desired cut system $C'$, with $d_2 = d_2(\xi(S)-3)$.
Note that by assumption $d_2(0)=3d_1$, and that $C_0$ satisfies
condition (**) vacuously.

We obtain $C_1$ from $C_0$ by letting
$\nslices{h}{C_1} = \nslices{h}{C_0}$ for all $h\ne g_H$, 
and removing slices on the main geodesic $g_H$:
If $\nslices{g_H}{C_0}$ has at most $\beta(r_0)+1$ slices,
set $\nslices{g_H}{C_1} = \emptyset$. In this case $|g_H|$ is at most
$3d_1(\beta(r_0)+2)$.

If $\nslices{g_H}{C_0}$ has at least $\beta(r_0)+2$ slices, 
we partition it into a 
sequence of ``intervals'' $\{J_i\}_{i\in\II}$,
indexed by an interval $\II\subset \Z$ containing $0$; 
by this we mean that, using the cut order $\cprec$ on slices, 
each $J_i$ contains all slices of $\nslices{g_H}{C_0}$
between $\min J_i$ and $\max J_i$, 
and that $\min J_{i+1}$ is the successor to $\max J_i$.
Furthermore, we may do this
so that $J_0$ is a singleton, and for $i\ne 0$ $J_i$ has size 
$\beta(r_0)+1$, except for the largest positive $i$ and smallest
negative $i$ (if any), for which $J_i$ has between $\beta(r_0)+1$ and 
$2\beta(r_0)+1$ elements.  Note that 
$\II$ is infinite if $g_H$ is, for example in the
doubly degenerate case. If $g_H$ and hence $\nslices{g_H}{C_0}$ are
finite, the condition on sizes of $J_i$ is easily arranged by
elementary arithmetic.

Let $J_0=\{c_0\}$. Proceeding inductively, we select some $c_i$ in $J_i$ for
each positive $i$ so that $f'(E_{c_i})$ is disjoint from $f'(E_{c_{i-1}})$.
Our constraints on the sizes of $J_i$, together with observation (*), 
guarantees that this choice is always
possible. Similarly for $i<0$ we choose $c_i$ such that 
$f'(E_{c_i})$ is disjoint from $f'(E_{c_{i+1}})$.

In this case we define $C_1$ so that $\nslices{g_H}{C_1} =
\{c_i\}_{i\in\II}$.
The spacing upper bound in this case is at most 
$d_2(1)\equiv(3\beta(r_0)+1)3d_1$.  We
note that $C_1 $ has the property that for two successive $c,c'$ in
$\nslices{g_H}{ C_1}$, $f'(E_c)$ and $f'(E_{c'})$ are disjoint. Since 
$\hhat F_c$ and $\hhat F_{c'}$ are in the boundary of the same region if
and only if they are consecutive, this establishes the inductive
hypothesis for $C_1$.

We proceed by induction. We will construct $C_{j+1}$ from $C_j$ by
applying a similar thinning process to every geodesic $h$ with
complexity $\xi(h)=\xi(S)-j$, and setting
$\nslices{m}{C_{j+1}} = \nslices{m}{C_j}$ for every $m$ with
$\xi(m)\ne\xi(S)-j$. 
For any address $\alpha$,
the total number of blocks in $\XX_j(\alpha)$ is bounded by
$b_j$ depending on $d_2(j)$, by Lemma \ref{Bound regions}, 
and hence there is a bound $s_j$ on the number of 
cut surfaces on the boundary of $\XX_j(\alpha)$
(certainly $s_j\le 4b_j$, since every block has at most
four gluing boundaries, but a better bound probably holds).

Enumerate the geodesics with complexity $\xi(S)-j$ as
$h_1,h_2,\ldots$, and thin them successively: 
At the $k$th stage we have already, 
for each $h_m$ with $m<k$, thinned $\nslices{h_m}{C_j}$
to obtain $\nslices{h_m}{C_{j+1}}$. 
If $h_k$ contains fewer than
$2s_j \beta(r_0) + 2$ slices, remove them all, obtaining 
$\nslices{h_k}{C_{j+1}}=\emptyset$. 
Otherwise, recall from Lemma \ref{product minus products}
that there is a unique address
$\alpha$ such that $h_k$ is an inner
boundary geodesic for $\XX_j(\alpha)$. 
Let $Q_1$ be the set of all slices $c\in C_j$
with $F_c$ in
$\boundary \XX_j(\alpha)$ and
$\xi(c)>\xi(h_k)$. 
Let $Q_2$ be the set of all slices that arise
as first and last slices in
$\nslices{h_m}{C_{j+1}}$ for each $m<k$ such that
$h_m$ is an inner boundary geodesic for $\XX_j(\alpha)$.
Thus each element of $Q_2$ is the ``replacement'' for a boundary
surface of $\XX_j(\alpha)$ of complexity $\xi(h_k)$ that may have been
removed by the thinning process. 
The bound on the number of boundary surfaces of $\XX_j(\alpha)$ implies
that the union $Q_1\union Q_2$ has at most $s_j-1$ elements. 

Partition $\nslices{h_k}{C_j}$ into a (possibly infinite)
sequence of consecutive, contiguous subsets $\{J_i\}_{i\in\II}$, such that 
the first and the last (if they exist) have length at least
$s_j\beta(r_0)+1$ and at most $(s_j+1)\beta(r_0) + 1$,
and the rest have length  $2\beta(r_0)+1$.

Now let $\nslices{h_k}{C_{j+1}}$ be the union of one cut from each
$J_i$, selected as follows: Supposing that there is a first $J_{i_p}$,
choose a cut $c_{i_p}\in J_{i_p}$ such that $f'(E_{c_{i_p}})$ is
disjoint from $f'(E_b)$ for each $b\in Q_1\union Q_2$.
If there is also a last $J_{i_q}$, choose $c_{i_q}$ such that 
$f'(E_{c_{i_q}})$ is
disjoint from $f'(E_b)$ for each $b\in Q_1\union Q_2\union
\{c_{i_p}\}$. 
Now for each $ i_p < i < i_q-1$ we successively choose
$c_i\in J_i$ so that $f'(E_{c_i})$ is disjoint from 
$f'(E_{c_{i-1}})$. If $i=i_q-1$, we choose $c_i$ so that $f'(E_{c_i})$
is disjoint from $f'(E_{c_{i-1}})$ and $f'(E_{c_{i_q}})$.
Note that all these selections are possible by the choice of sizes of
the $J_i$, the bound on the size of $Q_1\union Q_2$, and observation (*).
If there is a last $J_i$ but not a first, we proceed similarly but in
the opposite direction. (There must be either a last or a first since
$g_H$ is the only geodesic in $H$ that can be biinfinite).

We can then set $d_2(j+1)\equiv d_2(j)(2(s_j+1)\beta(r_0) + 1)$
to be the upper spacing bound for $C_{j+1}$.

To verify that $C_{j+1}$ satisfies the condition (**), 
consider $\XX_{j+1}(\alpha)$ for any address $\alpha=\langle d,d'\rangle$ occuring in
$C_{j+1}$.  
Let $h$ be $g_\alpha$ (as in the proof of Lemma \ref{Bound regions}), 
and denote $\xi(\alpha) = \xi(h)$. 

If $\xi(\alpha) < \xi(S)-j$ then
there is nothing to check since the boundary surfaces of
$\XX_{j+1}(\alpha)$ all have complexity less than $\xi(S)-j$. 

If $\xi(\alpha) = \xi(S) - j$ then only the two outer boundary
surfaces of $\XX_{j+1}(\alpha)$, namely $F_d$ and $F_{d'}$,  have
complexity greater than $\xi(S)-(j+1)$.  
Since in this case $h$ 
participated in the thinning step we just completed, and $d$ and $d'$
are successive slices in $\nslices{h}{C_{j+1}}$, 
we have $f'(E_d)$ disjoint from $f'(E_{d'})$. 

If $\xi(\alpha) > \xi(S) - j$, then the address $\alpha$ is also an
address of $C_j$, since $h$ was
not thinned in the construction of $C_{j+1}$. The outer boundary
surfaces of $\XX_j(\alpha)$ and $\XX_{j+1}(\alpha)$ are the same, namely
$F_d$ and $F_{d'}$. 
Now consider any inner boundary geodesic $m$ for $\XX_{j+1}(\alpha)$.
If $\xi(m) > \xi(S)-j$ then $m$ was not thinned in this step, and
hence $\nslices{m}{C_{j+1}} = \nslices{m}{C_j}$.
If $\xi(m) = \xi(S)-j$ then
$\nslices{m}{C_{j+1}} \subset \nslices{m}{C_j}$. 
In either case,  $m$ is an inner boundary geodesic for $\XX_j(\alpha)$
as well, since any address pair of $C_j$ in which slices on $m$ are nested
must have $\xi>\xi(m)$ and hence was not removed in this step.  

Now for any  boundary geodesics $m$ and $m'$ of $\XX_{j+1}(\alpha)$
with complexities at least $\xi(S)-j$, 
and slices $c$ on $m$ and $c'$ on $m'$ corresponding to boundary surfaces
of $\XX_{j+1}(\alpha)$, we must show $f'(E_c)$ and $f'(E_{c'})$ are
disjoint. If $\xi(m)>\xi(S)-j$ and $\xi(m')>\xi(S)-j$ then $c$ and
$c'$ correspond to boundary surfaces of $\XX_{j}(\alpha)$, and we have
disjointness by induction. If $\xi(m)=\xi(S)-j$ and $\xi(m')>\xi(S)-j$
then, when the slices on $m$ are thinned to yield
$\nslices{m}{C_{j+1}}$, the slice $c'$ is in $Q_1$, and by the
construction we have disjointness. If $\xi(m)=\xi(S)-j$ and
$\xi(m')=\xi(S)-j$, then if $m=m'$ then $c$ and $c'$ are the first and
last slices of $\nslices{m}{C_{j+1}}$, so again the construction 
makes $f'(E_c)$ and $f'(E_{c'})$ disjoint. Finally if $m\ne m'$
we may suppose that $m'$ is thinned before $m$, and then at the point
that $m$ is thinned we have $c'$ as one of the slices in $Q_2$, so
again the construction gives us disjointness. 

This gives the disjointness property for all boundary surfaces of
$\XX_{j+1}(\alpha)$ of complexity at least $\xi(S)-j$, which establishes
property (**) for $C_{j+1}$.
\end{proof}

\subsection{Preserving order of embeddings}
\label{preserving order of embeddings}

Let $C$ be a cut system with spacing lower bound at least $d_0$, and
such that $|\omega(c)|> \kotal$ for each annular $c\in C$. 
Let $f_c$ and $G_c$ be as in \S\ref{embedding an individual cut}.
The following lemma states that, if the spacing of $C$ is large enough then for
slices with overlapping domains {\em and with disjoint $f_c$-images},
topological order in the image is equivalent to the cut order $\cprec$.

\begin{lemma}{preserve cut order}
There exists a $d_1\ge d_0$ such that, if $C$ is a cut system with spacing
lower bound of $d_1$ and $|\omega(c)|>\kotal$ for every annular slice,  
$c$ and $c'$ are  two slices in $C$ such that $\check D(c)
\intersect \check D(c') \ne 
\emptyset$ 
and  $f_c$ and $f_{c'}$ are  disjoint, then
$$
c\cprec c' \implies f_{c} \topprec f_{c'}.
$$
\end{lemma}

\begin{proof}
Consider first the case where both $c$ and $c'$ are nonannular. 
Since $G_c$ and $G_{c'}$ are $(K_1,\hat\ep)$ uniform homotopies,
$f_c$ and $f_{c'}$ avoid $\MT(\gamma)$ whenever $\gamma$ has length
less than $\hat\ep$,
and the homotopies $G_c$ and $G_{c'}$
stay out of the $\hat\ep$-Margulis tubes. 

Set $k_1 = \max((k+b_1)/b_2,\LL(\hat\ep))$ where $b_1$ and $b_2$ are the constants in
Lemma \ref{big h big omega}, and $\LL$ is as in Lemma \ref{big h short
  l}.
Lemmas \ref{big h big omega} and \ref{big h short l} guarantee that
if $w$ is a component of $\boundary W$
with $W$ supporting a geodesic in $H_\nu$ with length greater than
$k_1$, then  $|\omega(w)|\ge k$ and $\ell_\rho(w)\le\hat\ep$.
Let $k_2\ge k_1$ be the constant produced by Theorem \ref{product
region} for this $k_1$, and for $Q=\Kone$.
When we obtain (within the domain of a slice) a subdomain with
geodesic of length at least $k_1$, we will get deep tubes which we can
use to control the topological ordering (Case 1a). When this does not
happen (Case 2a), we will obtain $(k_1,k_2)$-thick segments, and
apply Theorem \ref{product region} to get  geometric 
product regions which we can again use to control the ordering. 
We choose our new lower spacing bound to be $d_1 = \max(d_0,k_2+14))$.

\subsubsection*{Case 1:}
$D(c) = D(c')$. Thus the slices have a common bottom geodesic
$h$, and the base simplices satisfy $v_c < v_{c'}$.

The idea now is to find an intermediate subset in $\modl$ between
$\hhat F_c[\kotal]$ and $\hhat F_{c'}[\kotal]$, whose image will
separate their images in $N$, and force them to be in the
correct order. This separator will either be a Margulis tube with
large coefficient $\omega$, in case 1a below, or a ``product region'' 
isotopic to $D(c) \times [0,1]$ in case 1b. 

\subsubsection*{Case 1a:}
Suppose that there is some geodesic $m$ with $D(m)\subset D(h)$,
$|m| > k_1$  and such that
$\phi_h(D(m))$ is at least 5 forward of 
$v_c$ and at least 5 behind $v_{c'}$. 
There is at least one
boundary component $w$ of $D(m)$ which is nonperipheral in $D(h)$. 
Let $a$ be an annular slice with domain $\collar(w)$. Then 
$\{c,c',a\}$ satisfies the conditions of a cut system (with upper
spacing bound $d_2=\infty$). 
The footprint $\phi_h(w)$ contains $\phi_h(D(m))$, so by our choice of $m$
and the fact that footprints have diameter at most 2, we know that
$v_c$ is at least 3 behind $\min \phi_h(w)$ and
$v_c'$ is at least 3 ahead of $\max \phi_h(w)$. 
This implies that $c \wallto a \towall c'$, hence
$c\cprec a \cprec c'$.

Now Proposition \ref{topprec and cprec} implies that 
$$
\hhat F_c \topprec U(w) 
$$
and 
$$
U(w) \topprec \hhat F_{c'}.
$$

By our choice of $k_1$, $|m|>k_1$ implies that $\ell_\rho(w) <
\hat\ep$, and that 
$|\omega(w)|>k$ which implies
that $f(U(w))=\MT(w)$ and
that $f(\modl\setminus U(w))=\hhat C_N \setminus \MT(w)$.
Therefore, 
$f|_{\hhat F_c}$ is homotopic to $-\infty$ in the complement
of $ \MT(w)$, 
and $f|_{\hhat F_{c'}}$ is homotopic to $+\infty$ in the complement
of $\MT(w).$  (We make sense of this in the case when
$\boundary\modl\ne\emptyset$ just as in the proof of Lemma \ref{embed
  one cut}.)
Now since $G_c$ and $G_c'$ 
miss $\MT_{\hat \ep}(w)$, we find that also $f_c$ 
is homotopic to $-\infty$ in the complement
of $ \MT(w)$, 
and $f_{c'}$ is homotopic to $+\infty$ in the complement
of $\MT(w)$.
Let $\bar G_c$ and $\bar f_c$ be the extensions of $G_c$ and $f_c$
to $\hhat F_c$  given by Corollary \ref{fc on annuli}. Since they
differ from $G_c$ and $f_c$ only in tubes associated to vertices of
$c$, which are all disjoint from $U(w)$, we may conclude that $\bar f_c$
is homotopic to $-\infty$ in the complement of $\MT(w)$. Define $\bar
f_{c'}$ similarly and note that it is homotopic to $+\infty$ in the
complement of $\MT(w)$. 
Now since $\bar f_c$ and $\bar f_{c'}$ are embedded
surfaces anchored on Margulis tubes which are unknotted and unlinked
by Otal's theorem, we may apply Lemma \ref{embedded ordered}
to conclude that
$$
\bar f_c  \topprec \MT(w)$$
and 
$$\MT(w) \topprec \bar f_{c'}.$$
Now apply 
Lemma \ref{simple transitivity}, with
$R_1=\bar f_c$, $R_2=\bar f_{c'}$, $\VV = \overline\MT(\boundary D(h))$, 
and $Q=\overline\MT(w)$, to conclude
$$
\bar f_c \topprec \bar f_{c'}.
$$
Thus by Lemma \ref{subsurface inherits}
$$
 f_c \topprec  f_{c'}
$$

\subsubsection*{Case 1b:}  If Case 1a does not hold, then
for every geodesic $m$ with $D(m)\subset
D(h)$ such that 
$\phi_h(D(m))$ is at least 5 forward of 
$v_c$ and at least 5 behind $v_{c'}$, we must have  $|m| \le k_1$.
Let $\gamma$ be the 
subsegment of $[v_c,v_c']$ which excludes
7-neighborhoods of the endpoints. Then since $d_1 \ge k_2 + 14,$
$\gamma$  satisfies the $(k_1,k_2)$-thick condition of Theorem
\ref{product region}. Thus, Theorem
\ref{product region} provides 
slices 
$
\tau_{-2},\tau_{-1},\tau_0,\tau_1,\tau_2
$
with bottom
geodesic $h$ and bottom simplices in $\gamma$ satisfying
$$
v_{\tau_{-2}} < v_{\tau_{-1}} < v_{\tau_{0}} < v_{\tau_{1}} <
v_{\tau_{2}},
$$ 
with spacing of at least 5 between successive simplices, 
so that $f$ can be deformed, by a homotopy supported on
$\BB_2 = \BB(\tau_{-2},\tau_2)$, to an $L$-Lipschitz map $f'$ such that
$f'$ is an orientation-preserving embedding on $\BB_1 =
\BB(\tau_{-1},\tau_1)$, and $f'$ takes $\modl \setminus \BB_1$ to 
$N\setminus f'(\BB_1)$.

Since $v_{\tau_{-2}}$ and $v_{\tau_2}$ are at least 5 away
from  $v_c$ and $v_c'$ we may conclude 
that $\{c,\tau_{-2},\ldots,\tau_2,c'\}$ form a cut system (again with $d_2=\infty$), and
that $c\cprec \tau_{-2}$ 
and $\tau_2 \cprec c'$. Proposition \ref{topprec and cprec} now implies
that $\hhat F_c \topprec \hhat F_{\tau_{-2}}$ and
$\hhat F_{\tau_{2}} \topprec \hhat F_c' $.

This implies that $\hhat F_c$  can be pushed to $-\infty$ in $\modl$
in the complement of $\BB_2$.
Applying $f'$, we 
find that $f'|_{\hhat F_{c}}=f|_{\hhat F_{c}}$ may be pushed to $-\infty$ in
$N$ in the complement of $f'(\BB_1)$.

Since we invoked Theorem \ref{product region} with $Q=\Kone$, 
part (3) of that theorem tells us that 
$f'(\BB_1)$ contains a $\Kone$ neighborhood of $f'(\hhat F_{\tau_0})$
in $N\setminus \MT(\boundary D(h))$, and since the tracks of 
the homotopy $G_c$ have length at most $\Kone$, we may conclude
that $G_c$ avoids $f'(\hhat F_{\tau_0})$.

Again let $\bar G_c$ and $\bar f_c$ be the extensions of $G_c$ and $f_c$
to $\hhat F_c$ given by Corollary \ref{fc on annuli}. Each annulus
$A$ of $\hhat F_c\setminus \hhat F_c[\kotal]$ corresponds to an
annular slice $a$ such that, 
similarly to the argument in case 1a,
$\{a,\tau_{-2}\}$ form a cut system with $d_2=\infty$ where
$a\cprec \tau_{-2}$. Thus by Proposition \ref{topprec and cprec}
we have $U(a) \topprec \hhat F_{\tau_{-2}}$ and it follows
that $U(a)$
lies outside $\BB_2$. Hence, its image tube $\MT(A)$ lies outside
$f'(\BB_1)$. 

It follows that the extended homotopy $\bar G_c$ avoids
$f'(F_{\tau_0}).$ Thus
$\bar f_c$ can be pushed to $-\infty$ in the complement of 
$f'(\hhat F_{\tau_0}).$

Since $f'|_{\hhat F_{\tau_0}}$  and $\bar f_c$ are disjoint homotopic
embeddings 
anchored on the tubes of $\MT(\boundary D(h))$, 
part (1) of Lemma \ref{simple transitivity} implies they are
$\topprec$-ordered, and so the homotopy of $\bar f_c$ to $-\infty$
tells us that 
$$
\bar f_c \topprec  f'|_{ F_{\tau_0}}.
$$
Arguing similarly with $c'$, we obtain
$$
  f'|_{ F_{\tau_0}} \topprec \bar f_{c'}.
$$

Now we apply part (2) of Lemma \ref{simple transitivity} to conclude that 
$$ \bar f_c \topprec \bar f_{c'}.$$
It follows by Lemma \ref{subsurface inherits} that 
$$ f_c \topprec f_{c'}.$$

\subsubsection*{Case 2:} $D(c)$ and $D(c')$ intersect but are not
equal. In this case we will obtain the correct order by looking at the
tubes on the boundaries of $D(c)$ and $D(c')$.  

Since $c\cprec c'$, we have $\hhat F_c \topprec \hhat F_{c'}$ by
Proposition \ref{topprec and cprec}. Thus, for any component $\gamma'$
of $\boundary D(c')$ which overlaps $\hhat
F_c$, we can deform $\hhat F_c$ 
to $-\infty$ in the complement of $U(\gamma')$.
Since $d_1 > k_2\ge k_1$, we have
$|\omega(\gamma')| \ge k$.
It follows from the properties of the model map
that  $f|_{\hhat F_c}$ and $f|_{\hhat F_c[\kotal]}$ can
be deformed to $-\infty$ in the complement of
$\MT(\gamma')$. 
The choice of $k_1$ also tells us that
$\ell_\rho(\gamma') < \hat\ep$ and so 
the homotopy $G_c$ avoids the core of Margulis tubes $\MT(\gamma')$.
We can conclude that
the embedding $f_c$ can also be deformed to $-\infty$ in the
complement of $\MT(\gamma')$. 

Now the extended homotopy $\bar G_c$ (from Corollary \ref{fc on annuli})
takes each annulus $A$ of $\hhat F_c \setminus \hhat F_c[\kotal]$
to $\MT(A)$, and hence is still disjoint from
$\MT(\gamma')$. Thus 
$\bar f_c$ can also be deformed to $-\infty$ in the complement
of $\MT(\gamma')$,  and Lemma  \ref{embedded ordered} implies that
$$
\bar f_c \topprec \MT(\gamma').$$
(Since, $f_c$  and $f_{c'}$ are $\epotal$-anchored, and $\ell_\rho(\gamma')\le \hat\ep<\epotal$,
Theorem \ref{otal} implies that $\MT(\boundary D(c))\cup\MT(\boundary D(c'))\cup \MT(\gamma')$ is an unknotted and unlinked collection of solid tori.)
Similarly if $\gamma$ is a component of $\boundary D(c)$ which
intersects $D(c')$, we find that
$$
\MT(\gamma) \topprec \bar f_{c'}.$$
Now applying Lemma \ref{boundary ordering}, we conclude that
$$
\bar f_c \topprec \bar f_{c'}.
$$
Again by Lemma \ref{subsurface inherits} we conclude
$$
f_c \topprec f_{c'}.
$$

\subsubsection*{Annular cuts}
It remains to consider the case that at least one of $c$ and $c'$ are
annular.  
Suppose that both are. Since $\check D(c)=D(c)$ and $\check
D(c')=D(c')$ intersect but are not the same, 
Proposition \ref{topprec and cprec} implies that 
$U(c) \topprec U(c')$. 
Thus $U(c)$ can be pushed to $-\infty$ in $\modl \setminus U(c')$ and
$U(c')$ can be pushed to $+\infty$ in $\modl \setminus U(c)$. 
As before, we use the fact that
$f$ takes $\modl\setminus U(c)$ to $N\setminus \MT(c)$, and 
similarly for $c'$, to conclude that $\MT(c)$ and $\MT(c')$ can be 
pushed to $-\infty$ and $+\infty$, respectively, in the complement of
each other. It follows (Lemma \ref{basic ordering properties}) that
$\MT(c) \topprec \MT(c')$ or equivalently $f_c \topprec f_{c'}$. 

Suppose that $c$ is nonannular but $c'$ is annular. 
Since $c\cprec c'$ and the domains overlap, 
we may apply Proposition \ref{topprec and cprec}
to conclude that $\hhat F_{c} \topprec U(c')$.
Now again since $f$ takes $\modl \setminus U(c)$ to $N\setminus
\MT(c)$, we may conclude, using the same argument as in Case 1a above,
that $f_c \topprec \MT(c')$, or equivalently $f_c \topprec
f_{c'}$. The case where $c$ is annular is similar. 
\end{proof}

\subsection{Controlling complementary regions}
\label{region control}

In this section and Section \ref{tube control}, we will assume that
  $\rho$ is {\em 
  doubly degenerate}, i.e. that there are no nonperipheral parabolics
  or geometrically finite ends, $N=\hhat C_N$, and the model has no
  boundary blocks. The difference between this and the general case
  essentially involves taking care with notation and boundary
  behavior, and 
in Section \ref{hybrid} we will explain how to address these issues.

In outline, the argument in \S\ref{region control} is the following. 
The preceding sections give us a cut system $C'$ which cuts up the
model into complementary regions of controlled size, and with the
property that the model map, after adjustment, embeds the cut surfaces
and tubes disjointly, and in an order preserving way. This makes it
possible to apply the scaffold machinery of Section \ref{knotting}
to conclude that the map on each complementary region of $C'$ can be
replaced (after proper homotopy)  by an embedding as well. Such an embedding $\Phi$ does not
come with any uniform geometric control. On the other hand, the
presence of the Lipschitz model map allows us to obtain a map  $\Psi$
for each region which agrees with $\Phi$ on the boundary and admits
uniform Lipschitz bounds. To obtain from this a bilipschitz embedding,
we must make a geometric limit argument in which we argue by
contradiction as usual and must take some
care to be able to use both the topological properties of $\Phi$ and
the geometric properties of $\Psi$ in a limiting picture. 
We then put all the complementary regions together to get a locally
bilipschitz map, of the correct degree and homotopy class, on all of 
$\modl[\kotal]$. 

In \S\ref{tube control} we extend the bilipschitz control to the
remaining Margulis tubes $\UU[\kotal]$. This involves extending bilipschitz maps from the
boundaries of hyperbolic tubes to their interiors in a uniform way,
which is slightly trickier than one might at first suppose, but not
enormously difficult.

\medskip

From now on, assume $d_1$ has been chosen to be at least as large
as the constant $d_1$ given by Lemma \ref{preserve cut order}, {\em
and} the constant $d_1$ given by Proposition \ref{modl[k] regions}
when $k=\kotal$. 
Let $C$ be a $(d_1,3d_1)$  cut system, which exists by Lemma \ref{cut systems exist},
and which furthermore satisfies the condition that its annular cuts
correspond exactly to those curves $a$ such that $|\omega(a)|>\kotal$. 
Let $C'$ be the $(d_1,d_2)$ cut system obtained by applying
Lemma \ref{Thin for disjointness} to $C$. Note that the annular slices
of $C'$ are the same as those of $C$. 

Let $W$ be (the closure of) a complementary region of the union of non-annular
$\hhat F_c$ and $U(c)$ for all slices $c$ in $C'$.
Note that, 
because of our choice of annular slices, 
$W$ is also (the closure of) a complementary region of the
union of nonannular surfaces $\{\hhat F_c[\kotal]: c\in
C'\}$ and solid tori $\UU[\kotal]$.
That is, $int(W)$ is the closure of a connected component of $\WW_{\kotal}$
as defined in  \S\ref{filled regions}. By Proposition \ref{modl[k]
  regions},  $int(W)\intersect \modl[0]$ is a component of $\WW_0$, and in
particular  
every block in it has the same address by Lemma \ref{region address}. The
number of blocks in $W$ is uniformly bounded by Lemma \ref{Bound regions}.

Let $\Sigma$ be the scaffold in $\modl$ whose surfaces $\FF_\Sigma$ 
are components of the cut surfaces $\hhat F_c[\kotal]$  associated to non-annular slices $c\in C'$ that meet $\boundary W$, and whose 
solid tori $\VV_\Sigma$ are the closures of $\UU_\Sigma$, which are 
those components of 
$\UU[\kotal]$ whose closures meet $\boundary W$. By construction, if $U\in \UU_\Sigma$
then $U=U(c)$ for an annular slice $c\in C$.

Lemma \ref{consistency} and  Proposition \ref{topprec and cprec} imply that
$\topprec|_\Sigma$ satisfies the overlap condition, 
and by Lemma \ref{Topological Partial Order}, 
the transitive closure of $\topprec|_\Sigma$
is a partial order. Hence $\Sigma$ is combinatorially straight.

We want to consider $f'|_\Sigma$ as a good scaffold map. 
The first step is to 
identify $\modl$ with $N$ by an orientation-preserving homeomorphism
in the homotopy class 
of $f'$, so that from now on we may consider $f'$ to be homotopic to
the identity. 
By Lemma \ref{Thin for disjointness}, $f'$ is an embedding on
$\FF_\Sigma$, and the images of components of $\FF_\Sigma$ are all disjoint.
$f'(\VV_\Sigma)$ is a subcollection of the closed Margulis tubes
$\overline\MT[\kotal]$ which we denote $\overline\MT_\Sigma$, 
and is unknotted and unlinked by Otal's theorem. Hence $f'(\Sigma) =
f'(\FF_\Sigma) \union \overline\MT_\Sigma$ is a scaffold. 

Finally, Lemma \ref{preserve cut order} tells us that
$f'|_\Sigma$ is order preserving. To see this, let $p$ and $q$ be two
overlapping pieces of $\Sigma$ and let us show that
$p\topprec q \iff f'(p) \topprec f'(q)$. 
$p$ and $q$ are components of $\hhat F_c[\kotal]$ and
$\hhat F_{c'}[\kotal]$ for two slices $c,c'\in C'$, respectively
(where if $p$ or $q$ is a tube then the corresponding slice $c$ or $c'$ is
annular and $\hhat F_c[\kotal]=U(c)$ or $\hhat F_{c'}[\kotal]=U(c')$).
The overlap implies that $\check D(c)$ and $\check D(c')$ overlap, and
hence $c$ and $c'$ are $\cprec$-ordered by Lemma \ref{consistency}. 
If $c\cprec c'$ then 
$\hhat F_c \topprec \hhat F_{c'}$ by Proposition \ref{topprec and cprec}, 
and $f_c \topprec f_{c'}$ by 
Lemma \ref{preserve cut order}. For the components $p$ and $q$ this
implies
$p\topprec q$ and $f'(p)\topprec f'(q)$. 
If $c'\cprec c$ then the opposite orders hold in both the model and
the image. Therefore $f'|_\Sigma$ is order preserving.

This establishes all the
properties of Definition \ref{def good map}, 
and hence $f'|_\Sigma$ is a good scaffold map. 

By the properties of the model map, 
we also know that $f'(\modl\setminus \UU_\Sigma)$ is contained in
$N\setminus \MT_\Sigma$, and that $f'$ is proper and has degree 1. 
We can therefore 
apply Theorem \lref{Scaffold Extension} to find a 
homeomorphism of pairs
$$f'': (\modl,\VV_\Sigma)\to (N,\overline\MT_\Sigma)$$
which agrees with $f'$ on $\FF_\Sigma$ and
is homotopic to it, rel $\FF_\Sigma$, on each component of
$\VV_\Sigma$ (through proper maps to the corresponding component of
$\overline\MT_\Sigma$). 

We can now use the existence of $f''$ to obtain maps with geometric
control. We will find maps $\Phi$ and $\Psi$ from a neighborhood of 
$W$ to $N$ homotopic to $f''|_W$,
such that:
\begin{itemize}
\item $\Phi$ is an embedding, agrees  with $f''$ on
  $\FF_\Sigma$, is isotopic to $f''$ on $\boundary\VV_\Sigma$ rel
  $\FF_\Sigma$,  
  satisfies a uniform bilipschitz bound on a uniform bicollar of
  $\boundary W$, and respects the horizontal foliations on
  $\boundary\VV_\Sigma $ and $\boundary \MT_\Sigma$. 
\item $\Psi$ agrees with $\Phi$ on $\boundary W$, satisfies a uniform
  bilipschitz bound on a uniform bicollar of $\boundary W$, 
  and is uniformly Lipschitz in $W$.
\end{itemize}

Here a ``uniform bound'' is a bound independent of any of the data
except the topological type of $S$. A uniform bicollar
is the image of a piecewise-smooth embedding of $\boundary
W\times[-1,1]$ into $N$ with 
uniform bilipschitz bounds, so that $\boundary W\times\{0\}$
maps to $\boundary W$ and $\boundary W\times[0,1]$ maps into $W$.
Recall that the horizontal foliation on $\boundary \VV_\Sigma$ is the
foliation by Euclidean geodesic circles homotopic to the cores of
the constituent annuli, and the geodesic circles homotopic to their
images form the horizontal foliation of 
$\boundary\MT_\Sigma$. 

We remark that $\Phi$ is an embedding but not Lipschitz, whereas
$\Psi$ is Lipschitz but not an embedding. Converting these two maps
into a uniformly bilipschitz embedding will be our goal after
constructing them. 

\subsubsection*{Construction of $\Phi$}
To construct $\Phi$ from $f''$, we begin with $\boundary \VV_\Sigma$. 
Let $V$ denote a component of $\VV_\Sigma$ and $\MT_V$ its image under
$f''$. We claim that  $f''|_{\boundary V}$ is homotopic, through maps
$\boundary V \to \boundary \MT_V$, 
to a uniformly bilipschitz homeomorphism, 
where the homotopy is constant on $\FF_\Sigma \intersect \boundary V$, 

Consider first a component annulus $A$  of $\boundary V\setminus
\FF_\Sigma$. $f''|_A$ is an embedding into $\boundary \MT_V$ which is
homotopic to $f'|_A$ rel boundary. The height of $A$ in $\modl$ is
uniformly bounded since $W$ consists of boundedly many blocks by 
Lemma \ref{Bound regions}. Since $f'$ is uniformly Lipschitz, this
bounds the height of $f''(A)$ from above.
Since $f'$ on $\boundary A$ is a
(uniformly) bilipschitz bicollared embedding, the height of
$f''(A)$ is also uniformly 
bounded below. We conclude that $f''|_A$ is isotopic rel $\boundary
A$ to a bilipschitz embedding  with uniform constant. Since $f''$
already takes $\boundary \FF_\Sigma$ to geodesics in
$\boundary\MT_\Sigma$ by Theorem \ref{Relative embeddability}, 
this bilipschitz embedding can be chosen to respect the horizontal
foliations. We let $\Phi|_A$
be this embedding. Piecing
together over all the components of $\boundary V\intersect \boundary W
\setminus \FF_\Sigma$, we obtain a map which is an
embedding into $\boundary\MT_V$, because 
$\Phi(A)=f''(A)$ for each component $A$, and $f''$ is a homeomorphism.

Now consider the possibility that $\boundary U$ does not meet
$\FF_\Sigma$. We claim that the Euclidean tori $\boundary U$ and 
$\boundary \MT_U$ admit uniformly bilipschitz affine
identifications  with the standard torus $\R^2/\Z^2$.
For $\boundary U$ this follows because it is
composed of a bounded number of standard annuli
(again because of the bound on the size of $W$).
Since $f':\boundary U \to \boundary \MT_U$ is a Lipschitz
map that is homotopic to a homeomorphism (namely $f''|_{\boundary
  U}$),
the diameter of $\boundary\MT_U$ is uniformly bounded from above, and on
the other hand its area is uniformly bounded from below
since it is an $\epone$-Margulis tube
boundary (see e.g. \cite[Lemma 6.3]{minsky:torus}).
It follows that $\boundary\MT_U$ is also uniformly bilipschitz
equivalent to the standard torus.
These identifications conjugate $f'|_{\boundary U}$ to a uniformly
Lipschitz self-map of the standard torus which is homotopic to a
homeomorphism. It is now elementary to check that such a map can be
deformed to a uniformly bilipschitz affine map. 
In fact the homotopy can be chosen so that the tracks of all points
are uniformly bounded.

Thus we have defined $\Phi$ on $\boundary W$, so that it is
isotopic to $f''|_{\boundary W}$ through maps taking $\boundary W$ to 
$f''(\boundary W)$, and constant on $\FF_\Sigma$.

In order to extend the map to $W$, we observe first that $\boundary W$
has a uniform bicollar in $\modl$ (by the explicit construction of the model
manifold), and next that $\Phi(\boundary W) = f''(\boundary W)$ also
has a uniform bicollar in $N$. To see the latter, note that the cut
surface images $f''(\FF_\Sigma)$ have uniform bicollars which are just
the images by $f'$ of the collars $\{E_c\}$ in $\modl$ given by Lemma
\ref{Thin for 
  disjointness}. The boundary tori $\boundary \MT_\Sigma$ have uniform
bicollars because of the choice of $\epone$.
These collars can be fitted together to obtain a uniform bicollar of all
of $\Phi(\boundary W)$ because the pieces of $\Phi(\boundary W)$ fit
together at angles that are bounded away from 0 (due to the
uniformity of $f'$ on the bicollars $E_c$).
By standard methods we may now use the isotopy between $f''$
and $\Phi$ on $\boundary W$ to extend $\Phi$ to a map on $W$ which is
an embedding isotopic to $f''$, with uniform bilipschitz bounds on a
uniform subcollar of the boundary.  

\subsubsection*{Construction of $\Psi$}
We observe that the homotopy from $f'$ to $\Phi$ on $\boundary W$ can
be made to have uniformly bounded tracks, simply by taking the
straight-line homotopy in the Euclidean metric on $\boundary
\UU_\Sigma$. First note that the homotopy is constant except on
$\boundary W\intersect \boundary\UU_\Sigma$.
Let $U$ be a component of $\UU_\Sigma$.
If $\boundary U$ does not meet $\FF_\Sigma$ it has
bounded geometry, and the boundedness of tracks was already noted
above. If $\boundary U$ does meet $\FF_\Sigma$, the homotopy is constant 
by construction on a subset $X$ of $\boundary U$
which is a 
uniformly bounded distance from any point in $\boundary U$, and are
homotopic rel $X$. If $y\in
\boundary U$ let $\alpha$ denote a shortest arc from $y$ to $X$. The
union of $f'(\alpha)$ and $\Phi(\alpha)$ has uniformly bounded length since
both maps are uniformly Lipschitz, and this serves to bound the
shortest homotopy from $f'$ to $\Phi$.

Now let $\Xi$ be a uniform collar of $\boundary W$,  
such that there is a uniformly bilipschitz homeomorphism $h:W\setminus
\Xi \to W$  isotopic to the inclusion
(this is possible because the geometry of
$W$ is uniformly bounded). Define $\Psi|_{W\setminus \Xi} = f'
\circ h$, then extend $\Psi$ to $\Xi$ using the bounded-track homotopy
between $f'|_{\boundary W}$ and $\Phi|_{\boundary W}$. This map agrees
with $\Phi$ on 
$\boundary W$, and satisfies a
uniform Lipschitz bound.  Using the uniform collar structure for
$\boundary W$ and $\Phi(\boundary W)$, as in the construction of
$\Phi$, we can arrange for $\Psi$ to also satisfy uniform bilipschitz
bounds in a uniform collar of the boundary.

\subsubsection*{Uniformity via geometric limits}
We now have a uniform bilipschitz embedding of
$\boundary W$  which extends, by $\Phi$, to an embedding without
geometric control, and by $\Psi$ to a uniformly Lipschitz map which
may not be an embedding.
We claim next that $\Phi|_{\boundary W}=\Psi|_{\boundary W}$ can be
extended to an embedding of $W$ in $N$ with {\em uniform bilipschitz
constant.} 

If this is false, then there is
a sequence of examples $\{(M_{\nu_n}, W_n, N_n)\}$ 
where the best bilipschitz constant goes to infinity (we index 
our maps as $\Phi_n$, $\Psi_n,$ etc.). We shall reach a
contradiction by extracting a geometric limit. 

As before, 
$W_n$ contain a bounded number of blocks. Since the tubes 
whose interiors meet $W_n$ must 
have bounded coefficient $|\omega|< \kotal$, 
we may assume, after
restricting to a subsequence, that they have the same combinatorial
structure and tube coefficients. 
After applying a sequence of homeomorphisms to the model
manifolds we may assume that the $W_n$'s are all equal to a fixed $W$.
Choose a basepoint $x\in W$ and an orthornormal baseframe $\hat x$ for $T_xW$ and
let $y_n = \Psi_n(x)$ and $\hat y_n=d\Psi_n(\hat x)$.
After taking subsequences we 
may  assume that $\{(N_n,\hat y_n)\}$ converges geometrically to $(N_\infty,\hat y_\infty)$,
and that $\{\Psi_n\}$ converges geometrically  to a map $\Psi_\infty:W\to N_\infty$ (the latter
because of the uniform Lipschitz bounds on $\Psi_n$).

Because $\Psi_n|_{\boundary W}$ are uniformly bicollared embeddings, their 
limit $\Psi_\infty|_{\boundary W}$ is an embedding.
  
Since $\Psi_\infty(W)$ is a compact 3-chain with boundary
$\Psi_\infty(\boundary W)$, we know that $\Psi_\infty(\boundary W)$
bounds some compact region $W'_\infty\subset N_\infty$. Similarly let
$W'_n$ be the compact region bounded by $\Psi_n(\boundary W)$ (note
that  $W'_n = \Phi_n(W)$).
  
By definition of geometric convergence, given $R$ and $n$ large enough
there is a map $h_n:\NN_R(y_n) \to N_\infty$ which is an embedding
with bilipschitz constant going to 1, and taking the baseframe $\hat y_n$ to
$\hat y_\infty$. Geometric convergence of the maps means, taking $R$ larger than the
diameter of $W'_\infty$, that $h_n\circ \Psi_n$ converge pointwise
to $\Psi_\infty$ on $W$.

In fact we can arrange things so that eventually $h_n\circ\Psi_n =
\Psi_\infty$ on the boundary: 
note that $\Psi_\infty(\boundary W)$ is composed of finitely many pieces
(images of cut surfaces and annuli in Margulis tubes) which are
$C^2$-embedded, and meet transversely along boundary circles. Thus 
it has a collar neighborhood which is smoothly foliated by intervals
which $\Psi_\infty(\boundary W)$ intersects transversely. Since the
convergence of $h_n\circ\Psi_n(\boundary W)$ is $C^2$ on each cut
surface and annulus piece, they are eventually transverse to this
foliation too, and hence 
after adjusting $h_n$ by small isotopies of this collar neighborhood
we may assume that
$h_n\circ \Psi_n = \Psi_\infty$ on $\boundary W$.
With this
adjustment, we have $h_n(W_n') = W_\infty'$, with $h_n$ still
satisfying a uniform bilipschitz bound.

Now given (large enough) $m$ we note that the embedding $\Phi_m:W\to W_m'$
can be  assumed to be bilipschitz with {\em some} constant depending on $m$. 
Fix a value of $m$, and let $g_m = h_m\circ \Phi_m$.
This is a $K_m$-bilipschitz embedding of $W$ to $W'_\infty$, for some
$K_m$,
which restricts to $\Psi_\infty $ on the boundary. Finally, let $g_n =
h_n^{-1}\circ g_m$. Fixing $m$ and letting $n$ vary, we have a
{\em uniformly} bilipschitz sequence of embeddings taking $W$ to the region
$W_n'$ bounded by $\Psi_n(\partial W)$
and restricting to $\Psi_n$ on the boundary. This contradicts our
choice of sequence.  

With this contradiction we therefore conclude that in fact there is a
uniformly bilipschitz extension of $\Phi|_{\boundary W}$ to $W$, as desired. 
Denote this map by $\Theta_W : W \to N$. 

\subsubsection*{Degree of the map}
We claim that $\Theta_W$ maps with degree 1 onto its image. 
Consider first the case that $\VV_\Sigma$ is non-empty, and let $A$
be the intersection of $\boundary \VV_\Sigma$ with $\boundary W$. 
The map $f''$, since it is globally defined and of degree 1, must map
$A$ with degree 1 (and homeomorphically) to its image in $\boundary
\MT_\Sigma$. Since $\Theta_W$ is isotopic to $f''$ on $A$, it also must
map with degree 1. Any embedding of oriented manifolds $g:X\to Y$ which maps a
nonempty subset of $\boundary X$ with degree 1 to its image in
$\boundary Y$ must have degree 1 to its image in $Y$. Applying this to $\Theta_W:W\to
N\setminus\MT_\Sigma$, we conclude that $\Theta_W$ has degree 1 to its image.

Now if $\VV_\Sigma$ is empty, $W$ only meets components of
$\FF_\Sigma$, and hence these components must have no nonperipheral
boundary. Thus $W$ is the region between two slices $\hhat F_c[\kotal]$
and $\hhat F_{c'}[\kotal]$
with domain equal to all of $S$. Assume that $\hhat F_c \topprec \hat
F_{c'}$. If $\Theta_W$ does not have degree 1, it must switch the order of the
boundaries, that is $\Theta_W(\hhat F_{c'}) \topprec \Theta_W(\hat
F_c)$. Since $\Theta_W$ is equal to $f'$ on these  surfaces, this
contradicts Lemma \ref{preserve cut order}.

\subsubsection*{Putting together the maps}
The embeddings $\Theta_W$ can be pieced together over all regions $W$
to yield a global map $F: \modl[\kotal] \to  N$. This is because different
regions meet only along the cut surfaces $\hhat F_c[\kotal]$ ($c\in C'$),
and on these 
each $\Theta_W$ is equal to the original $f'$.

For each tube $U$ in $\UU[\kotal]$, $F|_{\boundary U}$ is homotopic to
$f'|{\boundary U}$ through maps to $\boundary \MT_U$: this is because
for each region $W$
the homotopy from $\Theta_W$ to $f'$ is constant on the boundary
circles of $\boundary W\intersect \boundary U$, so the homotopies can
be pieced together. Thus since $f'$ was defined on $U$ we can extend
$F$ to $U$ (without any geometric control at this point).
The resulting map $F:\modl\to\hhat C_N$ takes $\UU[\kotal]$ to $\MT[\kotal]$, 
and $\modl[\kotal]$ to $N\setminus
\MT[\kotal]$.

We next check that $F$ is in the right
homotopy class. This is not automatic, because in the
geometric limiting step that produced the maps $\Theta_W$, we did not
keep track of homotopy class. However, we note that $F$ agrees with $f'$ on
each of the cut surfaces, and is homotopic to $f'$ on the union of the
cut surfaces with the tube boundaries $\boundary\UU$.

Since we are in the doubly degenerate case $g_H$ is infinite, 
and hence there exists a slice $c\in C'$ with
$D(c) = S$.  We then have a cut surface $F_c$ which projects to $S$
minus the collar of a pants decomposition. The missing annuli can be
found on the boundaries of the tubes adjacent to $\boundary
F_c$. Adjoining these to $F_c$, we find a surface $S'\subset \modl$
which projects to all of $S$. Thus $\modl$ is homotopy-equivalent to
$S'$, and since $F|_{S'}$ is homotopic to $f'|_{S'}$ we conclude that
$F$ is homotopic to $f'$. 

Note that $F$ is a proper map, since the cut surfaces and tubes cannot
accumulate in $N$, and the diameters of images of the regions $W$ are
uniformly bounded. Thus $F$ has a well-defined degree. 
Since each $\Theta_W$ has degree 1 to its image, $F$ has positive
degree. Since it is a homotopy equivalence, the degree must be 1. 
The restriction to $\modl[\kotal]$ is then a uniformly bilipschitz
(with respect to path metrics) orientation-preserving homeomorphism 
to  $N\setminus \MT[\kotal]$.

\subsection{Control of Margulis tubes}
\label{tube control}

It remains to adjust $F$ on the tubes $\UU[\kotal]$ so that it is a
global bilipschitz homeomorphism. 

If $\MT$ is a hyperbolic tube with marking $(\alpha,\mu)$ (where $\mu$
is a meridian and $\alpha$ represents the core curve; see 
\S\ref{model definitions} and \cite{minsky:ELCI})
we let the {\em $\alpha$-foliation} be the foliation
of $\boundary\MT$ whose leaves are Euclidean geodesics in the homotopy
class of $\alpha$.

\begin{lemma}{Tube extension}
Let $\MT_1$ and $\MT_2$ be hyperbolic $\epone$-Margulis tubes with markings
$(\alpha_1,\mu_1)$ and $(\alpha_2,\mu_2)$ (where $\mu_i$ are meridians
and $\alpha_i$ are representatives of the core curve),
and let $h:\boundary\MT_1 \to \boundary\MT_2$ be a
marking-preserving $K$-bilipschitz homeomorphism which takes the
$\alpha_1$-foliation of $\boundary \MT_1$ to the $\alpha_2$-foliation
of $\boundary\MT_2$. Suppose that the radii of the tubes are at least
$a>0$, and that the length of $\alpha_1$ is at most $a'.$

Then $h$ can be extended to a $K'$-bilipschitz homeomorphism $\hhat
h:\MT_1\to\MT_2$, where $K'$ depends on $K, a$ and $a'$. 
\end{lemma}

\begin{proof}
It will be convenient to recall 
Fermi coordinates $(z,r,\theta)$ around a geodesic, where $z$ denotes
length along the geodesic and $(r,\theta)$ are polar coordinates in
orthogonal planes. The hyperbolic metric is given by
\begin{equation}\label{fermi metric}
\cosh^2r dz^2 + dr^2 + \sinh^2 r d\theta^2.
\end{equation}
This metric descends to any hyperbolic tube quotient (where the geodesic
$(z,0,0)$ descends to the core) in the usual way.

We begin by extending $h$ to all but bounded neighborhoods of the
cores of the tubes. 

Let $r_i\ge a$ be the radius of $\MT_i$, and $m_i$ the length of its
meridian. Because $h$ is marking-preserving and $K$-bilipschitz we have $m_1/m_2 \in
[1/K,K]$, and hence 
\begin{equation}\label{sinh r bound}
\sinh r_1/\sinh r_2 \in [1/K,K]
\end{equation}
since $\sinh r_i = m_i$ using (\ref{fermi metric}).

By hypothesis, $r_1,r_2>a$.
Letting $\MT_i(r)$ denote the $r$-neighborhood of the core in $\MT_i$, 
we extend $h$ to a map 
$$
h_1: \MT_1\setminus \MT_1(a/2) \to \MT_2\setminus \MT_2(a/2)
$$
using the foliations $\RR_i$ of $\MT_i$ minus its core by geodesics
perpendicular to the core. 
More precisely, choose an increasing $K'$-bilipschitz homeomorphism
$s:[a/2,r_1] \to [a/2,r_2]$ satisfying 
$\sinh s(r)/\sinh r \in [1/K',K']$, where $K'$ depends on $K$ and $a$
(one can easily do this with an affine map $s$, using a comparison of
$\sinh(x)$ to $e^x/2$).
Let $h_1$ be the unique extension of $h$ which takes $\RR_1$ to $\RR_2$
and takes $\boundary \MT_1(r)$ to $\boundary \MT_2(s(r))$.
The projection $\boundary \MT_i(r)\to \boundary \MT_i(a/2)$ along the
foliation $\RR_i$ is affine and contracts in each direction by a
factor between $\cosh(r)/\cosh(a/2)$ and $\sinh(r)/\sinh(a/2)$
(using (\ref{fermi metric})). Thus, the
properties of $s$ imply that the extension is bilipschitz.

It remains to extend $h_1$ to $h_2:\MT_1(a/2)\to\MT_2(a/2)$.  The
restriction of $h_1$ to $\boundary \MT_1(a/2)$ is bilipschitz with
constant $K''(K,a)$, and we note that $\boundary \MT_1(a/2)$ is a
torus with {\em bounded diameter}. This is true because both
generators in the boundary markings are bounded at radius $a/2$:
$\alpha_1$ is bounded by $a'$ by hypothesis, and the meridian length
at radius $a/2$ is bounded automatically by $2\pi\sinh a/2$, via
(\ref{fermi metric}).

We then use the following lemma:

\begin{lemma}{lipschitz isotopy}
Let $T$ be a Euclidean torus of diameter at most 1. 
Let $f:T\to T$ be a $K$-bilipschitz homeomorphism homotopic to the
identity, which preserves a linear foliation on $T$. 
Then there exists a map 
$$
F: T\times[0,1] \to T\times[0,1]$$
such that $F(\cdot,0) = id$ and $F(\cdot,1) = f$, and $F$ is
$K'$-bilipschitz for $K'$ depending only on $K$.
\end{lemma}

{\em Remark:} One would expect that the condition of preserving a
linear foliation is not necessary in this lemma. However this seems to
be a nontrivial matter. Luukkainen 
\cite{luukkainen:isotopy} has proven such a ``bilipschitz isotopy''
lemma when $f$ is a self-map of $\R^n$ with a bound on
$d(x,f(x))$ for $x\in\R^n$, building on work of Sullivan, Tukia, and V\"ais\"al\"a
\cite{tukia-vaisala:extension,vaisala:concordance,sullivan:bilipschitz}. 
One could try to obtain the result for the torus by considering the
universal cover, but getting equivariance for the isotopy with control
of the bilipschitz constant seems to be difficult.

At any rate with our added condition the proof is elementary:

\begin{proof}[Proof of Lemma \ref{lipschitz isotopy}]
Consider first this one-dimensional version: Let $h:\R\to
\R$ be a $K$-bilipschitz homeomorphism satisfying also
$|h(s)-s| < C$ for all $s\in\R$. The map 
$$
H(s,t) = ((1-t) h(s) + ts, t)
$$
is then a homeomorphism from $\R\times[0,1]$ to itself satisfying
a $K'$-bilipschitz bound (where $K'$ depends on $K$ and $C$) and such
that $H(\cdot,0) = h$ and $H(\cdot,1)=id$.

Now given our map $f$, let $F$ be a lift of $f$ to $\R^2$. Since $f$
is homotopic to the identity and $diam(T) \le 1$, $F$ can be chosen 
so that $|F(p) - p| \le K+2$ for all $p\in\R^2$. $F$ preserves a
foliation which we can assume is the horizontal foliation, so we can
express it as
$$
F(x,y) = (\xi(x,y), \eta(y))
$$
with $\eta:\R\to\R$ $K$-bilipschitz, and $\xi(x,y)$ $K$-bilipschitz in
$x$ for each $y$, and $K$-Lipschitz in $y$ for each $x$.

Now after applying the one-dimensional case to $\eta$ we may assume 
$\eta(y) = y$, and applying it again to $\xi(x,y)$ for each fixed $y$,
we have our desired bilipschitz isotopy. 

Since this construction is evidently invariant under isometries of
$\R^2$, it can be projected back to the torus $T$.
\end{proof}

Using this lemma we can extend $h_1$ to a $K'''$-bilipschitz homeomorphism
from the collar $\MT_1(a/2)\setminus \MT_1(a/4)$
to $\MT_2(a/2)\setminus \MT_2(a/4)$, 
so that on the inner boundary it 
is an {\em affine map} in the Euclidean metric.
We can then extend the map, again using the radial foliation, to the
rest of the solid torus. The bilipschitz control in this last step
follows from a simple calculation in the Fermi coordinates (\ref{fermi
  metric}), and 
depends on the fact that the map  on $\boundary\MT_1(a/4)$ is affine. 
It does not hold for a general bilipschitz boundary map; this was
the reason we needed to apply Lemma \ref{lipschitz isotopy}.
\end{proof}

Our model map, restricted to the boundary of each model tube,
satisfies the conditions of Lemma \ref{Tube extension}.  (Note that
the condition of preserving a linear foliation was supplied in the
construction, which respected the horizontal foliations on model tube
boundaries and their images. The length bound on the generator
$\alpha_1$ also follows from the properties of the model.) Thus we have
the desired bilipschitz extension.

\medskip

The resulting map is now a locally bilipschitz homeomorphism from
$\modl$ to $\hhat C_N$ (which in the doubly degenerate case is all of
$N$). Thus it is globally bilipschitz, and the Bilipschitz Model
Theorem is established in the doubly degenerate case.

\subsection{The mixed-end case}
\label{hybrid}

We will now consider the case of a Kleinian surface group
that is not necessarily doubly degenerate.

The boundary blocks of $\modl$, as described in \S\ref{model
  definitions},  have {\em outer boundaries} which are the
  boundary components of $\modl$. 
These outer boundaries behave essentially like cuts in a cut
  system. In particular in the proofs of Lemmas \ref{WB transitive}
  and \ref{block slice order} we observe that their
topological ordering properties in $\modl\subset \hhat S\times\R$ are
  as we would expect -- i.e. an outer boundary associated to a top
  boundary block lies above all overlapping cut surfaces, and vice versa for a
  bottom boundary block. 
The set $\XX(\varnothing)$ of blocks with address $\add\varnothing$ is
  nonempty in the case with   boundary, and in fact contains all of
  the boundary blocks (see \S\ref{address region structure}).

Theorem \ref{Relative embeddability}
provides us with uniform collars for the cut surfaces that lie in 
$\hhat C_N$, at a distance of at least $a$ from the boundary, where
$a$ is a uniform constant. 
The original model map $f:\modl\to \hhat C_N$ is already
$K$-bilipschitz on the boundaries. Because each boundary component has
a uniform collar in $\modl$ and in $\hhat C_N$, we may adjust the map 
to satisfy a uniform bilipschitz bound in these collars. We may assume
that the uniformly embedded collar obtained in $\hhat C_N$ is within an
$a$-neighborhood of the boundary.
Thus the collars of the cut surfaces are disjoint from 
the boundary collars. 

This tells us that the topological ordering of overlapping cut surfaces and
boundary surfaces is preserved by the adjusted model map $f'$
(generalizing Lemma \ref{preserve cut order}). 

The argument in \S\ref{region control} controlling the map on
complementary regions requires a few remarks.
The complementary regions contained in $\XX(\varnothing)$
will have outer boundary components in their boundary, so these should
be taken as components of $\FF_\Sigma$ for the scaffold $\Sigma$. 
The map $f''$ should take $(\modl, \VV_\Sigma)$ to $(\hhat
C_N,\overline\MT_\Sigma)$, and again the appeal to Theorem \lref{Scaffold
  Extension} is by way of first identifying the interiors of both
manifolds with $S\times\R$. The homotopy from $f'$ to $f''$ can be
assumed to be constant on the uniform collars of the outer
boundaries. The same holds for the construction of $\Phi$ and $\Psi$,
on the regions contained in $\XX(\varnothing)$. 
In the geometric limit step, when $W_n$ contain boundary blocks we
cannot assume that they are all identical after passage to a
subsequence. Boundary blocks do have finitely many combinatorial types,
so we may assume that these are constant on a subsequence. 
The geometry of a block can degenerate: the curves of $\I(H)$ or $\T(H)$
supported on the block can have lengths going to zero. The geometric
limit of a sequence of such blocks can be described as a union of
blocks based on smaller subsurfaces, where the curves whose lengths
vanish give rise to parabolic tubes in the limit. 
Thus we may assume that the $W_n$ minus these tubes 
are eventually combinatorially equivalent to a
fixed $W$ and geometrically converge to it.
This suffices to make the
argument work. 

Section \ref{tube control} on the extension of the bilipschitz map to
Margulis tubes goes through without change, noting that in the
general case there may be parabolic tubes that are not associated to
$\boundary S$, but that extension to these is no harder. Thus, we
obtain a bilipschitz homeomorphism of degree 1
$$
F:\modl \to \hhat C_N.
$$

Checking that $F$ is homotopic to $f$ is again done by
exhibiting a surface $S'$ in $\modl$ which projects to $S$ and
on which $F$ is known to be homotopic to $f$. In the general case
there may not be a single slice $c$ in the cut system with $D(c)=S$;
however we can piece $S'$ together from slices and outer boundaries. 
Let $P_+$ denote the annuli corresponding to parabolics facing the top
of the compact core, as in \S\ref{ends}. A component $Z$ of
$S\setminus P_+$ is either associated to a top outer boundary of
$\modl$, or supports a filling lamination component of $\nu_+$, and
hence a forward-infinite geodesic in $H$. In the latter case there is
a cut  $c_Z$ with domain $Z$ and an associated surface $\hhat
F_{c_Z}$. The union of these boundary surfaces and cut surfaces,
joined together with annuli along the parabolic model tubes associated
to $P_+$, 
gives our desired surface $S'$, and the argument goes through
as in the doubly degenerate case. 
This gives the desired map from $\modl$ to $\hhat C_N$. 
Since the map has not changed on $\boundary\modl$, we can use the
same extension to the exterior $E_\nu$ as given in Theorem
\ref{extended model}, 
so that we obtain the desired map from $\bME_\nu$ to $\bar N$. 
This completes the proof of the Bilipschitz Model Theorem.

\section{Proof of the ending lamination theorem}
\label{ELT proof}

The proof of the Ending Lamination Theorem is now an application of
the Bilipschitz Model Theorem and Sullivan's Rigidity Theorem. We give
the argument first in the surface group case. The general case
requires a bit more care in analyzing the covers associated to the
boundary components of the relative compact core. 

Before proceeding we give a more careful statement of the theorem. 

\state{Ending Lamination Theorem for Incompressible Ends.}{%
Let $G$ be a finitely-generated torsion-free nonabelian group. Let
$\rho_1, \rho_2:G\to \PSL 2(\C)$ be discrete, faithful 
representations whose quotient manifolds $N_{\rho_i}$ have relative
compact cores $(K_i,P_i)$ with $\boundary_0K_i$ incompressible.
If there is a homeomorphism $\phi:(K_1,P_1)\to
(K_2,P_2)$ such that $\phi_*\circ\rho_1 $ is conjugate to $ \rho_2$,
and such that $\phi$ takes the end invariant of $\rho_1$ on each
component of $\boundary_0K_1$ to the end invariant of $\rho_2$ on
its image,  then
there is an isometry $N_{\rho_1} \to N_{\rho_2}$ in the homotopy
class determined by $f$. 
}

\begin{proof} 
First consider the case that $\rho_1$ and $\rho_2$ are Kleinian
surface groups.
In this case we obtain a single model manifold
$\ME_\nu$ from the common end invariants $\nu$ of $\rho_1$ and
$\rho_2$, and the Bilipschitz Model Theorem gives us bilipschitz 
homeomorphisms $F_i : \ME_\nu \to N_{\rho_i}$ 
in the homotopy classes determined by $\rho_1$ and
  $\rho_2$, respectively. We also obtain
extensions  $\bar F_i : \bME_\nu \to \bar
  N_{\rho_i}$, which are homeomorphisms and map
  $\boundary_\infty\ME_\nu$ conformally to $\boundary_\infty
  N_{\rho_i}$. 
The composition $F_2\circ F_1^{-1}$ is therefore in the homotopy
  class of $f$, and lifts to a $K$-bilipschitz
  homeomorphism of $\Hyp^3$ that conjugates $\rho_1$ to $\rho_2$.
Up to possibly conjugating $\rho_2$ by an orientation-reversing
isometry, we may assume this homeomorphism is orientation-preserving.
 It therefore extends to a quasiconformal homeomorphism of $\hhat\C$
  (Mostow \cite{mostow:hyperbolic}), which
  is conformal from the domain of discontinuity of
  $\rho_1$ to that of $\rho_2$.

  Sullivan's
  Rigidity Theorem \cite{sullivan:rigidity} now implies that this map is
  in fact conformal on the whole sphere at
  infinity, and it follows that it is homotopic to an isometry on the
  interior.

\medskip

Before proceeding to the general case, we need the following
corollary of the surface group case, which treats the case where
geometrically finite end invariants do not match.
\begin{lemma}{bilipschitz corollary}
Let $\rho_1,\rho_2:\pi_1(S)\to\PSL 2(\C)$ be Kleinian surface groups
with relative compact cores $(K_1,P_1)$ and $(K_2,P_2)$, and a homeomorphism
$\phi:(K_1,P_1)\to (K_2,P_2)$ such that $\rho_2$ is conjugate to 
$\phi_*\circ\rho_1$. Suppose that each end $\EE$ of 
$N^0_{\rho_1}$ is geometrically infinite if and only if $\phi(\EE)$ is
geometrically infinite, in which case $\nu(\phi(\EE)) = \phi(\nu(\EE))$.
Then $\phi$ extends to a bilipschitz homeomorphism from $N_{\rho_1}$ to $N_{\rho_2}$. 
\end{lemma}

\begin{proof}
One may use the Measurable Riemann Mapping Theorem \cite{ahlfors-bers}
to construct
a quasiconformal map $\psi:\hhat\C\to\hhat\C$ such that $\rho_3=\psi\circ\rho_1\circ \psi^{-1}$ is
a Kleinian surface group and there exists a conformal map from
$\boundary_\infty N_{\rho_3}$
to $\boundary_\infty{N_{\rho_2}}$ in the homotopy class of $\rho_2\circ\rho_3^{-1}$.
The map $\psi$ extends equivariantly to a bilipschitz homeomorphism of
$\Hyp^3$ 
which
descends to a bilipschitz homeomorphism $F:N_{\rho_1}\to N_{\rho_3}$
(see e.g. Douady-Earle \cite{douady-earle}).
One may deform $F$ so that $(K_3,P_3)=(F(K_1),F(P_1))$ is a relative
compact core for $N_{\rho_3}^0$ 
Since $F$ is
bilipschitz, $F|_{K_1}$ preserves the end invariants of the geometrically infinite ends
of $N_{\rho_1}^0$. Thus, $\phi\circ (F|_{K_1})^{-1}:(K_3,P_3)\to (K_2,P_2)$ is
a homeomorphism preserving all the end invariants. The surface group
case of the Ending Lamination Theorem now implies that there exists an isometry
$I:N_{\rho_3}\to N_{\rho_1}$ in the homotopy class of
$\rho_1\circ\rho_3^{-1}$.  Then $I\circ F$ may deformed on a
neighborhood of $K_1$ to yield
the desired bilipschitz homeomorphism.
\end{proof}

\medskip

We now proceed to the proof of the general case.

Let $R_1$ be a component of $ \bdry_0 K_1$ and let $R_2 = \phi (R_1)$ be
its homeomorpic image in $\bdry_0 K_2$.  Letting $N_{R_i}$ be the
surface-group cover of
$N_{\rho_i}$ associated to  $R_i$, the lift of $\phi|_{R_1}$ extends
to an orientation-preserving homeomorphism from $N^0_{R_1}$ to $N^0_{R_2}$. Our next step will be
to replace this with a bilipschitz homeomorphism, and in order to do
this we will examine the end invariants of $N_{R_i}$ and apply the
Bilipschitz Model Theorem. Once this is
done, we will apply it to obtain bilipschitz maps from neighborhoods of
each end of $N^0_{\rho_1}$ to neighborhoods of the corresponding ends
of $N^0_{\rho_2}$. Extending across the remaining compact core and the
cusps we will obtain the desired bilipschitz map from $N_{\rho_1}$ to
$N_{\rho_2}$ and  finish the proof as before. 

\medskip

We first construct a compact submanifold $J_1$ of $K_1$ which lifts
to a relative compact core $\hat J_1$ of $N_{R_1}^0$. Let $\{\alpha_i\}$ be the  collection
of simple closed curves on $R_1$ which are homotopic into $P_1$.
One may show that the $\{\alpha_i\}$ are disjoint and use the Annulus
Theorem to construct a disjoint collection $\{A_i\}$ such that $A_i$ joins
$\alpha_i$ to a simple closed curve in $P_1$. Let $J_1$ be a regular
neighborhood in $K_1$ of the 2-complex $R_1\cup \bigcup A_i$. (This is a special
case of the construction of a refined relative compression body neighborhood
of a relative boundary component of a pared manifold from  \cite{canary-mccullough:refined}.)

Let $J_2=\phi(J_1)$. Then $J_2$ lifts to a relative compact core $\hat
J_2$ of $N_{R_2}$ and $\phi|_{J_1}$
lifts to a homeomorphism
$$\hat \phi_{R_1}:\hat J_1 \to \hat J_2.$$

Let  $\hat F_1$ be a component of $\partial_0\hat J_1$ which
faces an  end  $\EE_1$ of $N_{R_1}^0$.   We claim that 
\begin{itemize}
\item[(*)] $\EE_1$ is either geometrically finite, or has a
  neighborhood which maps isometrically to a neighborhood of an end of
  $N_{\rho_1}^0$.
\end{itemize}

Let $X_1$ be the component of $K_1-J_1$ bounded by the image $F_1$ of
$\hat F_1$ in $J_1$. If $\EE_1$ is geometrially infinite, the 
Covering Theorem 
\cite{wpt:notes,canary:covering} implies that a neighborhood of $\EE_1$ projects to
$N_{\rho_1}$ by a finite-to-one covering (otherwise, $N_{\rho_1}$
would have a finite cover fibering over the circle). This then implies
that $\pi_1(F_1)$ has finite index in $\pi_1(X_1)$,
which implies (see \cite[Theorem 10.5]{hempel}),
that $(X_1,X_1\cap P)$ is an interval bundle pair.

If the interval bundle is trivial, then $F_1$ is parallel to a
component of $\boundary_0K_1$, which implies that $\EE_1$ has a
neighborhood mapping injectively to the end associated to this
component. This establishes (*) in this case. 

If the interval bundle is twisted, we consider the cover $N_{X_1}$
associated to $\pi_1(X_1)$, which is double-covered by $N_{F_1}$, the
cover associated $\pi_1(F_1)$. A neighborhood of the end $\EE_1$ lifts
isometrically to an end of $N_{F_1}$, which descends isometrically to
the (unique) end of $N^0_{X_1}$. Note that $F_1$ is isotopic in $J_1$ to a subsurface of
$R_1$. If this is a proper subsurface then $\pi_1(X_1)$ is infinite
index in $\pi_1(K_1)$, which implies that a neighborhood of the end of
$N_{X_1}$ maps with infinite degree. But then the same is true for
$\EE_1$, a contradiction. We conclude that $F_1$ is parallel to $R_1$
in $J_1$, and hence that a neighborhood of $\EE_1$ maps injectively to
the end associated to $R_1$ in $N_{\rho_1}^0$, again establishing (*).

\medskip

Claim (*) together with the fact that $\phi$ preserves end invariants
implies that, for each end $\EE$ of $N^0_{R_1}$, if $\EE$ is geometrically infinite
then  $\phi_{R_1}(\nu(\EE)) = \nu(\phi_{R_1}(\EE))$. The invariants of
geometrically finite ends, which are points in Teichm\"uller spaces,
may differ.  We may therefore apply Lemma \ref{bilipschitz corollary} to extend
$\phi_{R_1}$ to a bilipschitz homeomorphism from $N^0_{R_1}$ to $N^0_{R_2}$.

A restriction of this homeomorphism to 
a sufficiently small neighborhood of the end of $R_1$ in $N_{R_1}$
descends to the corresponding end of $N^0_{\rho_1}$.
We can assemble these homeomorphisms for all of the ends of
$N^0_{\rho_1}$, extend across the remaining compact subset of
$N^0_{\rho_1}$, and then radially across the cusps as before to
produce a bilipschitz map from $N_{\rho_1}$ to $N_{\rho_2}$ in the
homotopy class of $\phi$. 

As in the surface case, an application of Sullivan's Rigidity Theorem
finishes the proof.
\end{proof}

\subsubsection*{Remarks on the proof}
Lemma \ref{bilipschitz corollary} produces a bilipschitz
homeomorphism between quotients of Kleinian surface groups whose relative compact
cores are homeomorphic and whose  corresponding geometrically
infinite ends have matching laminations.  The same result holds in the
general incompressible-ends case, with the same proof, once the Ending
Lamination Theorem is established in that setting. 

A direct argument providing bilipschitz comparisons between ends with
corresponding end-invariants would be perhaps logically
preferable. Moreover it would allow us to simplify the final argument
of the Ending Lamination Theorem as well, removing the need for the
analysis of the covers $N_R$. 
Such an approach, albeit straightforward, requires some
additional technical tools: either combinatorial arguments using
hierarchies, or an application of the drilling technology of
\cite{brock-bromberg:density} to isolate ends geometrically.  We have
chosen a more indirect method using available tools (quasiconformal
deformation theory in Lemma \ref{bilipschitz corollary}, and the
Covering Theorem in the Ending Lamination Theorem)
for brevity.

\section{Corollaries}
\label{corollaries}

In this section we give proofs of the corollaries mentioned in
the introduction and  of the Length Bound Theorem from Section \ref{length estimates}. 

Our first corollary is the resolution of the Bers-Sullivan-Thurson
Density conjecture in the setting of {\em pared manifolds} with
incompressible boundary. A pared manifold is a pair $(M,P)$ where $M$
is a compact irreducible 3-manifold and $P$ is a submanifold of
$\boundary M$ consisting of incompressible annuli and tori, such that 
every noncyclic abelian subgroup of $\pi_1(M)$ is conjugate into the
fundamental group of a component of $P$, and every $\pi_1$-injective map of an annulus
$\phi:(S^1\times I,S^1\times\boundary I)\to (M,P)$ is homotopic, as a
map of pairs, to a map with image in $P$. We note that a relative
compact core $(K,P)$ of a hyperbolic 3-manifold is always a pared
manifold. 

We define $AH(M,P)$ to be the space of conjugacy classes of discrete,
faithful representations $\rho:\pi_1(M)\to \PSL 2(\C)$ such that every
conjugacy class represented by a curve in $P$ is mapped to parabolics.
We endow $AH(M,P)$ with the {\em algebraic topology}, which is just
the topology inherited from the representation variety of $\pi_1(M)$.

\begin{corollary}{}{{\em (Density Theorem for Incompressible Ends)}}
Let $(M,P)$ be a pared manifold with non-abelian fundamental group,
such that $\boundary M \setminus P$ is incompressible. Then 
$$\overline{\interior(AH(M,P))} = AH(M,P).$$ 
\end{corollary}

\begin{proof}{}
  Results of Sullivan \cite{sullivan:QCDII} and Marden \cite{marden:geometry} imply that
  $\interior(AH(M,P))$ consists exactly of those representations that
  are geometrically finite, and send {\em only} elements represented
  by curves in $P$ to parabolics. Ohshika \cite{ohshika:ending-lams}
  used convergence results of Thurson
  \cite{wpt:II,wpt:III}  to prove that every collection
  of end invariants that occurs for points in $AH(M,P)$ arises as the
  end invariants of a limit of elements of this type. 

  The Ending Lamination Theorem asserts that elements of $AH(M,P)$ are
  determined by their end invariants. 
  The Density Theorem follows.
\end{proof}

The proof of our rigidity theorem is somewhat more involved
as we must observe that a topological conjugacy can detect
the (marked) homeomorphism type of the relative compact core and
the ending invariants.

\begin{corollary}{}{{\em (Rigidity Theorem)}}
Let $G$ be a finitely generated, torsion-free, non-abelian group.
 If $\rho$ and $\rho'$ are two discrete faithful representations of
  $G$ into $\PSL 2(\C)$ that are conjugate by an
  orientation-preserving homeomorphism $\phi$ of $\hhat \C$ and $N_\rho^0$
  has incompressible ends, then $\rho$ and $\rho'$ are
  quasiconformally conjugate. Moreover, if $\phi$ is conformal on
  $\Omega(\rho_1)$, then $\phi$ is conformal.
\end{corollary}

\begin{proof}{}
We first reduce to the case where $\phi$ is conformal on
$\Omega(\rho_1)$. 
Since $\phi(\Omega(\rho_1))=\Omega(\rho_2)$, $\phi$ induces a 
homeomorphism between $\partial_\infty N_1$ and $\partial_\infty N_2$,
where $N_i=\Hyp^3/\rho_i(\pi_1(S))$.
Ahlfors' Finiteness Theorem \cite{ahlfors:finitegen}
assures us that $\partial_\infty N_i$ is
a Riemann surface of finite type, so we may deform
$\phi$ so that it is quasiconformal on $\Omega(\rho_1)$.
One may use the Measurable Riemann Mapping Theorem \cite{ahlfors-bers}
to construct a quasiconformal map $\psi:\hhat \C \to \hhat\C$ such
that $\rho_2'=\psi\circ\rho_2\circ \psi^{-1}$ is a Kleinian surface group
and $\psi\circ\phi$ is conformal on $\Omega(\rho_1)$.

For the remainder of the argument we will assume that $\phi$
is conformal on $\Omega(\rho_1)$.
Let $(K_1,Q_1)$ be a relative compact core for $N_1^0$ and
let $(K_2,Q_2)$ be a relative compact core for $N_2^0$.
Since $\phi$ identifies $\rho_1(\pi_1(S))$ with $\rho_2(\pi_1(S))$,
it induces a homotopy equivalence $\bar\phi$ from $K_1$ to $K_2$.
Recall that $\rho_i(g)$ is
parabolic if and only if it has exactly one fixed point in $\hhat \C$.
Therefore, $\rho_1(g)$ is parabolic if and only if $\rho_2(g)$ is
parabolic. Thus, $\bar\phi(Q_1)$ is homotopic to $Q_2$.

Let $G_i$ be the union of the components of $\partial K_i-Q_i$ which
are associated to geometrically finite ends of $N_i$. One may identify
$G_i$ with $\partial_\infty N_i$ and assume that $\bar\phi$
is a conformal homeomorphism from $G_1$ to $G_2$.

Let $R$ be a component of $\boundary K_1 - Q_1$ associated to a
geometrically infinite end, with ending lamination $\lambda$.
The restriction $\sigma_1=\rho_1|_{\pi_1(R)}$ is a Kleinian surface group with
ending lamination $\lambda$. We claim that the same holds for
$\sigma_2 = \rho_2|_{\pi_1(R)}$.

For an element $\gamma\in \pi_1(R)$ and $i=1$ or $2$, let $d_i(\gamma)$ be
the maximal distance between the fixed points of $\sigma_i(g)$ in the
ball model of $\Hyp^3$, where $g$ runs over the conjugacy class of
$\gamma$. Given a sequence $(\gamma_k)$ of elements of $G$ with
non-parabolic images, we note
that the geodesic representatives of $\gamma_k$ in $N_{\sigma_i}$ leave
every compact set if and only if $d_i(\gamma_k)\to 0$. 

Because $\lambda$ is an ending lamination of
 $\sigma_1$, there exists a sequence $\gamma_k$ of elements
represented by simple closed curves converging to $\lambda$ in
$\PML(R)$, whose geodesic representative leave every compact set in
$N^0_{\sigma_1}$. Hence $d_1(\gamma_k) \to 0$. Since $\phi$ is a
homeomorphism conjugating $\sigma_1$ to $\sigma_2$, we conclude that
$d_2(\gamma_k)\to 0$ as well. It follows that $\lambda$ is an ending
lamination of $\sigma_2$. Therefore there is a geometrically infinite
end $\EE$ of $N^0_{\sigma_2}$, with base surface $R$. The Covering
Theorem \cite{wpt:notes,canary:covering} implies that the projection of $\EE$ to $N_{\rho_2}$ is
finite-to-one, and a neighborhood of $\EE$ maps to a neighborhood of
an end of $N^0_{\rho_2}$.

From this it follows that $\bar\phi$ can be chosen to take $R$
properly to a component of $\boundary K_2 \setminus Q_2$. 
We do this for all the geometrically
infinite ends. Now $\bar\phi$ maps $\boundary K_1 \setminus Q_1$
properly to $\boundary K_2 \setminus Q_2$, and in particular maps
$\boundary Q_1$ to $\boundary Q_2$. Since each $(K_i,Q_i)$ is a pared
manifold, the map on $Q_1$ can be deformed rel boundary into $Q_2$.
By Johannson's version of Waldhausen's
Theorem (Proposition 3.4 in \cite{johannson:3mfds}, see also the
discussion in \S2.5 of \cite{canary-mccullough}), $\bar\phi$ may be
deformed to a homeomorphism of pared manifolds. 

Moreover our argument has shown that $\bar\phi$ takes the end
invariants of $N_{\rho_1}$ to those of $N_{\rho_2}$. We may therefore
apply the Ending Lamination Theorem to conclude that there is an
isometry $F:N_{\rho_1}\to N_{\rho_2}$ in the homotopy class of
$\bar\phi$. 

Let $\phi':\hhat\C\to\hhat\C$ be the map
which is the extension of the lift of $F$ to
$\Hyp^3$. 
Then $\phi'$ is either conformal or anti-conformal and
conjugates $\rho_1$ to $\rho_2$.
Notice that since fixed points of elements of $\rho_1(G)$
are dense in $\Lambda(\rho_1)$, $\phi$ and $\phi'$ agree on $\Lambda(\rho_1)$.
Since our initial map $\phi$ was conformal on $\Omega(\rho_1)$, $\phi$ and $\phi'$
must agree on $\Omega(\rho_1)$ and hence on $\hhat\C$. Therefore,
since $\phi$ is orientation-preserving and $\phi'$ is either conformal
or anti-conformal, $\phi=\phi'$ is conformal.
\end{proof}

We next turn our attention to: 

\begin{corollary}{}{{\em (Volume Growth Theorem)}}
If $N$ is the quotient of a Kleinian surface group $\rho\in\DD(S)$,
then for any $x$ in the $\ep_1$-thick part of the convex core $C_N$
and $r\ge 1$ we have
$${\rm volume}\left( B^{thick}_r(x)\right) \le c_1r^{d(S)},$$
where $c_1$ depends only on the topological type of $S$. 

In general, if $N$ is a complete hyperbolic 3-manifold with 
relative compact core $K$ so that $N^0$ has incompressible ends, 
we have
$${\rm volume}\left( B^{thick}_r(x)\right) \le c_1r^{d(\bdry_0 K)} + c_2,$$
where $c_1$ depends only on
the topological type of $\bdry_0 K$, and
$c_2$ depends on the hyperbolic structure of $N$.
\end{corollary}
Recall that, for connected $S$, $d(S) = -\chi(S)$ when $\text{genus}(S)>0$
and  $d(S) = -\chi(S)-1$ when $\text{genus}(S)=0$, and for
disconnected $S$ $d(S)$ is the maximum over its components. Recall
also that $B^{thick}_r(x)$ denotes the $r$-neighborhood of $x$ in the
path metric of the $\ep_1$ thick part of $C_N$.

\begin{proof}
We first consider the surface group case. 

We can replace $C_N$ by $\hhat C_N$, which contains it. 
The $\ep_1$-thick part of $\hhat C_N$ is almost the same as
$\hhat C_N\setminus\MT[k]$; the latter may include some
$\ep_1$-Margulis tubes with $\omega$ coefficients bounded (in terms of
$k$ and the bounds of the model theorem). Since all 
such tubes have uniformly bounded diameters and volumes, and uniformly
large disjoint regular neighborhoods, it suffices
to prove the theorem for $\hhat C_N\setminus\MT[k]$. Now since
the Bilipschitz Model Theorem gives a uniformly bilipschitz homeomorphism of
$\modl[k]$ to $\hhat C_N\setminus \MT[k]$, it suffices to prove the
theorem for $\modl[k]$. Finally, this is equivalent to proving the theorem
for $\modl[0]$, again because the difference consists of tubes with
bounded diameters and volumes, and uniform separation.  This is what we will do. 

\medskip

Fix a cut system $C$, and 
recall from \S\ref{address region structure} the definition of the
product regions $\BB(h)\subset \modl$ where $h\in H$ and $C|_h$ is
nonempty.
Each $\BB(h)$ is isotopic to $D(h)\times I$ for an interval $I$, and is defined
as the region between the first and last slices $a_h$ and $z_h$ in $C|_h$
(note that $a_h$ or $z_h$ could fail to exist if $h$ is infinite in the backward
or forward direction, in which case $\BB(h)$ is defined accordingly).

Define also $\BB_0(h) = \BB(h)\intersect \modl[0]$. 
For $x\in \modl[0]$, 
let $\NN_r(x)$ denote the $r$-neighborhood of $x$ with respect to the
path metric in $\modl[0]$. 

We shall prove the following statement by induction on $d(D(h))$:

\begin{itemize}
\item[(*)] For any $h\in H$ with $C|_h\ne\emptyset$ and $d(D(h))\ge 1$, 
given $x\in\BB_0(h)$  the volume of  
$\BB_0(h)\intersect \NN_r(x)$ is bounded by 
$cr^{d(D(h))}$, where $c$ depends only on $d(D(h))$. 
\end{itemize}

Note that the boundary of $\BB(h)$ consists of 
the bottom and top slice surfaces $\hhat F_{a_h}$ and $\hhat F_{z_h}$
(these could be empty if $h$ is infinite)
together with
tube-boundary annuli associated to $\boundary D(h)$. Thus the frontier
of $\BB_0(h)$ in $\modl[0]$ is just the surfaces $ F_{a_h}$ and
$F_{a_z}$, each of which have at most $-\chi(S)$ components, 
with uniformly bounded diameter.
Since $\BB_0(h)\intersect \NN_r(x)$ is contained in the union of
$r$-neighborhoods, in the path metric of $\BB_0(h)$, of $x$ and
the frontier of $\BB_0(h)$, this implies that
proving (*) for
$\BB_0(h)\intersect \NN_r(x)$ is equivalent to proving it for
the $r$-neighborhood of a point $y\in\BB_0(h)$ in the path metric of 
$\BB_0(h)$, which we can denote $\NN_{r,h}(y)$.
We proceed to do this. 

When $d(D(h))=1$, $D(h)$ can only be a one-holed torus or a 4-holed
sphere (this is the place where the variation in the definition of $d$
for genus 0 comes into play). In this case, 
each block in $\BB(h)$ is associated
with exactly one edge of the geodesic $h$ and is isotopic to
a sub-product region $D(h)\times J$ -- in particular it
separates $\BB(h)$. It follows immediately that 
$$
vol(\NN_{r,h}(y)) \le 2 v_0  r/r_0
$$
where $v_0$ is an upper bound on the volume of a block and $r_0$ a
lower bound on the distance between the top and bottom boundaries of
a block. This establishes (*) for $d(D(h))=1$. 

Now for $d(D(h))=u >1$, consider the set $E$ of slices $e\in C$ 
such that $F_e \subset \BB_0(h)$, and $d(D(e)) = u$.
It follows from the definition of the function $d$ that for each $e\in E$, 
$D(e)$ is equal to $D(h)$ minus a (possibly empty) disjoint union of
annuli. 

We extend $\hhat F_e$ across all the tubes associated to annuli of
$D(h)\setminus D(e)$, obtaining a surface 
$\ddot F_e$ 
isotopic to $D(h)\times\{t\}$ in the product
structure of $\BB(h)$.  Since these surfaces are disjoint, they are topologically
ordered in $\BB(h)$, so we can number them 
$\cdots \ddot F_{e_i} \topprec \ddot F_{e_{i+1}} \cdots$
and let $P_i$ denote the product region between 
$\ddot F_{e_i}$ and $\ddot F_{e_{i+1}}$.

Each $\ddot F_{e_j}$ separates $P_i$ from $P_l$ for $i<j\le l$. Since 
$\ddot F_e \setminus F_e$ is a union of annuli in the tubes of
$\modl$, it follows that $F_{e_j}$ also separates
$P_i[0] = P_i\intersect \modl[0]$ from
$P_l[0] = P_l\intersect \modl[0]$.

Note that $r_0$ is a lower bound on the distance in $P_i[0]$
between $F_{e_i}$ and $F_{e_{i+1}}$, since the slices cannot have any
pieces in common (Lemma \ref{cuts dont meet}) and hence must be
separated by a layer 
at least one block thick. It follows that $\NN_{r,h}(y)$ can meet 
at most $2r/r_0$ different regions $P_i[0]$.

It remains to estimate the volume of $P_i[0]\intersect \NN_{r,h}(y)$.
Inside $P_i$ there may be block regions $\BB(m)$, where $m\in H$
such that, necessarily, $C|_m\ne\emptyset$ and $d(D(m))<u$. 
Any slice surface $\hhat F_b$ for $b\in C$ which meets $int(P_i)$ must
in fact be contained in one of these $\BB(m)$'s, and hence
the complement in $P_i$ of all such $\BB(m)$ is contained in a
complementary region $W$ of $C$. By Lemma \ref{Bound regions} and
Proposition \ref{modl[k] regions}, $W$ has 
a uniformly bounded number of blocks, which 
bounds its volume by  
a constant $v_1$. It also implies that the 
number of $\BB(m)$ contained in $P_i$ and adjacent to $W$ (which we
will call {\em outermost} $\BB(m)$) is bounded by some $n_1$. 
For each outermost $\BB_0(m)$ we have, by induction, 
$$
vol(\NN_{r,m}(y)) \le c r^{d(D(m))} \le cr^{u-1}
$$
for any $y\in \BB_0(m)$. Since the frontier of $\BB_0(m)$ in $P_i$
consists of a top and a bottom surface, each of bounded diameter, 
we may conclude that a bound of the form $c'r^{u-1}$ holds for
$\BB_0(m)\intersect \NN_{r,h}(y)$, whether $y\in\BB_0(m)$ or not.
Summing these over all outermost $\BB(m)$ in 
$P_i$ and including the rest of $P_i$ we have a bound
$$
vol(P_i\intersect\NN_{r,h}(y)) \le v_1 + n_1c'r^{u-1} \le c''r^{u-1}
$$
for some $c''$ depending only on $d(D(h))$. 
Now summing this over all the (at most $2r/r_0$) $P_i$ that meet the
$r$-neighborhood of $y$ gives
$$
vol(\NN_{r,h}(y)) \le c''' r^u
$$
which establishes (*). 
(Note that the union of closures of $P_i$ fills up $\BB(h)$ except possibly
if $h=g_H$ and the top or bottom boundary of $\modl$ is nonempty; at any
rate the exterior of $\union \bar P_i$ is contained in a single
address region, and has bounded volume by Lemma \ref{Bound regions}.)

Now applying (*) with $h=g_H$ yields the desired growth estimate for
$\modl[0]$.

\medskip

We now consider the general case. 

Given a relative compact core $K\subset C_N$ for $N^0$, for each component $R$ of
$\boundary_0 K$ the component $U_R$ of $C_N
\setminus K$ bounded by $R$
lifts isometrically to $\hat U_R$ in the surface-group cover $N_R$. We select
$K$ large enough so that $\hat U_R$ lies in the convex hull of $N_R$, 
for each $R$. Throughout this proof let $X^\epone$ denote the
$\epone$-thick part of $X$. 

Let $d^{thick}$ be the path-metric in $C^\epone_N$
and let $\hat d^{thick}_R$ be the metric on $U^\epone_R$ inherited
from  the restriction to $\hat U^\epone_R$ of the
path-metric in $N^\epone_R$. We have
$d^{thick}\le \hat d^{thick}_R$. Let $\delta_R$ denote the 
$\hat d^{thick}_R$-diameter of $R$. Let $K'$ denote the
$(\max_R\delta_R)$-neighborhood of $K$, in the $d^{thick}$ metric.

We first note that, for $x,y\in U^\epone_R \setminus K'$, 
$\hat d^{thick}_R(x,y) \le 2 d^{thick}(x,y)$. This is because any path
connecting the lifts of $x$ and $y$ to the cover $\hat U^\epone_R$ can only
exit $\hat U^\epone_R$ through the lift of $R$, where the length savings is
at most $\delta_R$; but in that case the length of the path is at
least $2\delta_R$.

The volume bound for $C^\epone_N$ now follows from the surface-group case
applied to each $U^\epone_R$, where a multiplicative error of at most 2 is
introduced by the ratio of $d^{thick}$ to $\hat d^{thick}_R$, and an
additive error of at most $vol(K')$ is introduced by the volume of the
core. 
\end{proof}

Finally, we recall and prove the Length Bound Theorem. 

\state{Length Bound Theorem.}{%
There exist $\bar \ep>0$ and $c>0$ depending only on $S$, such that
the following holds:

Let $\rho:\pi_1(S)\to \PSL 2(\C)$ be a Kleinian surface group and
$v$ a vertex of $\CC(S)$, and let $H_{\nu_\rho}$ be an associated
hierarchy.
\begin{enumerate}
\item If $\ell_\rho(v) < \bar\ep$ then $v$  appears in $H_{\nu_\rho}$.

\item  
If $v$ appears in $H_{\nu_\rho}$ then
$$
d_{\Hyp^2}\left(\omega(v),\frac{2\pi i}{\lambda_\rho(v)}\right)
\le c
$$
\end{enumerate}
}

\begin{proof}
The Short Curve Theorem (\S\ref{length estimates}) already contains 
part (1) of the Length Bound Theorem. It remains to prove part (2).
For simplicity
we suppress $v$ in the proof, writing $\omega$, $\lambda$, etc. 

Suppose first that $|\omega| \in (k,\infty)$ where $k$ is the
constant in the 
Bilipschitz Model Theorem. Then the tube $U=U(v)$ is in $\UU[k]$ and
the model map
$F$ takes $U$ to the correponding Margulis tube $\MT(v)$ by a
$K$-bilipschitz map.
Letting $\omega_\MT$ denote the Teichm\"uller parameter of $\MT(v)$
with respect to the marking induced by the model map, 
it follows that 
\begin{equation}\label{omega omegaT}
d_{\Hyp^2}(\omega,\omega_\MT) \le \log K
\end{equation}
where $\Hyp^2$, the upper half plane, is identified with the
Teichm\"uller space of the torus, and $d_{\Hyp^2}$ is the
Teichm\"uller metric and the Poincar\'e metric. 

Equation (\ref{lambda tube}), which comes from                   
\S3.2 of \cite{minsky:ELCI}, says that
\begin{equation}\label{lambda lambdaT}
\lambda' \equiv h_r(2\pi i/\omega_\MT)
\end{equation}
is equal to $\lambda$ modulo $2\pi i$, 
where $r$ is the radius of $\MT(v)$ and 
where $h_r(z) = \Re z \tanh r + i\Im z$.

The radius of $U(v)$, by (\ref{r tube}) in Section \ref{model definitions}, is given by
$$
r_U = \sinh^{-1}(\epone|\omega|/2\pi).
$$
The $K$-bilipschitz map between $U(v)$ and $\MT(v)$ guarantees that
$ r \ge r_U/K $. Thus since $|\omega| > k$ we obtain a uniform positive lower
bound on $r$. Now, noting that 
that $h_r$ preserves the right half plane $\Hyp'=\{z:\Re z>0\}$
and, letting $d_{\Hyp'}$   denote the Poincar\'e metric on
$\Hyp'$, it is easy to check that
$d_{\Hyp'}(z,h_r(z)) $ is uniformly bounded above, and this gives an
upper bound of the form
\begin{equation}\label{omegaT lambdaprime}
d_{\Hyp'}(\frac{2\pi i}{\omega_\MT}, \lambda') \le C_1.
\end{equation}
The lower bound on $r$ also implies that, as $|\omega|\to\infty$,
$|\lambda'|\to 0$, and in particular there is some $k_2\ge k$ such that, when
$|\omega| > k_2$, $|\Im \lambda'|<\pi $.
Recalling from Section \ref{length estimates} that the
  imaginary part of $\lambda$ was normalized to lie in $(-\pi,\pi]$,
    we conclude that   $\lambda = \lambda'$ whenever $|\omega| >
    k_2$.

Now since the map $z \to 2\pi
i/z$ is an isometry in the Poincar\'e metric from $\Hyp'$ to
$\Hyp^2$, we obtain from (\ref{omegaT lambdaprime}) a uniform bound 
on $d_{\Hyp^2}(\omega_\MT,2\pi i/\lambda')$. 
Together with 
(\ref{omega omegaT}), we have the desired bound on 
on $d_{\Hyp^2}(\omega,2\pi i/\lambda)$, whenever $|\omega| > k_2$.

If $|\omega| \le k_2$ then we have uniform lower and upper bounds on
$\Re \lambda$ by the Short Curve Theorem (\S\ref{bilip model
  intro}), and on $\Im \lambda$ by definition. Moreover we know that
$\Im \omega \ge 1$ always. This constrains both $\omega$ and $\lambda$
to compact sets, the estimate is immediate. 

\end{proof}

\subsection*{Erratum}

We use this opportunity to point out a small error in 
Minsky \cite{minsky:ELCI}. A central result in \cite{minsky:ELCI} is Theorem
7.1, which gives projection bounds for the short-curve set in a
Kleinian surface group. This theorem contains two statements, of which
the first contains the error and the second is the one which is
actually applied in the paper. The corrected version of the theorem is
as follows:

\state{Theorem 7.1 of \cite{minsky:ELCI}}{%
Fix a surface $S$. 
There exists $L_1\ge L_0$ such that for every 
$L\ge L_1$ there exist $B,D_2>0$ such that, 
given $\rho\in\DD(S)$, a hierarchy $H=H_{\nu(\rho)}$,
and an essential  subsurface $Y$ in $S$ with $\xi(Y)\ne 2,3$, the set
$$
\pi_Y(\CC(\rho,L))
$$
is $B$-quasiconvex in $\AAA(Y)$.  Furthermore, when $\xi(Y)\ne 3$,
\begin{equation*}
d_Y(v,\Pi_{\rho,L}(v)) \le D_2  \eqno{(*)}
\end{equation*}
for every vertex $v$ appearing in $H$ such that the left-hand
side is defined.
}

Here, $\CC(\rho,L)$
denotes the subset of $\CC(S)$ consisting of curves $\alpha$ with
$\ell_\rho(\alpha) \le L$. $\Pi_{\rho,L}$ is a partially-defined map
from $\CC(S)$ to $\CC(\rho,L)$. The second part, the inequality $(*)$,
is unchanged from the original. In the first part the condition that
$\xi(Y)\ne 2$ -- namely that $Y$ not be an annulus -- has been
added. Indeed, the quasiconvexity property fails to hold when $Y$ is
an annulus. 

The second statement implies the
quasiconvexity statement for nonannular $Y$, but not in the case of annuli.
The quasiconvexity property is not used anywhere in
\cite{minsky:ELCI}.

\bibliographystyle{hamsplain}
\bibliography{math}

\providecommand{\bysame}{\leavevmode\hbox to3em{\hrulefill}\thinspace}
\begin{thebibliography}{10}

\bibitem{agol:volume}
I.~Agol, \emph{Volume change under drilling}, Geometry \& Topology \textbf{6}
  (2002), 905--916.

\bibitem{agol:tame}
\bysame, \emph{Tameness of hyperbolic 3-manifolds}, Preprint, 2004,
  \mbox{arXiv:math.GT/0405568}.

\bibitem{ahlfors-bers}
L.~Ahlfors and L.~Bers, \emph{Riemann's mapping theorem for variable metrics.},
  Ann. of Math. \textbf{72} (1960), 385--404.

\bibitem{ahlfors:finitegen}
Lars~V. Ahlfors, \emph{Finitely generated {K}leinian groups}, Amer. J. Math.
  \textbf{86} (1964), 413--429.

\bibitem{anderson-canary:pages}
J.~Anderson and R.~Canary, \emph{Algebraic limits of {Kleinian} groups which
  rearrange the pages of a book}, Invent. Math. \textbf{126} (1996), 205--214.

\bibitem{anderson-canary:cores}
\bysame, \emph{Cores of hyperbolic 3-manifolds and limits of {Kleinian}
  groups}, Amer. J. Math. \textbf{118} (1996), 745--779.

\bibitem{anderson-canary-culler-shalen}
J.~Anderson, R.~Canary, M.~Culler, and P.~Shalen, \emph{Free {K}leinian groups
  and volumes of hyperbolic $3$-manifolds}, J. Differential Geom. \textbf{43}
  (1996), no.~4, 738--782.

\bibitem{anderson-canary-mccullough:bumping}
J.~Anderson, R.~Canary, and D.~McCullough, \emph{The topology of deformation
  spaces of {K}leinian groups}, Ann. of Math. (2) \textbf{152} (2000), no.~3,
  693--741.

\bibitem{anderson:intersection}
J.~W. Anderson, \emph{Intersections of topologically tame subgroups of
  {Kleinian} groups}, Journal d'Analyse \textbf{65} (1995), 77--94.

\bibitem{behrstock-kleiner-minsky-mosher:qirigid}
J.~Behrstock, B.~Kleiner, Y.~Minsky, and L.~Mosher, \emph{Geometry and rigidity
  of mapping class groups}, E-print arXiv:0801.2006, 2008.

\bibitem{benedetti-petronio}
R.~Benedetti and C.~Petronio, \emph{Lectures on hyperbolic geometry},
  Springer-Verlag Universitext, 1992.

\bibitem{bers:spaces}
L.~Bers, \emph{Spaces of {K}leinian groups}, Maryland conference in Several
  Complex Variables {I}, Springer-Verlag Lecture Notes in Math, No. 155, 1970,
  pp.~9--34.

\bibitem{bonahon}
F.~Bonahon, \emph{Bouts des vari\'et\'es hyperboliques de dimension 3}, Ann. of
  Math. \textbf{124} (1986), 71--158.

\bibitem{bowditch:elcgen}
B.~Bowditch, \emph{End invariants of hyperbolic 3-manifolds}, Preprint,
  Southampton, 2005.

\bibitem{bowditch:complex}
\bysame, \emph{Intersection numbers and the hyperbolicity of the curve
  complex}, J. Reine Angew. Math. \textbf{598} (2006), 105--129.

\bibitem{bowditch:bands}
\bysame, \emph{Systems of bands in hyperbolic 3-manifolds}, Pacific J. Math.
  \textbf{232} (2007), 1--42.

\bibitem{brock-bromberg:density}
J.~Brock and K.~Bromberg, \emph{{On the density of geometrically finite
  Kleinian groups.}}, Acta. Math. \textbf{192} (2004), 33--93.

\bibitem{brock-canary-minsky:ELCIII}
J.~Brock, R.~Canary, and Y.~Minsky, \emph{The classification of
  finitely-generated {Kleinian} groups}, in preparation.

\bibitem{bromberg:density}
K.~Bromberg, \emph{Projective structures with degenerate holonomy and the
  {Bers} density conjecture}, Ann. of Math. \textbf{166} (2007), 77--93.

\bibitem{calegari-gabai:tame}
D.~Calegari and D.~Gabai, \emph{Shrinkwrapping and the taming of hyperbolic
  3-manifolds}, J. Amer. Math. Soc. \textbf{19} (2006), no.~2, 385--446
  (electronic).

\bibitem{canary-mccullough:refined}
R.~Canary and D.~McCullough, \emph{The refined relative compression body
  neighborhood}, Preprint, http://www.math.ou.edu/\textasciitilde
  dmccullough/research/preprints.html, 2000.

\bibitem{canary-mccullough}
\bysame, \emph{Homotopy equivalences of 3-manifolds and deformation theory of
  {Kleinian} groups}, Mem. Amer. Math. Soc., vol. 172, AMS, 2004.

\bibitem{canary-minsky:tamelimits}
R.~Canary and Y.~Minsky, \emph{On limits of tame hyperbolic 3-manifolds}, J.
  Differential Geom. \textbf{43} (1996), 1--41.

\bibitem{canary:ends}
R.~D. Canary, \emph{Ends of hyperbolic 3-manifolds}, J. Amer. Math. Soc.
  \textbf{6} (1993), 1--35.

\bibitem{canary:covering}
\bysame, \emph{A covering theorem for hyperbolic 3-manifolds and its
  applications}, Topology \textbf{35} (1996), 751--778.

\bibitem{ceg}
R.~D. Canary, D.~B.~A. Epstein, and P.~Green, \emph{Notes on notes of
  {T}hurston}, Analytical and Geometric Aspects of Hyperbolic Space, Cambridge
  University Press, 1987, London Math. Soc. Lecture Notes Series no. 111,
  pp.~3--92.

\bibitem{douady-earle}
A.~Douady and C.J. Earle, \emph{Conformally natural extension of homeomorphisms
  of the circle}, Acta Math. \textbf{157} (1986), 23--48.

\bibitem{epstein-marden}
D.~B.~A. Epstein and A.~Marden, \emph{Convex hulls in hyperbolic space, a
  theorem of {Sullivan}, and measured pleated surfaces}, Analytical and
  Geometric Aspects of Hyperbolic Space, Cambridge University Press, 1987,
  London Math. Soc. Lecture Notes Series no. 111, pp.~113--254.

\bibitem{evans:tamenesspersists}
R.~Evans, \emph{Tameness persists in weakly type-preserving limits.}, Amer. J.
  Math. \textbf{126} (2004), no.~4, 713--737.

\bibitem{freedman-hass-scott:area}
M.~Freedman, J.~Hass, and P.~Scott, \emph{Least area incompressible surfaces in
  $3$-manifolds}, Invent. Math. \textbf{71} (1983), no.~3, 609--642.

\bibitem{hamenstadt:boundary}
U.~Hamenst{\"a}dt, \emph{Train tracks and the {Gromov} boundary of the complex
  of curves}, Spaces of Kleinian Groups (Y.~Minsky, M.~Sakuma, and C.~Series,
  eds.), LMS Lecture Note Series, no. 329, Cambridge U. Press, 2006,
  pp.~187--207.

\bibitem{harvey:boundary}
W.~J. Harvey, \emph{Boundary structure of the modular group}, Riemann Surfaces
  and Related Topics: Proceedings of the 1978 Stony Brook Conference (I.~Kra
  and B.~Maskit, eds.), Ann. of Math. Stud. 97, Princeton, 1981.

\bibitem{harvey:modular}
\bysame, \emph{Modular groups and representation spaces}, Geometry of group
  representations (Boulder, CO, 1987), Amer. Math. Soc., 1988, pp.~205--214.

\bibitem{hempel}
J.~Hempel, \emph{3-manifolds}, Annals of Math. Studies no. 86, Princeton
  University Press, 1976.

\bibitem{jaco-rubinstein}
W.~Jaco and H.~Rubinstein, \emph{Pl minimal surfaces in 3-manifolds}, J.
  Differential Geom. \textbf{27} (1988), 493--524.

\bibitem{jaco}
W.~H. Jaco, \emph{Lectures on three-manifold topology}, CBMS Regional
  Conference Series no. 43, Amer. Math. Soc., 1980.

\bibitem{jaco-shalen}
W.~H. Jaco and P.~B. Shalen, \emph{Seifert fibered spaces in 3-manifolds},
  vol.~21, Memoirs of the Amer. Math. Soc., no. 220, A.M.S., 1979.

\bibitem{johannson:3mfds}
K.~Johannson, \emph{Topology and combinatorics of 3-manifolds}, Lecture Notes
  in Mathematics, vol. 1599, Springer-Verlag, Berlin, 1995.

\bibitem{jorgensen-marden:convergence}
T.~J{\o }rgensen and A.~Marden, \emph{Alebraic and geometric convergence of
  {K}leinian groups}, Math. Scand. \textbf{66} (1990), 47--72.

\bibitem{klarreich:boundary}
E.~Klarreich, \emph{The boundary at infinity of the curve complex and the
  relative {Teichm\"uller} space}, Preprint,
  http://www.nasw.org/users/klarreich/research.htm.

\bibitem{kra:spaces}
I.~Kra, \emph{On spaces of {K}leinian groups}, Comment. Math. Helv. \textbf{47}
  (1972), 53--69.

\bibitem{luukkainen:isotopy}
J.~Luukkainen, \emph{Bi-{L}ipschitz concordance implies bi-{L}ipschitz
  isotopy}, Monatsh. Math. \textbf{111} (1991), no.~1, 35--46.

\bibitem{marden:geometry}
A.~Marden, \emph{The geometry of finitely generated {K}leinian groups}, Ann. of
  Math. \textbf{99} (1974), 383--462.

\bibitem{maskit:self}
B.~Maskit, \emph{Self-maps of {K}leinian groups}, Amer. J. of Math. \textbf{93}
  (1971), 840--856.

\bibitem{masur-minsky:complex1}
H.~A. Masur and Y.~Minsky, \emph{Geometry of the complex of curves {I}:
  Hyperbolicity}, Invent. Math. \textbf{138} (1999), 103--149.

\bibitem{masur-minsky:complex2}
\bysame, \emph{Geometry of the complex of curves {II}: Hierarchical structure},
  Geom. Funct. Anal. \textbf{10} (2000), 902--974.

\bibitem{mcmullen:grafting}
C.~McMullen, \emph{Complex earthquakes and {T}eichm\"uller theory}, J. Amer.
  Math. Soc. \textbf{11} (1998), no.~2, 283--320.

\bibitem{minsky:knoxville}
Y.~Minsky, \emph{On {T}hurston's ending lamination conjecture}, Proceedings of
  Low-Dimensional Topology, May 18-23, 1992, International Press, 1994.

\bibitem{minsky:torus}
\bysame, \emph{The classification of punctured-torus groups}, Annals of Math.
  \textbf{149} (1999), 559--626.

\bibitem{minsky:kgcc}
\bysame, \emph{Kleinian groups and the complex of curves}, Geometry and
  Topology \textbf{4} (2000), 117--148.

\bibitem{minsky:boundgeom}
\bysame, \emph{Bounded geometry in {Kleinian} groups}, Invent. Math.
  \textbf{146} (2001), 143--192, \mbox{arXiv:math.GT/0105078}.

\bibitem{minsky:warwick}
\bysame, \emph{Combinatorial and geometrical aspects of hyperbolic
  3-manifolds}, {Kleinian Groups and Hyperbolic 3-Manifolds} (V.~Markovic
  Y.~Komori and C.~Series, eds.), London Math. Soc. Lec. Notes, vol. 299,
  Cambridge Univ. Press, 2003, pp.~3--40.

\bibitem{minsky:cdm}
\bysame, \emph{End invariants and the classification of hyperbolic
  3-manifolds}, Current developments in mathematics, 2002, Int. Press,
  Somerville, MA, 2003, (On web at
  http://www.math.yale.edu/users/yair/research/cdm.ps), pp.~181--217.

\bibitem{minsky:ELCI}
\bysame, \emph{The classification of {Kleinian} surface groups {I}: models and
  bounds}, Annals of Math. \textbf{171} (2010), 1--107.

\bibitem{mostow:hyperbolic}
G.~D. Mostow, \emph{Quasiconformal mappings in $n$-space and the rigidity of
  hyperbolic space forms}, Publ. I.H.E.S. \textbf{34} (1968), 53--104.

\bibitem{namazi:thesis}
H.~Namazi, \emph{Heegaard splittings and hyperbolic geometry}, ProQuest LLC,
  Ann Arbor, MI, 2005, Thesis (Ph.D.)--State University of New York at Stony
  Brook.

\bibitem{namazi-souto:uniform-inj}
H.~Namazi and J.~Souto, \emph{{Revisiting Thurston's uniform injectivity
  theorem}}, in preparation.

\bibitem{namazi-souto:density}
\bysame, \emph{{Non-realizability and ending laminations - Proof of the Density
  Conjecture}}, Preprint, 2010.

\bibitem{ohshika:density}
K.~Ohshika, \emph{Realising end invariants by limits of minimally parabolic,
  geometrically finite groups}, \mbox{arXiv:math/0504546}.

\bibitem{ohshika:ending-lams}
\bysame, \emph{Ending laminations and boundaries for deformation spaces of
  {K}leinian groups}, J. London Math. Soc. \textbf{42} (1990), 111--121.

\bibitem{ohshika:constructing-limits}
\bysame, \emph{Constructing geometrically infinite groups on boundaries of
  deformation spaces}, J. Math. Soc. Japan \textbf{61} (2009), no.~4,
  1261--1291.

\bibitem{miyachi-ohshika:tamebg}
K.~Ohshika and H.~Miyachi, \emph{On topologically tame {K}leinian groups with
  bounded geometry}, Spaces of {K}leinian groups, London Math. Soc. Lecture
  Note Ser., vol. 329, Cambridge Univ. Press, Cambridge, 2006, pp.~29--48.

\bibitem{otal:knotting}
J.-P. Otal, \emph{Sur le nouage des geodesiques dans les varietes
  hyperboliques.}, C. R. Acad. Sci. Paris S\`er. I Math. \textbf{320} (1995),
  no.~7, 847--852.

\bibitem{otal:knotting2}
\bysame, \emph{Les g\'eod\'esiques ferm\'ees d'une vari\'et\'e hyperbolique en
  tant que n\oe uds}, Kleinian groups and hyperbolic 3-manifolds (Warwick,
  2001), London Math. Soc. Lecture Note Ser., vol. 299, Cambridge Univ. Press,
  Cambridge, 2003, pp.~95--104.

\bibitem{prasad}
G.~Prasad, \emph{Strong rigidity of {$\bf Q$}-rank $1$ lattices}, Invent. Math.
  \textbf{21} (1973), 255--286.

\bibitem{rees:elc}
M.~Rees, \emph{The ending laminations theorem direct from {Teichm\"uller}
  geodesics}, E-print Arxiv:math.GT/0404007, 2004.

\bibitem{rourke-sanderson}
C.~P. Rourke and B.~J. Sanderson, \emph{Introduction to piecewise-linear
  topology}, Springer-Verlag, New York, 1972, Ergebnisse der Mathematik und
  ihrer Grenzgebiete, Band 69.

\bibitem{soma:function}
T.~Soma, \emph{Function groups in {Kleinian} groups}, Math. Annalen
  \textbf{292} (1992), 181--190.

\bibitem{soma:elc}
\bysame, \emph{Geometric approach to {Ending Lamination Conjecture}},
  arXiv:0801.4236, 2008.

\bibitem{souto:knotting}
J.~Souto, \emph{Short curves in hyperbolic manifolds are not knotted},
  preprint, 2004.

\bibitem{sullivan:bilipschitz}
D.~Sullivan, \emph{Hyperbolic geometry and homeomorphisms}, Geometric topology
  (Proc. Georgia Topology Conf., Athens, Ga., 1977), Academic Press, New York,
  1979, pp.~543--555.

\bibitem{sullivan:rigidity}
D.~Sullivan, \emph{On the ergodic theory at infinity of an arbitrary discrete
  group of hyperbolic motions}, Riemann Surfaces and Related Topics:
  Proceedings of the 1978 Stony Brook Conference, Ann. of Math. Stud. 97,
  Princeton, 1981.

\bibitem{sullivan:QCDII}
\bysame, \emph{Quasiconformal homeomorphisms and dynamics {II}: {S}tructural
  stability implies hyperbolicity for {K}leinian groups}, Acta Math.
  \textbf{155} (1985), 243--260.

\bibitem{susskind:limitsetints}
P.~Susskind, \emph{Kleinian groups with intersecting limit sets}, J. d'Analyse
  Math. \textbf{52} (1989), 26--38.

\bibitem{wpt:II}
W.~Thurston, \emph{Hyperbolic structures on 3-manifolds, {II}: surface groups
  and manifolds which fiber over the circle}, preprint,
  \mbox{arXiv:math.GT/9801045}.

\bibitem{wpt:III}
\bysame, \emph{Hyperbolic structures on 3-manifolds, {III}: deformations of
  3-manifolds with incompressible boundary}, preprint,
  \mbox{arXiv:math.GT/9801058}.

\bibitem{wpt:notes}
\bysame, \emph{The geometry and topology of 3-manifolds}, Princeton University
  Lecture Notes, online at http://www.msri.org/publications/books/gt3m, 1982.

\bibitem{wpt:bull}
\bysame, \emph{Three dimensional manifolds, {K}leinian groups and hyperbolic
  geometry}, Bull. Amer. Math. Soc. \textbf{6} (1982), 357--381.

\bibitem{wpt:I}
\bysame, \emph{Hyperbolic structures on 3-manifolds, {I}: deformation of
  acylindrical manifolds}, Ann. of Math. \textbf{124} (1986), 203--246.

\bibitem{wpt:textbook}
\bysame, \emph{{Three-Dimensional Geometry and Topology}}, Princeton University
  Press, 1997, (S. Levy, ed.).

\bibitem{tukia-vaisala:extension}
P.~Tukia and J.~V{\"a}is{\"a}l{\"a}, \emph{Lipschitz and quasiconformal
  approximation and extension}, Ann. Acad. Sci. Fenn. Ser. A I Math. \textbf{6}
  (1981), no.~2, 303--342 (1982).

\bibitem{vaisala:concordance}
J.~V{\"a}is{\"a}l{\"a}, \emph{Quasiconformal concordance}, Monatsh. Math.
  \textbf{107} (1989), no.~2, 155--168.

\bibitem{waldhausen}
F.~Waldhausen, \emph{On irreducible 3-manifolds which are sufficiently large},
  Ann. of Math. \textbf{87} (1968), 56--88.

\end{thebibliography}

\end{document}